# The ℝ-invariant solutions to the Kapustin-Witten equations on $(0,\infty) \times \mathbb{R}^2 \times \mathbb{R}$ with generalized Nahm-pole asymptotics

Clifford Henry Taubes[†]

Department of Mathematics
Harvard University
Cambridge, MA 02138

chtaubes@math.harvard.edu

ABSTRACT: This paper supplies a new characterization of the Kapustin-Witten equation solutions on $(0,\infty) \times \mathbb{R}^2 \times \mathbb{R}$ that play a key role in Edward Witten's program to obtain the Jones polynomial knot invariants using solutions to the same equation on $(0,\infty) \times S^3$.

[†]Supported in part by the National Science Foundation

# 1. Introduction

Two theorems in this paper supply novel characterizations of the model solutions to the Kapustin-Witten equations on $(0,\infty)\times\mathbb{R}^3$ that Edward Witten uses for his conjecture concerning these equations and the Jones polynomial knot invariant (see [W1-4], [GW] and [MW1] for more on this conjecture).

By way of motivation: The solutions of interest here are local models for the singular behavior of solutions of the analogous equations on the product of $(0,\infty)$ with a compact, oriented 3-dimensional Riemannian manifold. The latter appear in Witten's conjecture to the effect that a certain algebraic count of solutions with prescribed asymptotics involving a knot or link in the 3-manifold is an invariant of the knot which is the knot's Jones polynomial when the 3-manifold is $S^3$. Seminal papers of Mazzeo and Witten [MW1, 2] prove that this counting problem fits into a standard gauge theory package for computing low dimensional differential topology invariants. This is to say that a neighborhood of the automorphism group orbit of any given solution is homeomorphic to the zero locus of a smooth map from a ball about the origin in one finite dimensional Euclidean space to another finite dimensional Euclidean space; and that these Euclidan spaces are the respective kernel and cokernel of a Fredholm operator that is determined by the solution. A certain canonical section of the determinant line bundle for the Fredholm operator should then determine the algebraic weight of the given solution (Siqi He proved that this happens for the case for the unknot [H]). Granted that there is a well defined algebraic weight for each solution, then the issue with regards to defining a knot or link invariant is this: Is the count finite for any given depiction of the knot as an embedded loop, and for any given metric on the 3-manifold. And, if so, is the count constant along 1-parameter families embedded loops and 1-parameter families of metrics? This is a compactness issue: In what sense is the space of solutions (modulo bundle automorphisms) compact?

The analysis of the non-compactness (as approached by this author) requires a priori control of the asymptotics of solutions with respect to the 0 and $\infty$ limits of the $(0,\infty)$ factor. The asymptotics in the case of the unknot was studied in [T]. Suitable control in the case where a knot or link specifies the asymptotic conditions requires a priori control of certain scaling limits of solutions; and whether there is or isn't sufficient control leads, in turn, to the question that is studied here on $(0,\infty)\times\mathbb{R}^3$ which is this: Are the model solutions from [W1] the *unique* scaling limits? (Siqi He and Rafe Mazzeo have a uniqueness theorem for the model solutions (see [HM1]; and also [HM2]). Unfortunately, their theorem does not appear to be strong enough for the author's hypothetical application to the compactness question.)

**a) The Kapustin-Witten equations**

Suppose for the moment that X denotes a smooth, oriented, Riemannian 4-manifold and that $P \to X$ denotes a principal $SU(2)$ bundle. The Kapustin-Witten equations on X comprise a system of equations for a pair consisting of a connection on P (to be denoted by A) and 1-form on X with values in the adjoint bundle to P (this ad(P) valued 1-form is denoted by $\mathfrak{a}$). To write these equations, let $F_A$ denote curvature 2-form of a given connection A (it is an ad(P) valued 2-form); and let $D_A$ denote the exterior covariant as defined by the same connection A. Meanwhile, let $*_X$ denote the metric's Hodge dual operator. Any given pair $(A, \mathfrak{a})$ obeys the Kapustin-Witten equations when

$$F_A - \mathfrak{a} \wedge \mathfrak{a} = *_X D_A \mathfrak{a} \quad and \quad D_A *_X \mathfrak{a} = 0 \ .$$
(1.1)

The manifold X in this paper will be $(0, \infty) \times \mathbb{R}^3$ and the metric will be the Euclidean one that is defined using a chosen coordinate for the $(0, \infty)$ factor (to be denoted by t) and chosen coordinates $(x_1, x_2, x_3)$ for the $\mathbb{R}^3$ factor. The corresponding volume 4-form $dt \wedge dx_1 \wedge dx_2 \wedge dx_3$ defines the orientation. The bundle P here is necessarily isomorphic to the product principle $SU(2)$ bundle.

The solutions of interest are first constrained so that the dt component of $\mathfrak{a}$ is everywhere zero. This constraint is implicit in all of what follows. In this regard, the respective $dx_1$, $dx_2$ and $dx_3$ components of $\mathfrak{a}$ are denoted by $\mathfrak{a}_1$, $\mathfrak{a}_2$ and $\mathfrak{a}_3$, each being a section of ad(P). Thus $\mathfrak{a}$ can be written as $\mathfrak{a}_1 dx_1 + \mathfrak{a}_2 dx_2 + \mathfrak{a}_3 dx_3$.

Two additional global constraints can be written schematically as

$$F_{A3} \equiv 0 \quad and \quad \nabla_{A3} \mathfrak{a} \equiv 0 \ .$$
(1.2)

Here, $F_{A3}$ denotes the ad(P) valued 1-form that results when $F_A$ is paired with the tangent vector $\frac{\partial}{\partial x_3}$. Meanwhile, $\nabla_{A3}$ denotes the directional covariant derivative defined by A along this same tangent vector. (The respective directional covariant derivatives defined by A along the vector fields $\frac{\partial}{\partial x_1}$, $\frac{\partial}{\partial x_2}$ and $\frac{\partial}{\partial t}$ are denoted by $\nabla_{A1}$, $\nabla_{A2}$ and $\nabla_{At}$.) A pair that obeys the contraints in (1.2) is invariant with respect to an $\mathbb{R}$ action on the set of solutions to (1.1) that is induced by the $\mathbb{R}$ action on $(0, \infty) \times \mathbb{R}^3$ via constant translations of the coordinate $x_3$. (This action is generated by $\frac{\partial}{\partial x_3}$). This implies in particular that all norms and Aut(P)-invariant inner products are independent of the variable $x_3$. With the preceding in mind, the Euclidean space $\mathbb{R}^3$ is henceforth written as $\mathbb{R}^2 \times \mathbb{R}$ with $(x_1, x_2)$ being the coordinates on the $\mathbb{R}^2$ factor and with $x_3$ being the coordinate on the $\mathbb{R}$ factor. This is done to highlight the distinguished $x_3$ coordinate.

If $F_{A3}$ is identically zero, then $F_A$ can be written as

$$F_{A3} = E_{A1} dt \wedge dx_1 + E_{A2} dt \wedge dx_2 + B_{A3} dx_1 \wedge dx_2$$

(1.3)

with $E_{A1}$, $E_{A2}$ and $B_{A3}$ denoting respective sections of ad(P). (Each is annihilated by $\nabla_{A3}$ too.) Granted (1.2), then the equations in (1.1) can be written using these sections and with the help of $\varphi = \mathfrak{a}_1 - i\mathfrak{a}_2$ and $\varphi^* = \mathfrak{a}_1 + i\mathfrak{a}_2$ as follows:

- $\nabla_{At}\varphi = i[\mathfrak{a}_3, \varphi]$ .
- $(\nabla_{A1} + i\nabla_{A2})\varphi = 0$ .
- $\nabla_{At}\mathfrak{a}_3 = B_{A3} + \frac{i}{2}[\varphi, \varphi^*]$ .
- $E_{A1} = \nabla_{A2}\mathfrak{a}_3$
- $E_{A2} = -\nabla_{A1}\mathfrak{a}_3$

(1.4)

The model solutions to (1.4) are described momentarily as are other less useful solutions.

**b) Asymptotic constraints**

Of concern in what follows are solutions $(A, \mathfrak{a})$ to (1.2) that obey certain constraints on their behavior as $t \to 0$ or as $t \to \infty$ or as $|x_1|^2 + |x_2|^2 \to \infty$. Two different sets of constraints will be considered (the first set imposes some algebraic conditions too).

CONSTRAINT SET 1: *A pair $(A, \mathfrak{a})$ is described by this constraint set when*
- *The function $t|\mathfrak{a}|$ is uniformly bounded on $(0, \infty) \times \mathbb{R}^2 \times \mathbb{R}$.*
- *There exists $t_0 > 0$ and $\varepsilon > 0$ such that $t|\mathfrak{a}_3| > \varepsilon$ on $(0, t_0) \times \mathbb{R}^2 \times \mathbb{R}$.*
- *$\varphi$ is not identically zero.*
- *$\langle \mathfrak{a}_3 \varphi \rangle$ is identically zero.*

By way of notation for the second constraint set, let $D_R$ (for a given $R > 0$) denote the radius R disk about the origin in $\mathbb{R}^2$.

CONSTRAINT SET 2: *A pair $(A, \mathfrak{a})$ is described by this constraint set when*
- *The function $t|\mathfrak{a}|$ is uniformly bounded on $(0, \infty) \times \mathbb{R}^2 \times \mathbb{R}$.*
- *There exist positive numbers $t_0$, $\varepsilon$ and $r > 0$ such that $t|\mathfrak{a}| > \varepsilon$ where both $t < t_0$ and $(x_1^2 + x_2^2)^{1/2} \geq rt$ in $(0, \infty) \times \mathbb{R}^2 \times \mathbb{R}$.*
- *The function $|F_A|^2$ (which is $|B_{A3}|^2 + |E_{A1}|^2 + |E_{A2}|^2$) has finite integral on any given $R > 0$ version of the domain $(0, \infty) \times (\mathbb{R}^2 - D_R) \times \{0\}$; and that integral is bounded by an R-independent multiple of $\frac{1}{R}$.*

The model solutions that are described momentarily obey both constraint sets.

It is important to note that both constraint sets (and the constraints in (1.2)) are invariant under the action of Aut(P). This is to say that if $(A, \mathfrak{a})$ obeys (1.2) and either constraint set, then so does its pull-back under any automorphism of P.

Just as important: The constraint sets and (1.2) are scale invariant in the following sense: Any given positive number (call it $\lambda$) defines a corresponding coordinate rescaling diffeomorphism of $(0, \infty) \times \mathbb{R}^2 \times \mathbb{R}$ by the rule whereby $(t, x_1, x_2, x_3) \to (\lambda t, \lambda x_1, \lambda x_2, \lambda x_3)$. If $(A, \mathfrak{a})$ obeys (1.2) and (1.4) and if it is described by either constraint set, then the pull-back pair via such a scaling diffeomorphism also obeys (1.2) and (1.4) and the relevant constraint set. (Any given bundle isomorphism between P and its pull-back by a scaling diffeomorphism can be used to identify the pull-back pair as a pair of connection on P and ad(P) valued 1-form.)

**c) Witten's model solutions**

As noted at the outset, certain model solutions to (1.4) play a central role in Witten's program (see [W1 and also [MW1]). These are depicted momentarily. The depiction uses the complex coordinate $z = x_1 + i x_2$ for the factor of $\mathbb{R}^2$ in $\mathbb{R}^3$. It is also useful to define functions $\Theta$ and $x$ with values in $(0, \infty)$ by setting

$$\sinh\Theta = \tfrac{t}{|z|} \quad and \quad x = (t^2 + |z|^2)^{1/2}.$$

(1.5)

Thus, $|z| = x \tfrac{1}{\cosh\Theta}$ and $t = x \tfrac{\sinh\Theta}{\cosh\Theta}$. Note that $\Theta \to 0$ at fixed $x$ corresponds to the limit where $t \ll |z|$, where as $\Theta \to \infty$ at fixed $x$ corresponds to the limit where $|z| \ll t$. Keep in mind for what follows that $\Theta$ does not extend over $(0, \infty) \times \{0\} \times \mathbb{R}$ in $(0, \infty) \times \mathbb{R}^2 \times \mathbb{R}$.

To continue with notation: The product connection on P is denoted by $\theta_0$. Let $(\sigma_1, \sigma_2, \sigma_3)$ denote an orthonormal basis for the Lie algebra of SU(2) chosen so that the following commutation rules are in play:

$$[\sigma_1, \sigma_2] = -2\sigma_3 \quad and \quad [\sigma_2, \sigma_3] = -2\sigma_1 \quad and \quad [\sigma_3, \sigma_1] = -2\sigma_2 .$$

(1.6)

This Lie algebra basis will be viewed as a $\nabla_{\theta_0}$-covariantly constant, orthonormal basis for the (product) bundle associated, SU(2) Lie algebra bundle over $(0, \infty) \times \mathbb{R}^2 \times \mathbb{R}$.

Witten's model solutions are labeled by a non-negative integer to be denoted by $m$. The integer $m$ version is equivalent via the action of Aut(P) (the action of a gauge transformation) to the following solution: The version of ad(P)-valued 1-form $\mathfrak{a}$ is written as $\mathfrak{a} = \mathfrak{a}_1 dx_1 + \mathfrak{a}_2 dx_2 + \mathfrak{a}_3 dx_3$ with

- $\mathfrak{a}_3 = -\tfrac{1}{2t} \tfrac{(m+1)\sinh(\Theta)}{\sinh((m+1)\Theta)} \tfrac{\cosh((m+1)\Theta)}{\cosh(\Theta)} \sigma_3$ ,

- $\mathfrak{a}_1 - i\mathfrak{a}_2 = -\frac{1}{2t} \frac{(m+1)\sinh(\Theta)}{\sinh((m+1)\Theta)} \frac{z^m}{|z|^m} (\sigma_1 - i\sigma_2)$ ,

(1.7)

The corresponding connection A is written as $A = \theta_0 + \hat{A}$ with $\hat{A}$ denoting the 1-form with values in $ad(P_0)$ that is given by the formula:

$$\hat{A} = \frac{(m+1)}{2} \left(1 - \frac{\sinh(\Theta)}{\cosh(\Theta)} \frac{\cosh((m+1)\Theta)}{\sinh((m+1)\Theta)}\right) \frac{1}{|z|^2} (z\,d\bar{z} - \bar{z}\,dz) \sigma_3$$

(1.8)

Although (1.7) and (1.8) are defined a priori only on $(0,\infty) \times (\mathbb{R}^2 - 0) \times \mathbb{R}$, both $\mathfrak{a}$ and $\hat{A}$ extend smoothly to the whole of $(0,\infty) \times \mathbb{R}^2 \times \mathbb{R}$ so as to define a pair of smooth ad(P) valued 1-forms. (Write $\mathfrak{a}$ and $\hat{A}$ where $\Theta$ is large in terms of $x_1$ and $x_2$ to see this.) The corresponding $E_A$ and $B_A$ parts of the curvature 2-form can be computed from A:

- $B_A = \frac{(m+1)}{2x^2} \frac{\sinh(\Theta)}{\cosh(\Theta)} \frac{\cosh((m+1)\Theta)}{\sinh((m+1)\Theta)} \left(1 - \frac{(m+1)\sinh(2\Theta)}{\sinh(2(m+1)\Theta)}\right) \sigma_3 \, dx_3$

- $E_A = -\frac{(m+1)}{2x^2} \frac{\cosh((m+1)\Theta)}{\sinh((m+1)\Theta)} \left(1 - \frac{(m+1)\sinh(2\Theta)}{\sinh(2(m+1)\Theta)}\right) \sigma_3 \frac{1}{x} \frac{i}{2} (z\,d\bar{z} - \bar{z}\,dz)$ .

(1.9)

An integer $m$ model pair is henceforth denoted by $(A^{(m)}, \mathfrak{a}^{(m)})$.

The $m = 0$ version (which is $(A^{(0)}, \mathfrak{a}^{(0)})$) is special, it is the Nahm pole solution from [W1] and [MW2]. In this case, $(A^{(0)}, \mathfrak{a}^{(0)})$ is Aut(P) equivalent to the solution with $A^{(0)} \equiv \theta_0$ and $\{\mathfrak{a}^{(0)}{}_a = -\frac{1}{2t} \sigma_a\}_{a=1,2,3}$. Note in particular that $A^{(0)}$ is flat, that $\mathfrak{a}^{(0)}$ is $\nabla_{A^0}$-covariantly constant and that its components are orthogonal with norm equal to $\frac{1}{2t}$.

Various salient properties of the $m \geq 0$ versions of $(A^{(m)}, \mathfrak{a}^{(m)})$ are noted in Section 2a. In particular, the properties listed there can be used to verify that each such pair is described by both CONSTRAINT 1 and CONSTRAINT 2. There are two other crucial properites to note; the first being that $\varphi^{(0)}$ is nowhere zero and $\varphi^{(m)}$ for $m > 0$ is zero only on the locus $(0,\infty) \times \{0\} \times \mathbb{R}$. The second crucial property is this: The coordinate scaling diffeomorphisms (the maps sending $(t, z, x_3)$ to $(\lambda t, \lambda z, \lambda x_3)$ with $\lambda$ being a positive number) fix each $m \geq 0$ version of $(A^{(m)}, \mathfrak{a}^{(m)})$ in the following sense: The pull-back of $(A^{(m)}, \mathfrak{a}^{(m)})$ by the scaling diffeomorphism is the same as $(A^{(m)}, \mathfrak{a}^{(m)})$ if parallel transport by $A^{(m)}$ is used to identify P with its corresponding pull-back.

**d) One or maybe two other families of solutions**

There is a second family of solution to (1.4) on $(0,\infty) \times \mathbb{R}^2 \times \mathbb{R}$ that plays a tangential role in the upcoming story. These are the *Nahm pole imposter* solutions. The Nahm pole imposters are parameterized by the complement of the origin in the closed disk $\bar{D} = \{w \in \mathbb{C}: |w| \leq 1\}$. The $w = 0$ point in $\bar{D}$ parametrizes the Nahm pole solution which is the $m = 0$ version of the solution that is depicted in (1.7) and (1.8). This $\bar{D}$- parametrization of solutions is such that each point in the interior of $\bar{D}$ corresponds to a unique automorphism equivalence class of solutions, where as all of the points on the boundary of $\bar{D}$ parametrize the same equivalence class. For $w \in \bar{D}$, the corresponding Nahm-pole imposter is Aut(P) equivalent to a pair that is defined as follows:

- $\mathfrak{a}_3 = -\frac{1}{2t}\sigma_3$.
- $\mathfrak{a}_1 - i\mathfrak{a}_2 = -\frac{1}{2t}(1-|w|^2)^{1/2}(\sigma_1 - i\sigma_2)$.
- $A = \theta_0 + \mathfrak{b}$ *where*
  a) $\mathfrak{b}_t = \mathfrak{b}_3 = 0$.
  b) $\mathfrak{b}_1 + i\mathfrak{b}_2 = -\frac{1}{2t}w(\sigma_1 - i\sigma_2)$.

(1.10)

These solutions are described by CONSTRAINT SET 1 when $|w| < 1$; the $|w| = 1$ version violates the third bullet of CONSTRAINT SET 1. But all of these Nahm pole imposters violate the third bullet in CONSTRAINT SET 2. This fact distinguishes these Nahm pole imposters from the honest Nahm pole solution and from the $m \geq 1$ model solutions.

There is another, and purely hypothetical family of solutions to (1.4) that play an important role in the upcoming theorems should they exist. For lack of a better term, these hypothetical solutions are called *Kapustin-Witten instantons*. These solutions obey the conditions in CONSTRAINT SET 2. Distinguishing features are:

- *No Kapustin-Witten instanton is invariant under the coordinate rescaling diffeomorphism.*
- *A Kapustin-Witten instanton solution converges in a certain strong sense as* $t \to \infty$ *on* $(0,\infty) \times \mathbb{R}^2 \times \mathbb{R}$ *(with* $\frac{|z|}{t}$ *bounded) to some* $m \geq 0$ *model solution; and it converges as* $t \to 0$ *(also with* $\frac{|z|}{t}$ *bounded) in a slightly less strong sense to some* $m+2p$ *model solution with* p *being a positive integer*.
- *A Kapustin-Witten instanton converges as* $\frac{|z|}{t} \to \infty$ *to the Nahm pole model solution*.

(1.11)

More is said about these hypothetical solutions and their small and large t behavior and their large $\frac{|z|}{t}$ behavior in subsequent parts of this paper. I hope to address the existence question elsewhere. (Edward Witten and Davide Gaiotto also conjecture that these instantons exist; Edward Witten wrote that this conjecture follows readily from the observations in Section 2 of their joint paper [GW].)

One point to note is that the components of the moduli space of Kapustin-Witten instantons (the set of solutions modulo the action of Aut(P)) should have dimension at least 2 because there is a free circle action on the moduli spaces (if non-empty) which commutes with the $\mathbb{R}$-action given by coordinate rescalings. (The moduli spaces of Kapustin-Witten instantons that interpolate between the integer $m+2p$ model solution at small t and the integer $m$ model solution at large t should have dimension at least $2p$.)

**e) Characterizations of the model solutions**

The two theorems that follow momentarily are the central results of this paper. The first theorem implicitly refers to a crucial implication of the first bullet in (1.4) which is this: If $(A,\mathfrak{a})$ obeys (1.4), then the zero locus of $\varphi$ in $(0,\infty)\times\mathbb{R}^2\times\mathbb{R}$ has the form $(0,\infty)\times Z\times\mathbb{R}$ with Z being a subset of $\mathbb{R}^2$. (The norm of $\frac{\partial}{\partial t}|\varphi|^2$ is bounded by $4|\mathfrak{a}_3||\varphi|^2$ by virtue of that first bullet. This leads in turn to a $4|\mathfrak{a}_3|$ bound for the t-derivative of $\ln|\varphi|^2$; and that bound implies that $\ln|\varphi|^2$ is finite at any given $(t, x_1, x_2, x_3)$ iff and only if it is finite at any given $(t´, x_1, x_2, x_3)$.) All of the model solutions in Section 1d have $Z = \{0\}$. Here is the first theorem:

**Theorem 1**: *Suppose that $(A,\mathfrak{a})$ is a pair of connection on P and ad(P)-valued 1-form (with no dt component) that obeys (1.1) and is described by (1.2) and* CONSTRAINT SET 1. *If the corresponding version of Z is the origin in $\mathbb{R}^2$, then $(A,\mathfrak{a})$ is Aut(P) equivalent to some* $m > 0$ *model solution. If corresponding version of Z is empty, then $(A,\mathfrak{a})$ is a Nahm pole solution (which is the* $m = 0$ *model solution) or one of the Nahm pole imposters*.

Section 4 has a generalization of Theorem 1 asserting this: All solutions to (1.1) that are described by (1.2) and CONSTRAINT SET 1 and with the corresponding set Z being non-empty and finite are characterized up to the action of Aut(P) by the set Z and an associated set positive integer multiplicities assigned to each of its points. (A version of Theorem 1 and its generalization in Section 4 with stronger initial assumptions on the solution was proved by Siqi He and Rafe Mazzeo in [HM1] and also [HM2].)

The second theorem says in effect that CONSTRAINT SET 2 solutions to (1.2) are either model solutions from Section 1c or Kapustin-Witten instantons from Section 1d. By way of notation: The second theorem uses $P/\{\pm1\}$ to denote quotient of P by the action of multiplication by the central automorphisms $\{\pm1\}$. This is a priinicpal SO(3) bundle. Automorphisms of the latter also act on the space of solutions to (1.2) (preserving the CONSTRAINT 2 conditions) because any automorphism of $P/\{\pm1\}$ lifts locally to an automorphism of P with any two lifts differing by the action of the central element -1 which acts as the identity on the space of connections on P and sections of ad(P). The notation uses g*A and g*$\mathfrak{a}$ to denote the respective pull-backs of a connection (which is A) and 1-form valued section of ad(P) (which is $\mathfrak{a}$) by a given automorphism, g, of either P or $P/\{\pm1\}$.

**Theorem 2**: *Suppose that $(A, \mathfrak{a})$ is a pair of connection on P and ad(P)-valued 1-form (with no dt component) that obeys (1.1) and is described by (1.2) and* CONSTRAINT SET 2. *Then $(A, \mathfrak{a})$ is one of the $m \geq 0$ model solutions from Section 1c if and only if there exists $t_0 > 0$ such that $\langle \mathfrak{a}_3[\mathfrak{a}_1, \mathfrak{a}_2] \rangle \geq 0$ where $t < t_0$ on $(0, \infty) \times \mathbb{R}^2 \times \mathbb{R}$. If $(A, \mathfrak{a})$ is not a model solution from Section 1c, then there exist a non-negative integer m, a positive integer p and automorphisms $g_\infty$ and $g_0$ of P defined on the respective $t > 1$ and $t < 1$ parts of $(0, \infty) \times \mathbb{R}^2 \times \mathbb{R}$, and an automorphism $g_>$ of $P/\{\pm 1\}$ defined in the $|z| > t$ part; and these are such that the following are true:*

- $\lim_{t \to \infty} (t|g_\infty^* A - A^{(m)}| + t|g_\infty^* \mathfrak{a} - \mathfrak{a}^{(m)}| + |\ln(\frac{|\varphi|}{|\varphi^{(m)}|})|) = 0.$
- $\lim_{t \to 0} (t|g_0^* A - A^{(m+2p)}| + t|g_0^* \mathfrak{a} - \mathfrak{a}^{(m+2p)}|) = 0$
- $\lim_{|z|/t \to \infty} (t|g_>^* A - A^{(0)}| + t|g_>^* \mathfrak{a} - \mathfrak{a}^{(0)}|) = 0.$

*(The limits in the top two bullets are uniform with respect to the $\mathbb{R}^2$ coordinate z, and the limit in the third bullet is uniform with respect to $x = (t^2 + |z|^2)^{1/2}$.)*

    With regards to the constraint assumptions in Theorems 1 and 2: There are solutions to (1.4) that violate one or more of the bulleted items in either constraint set. With regards to the first constraint in either constraint set: Given a positive number to be denoted by $c$, there is the Nahm pole like solution with $\mathfrak{a}_3 = -\frac{c}{4} \frac{1+e^{-ct}}{1-e^{-ct}} \sigma_3$ and $\varphi = \frac{c}{\sqrt{2}} \frac{e^{-ct/2}}{1-e^{-ct}} (\sigma_1 - i\sigma_2)$, and with A being $\theta_0$. (There are also solutions like this that look like Nahm pole imposters.) With regards to the second constraint in either constraint set: Fix a time $\tau > 0$ and pull-back a model solution (or a Nahm pole imposter) by the re-imbedding of $(0, \infty) \times \mathbb{R}^2 \times \mathbb{R}$ into itself that sends t to $t + \tau$. The resulting pull-back of a model solution obeys all but the second bullet in either constraint set. (The solution has a smooth extension to the domain $(-\tau, \infty) \times \mathbb{R}^2 \times \mathbb{R}$.) As noted before, the Nahm pole imposters are all rescaling invariant and they all violate the third bullet of CONSTRAINT SET 2. There are other constructions: For an Abelian family with members fixed by all rescaling diffeomorphism, take a positive number, $r$, and then set

$$\mathfrak{a}_3 = -\frac{r}{x} \sigma_3 \text{ and } \varphi \equiv 0 \text{ and } A = \theta_0 + r(1 - \frac{\sinh\Theta}{\cosh\Theta}) \frac{1}{|z|^2} (z d\bar{z} - \bar{z} dz) \sigma_3.$$

(1.12)

These run afoul of the second bullet constraint in both constraint sets (and the third bullet in the first constraint set).

    With regards to the sign of $\langle \mathfrak{a}_3[\mathfrak{a}_1, \mathfrak{a}_2] \rangle$ in Theorem 2: In the case of the Nahm pole solution, this function is everywhere positive; and it is positive except on the $z = 0$ locus for the integer $m > 0$ model solutions. (The function has to vanish at $z = 0$ because $\mathfrak{a}_1$ and $\mathfrak{a}_2$ vanish there.) A hypothetical Kapustin-Witten instanton that interpolates between an

integer *m* model solution as t → ∞ and an integer *m*+2*p* model solution as t → 0 will have positive $\langle \mathfrak{a}_3[\mathfrak{a}_1,\mathfrak{a}_2]\rangle$ for sufficiently small t on the domain where |z| is greater than $\mathcal{O}(t^{(2p+m+1)/p})$; and then negative $\langle \mathfrak{a}_3[\mathfrak{a}_1,\mathfrak{a}_2]\rangle$ for smaller |z| except for the z = 0 locus where it also has to vanish in the case when *m* > 0.

With regards to the $C^0$ nature of the limits in Theorem 2: The convergence is stronger than this. See the Propositions 8.1 and 8.2. See the Appendix for even more about the very small t behavior of (hypothetical) Kapustin-Witten instantons.

**f) The rest of this paper**

This paper has nine sections in all (including this introduction) plus a short appendix. Theorem 1 is proved in Section 4 and Theorem 2 is proved in Section 9. The intervening sections supply the tools for these proofs. Here is the table of contents:

1. Introduction.
2. Some preliminary observations.
3. The general structure of a solution.
4. When $\langle \mathfrak{a}_3\varphi\rangle \equiv 0$.
5. CONSTRAINT SET 2 solutions: The norm of |φ|.
6. Norms of β and $\mathfrak{b}$.
7. The natural logarithm of |φ|.
8. The t → ∞ and t → 0 limits of (A, $\mathfrak{a}$).
9. Proof of Theorem 9.1
Appendix: The small t behavior of Kapustin-Witten instantons.
References.

**g) Conventions**

The first set of conventions concern the Lie algebra of the group SU(2). This is vector space is denoted in what follows by $\mathfrak{S}$; and it is identified implicitly with the vector space of 2 × 2, anti-hermitian matrices. The commutator of elements in $\mathfrak{S}$ is denoted by [·, ·]. If σ is any given 2×2, ℂ-valued matrix, the notation $\langle\sigma\rangle$ means this:

$$\langle\sigma\rangle = -\tfrac{1}{2}\text{trace}(\sigma) .$$

(1.13)

The inner product of two elements σ, τ ∈ $\mathfrak{S}$ is taken to be $\langle\sigma\tau\rangle$. A basis for $\mathfrak{S}$ obeying (1.6) is orthonormal with respect to this inner product. (Note in particular the minus signs in (1.6). This is a non-standard choice for a basis!) This inner product on $\mathfrak{S}$ defines a fiberwise inner product on ad(P). The inner product between elements σ and τ in the fiber of ad(P) over any given point is denoted by $\langle\sigma\tau\rangle$ also.

Inner products on vectors, differential forms and higher rank tensors on $(0,\infty)\times\mathbb{R}^2\times\mathbb{R}$ are defined using the Euclidean metric that is defined so that dt, $dx_1$, $dx_2$ and $dx_3$ are orthonormal. The inner product between tensor $\mathfrak{t}$ and $\mathfrak{h}$ of any given rank is denoted by $\langle\mathfrak{t},\mathfrak{h}\rangle$. This Euclidean inner product and the previously defined fiber inner product for ad(P) gives an inner product on ad(P) valued tensors of any rank. If $\mathfrak{t}$ and $\mathfrak{h}$ are any two such tensors, their inner product is also denoted by $\langle\mathfrak{t},\mathfrak{h}\rangle$.

With regards to covariant derivatives: Supposing that $\mathfrak{t}$ is an tensor, its covariant derivative (using the Euclidean metric's flat connection) is denoted by $\nabla\mathfrak{t}$. Supposing that $\mathfrak{t}$ is an ad(P) valued tensor, the covariant derivative using the connnection A and the Euclidean metric's connection is denoted by $\nabla_A$. (It has components $\nabla_{At}$, $\nabla_{A1}$, $\nabla_{A2}$ and $\nabla_{A3}$ with respect to the basis dt, $dx_1$, $dx_2$ and $dx_3$. But $\nabla_{A3}$ will act as zero on every tensor and ad(P) valued tensor of interest.) The covariant derivative along the constant t and slices of $(0,\infty)\times\mathbb{R}^2\times\mathbb{R}$ is denoted by $\nabla_A^\perp$ (its components are $\nabla_{A1}$, $\nabla_{A2}$ and $\nabla_{A3}$).

Continuing with conventions: A convention used throughout this article concerns the choice of a standard 'cut-off' function on $\mathbb{R}$. This is a non-increasing function that is equal to 1 on $(-\infty, \frac{1}{4}]$ and equal to 0 on $[\frac{3}{4}, \infty)$. Such a function should be chosen now and then used throughout. The chosen function is denoted by $\chi$. Since the precise choice has no bearing on subsequent discussions, please choose your favorite. All 'bump' functions and cut-off functions are (sometimes implicitly) made from $\chi$ by rescalings (using $\chi(rt)$ in lieu of $\chi$ with $r$ a positive number) and taking suitable products of rescalings of $\chi$ and/or $(1-\chi)$. This insures that their derivatives to any given order have uniform norm bounds after accounting for the rescalings.

The next convention concerns the notation for numbers that are independent of all relevant quantities in any given inequality (relevant quantities can be coordinates such as t, or the choice of $(A,\mathfrak{a})$, etc.). Such a number is denoted in what follows by $c_0$. It is always greater than 1 and it can be assumed to increase between successive appearances.

Here is a final convention: The set of positive integers, the set $\{1, 2, \ldots\}$, is denoted in this paper by $\mathbb{N}$.

**h) Acknowledgements**

Many thanks to Edward Witten and Davide Gaiotto for sharing their thoughts about the hypothetical Kapustin-Witten instantons.

## 2. Some preliminary observations

The first subsection below summarizes some of the salient properties of the model solutions from (1.7) and (1.8). The remaining subsections make observations about general properties of solutions and sequences of solutions to (1.1) and (1.4).

### a) Properties of the model solutions

The salient properties of the $m > 0$ model solutions are as follows:

- At any fixed $t > 0$ (and $x_3 \in \mathbb{R}$), the function $|\varphi^{(m)}|$ is zero only at the origin; and it is a monotonically increasing function of the distance to the origin in $\mathbb{R}^2$ with limit $\frac{1}{\sqrt{2t}}$ as the distance limits to $\infty$.
- $\mathfrak{a}_3$ can be written as $\alpha \sigma_3$ with $\alpha$ at any fixed $t > 0$ (and $x_3 \in \mathbb{R}$) being negative, and being a monotonically decreasing function of the distance to the origin in $\mathbb{R}^2$ with $\alpha = -\frac{(m+1)}{2t}$ at the origin in $\mathbb{R}^2$ and with limit $-\frac{1}{2t}$ as the distance limits to $\infty$.
- The norms of the curvature (which is $F_A$) and the A-covariant derivative of $\mathfrak{a}$ (which is $\nabla_A \mathfrak{a}$) are bounded by $\frac{c_m}{x^2}$ with $c_m$ being a constant.

(2.1)

These bulleted facts have the following implication which is good to keep in mind: At any fixed $t > 0$ (and $x_3 \in \mathbb{R}$), any $m > 0$ version of $(A^{(m)}, \mathfrak{a}^{(m)})$ looks more and more like the Nahm pole solution $(A^{(0)}, \mathfrak{a}^{(0)})$ as the distance to the origin gets ever larger. (This can be made precise using what is said in Section 2d.)

### b) A second order equation for the components of $\mathfrak{a}$

The first order system that is depicted in (1.1) leads to a second order Laplacian-like differential equation for $\mathfrak{a}$. In the context of this paper, the equation is this:

$$-(\nabla_{At}\nabla_{At} + \nabla_{A1}\nabla_{A1} + \nabla_{A2}\nabla_{A2})\mathfrak{a} + [\mathfrak{a}_t,[\mathfrak{a},\mathfrak{a}_t]] + [\mathfrak{a}_1,[\mathfrak{a},\mathfrak{a}_1]] + [\mathfrak{a}_2,[\mathfrak{a},\mathfrak{a}_2]] = 0 .$$

(2.1)

This is derived by taking the exterior covariant derivative of both sides of (1.1) and then employing (1.1) to write resulting curvature terms that arise when commuting covariant derivatives.

### c) Elliptic regularity

The proposition that follows directly makes a statement to the effect that a suitable bound for the pointwise norm of $\mathfrak{a}$ on a given ball in $(0,\infty) \times \mathbb{R}^2 \times \mathbb{R}$ determines an upper bound for the norms of $F_A$ and $\nabla_A \mathfrak{a}$ and their higher order, $\nabla_A$-covariant derivatives on the half-radius, concentric ball.

**Proposition 2.1**: *There exists $\kappa > 10$ with the following significance: Suppose that $r > 0$ and that $B_r$ is a ball in $(0,\infty) \times \mathbb{R}^2 \times \mathbb{R}$ with compact closure. Let $(A,\mathfrak{a})$ denote a solution to (1.1) on $B_r$ with the norm of $\mathfrak{a}$ obeying $|\mathfrak{a}| < \frac{1}{\kappa r}$.*
- *$|F_A|$ and $|\nabla_A \mathfrak{a}|$ are bounded by $\frac{\kappa}{r^2}$ on the concentric ball with radius $\frac{1}{2} r$.*

- *Let* n *denote a positive integer. Both* $|\nabla_A^{\otimes(n+1)}\mathfrak{a}|$ *and* $|\nabla_A^{\otimes n}F_A|$ *are bounded by* $\kappa_n \frac{1}{r^{n+2}}$ *on the concentric, radius* $\frac{1}{2}r$ *ball with* $\kappa_n$ *being constant that is independent of* $(A,\mathfrak{a})$ *and* $B_r$ *and* r.
- *There is an automorphism of* P *over* B *that pulls back* A *to have the form* $\theta_0 + \hat{A}$ *with* $\hat{A}$ *being an* ad(P) *valued 1-form on* B *whose norm on the concentric, radius* $\frac{1}{2}r$ *ball is bounded by* $\frac{\kappa}{r}$ *and whose* $\nabla_{\theta_0}$ *-covariant derivative to any given order (call it* n*) on this same concentric ball is bounded by* $\kappa_n \frac{1}{r^{n+1}}$.

**Proof of Proposition 2.1**: This proposition is an instance of what is asserted by Proposition 2.1 in [T]. The latter proof in the context of this paper amounts to a bootsrapping argument that starts by taking the inner product of both sides of (2.1) with $\mathfrak{a}$, multiplying the result by a bump function with support in the given ball, and then integrating by parts to bound the integral of $|\nabla_A \mathfrak{a}|^2$ in a smaller, concentric ball by the integral of $|\mathfrak{a}|^2$ over the larger ball. Then, note that a supremum norm bound for $|\mathfrak{a}|$ over the ball and a bound for the integral $|\nabla_A \mathfrak{a}|^2$ over the smaller ball leads via (1.1) to a bound for the integral of $|F_A|^2$ over the smaller ball too. This bound with Karen Uhlenbeck's work in [U] can then be used to get pointwise bounds and integral bounds for the norms of higher order covariant derivatives.

The following lemma is essentially an instance of the preceding proposition. Think about this lemma in conjunction with what is required by the first bullet of either CONSTRAINT SET 1 or CONSTRAINT SET 2.

**Lemma 2.2**: *Given* $\mathfrak{z} > 0$, *there exists* $\kappa_\mathfrak{z} > 1$ *with the following significance: Let* $(A,\mathfrak{a})$ *denote a solution to (1.4) on* $(0,\infty) \times \mathbb{R}^2 \times \mathbb{R}$ *with the norm of* $\mathfrak{a}$ *obeying* $|\mathfrak{a}| \leq \frac{\mathfrak{z}}{t}$.
- $|F_A|$ *and* $|\nabla_A \mathfrak{a}|$ *are bounded by* $\frac{\kappa_\mathfrak{z}}{t^2}$ *on* $(0,\infty) \times \mathbb{R}^2 \times \mathbb{R}$.
- *If* n *is a positive integer, then both* $|\nabla_A^{\otimes(n+1)}\mathfrak{a}|$ *and* $|\nabla_A^{\otimes n}F_A|$ *are bounded by* $\kappa_{\mathfrak{z},n} \frac{1}{t^{n+2}}$ *on* $(0,\infty) \times \mathbb{R}^2 \times \mathbb{R}$ *with* $\kappa_{\mathfrak{z},n}$ *being constant. It can depend on* $\mathfrak{z}$ *but it is otherwise independent of* $(A,\mathfrak{a})$.

*Proof of Lemma 2.2*: Let $\kappa_*$ denote Proposition 2.1's version of $\kappa$. Fix $t > 0$ and let r denote the number $\frac{1}{4\kappa_*\mathfrak{z}}t$. The assumptions are such that $|\mathfrak{a}| < \frac{1}{\kappa_* r}$ in the radius r ball centered at any point in $\{t\} \times \mathbb{R}^2 \times \mathbb{R}$. Granted this, then the bullets of the lemma follow directly from what is said by Proposition 2.1.

### d) Limits of sequences

Sequential compactness is an important tool that will be brought to bear in subsequent arguments. The lemma that follows directly supplies the relevant sequential compactness theorem for solutions to (2.1) on $(0,\infty) \times \mathbb{R}^2 \times \mathbb{R}$.

**Lemma 2.3**: *Fix $\mathfrak{z} > 1$ and let $\{(A^n, \mathfrak{a}^n)\}_{n \in \mathbb{N}}$ denote a sequence of solutions to (1.1) on $(0,\infty) \times \mathbb{R}^2 \times \mathbb{R}$ such that $|\mathfrak{a}^n| \leq \frac{\mathfrak{z}}{t}$ for each $n \in \mathbb{N}$. There is a subsequence $\Lambda \subset \mathbb{N}$, a corresponding sequence of automorphisms of P (denoted by $\{g_n\}_{n \in \Lambda}$) such that $\{(g_n^*A^n, g_n^*\mathfrak{a}^n)\}_{n \in \Lambda}$ converges in the $C^\infty$ topology on compact subsets of $(0,\infty) \times \mathbb{R}^2 \times \mathbb{R}$ to a solution to (1.1) on $(0,\infty) \times \mathbb{R}^2 \times \mathbb{R}$. This limit (denoted by $(A^\infty, \mathfrak{a}^\infty)$) obeys $|\mathfrak{a}^\infty| \leq \frac{\mathfrak{z}}{t}$ also. In addition, suppose that every $(A, \mathfrak{a})$ from the set $\{(A^n, \mathfrak{a}^n)\}_{n \in \mathbb{N}}$ obeys the same inequality from the following list:*

- *$|\mathfrak{a}_3| \geq \frac{\delta}{t}$ where $t \leq t_0$ with $\delta > 0$ and $t_0 > 0$ being independent of n.*
- *$|\varphi|_{(t=t_0, z=z_0)} \geq \frac{\delta}{t_0}$ with $\delta > 0$ and $t_0 > 0$, with $t_0$, $z_0$ and $\delta$ being independent of n.*
- *$\displaystyle\int_{(0,\infty) \times (\mathbb{R}^2 - D_R) \times \{0\}} |F_A|^2 < \frac{\Xi}{R}$ for some $R > 0$, with R and $\Xi$ being independent of n,*

*Then that same inequality is also obeyed when $(A, \mathfrak{a}) = (A^\infty, \mathfrak{a}^\infty)$.*

**Proof of Lemma 2.3**: Given any fixed pair of positive times $t_1 < t_2$, there is a bound on the norm $t|\mathfrak{a}^n|$ on $[t_1, t_2] \times \mathbb{R}^2 \times \mathbb{R}$ that is independent of the index n. This implies (via Lemma 2.2) that there are index n-independent bounds on $[t_1, t_2] \times \mathbb{R}^2 \times \mathbb{R}$ for the norms of $F_{A^n}$ and $\nabla_{A^n} \mathfrak{a}^n$ and their $\nabla_{A^n}$-covariant derivatives to any given order. These bounds (plus Karen Uhlenbeck's compactness theorem [U]) can be used to infer the existence of a subsequence $\Lambda \subset \mathbb{N}$ and a corresponding sequence of automorphisms of P (this is the sequence $\{g_n\}_{n \in \Lambda}$) such that $\{(g_n^*A^n, g_n^*\mathfrak{a}^n)\}_{n \in \Lambda}$ converges in the $C^\infty$ topology on compact subsets of $(0,\infty) \times \mathbb{R}^2 \times \mathbb{R}$ to a smooth pair of connection on P and ad(P)-valued 1-form. This is the limit pair $(A^\infty, \mathfrak{a}^\infty)$. By virtue of the $C^\infty$ convergence on compact sets, the pair $(A, \mathfrak{a})$ obeys the equations in (1.4); and it obeys $|\mathfrak{a}^\infty| \leq \frac{\mathfrak{z}}{t}$. For the same reason, it obeys any of the three bulleted inequalities that is obeyed by all of the $\{(A^n, \mathfrak{a}^n)\}_{n \in \mathbb{N}}$ pairs. (The verification of third bullet's inequality also uses the dominated convergence theorem.)

### d) Translations and rescalings of solutions

Sequences of solutions to (1.1) or (1.4) on $(0,\infty) \times \mathbb{R}^2 \times \mathbb{R}$ for use in Lemma 2.3 can be obtained from a fixed solution using a sequence of its translations and rescalings. These are the sequences that described momentarily. The next paragraph is a digression that comes first to make things more precise.

The group of rigid translations of the $\mathbb{R}^2$ factor of $(0,\infty)\times\mathbb{R}^2\times\mathbb{R}$ acts on the set of solutions to (1.4) (or (1.1)) because the pull-back by a translation of any given pair obeying (1.4) (or (1.1)) also obeys the (1.4) (or (1.1) as the case may be. (The pull-back of a pair $(A,\mathfrak{a})$ via a diffeomorphism is technically a pair consisting of a connection on the pull-back of the product principle bundle and 1-form with values in the adjoint bundle of the pull-back principle bundle. The convention used implicitly here is to identify the pull-back bundle with the original product principle bundle via parallel transport using the product connection $\theta_0$.) As noted previously, the group of coordinate rescalings of $(0,\infty)\times\mathbb{R}^2\times\mathbb{R}$ also acts on the set of solutions to (1.4) (or (1.1)). By way of a reminder, this group is isomorphic to $(0,\infty)$. The action of any given $\lambda\in(0,\infty)$ sends $(t, z, x_3)$ to the point $(\lambda t, \lambda z, \lambda x_3)$. The translation and rescaling actions together give an action of the semi-direct product of $(0,\infty)\times\mathbb{R}^2$ on the set of solutions of (1.4) (or 1.1).

Suppose now that $(A,\mathfrak{a})$ is a given solution to (1.1) on $(0,\infty)\times\mathbb{R}^2\times\mathbb{R}$ with the norm of $\mathfrak{a}$ obeying the bound $|\mathfrak{a}|\le\frac{3}{t}$. Let $\phi: (0,\infty)\times\mathbb{R}^2\times\mathbb{R} \to (0,\infty)\times\mathbb{R}^2\times\mathbb{R}$ denote a combination of translation along the $\mathbb{R}^2$ factor and coordinate rescaling. The norm of $\phi^*\mathfrak{a}$ also obeys the bound $|\phi^*\mathfrak{a}|\le\frac{3}{t}$. Therefore, if $\{\phi_n\}_{n\in\mathbb{N}}$ is any sequence of translations along the $\mathbb{R}^2$ factor of $(0,\infty)\times\mathbb{R}^2\times\mathbb{R}$ and rescalings of $(0,\infty)\times\mathbb{R}^2\times\mathbb{R}$, then the corresponding sequence $\{(A^n,\mathfrak{a}^n) = (\phi_n^*A, \phi_n^*\mathfrak{a})\}_{n\in\mathbb{N}}$ can be used as input for Lemma 2.3. In that event, Lemma 2.3 supplies a subsequence $\Lambda\subset\mathbb{N}$, a sequence $\{g_n\}_{n\in\mathbb{N}}$ of automorphisms of P and a solution, $(A^\infty,\mathfrak{a}^\infty)$, to (1.4) that is the limit (in the $C^\infty$ topology on compact subsets of $(0,\infty)\times\mathbb{R}^2\times\mathbb{R}$) of the sequence $\{(g_n^*A_n, g_n^*\mathfrak{a}_n)\}_{n\in\mathbb{N}}$. As noted by Lemma 2.3, the norm of $\mathfrak{a}^\infty$ also obeys the bound $|\mathfrak{a}^\infty|\le\frac{3}{t}$.

By way of an example: Suppose that $(A,\mathfrak{a})$ is one of the $m\ge 0$ model solutions from Section 1c. Let $\{\phi_n\}_{n\in\mathbb{N}}$ denote any sequence of translations along the $\mathbb{R}^2$ factor of $(0,\infty)\times\mathbb{R}^2\times\mathbb{R}$ and coordinate rescalings of $(0,\infty)\times\mathbb{R}^2\times\mathbb{R}$. For any given positive integer n, write $\phi_n$ as the diffeomorphism $(t, z, x_3) \to (\lambda_n t, \lambda_n(p-a_n), \lambda_n x_3)$ with p denoting here a point in $\mathbb{R}^2$. The sequence $\{(A^n,\mathfrak{a}^n)=(\phi_n^*A, \phi_n^*\mathfrak{a})\}_{n\in\mathbb{N}}$ can be used as input for Lemma 2.3 because of what is said by the top two bullets of (2.1). There are only two possibilities for the resulting limit $(A^\infty,\mathfrak{a}^\infty)$: It is a rigid translation of $(A,\mathfrak{a})$ along the $\mathbb{R}^2$ factor of $(0,\infty)\times\mathbb{R}^2\times\mathbb{R}$ if the sequence $\{a_n\}_{n\in\Lambda}$ converges; and it is the Nahm pole solution otherwise.

## 3. The general structure of a solution

The upcoming Proposition 3.1 makes preliminary observations about the structure of any given solution to (1.4). Much of what is said by this proposition can be found (in one form or another) either in [W1] or [MW1] or [HM1,2].

Proposition 3.1 refers to an eigenbundle decomposition of ad(P)$_\mathbb{C}$ (the complexification of ad(P)) that can be defined from any given unit length section of ad(P). Supposing that $\sigma$ denotes the given section, then the eigenbundles are defined with respect to the Hermitian endomorphism of ad(P)$_\mathbb{C}$ that sends any given $\eta$ from ad(P)$_\mathbb{C}$ to $[\frac{i}{2}\sigma, \eta]$. This endomorphism has eigenvalues 0, 1 and -1 at each point. The eigenvalue 0 subbundle is the span of $\sigma$. The eigenvalue +1 subbundle is denoted by $\mathcal{L}^+$; the eigenvalue -1 subbundle is $\mathcal{L}^-$ which is the Hermitian conjugate of $\mathcal{L}^+$. The bundles $\mathcal{L}^+$ and $\mathcal{L}^-$ are complex line bundles over $(0,\infty)\times\mathbb{R}^2\times\mathbb{R}$.

**Proposition 3.1**: *Suppose that* (A, $\mathfrak{a}$) *obeys bullets 1,2,4,and 5 of (1.4) on* $(0,\infty)\times\mathbb{R}^2\times\mathbb{R}$, *that* $t|\mathfrak{a}|$ *is a bounded function on* $(0,\infty)\times\mathbb{R}^2\times\mathbb{R}$, *and that* $\varphi \equiv \mathfrak{a}_1 - i\mathfrak{a}_2$ *is not identically zero*.
- *The sections* $\{\mathfrak{a}_1, \mathfrak{a}_2\}$ *of* ad(P$_0$) *are orthogonal (pointwise); and* $|\mathfrak{a}_1|^2 = |\mathfrak{a}_2|^2$. *(These identities say that* $\langle\varphi\varphi\rangle \equiv 0$.)
- *If* $\varphi$ *vanishes at a given point* (t, p, $x_3$) *in* $(0,\infty)\times\mathbb{R}^2\times\mathbb{R}$ *then* $\varphi$ *vanishes on the whole of* $(0,\infty)\times\{p\}\times\mathbb{R}$.
- *Fix* (t, $x_3$) *and view* $\varphi|_{(t,x_3)}$ *as a section of* ad(P)$_\mathbb{C}$ *on* $\mathbb{R}^2$ *by varying only* $x_1$ *and* $x_2$.
    a) *The set of zeros of* $\varphi|_{(t,x_3)}$ *on* $\mathbb{R}^2$ *is a discrete set*.
    b) *View* $\mathbb{R}^2$ *as* $\mathbb{C}$ *using the coordinate* $z = x_1 + ix_2$. *If* $p \in \mathbb{C}$ *is a zero of* $\varphi|_{(t,x_3)}$, *there is a positive integer* m *and a non-zero section of* ad(P)$_\mathbb{C}$ *defined near* p *(to be denoted by* $\eta$) *such that* $\varphi|_{(t,x_3)}$ *near* p *can be written as* $\eta(z-p)^m + \mathcal{O}(|z-p|^{m+1})$.
    c) *The integer* m *as defined above for any given zero of* $\varphi$ *is the same for all* (t, $x_3$).

*Moreover, there exists a unit length section of* ad(P) *(denoted by* $\sigma$*) with the following significance*:
- $\langle\sigma\varphi\rangle = 0$ *and* $[\frac{i}{2}\sigma, \varphi] = \varphi$; *so* $\varphi$ *is a section of the associated complex line bundle* $\mathcal{L}^+$.
- *Define* $\hat{A} = A - \frac{1}{4}[\sigma, \nabla_A\sigma]$. *Then*
    a) $\nabla_{\hat{A}}\sigma = 0$ *and* $\nabla_{\hat{A}}$ *defines a covariant derivative on sections of* $\mathcal{L}^+$.
    b) $(\nabla_{\hat{A}1} + i\nabla_{\hat{A}2})\varphi = 0$ *so* $\varphi$ *is a* $\nabla_{\hat{A}}$-*holomorphic section of* $\mathcal{L}^+$.
    c) *Fix* (t, $x_3$) *and view* $\varphi|_{(t,x_3)}$ *as a section of* $\mathcal{L}^+$ *on* $\mathbb{R}^2$ *by varying only* $x_1$ *and* $x_2$. *The local degree of* $\varphi|_{(t,x_3)}$ *at any given zero is positive and independent of* t *and* $x_3$.

Sections 3a and 3b contain the proof of this proposition. Section 3c uses what is said by the proposition to rewrite the equations that are depicted by (1.4).

**a) Proof of Proposition 3.1: The section $\varphi$**

This subsection proves the assertions in the the first three bullets of the proposition. With regards to the top bullet: It's assertion is equivalent to the assertion that $\langle\varphi\varphi\rangle$ is everywhere zero. To see that it is, note that its t-derivative is zero because of the top bullet in (1.4). Therefore, $\langle \varphi \varphi \rangle$ is constant on any half-line of the form $(0,\infty) \times \{(z, x_3)\}$ in $(0,\infty) \times \mathbb{R}^2 \times \mathbb{R}$. It must therefore be zero because $\lim_{t \to \infty} |\varphi| = 0$ (since $t|\varphi|$ is bounded).

With regards to the second bullet: By virtue of the first bullet in (1.4), the derivative of $|\varphi|$ obeys the differential inequality $\left|\frac{\partial}{\partial t}|\varphi|\right| \leq 2|\mathfrak{a}_3||\varphi|$. This has the folowing implications: Suppose that $(t, z, x_3)$ is a point in $(0,\infty) \times \mathbb{R}^2 \times \mathbb{R}$ where $|\varphi| \neq 0$ and that s is a point in $(t, \infty)$. Then $|\varphi|$ at $(s, z, x_3)$ is no smaller than $(\frac{t}{s})^{\mathfrak{z}} |\varphi|(t, z, x_3)$ with $\mathfrak{z}$ being the supremum of $t|\mathfrak{a}_3|$ along the line segment $[t, s] \times \{(z, x_3)\}$. By the same token, if s is in $(0, t)$, then $|\varphi|$ at $(s, z, x_3)$ is no smaller than $(\frac{s}{t})^{\mathfrak{z}} |\varphi|(t, z, x_3)$ with $\mathfrak{z}$ now being the supremum of $t|\mathfrak{a}_3|$ along the segment $[s, t] \times \{(z, x_3)\}$. By arguing the contrapositive, if $(t, z, x_3)$ are such that $\varphi(t, z, x_3) = 0$, then $\varphi(s, z, x_3) = 0$ for all $s \in (0, \infty)$ (and also for all $x_3$ since $\nabla_{A3}\varphi = 0$.)

With regards to Items a) and b) of the third bullet: The second bullet of (1.4) says in effect that $\varphi$'s restriction to any given constant $(t, x_3)$ slice of $(0, \infty) \times \mathbb{R}^2 \times \mathbb{R}$ (which is a copy of $\mathbb{R}^2$) defines a holomorphic section of $\text{ad}(P)_{\mathbb{C}}$ with respect to the holomorphic structure that is defined as follows: The holomorphic structure on $\mathbb{R}^2$ is the standard one whereby the complex coordinate z is $x_1 + ix_2$. With this understood, then the holomorphic structure for the bundle $\text{ad}(P)_{\mathbb{C}}$ is defined by specifying the associated $\bar{\partial}$-operator, that being $\frac{1}{2}(\nabla_{A1} + i\nabla_{A2})$. Items a) and b) of the third bullet of Proposition 3.1 follow immediately from the fact that $\varphi$ is holomorphic.

Item c) of the third bullet concerns the integer n that is defined by writing $\varphi|_{(t, x_3)}$ in the manner of Item b) near one of its zeros (thus as $\eta(z-p)^m + \mathcal{O}(|z-p|^{m+1})$ with $\eta \neq 0$). The fact that m is independent of the $x_3$ coordinate follows from two observations: The first is that the operator $\nabla_{A3}$ commutes with $\nabla_A^\perp$ (which is to say, both $\nabla_{A1}$ and $\nabla_{A2}$). This is because $F_{A3} \equiv 0$. The second is that $\nabla_{A3}$ annihilates $\varphi$ (which is because $\nabla_{A3}\varphi \equiv 0$). The fact that the integer n does not depend on the t-coordinate follows from the top bullet of (1.4). Note in this regard that the top bullet of (1.4) can be used to write any $k \geq 0$ version of $\nabla_{At}((\nabla_A^\perp)^{\otimes k}\varphi)$ schematically as

$$[\mathfrak{a}_3, (\nabla_A^\perp)^{\otimes k}\varphi] + \sum_{1 \leq j < m}[\mathfrak{F}^k, (\nabla_A^\perp)^{\otimes j}\varphi] + [(\nabla_A^\perp)^{\otimes k}\mathfrak{a}_3, \varphi]$$

(3.1)

with $\{\mathfrak{F}^j\}_{1 \leq j \leq k}$ denoting a certain set of $\text{ad}(P)$-valued tensors. Thus, $(\nabla_A^\perp)^{\otimes j}\varphi$ is zero for all $0 \leq j \leq k$ at a given $(t, z, x_3)$, if and only if this is the case for its $\nabla_{At}$ covariant derivative.

**b) Proof of Proposition 3.1: The section $\sigma$**

This subsection proves fourth and fifth bullet assertions of Proposition 3.1. To this end, it is necessary to first define $\sigma$; and to do that, write the zero locus of $\varphi$ as $(0,\infty) \times Z \times \mathbb{R}$ with Z being a discrete set in $\mathbb{R}^2$ with no accumulation points. The section σ on the complement of Z in any given constant $(t, x_3)$ slice of $(0,\infty) \times \mathbb{R}^2 \times \mathbb{R}$ is

$$\sigma = \tfrac{i}{2|\varphi|^2} [\varphi, \varphi^*]$$

(3.2)

with $\varphi^*$ defined here to be $\mathfrak{a}_1 + i\mathfrak{a}_2$. When written in terms of $\mathfrak{a}_1$ and $\mathfrak{a}_2$, the definition in (3.2) says that $\sigma = -\tfrac{1}{2|\mathfrak{a}_1|^2}[\mathfrak{a}_1, \mathfrak{a}_2]$. The section σ has norm 1 because $\mathfrak{a}_1$ and $\mathfrak{a}_2$ are pointwise orthogonal and have the same norm (which is the top bullet of Proposition 3.1).

To extend σ over Z, fix a point $p \in Z$ (but viewed as a point in $\mathbb{C}$) and fix a product structure for ad(P) near p. Then φ near p can be written using the complex coordinate $z = x_1 + ix_2$ as

$$\varphi = \eta (z-p)^m + \mathcal{O}(|z-q|^{m+1})$$

(3.3)

with m being a positive integer and with η being a non-zero element in the the complification of the Lie algebra of SU(2). This is because of the third bullet of Proposition 3.1. Note in particular that $\langle \eta \eta \rangle = 0$ which is because of the top bullet of Proposition 3.1). Use the depiction of φ in (3.3) to see that σ near q can be written as

$$\sigma = \tfrac{i}{2|\eta|^2} [\eta, \eta^*] + \mathcal{O}(|z-q|)$$

(3.4)

which shows that σ extends over Z on the constant $(t, x_3)$ slice of $(0,\infty) \times \mathbb{R}^2 \times \mathbb{R}$ as a Holder continuous section of ad(P). Slightly more work shows that the extension is smooth. (The argument for this also uses the top bullet of Proposition 3.1.) These extensions for each constant $(t, x_3)$ define a smooth section of ad(P) over $(0,\infty) \times \mathbb{R}^2 \times \mathbb{R}$ because η varies smoothly with $(t, x_3)$ by virtue of p being isolated from the other zeros of φ. (Or, one can argue using the fact that $\nabla_{A_3}\mathfrak{a} = 0$ to see that the A-covariant derivative of η in the $x_3$ direction is zero; and use the fact that and $\nabla_{A_t}\varphi = i[\mathfrak{a}_3, \varphi]$ to see that η varies smoothly in the t-direction.)

With regards to the fourth bullet of Proposition 3.1. The definition of σ implies directly that $\langle \sigma\varphi \rangle$ is zero. To prove that $[\tfrac{i}{2}\sigma, \varphi]$ is equal to φ, note first that $[\tfrac{i}{2}\sigma, \varphi]$ where $\varphi \neq 0$ can be always be written as

$$[\tfrac{i}{2}\sigma, \varphi] = \nu_+\varphi + \nu_-\varphi^*$$

(3.5)

with $\nu_+$ and $\nu_-$ being $\mathbb{C}$-valued functions. This is because $\mathfrak{a}_1$ and $\mathfrak{a}_2$ are orthogonal. To see that $\nu_-$ must be zero, take the inner product of both sides of (3.5) with φ. The left hand side of the result is zero because $\langle \varphi[\sigma, \varphi] \rangle$ is the same as $\langle \sigma[\varphi, \varphi] \rangle$. Meanwhile the right hand side of the resulting identity is $v_-|\varphi|^2$ because $\langle\varphi\varphi\rangle = 0$ and $\langle\varphi\varphi^*\rangle = |\varphi|^2$ (these restate the top bullet of Proposition 3.1). To see that $v_+ = 1$, take the inner product of both sides of (4.15) with $\varphi^*$ and use the definition of $\sigma$ to see that the left hand side of the result is $|\varphi|^2$ and the right hand side is $v_+|\varphi|^2$.

With regards to the fifth bullet: Note first that $\hat{A}$ is *defined* so that $\nabla_{\hat{A}}\sigma = 0$. And, granted that, then the endomorphism $[\frac{i}{2}\sigma,\cdot]$ of $\mathrm{ad}(P)_{\mathbb{C}}$ commutes with $\nabla_{\hat{A}}$; and so $\nabla_{\hat{A}}$ defines a covariant derivative on sections of $\mathcal{L}^+$. This is Item a) of the fifth bullet.

To prove Item b), introduce by way of notation $\mathfrak{b}$ to denote $\frac{1}{4}[\sigma, \nabla_A\sigma]$. An important point is that all components of this $\mathrm{ad}(P)$ valued 1-form are orthogonal to $\sigma$. The identity in the second bullet of (1.4) when written using $\hat{A}$ and $\mathfrak{b}$ says that

$$(\nabla_{\hat{A}1} + i\nabla_{\hat{A}2})\varphi + [\mathfrak{b}_1 + i\mathfrak{b}_2, \varphi] = 0 .$$

(3.6)

Now $(\nabla_{\hat{A}1}+i\nabla_{\hat{A}2})\varphi$ is a section of $\mathcal{L}^+$ because the $\hat{A}$-covariant derivatives map sections of $\mathcal{L}^+$ to sections of $\mathcal{L}^+$. Meanwhile, $[\mathfrak{b}_1+i\mathfrak{b}_2, \varphi]$ is proportional to $\sigma$ because $\varphi$ and $\mathfrak{b}_1$ and $\mathfrak{b}_2$ are all orthogonal to $\sigma$. Thus, both $(\nabla_{\hat{A}1}+i\nabla_{\hat{A}2})\varphi$ and $[\mathfrak{b}_1+i\mathfrak{b}_2, \varphi]$ must separately vanish. The vanishing of $(\nabla_{\hat{A}1}+i\nabla_{\hat{A}2})\varphi$ is the assertion of Item b) of Proposition 3.1's sixth bullet. For future reference: The commutator $[\mathfrak{b}_1+i\mathfrak{b}_2, \varphi]$ is zero if and only if $\mathfrak{b}_1+i\mathfrak{b}_2$ is a section of $\mathcal{L}^+$. (This is because $\varphi$ is a section of $\mathcal{L}^+$ and because $\mathfrak{b}_1$ and $\mathfrak{b}_2$ are orthogonal to $\sigma$.)

Item c) of the fifth bullet follows from Item c) of the proposition's third bullet.

**c) Rewriting the equation in (1.4)**

Suppose in what follows that $(A,\mathfrak{a})$ obeys (1.4), that $t|\mathfrak{a}|$ is a bounded function on $(0,\infty)\times\mathbb{R}^2\times\mathbb{R}$ and that $\varphi$ is not identically zero. Since Proposition 3.1 can be invoked granted these assumptions, there is a corresponding version of $\sigma$ and a connection $\hat{A}$ and complex line bundles $\mathcal{L}^+$ and $\mathcal{L}^-$. Each of the bulleted equations in (1.4) depicts an equality between sections of $\mathrm{ad}(P)$. As such, each can be orthogonally projected onto the span of $\sigma$ and likewise onto the span of $\mathcal{L}^+$. These projected equations are depicted momentarily. (There is also a projection to the span of $\mathcal{L}^-$, but that is the complex conjugate of the projection to the span of $\mathcal{L}^+$.) What follows directly sets the stage for the upcoming depiction of (1.4)

To set the stage for what is to come, introduce an $\mathbb{R}$-valued function $\alpha$ and a section $\beta$ of $\mathcal{L}^+$ by writing $\mathfrak{a}_3$ as

$$\mathfrak{a}_3 = \alpha\sigma + \beta + \beta^*$$

(3.7)

with $\beta^*$ denoting -1 times the Hermitian conjugate. (By way of notation: Supposing that $\eta$ is a section of $\mathcal{L}^+$, then $\eta^*$ is used to denote -1 times its Hermitian conjugate (this is to say that $\eta^* = -\eta^\dagger$). For example, $\varphi^* = \mathfrak{a}_1 + i\mathfrak{a}_2$.) By way of more notation: The $\mathrm{ad}(P)$ valued 1-form $\frac{1}{4}[\sigma, \nabla_A \sigma]$ is denoted henceforth by $\mathfrak{b}$. Thus, $A = \hat{A} + \mathfrak{b}$. Meanwhile, $\mathfrak{b}$ is written as $\mathfrak{b}_t dt + \mathfrak{b}_1 dx_1 + \mathfrak{b}_2 dx_2$ and then $\flat$ is introduced to denote $\frac{1}{2}(\mathfrak{b}_1 + i\mathfrak{b}_2)$; it is a section of $\text{ad}(P)_\mathbb{C}$ that is pointwise orthogonal to $\sigma$. (There is no $dx_3$ component of $\mathfrak{b}$. To see why, note first that $\nabla_{A3}\sigma$ is orthogonal to $\sigma$ in $\text{ad}(P)$ because $\sigma$ has norm 1. It is therefore determined by its $\mathcal{L}^+$ projection which is everywhere zero if $\langle \nabla_{A3}\sigma \varphi\rangle$ is everywhere zero (assuming that $\varphi$ is not identically zero). Meanwhile, $\langle \nabla_{A3}\sigma \varphi\rangle$ is equal to $-\langle \sigma \nabla_{A3}\varphi\rangle$ because $\varphi$ is orthogonal to $\sigma$ in $\text{ad}(P)$; and $\langle \sigma \nabla_{A3}\varphi\rangle$ is zero because $\nabla_{A3}\mathfrak{a} = 0$.)

The respective projections of the first and second bullets in (1.4) along $\sigma$ are algebraic constraints on $\mathfrak{b}_t$ and $\mathfrak{b}_1$ and $\mathfrak{b}_2$:

- $\mathfrak{b}_t = -i(\beta - \beta^*)$.
- $\flat$ is a section of $\mathcal{L}^+$.

(3.8)

The first bullet follows by $\langle \sigma \nabla_{At}\varphi\rangle$ as $-\langle \nabla_{At}\sigma \varphi\rangle$ and using the fact that $\frac{1}{2}[\sigma, \varphi] = \varphi$. The second bullet follows from (3.6). Meanwhile, the respective $\mathcal{L}^+$ projections of the top two bullets in (1.4) give the following:

- $\nabla_{\hat{A}t}\varphi = 2\alpha \varphi$.
- $(\nabla_{\hat{A}1} + i\nabla_{\hat{A}2})\varphi = 0$.

(3.9)

The projection along $\sigma$ of the equations in the last three bullets of (1.4) can be rewritten using (3.8) as

- $\langle \sigma B_{\hat{A}3}\rangle = \frac{\partial}{\partial t}\alpha - |\varphi|^2 - 4|\flat|^2 - 4|\beta|^2$.
- $\langle \sigma(E_{\hat{A}1} + iE_{\hat{A}2})\rangle = -i(\frac{\partial}{\partial x_1} + i\frac{\partial}{\partial x_2})\alpha$

(3.10)

And, a useful rewriting of the $\mathcal{L}^+$ projection of these equations (again using (3.8)) is:

- $\nabla_{\hat{A}t}\beta + 2\alpha \beta = -i(\nabla_{\hat{A}1} - i\nabla_{\hat{A}2})\flat$.
- $\nabla_{\hat{A}t}\flat - 2\alpha \flat = -i(\nabla_{\hat{A}1} + i\nabla_{\hat{A}2})\beta$.

(3.11)

For the record (and for future use), this rewriting of the equations in (1.4) uses depictions of the $\sigma$ projections of $B_{A3}$ and $E_{A1}$ and $E_{A2}$ that can be written as

- $\langle \sigma B_{A3}\rangle = \langle \sigma B_{\hat{A}3}\rangle + 4|\flat|^2$.
- $\langle \sigma(E_{A1} + iE_{A2})\rangle = \langle \sigma(E_{\hat{A}1} + iE_{\hat{A}2})\rangle - 4\langle \beta^* \flat\rangle$.

(3.12)

Meanwhile, the corresponding $\mathcal{L}^+$ projections can be written as:

- $B_{A3}{}^+ = -i(\nabla_{\hat{A}1} - i\nabla_{\hat{A}2})\flat$ .
- $E_{A1}{}^+ = i\nabla_{\hat{A}1}\beta + \nabla_{\hat{A}t}\flat$ .
- $E_{A2}{}^+ = i\nabla_{\hat{A}2}\beta - i\nabla_{\hat{A}t}\flat$ .

(3.13)

To put things in perspective, notice first that the respective bullets in (3.9) determine the behavior of φ as a function of t with constant $\mathbb{R}^2$ (and $\mathbb{R}$) coordinate, and then with constant t coordinate. These equations determine φ given only knowledge of α and Â because there is no explicit appearances of β and $\flat$. Meanwhile, (3.11) is a linear elliptic equation for β and $\flat$ given only knowledge of α and Â because there is no explicit appearances of φ. Then, (3.10) determines α and Â given knowledge of φ and β and $\flat$.

A parenthetical remark: There is a U(1) action on the set of solutions to (1.4) that does not factor through Aut(P) unless either φ is identically zero or β and $\flat$ are identically zero. This action is defined by the rule whereby any given unit norm complex number (call it $u$) has no effect on (σ,Â,φ) but sends (β,$\flat$) to ($u$β,$u\flat$).

## 4. When $\langle \mathfrak{a}_3\varphi \rangle \equiv 0$

Theorem 1 is an instance of Theorem 4.1 which is given momentarily and proved later in this section. The proof of Theorem 4.1 calls on various preliminary observations about $\langle \mathfrak{a}_3\varphi \rangle \equiv 0$ solutions that are established first. (These are the solutions with β ≡ 0.) To set the stage for Theorem 4.1, introduce by way of notation a set $\mathcal{M}_f$ whose members are pairs (A,$\mathfrak{a}$) that obey (1.4), are described by CONSTRAINT SET 1 and whose corresponding version of Z is finite. (Remember that the zero locus of φ has the form $(0,\infty) \times Z \times \mathbb{R}$ with Z being a countable set in $\mathbb{R}^2$ with no accumulation points if it isn't the whole of $\mathbb{R}^2$.)

**Theorem 4.1**: *Suppose that (A,$\mathfrak{a}$) and (A´,$\mathfrak{a}$´) are elements in $\mathcal{M}_f$ with Z being the same as Z´ and with the local degree of vanishing of φ´ at any given point being the same as that of φ at that point. Then (A´,$\mathfrak{a}$´) is the pull-back of (A,$\mathfrak{a}$) via an automorphism of* P.

Just to be sure about terminology: Having fixed (t, $x_3$), then φ can be viewed as a section of the complex line bundle $\mathcal{L}^+$ (see the fifth bullet of Proposition 3.1), and likewise φ´ can be viewed as a section of a primed version of $\mathcal{L}^+$. As sections of respective complex line bundles, φ and φ´ have well defined local vanishing degrees at each of their zeros.

This theorem is proved in Section 4d (Siqi He and Rafe Mazzeo in [HM1] state an analog of Theorem 4.1 that assumes more with regards to the t → 0 asymptotics.)

### a) The structure of (A, 𝔞) with ⟨𝔞φ⟩ ≡ 0

The lemma that follows adds to what is said by Proposition 3.1 when ⟨𝔞₃φ⟩ is assumed to vanish. The lemma refers to the section σ and the connection Â which are supplied by Proposition 3.1 when φ is not identically zero.

**Lemma 4.2**: *Suppose that (A, 𝔞) is a solution to (1.4) on $(0,\infty)\times\mathbb{R}^2\times\mathbb{R}$ with $t|𝔞|$ being a bounded function on $(0,\infty)\times\mathbb{R}^2\times\mathbb{R}$. Assume that φ is not identically zero but that ⟨𝔞₃φ⟩ (and hence β) is identically zero. Then A = Â and $\nabla_A \sigma \equiv 0$ unless (A, 𝔞) is either a Nahm pole imposter or the pull-back of a Nahm pole imposter by a reimbedding of $(0,\infty)\times\mathbb{R}^2\times\mathbb{R}$ into itself that sends t to t + τ with τ being a positive number.*

This lemma says in effect that $\frac{𝔞_3}{|𝔞_3|}$ is everywhere $\nabla_A$-covariantly constant when ⟨𝔞₃φ⟩ = 0 except for the case of Nahm pole imposters or their constant time translation pull-backs. (The Nahm pole imposters and their pull-backs are described in Section 1c).

*Proof of Lemma 4.2*: The assertion to the effect that A = Â follows from the definition in Proposition 3.1 of Â if $\nabla_A \sigma \equiv 0$. The proof of the latter claim (assuming that (A, 𝔞) is not the pull-back of a Nahm-pole imposter by a constant time translation) has five parts

*Part 1*: As noted previously, $\nabla_{A3}\sigma \equiv 0$. Meanwhile, $\nabla_{At}\sigma \equiv 0$ because that is $[𝔟_t, \sigma]$ and $𝔟_t$ is zero because $𝔟_t = -i(\beta - \beta^*)$ and $\beta = 0$. Note that the vanishing of $\nabla_{At}\sigma$ implies that $\nabla_{At}$ maps sections of $\mathcal{L}^+$ to sections of $\mathcal{L}^+$. In particular, $\nabla_{At}\varphi$ is $2\alpha\varphi$, a section of $\mathcal{L}^+$.

To consider $\nabla_{A1}\sigma$ and $\nabla_{A2}\sigma$, note that these are zero if and only if their respective commutators with σ are zero. Thus, both are zero if and only if $𝔟_1$ and $𝔟_2$ are both zero, which is if and only if $𝔟$ is zero. (Remember that $𝔟_1$ is $\frac{1}{4}[\sigma, \nabla_{A1}\sigma]$ and $𝔟_2$ is $\frac{1}{4}[\sigma, \nabla_{A2}\sigma]$ and $𝔟$ is $\frac{1}{2}(𝔟_1 + i𝔟_2)$.) The proof that $𝔟$ is zero (given the assumptions) uses the equations in (3.11) which say this when β is zero:

$$-i(\nabla_{\hat{A}1} - i\nabla_{\hat{A}2})𝔟 = 0 \quad \text{and} \quad \nabla_{\hat{A}t}𝔟 - 2\alpha 𝔟 = 0 .$$

(4.1)

Now $𝔟$ is a section of $\mathcal{L}^+$ (according to (3.8); and because the fibers of $\mathcal{L}^+$ are 1-dimensional, $𝔟$ can be written where $\varphi \neq 0$ as $\mu\varphi$ with $\mu$ being a $\mathbb{C}$-valued function. (It is defined on $(0,\infty)\times(\mathbb{R}^2 - Z)\times\mathbb{R}$.) The equations in (4.1) when written using $\mu$ are equivalent to these:

$$\tfrac{\partial}{\partial z}(\mu|\varphi|^2) = 0 \quad \text{and} \quad \tfrac{\partial}{\partial t}\mu = 0 .$$

(4.2)

The left most equation in (4.2) asserts that $\mu|\varphi|^2$ is an anti-holomorphic function on each constant $(t, x_3)$ slice of $(0,\infty) \times \mathbb{R}^2 \times \mathbb{R}$ and the right most equation says $\mu$ is independent of the t-coordinate of $(0,\infty) \times \mathbb{R}^2 \times \mathbb{R}$ (and note that it is independent of the $\mathbb{R}$ factor also)

*Part 2*: Invoke Lemma 2.2 to see that $\mu|\varphi|^2$ is a priori bounded by constant multiple of $\frac{1}{t^2}$. Indeed, the latter bound follows from Lemma 2.2 because $\langle\varphi^*(\mathfrak{b}_1+i\mathfrak{b}_2)\rangle$ can be written as $-\frac{i}{2}\langle\varphi^*(\nabla_{A1}+i\nabla_{A2})\sigma\rangle$ using the definitions of $\mathfrak{b}_1$ and $\mathfrak{b}_2$ as covariant derivatives of $\sigma$; and then $-\frac{i}{2}\langle\varphi^*(\nabla_{A1}+i\nabla_{A2})\sigma\rangle$ is the same as $\frac{i}{2}\langle(\nabla_{A1}+i\nabla_{A2})\varphi^*\sigma\rangle$ because $\varphi^*$ is orthogonal to $\sigma$. Since $\mu|\varphi|^2$ is anti-holomorphic on each constant $(t,x_3)$ slice and since it is uniformly bounded on each such slice, it must be constant on each $(t,x_3)$ slice.

Meanwhile: Keeping in mind that $\sigma$ is a differentiable section of ad(P) (it is smooth, in fact), it follows that $\mu|\varphi|$ is bounded at each zero of $\varphi$. (This is because $\mathfrak{b}_1$ and $\mathfrak{b}_2$ are bounded by multiples of $|\nabla_A\sigma|$.) Therefore, $\mu|\varphi|^2$ (which is constant on any given $(t,x_3)$ slice of $(0,\infty)\times\mathbb{R}^2\times\mathbb{R}$) must vanish at each zero of $\varphi$. Therefore, $\mu|\varphi|^2$ and hence $\mathfrak{b}_1$ and $\mathfrak{b}_2$ are everywhere zero if $\varphi$ has at least one zero on any constant $(t,x_3)$ slice of $(0,\infty)\times\mathbb{R}^2\times\mathbb{R}$. Thus, $\nabla_A\sigma \equiv 0$ if there are zeros of $\varphi$.

*Part 3*: It remains now to consider the case where $\mu|\varphi|^2$ is not identically zero. In this case, as noted in Part 2, it depends only on t. To see what this implies, use the right most equation in (4.2) and the top bullet equation in (3.9) to see that

$$\tfrac{\partial}{\partial t}(\mu|\varphi|^2) = 4\alpha\mu|\varphi|^2 .$$

(4.3)

If the assumption that $\mu|\varphi|^2$ depends only on t is enforced, then (4.3) needs $\alpha$ to depend only on t also. And, if $\alpha$ depends only on t, then $E_{\hat{A}}$ must vanish on $(0,\infty)\times\mathbb{R}^2\times\mathbb{R}$ (look at the second bullet equation in (3.10)). And, if $E_{\hat{A}} \equiv 0$, then the t-derivative of $B_{\hat{A}3}$ must vanish too because of a Bianchi identity. As explained in the next paragraphs, if $B_{\hat{A}3}$ is t-independent, then it must be zero so not to run afoul of the Lemma 2.2.

The claim that $B_{\hat{A}3}$ must vanish will follow from the top bullet identity in (3.12) and Lemma 2.2's bound for $|F_A|$ if the following is true: If $z \in \mathbb{R}^2 - Z$, then there is an increasing, unbounded sequences $\{t_n\}_{n\in\mathbb{N}} \subset (0,\infty)$ such that $\lim_{n\to\infty}|\mathfrak{b}|(t_n,z,0) = 0$. By virtue of the definition of $\mathfrak{b}$, this will happen if and only if $\lim_{n\to\infty}|\nabla_A\sigma| = 0$. The identity $\nabla_A\sigma = \frac{1}{4\alpha}[\sigma,[\nabla_A\mathfrak{a}_3,\sigma]]$ will be used to find a sequence with the latter property. Indeed, because of this identity, and because of Lemma 2.2 (which bounds $|\nabla_A\mathfrak{a}|$ by a constant multiple of $\frac{1}{t^2}$), it is enough to find a sequence $\{t_n\}_{n\in\mathbb{Z}}$ such that $\lim_{n\to\infty}t_n|\alpha| > 0$.

To see that there is an unbounded sequence of times with this property, introduce by way of notation $w$ to the function of $t$ in $(0, \infty)$ given by the integral of $\alpha(\cdot, z, x_3 = 0)$ from 1 to $t$. Use this definition with the top bullet in (3.9) to write

$$|\varphi|(t, z, x_3) = e^{2w} |\varphi|(1, z, x_3) .$$

(4.4)

If $T > 1$ and if $\alpha$ on the whole of $[T, \infty)$ is less negative than $-\frac{1}{2t}(1-\delta)$ for any given $\delta > 0$, then (4.4) will run afoul of the assumption that $t|\mathfrak{a}|$ is bounded. This is because the right hand side of (4.4) will be greater than $(\frac{T}{t})^{1-\delta}|\varphi|(T, z, x_3)$.

*Part 4*: With $B_{\hat{A}3}$ being zero, then the top bullet in (3.10) asserts the identity

$$\tfrac{\partial}{\partial t} \alpha = (1 + |\mu|^2)|\varphi|^2 .$$

(4.5)

Since $\alpha$ depends only on $t$, this last identity requires that $(1 + |\mu|^2)|\varphi|^2$ depend only on $t$ also. And, because of this, and because $\mu|\varphi|^2$ depends only on $t$, both $\mu$ and $|\varphi|^2$ must be be functions of only $t$ also (which implies that $\mu$ is constant since it is independent of $t$). Indeed, if both $(1+|\mu|^2)|\varphi|^2$ and $\mu|\varphi|^2$ depend only on $t$, then their ratio, $\frac{\mu}{1+|\mu|^2}$, depends only on $t$; and this can happen only in the event that $\mu$ depends only on $t$.)

With $B_{\hat{A}3}$ being zero, Item b) of the fifth bullet of Proposition 3.1 implies that $(\nabla_{\hat{A}1}{}^2 + \nabla_{\hat{A}2}{}^2)\varphi = 0$, and this implies in turn that $\tfrac{1}{2}(\tfrac{\partial^2}{\partial x_1^2} + \tfrac{\partial^2}{\partial x_1^2})|\varphi|^2 = |\nabla_{\hat{A}1}\varphi|^2 + |\nabla_{\hat{A}2}\varphi|^2$. Therefore, both $\nabla_{\hat{A}1}\varphi$ and $\nabla_{\hat{A}2}\varphi$ must be zero because $|\varphi|^2$ is a function only of $t$.

*Part 5*: To determine $\alpha$ and $|\varphi|$, differentiate (4.5) using the identity $\tfrac{\partial}{\partial t}|\varphi|^2 = 4\alpha|\varphi|^2$ and the fact that $\mu$ is constant. Then, use the identity in (4.5) to see that

$$\tfrac{\partial}{\partial t}(\tfrac{\partial}{\partial t}\alpha - 2\alpha^2) = 0 .$$

(4.6)

This is to say that

$$\tfrac{\partial}{\partial t}\alpha - 2\alpha^2 = z$$

(4.7)

with $z$ being constant. If $z \neq 0$, then no solution to (4.7) is compatible with the fact that $t|\mathfrak{a}|$ and $t^2|\nabla_A \mathfrak{a}|$ are a priori bounded (the latter fact is courtesy of Lemma 2.2). (The solutions with $z < 0$ are the ones that are described in Section 1c after the statement of Theorem 2 where $|\mathfrak{a}_3|$ has a positive $t \to \infty$ limit.) If $z$ is zero, then the solution is either a Nahm pole solution or a Nahm pole imposter, or one of their pull-backs by the $t \to t + \tau$ reimbedding of $(0,\infty) \times \mathbb{R}^2 \times \mathbb{R}$.

To elaborate: The solutions to the $z = 0$ version of (4.7) have the form $\alpha = -\frac{1}{2(t+\tau)}$ for $\tau \in \mathbb{R}$. The nonsingular ones on $(0, \infty)$ require $\tau > 0$. Then, (4.5) implies that $|\varphi| = \frac{1}{2(t+\tau)} \frac{1}{\sqrt{1+|\mu|^2}}$. Granted that $\hat{A}$ is flat, that $\nabla_{\hat{A}1}\varphi = \nabla_{\hat{A}2}\varphi = 0$, that $\mu$ is constant, and granted that $\alpha$ and $|\varphi|$ are as described in the preceding paragraph, then $(A,\mathfrak{a})$ is a Nahm pole solution when $\mu = 0$ and $\tau = 0$; and $(A,\mathfrak{a})$ when $\mu \neq 0$ and $\tau = 0$ is a $w = \frac{\mu}{1+|\mu|^2}$ Nahm pole imposter solution. When $\tau > 0$, then the solutions are the respective pull-backs of these by the $t \to t+\tau$ reimbedding.

**b) Solutions where $t|\mathfrak{a}_3|$ has no positive lower bound**

The proposition that follows directly limits the possiblilities for $(A,\mathfrak{a})$ when $t|\mathfrak{a}_3|$ has no positive lower bound and when $\varphi = \mathfrak{a}_1 - i\mathfrak{a}_2$ is not identically zero. This proposition asserts in part that the condition in the second bullet of CONSTRAINT SET 1 imply that $t|\mathfrak{a}_3|$ (which is $t|\alpha|$) has a positive lower bound on the whole of $(0,\infty) \times \mathbb{R}^2 \times \mathbb{R}$.

**Proposition 4.3**: *Suppose that $(A,\mathfrak{a})$ is a solution to (1.4) on $(0,\infty) \times \mathbb{R}^2 \times \mathbb{R}$ with $t|\mathfrak{a}|$ being a bounded function on $(0,\infty) \times \mathbb{R}^2 \times \mathbb{R}$ and with $\langle \mathfrak{a}_3 \varphi \rangle \equiv 0$. If $\varphi$ is not identically zero, and if there exists a positive $t_0$ such that $|\mathfrak{a}_3| > 0$ when $t < t_0$, then $|\mathfrak{a}_3| > 0$ on the whole of $(0,\infty) \times \mathbb{R}^2 \times \mathbb{R}$. Moreover, if there exists $t_0 > 0$ and any positive $\varepsilon$ such that $t|\mathfrak{a}_3| > \varepsilon$ where $t < t_0$, then $t|\mathfrak{a}_3|$ has a positive lower bound on the whole of $(0,\infty) \times \mathbb{R}^2 \times \mathbb{R}$.*

*Proof of Proposition 4.3*: The proof has ten parts. Parts 1-3 prove that $|\mathfrak{a}_3|$ can't be zero at any point. Parts 4-10 establish the uniform lower bound for $t|\mathfrak{a}_3|$ under the extra assumption about this function's behavior as $t \to 0$.

*Part 1*: If $\varphi$ is not identically zero, then Proposition 3.1 and Lemma 4.2 describe the pair $(A,\mathfrak{a})$. If $(A,\mathfrak{a})$ is a Nahm pole imposter or its pull-back by a $t \to t +\tau$ reimbedding, then the conclusions of Proposition 4.3 can be verified by inspection because $\mathfrak{a}_3$ in this case is $\mathfrak{a}_3 = -\frac{1}{2(t+\tau)}\sigma$ where $\tau$ is a nonnegative number. Assume henceforth that $(A, \mathfrak{a})$ is not one of these solutions. In this event, $\sigma$ is an A-covariantly constant, norm 1 section of $\mathrm{ad}P$. (Remember that $\mathfrak{a}_3 = \alpha\sigma$.)

Since $\mathfrak{a}_3 = \alpha\sigma$ and since $\sigma$ is covariantly constant, the $dx_3$ component of the equation in (2.1) when written in terms of $\alpha$ is

$$-(\tfrac{\partial^2}{\partial t^2} + \tfrac{\partial^2}{\partial x_1^2} + \tfrac{\partial^2}{\partial x_2^2})\alpha + 4|\varphi|^2 \alpha = 0 \ .$$

(4.8)

The latter equation has an associated maximum principle which implies the following: The function $\alpha$ can not take a non-negative local maximum, nor can it take on a non-positive local minimum. The preceding conclusion implies in particular that the locus where $\alpha$ has a particular sign (positive or negative) can not have components with compact closure.

*Part 2*: The assumption that $|\mathfrak{a}_3| > 0$ for $t < t_0$ means that $\alpha$ is either strictly positive or strictly negative for $t < t_0$. As explained in the subsequent paragraphs, the second bullet in (3.9) implies that $\alpha$ must be negative for $t < t_0$.

If $\alpha$ were strictly positive for $t < t_0$, then $\varphi$ would have a finite limit as $t \to 0$ for each fixed $(z, x_3) \in \mathbb{R}^2 \times \mathbb{R}$ because of the second bullet in (3.9). In addition, $|\varphi|$ would be decreasing towards this $t = 0$ limit as $t \to 0$. By virtue of $|\varphi|$ being bounded on each constant t slice (it is bounded by some constant multiple of $\tfrac{1}{t}$), the preceding observation implies in turn that $|\varphi|$ is uniformly bounded on $(0,\infty) \times \mathbb{R}^2 \times \mathbb{R}$. In particular, there exists some positive number $\mathfrak{z}$ such that $\tfrac{\mathfrak{z}}{t_0}$ is a global upper bound for $|\varphi|$.

With the preceding in mind, note that the $dx^1$ and $dx^2$ versions of (2.1) lead to an equation for $|\varphi|^2$ that has the form

$$-\tfrac{1}{2}(\tfrac{\partial^2}{\partial t^2} + \tfrac{\partial^2}{\partial x_1^2} + \tfrac{\partial^2}{\partial x_2^2})|\varphi|^2 + |\nabla_{At}\varphi|^2 + |\nabla_{A1}\varphi|^2 + |\nabla_{A2}\varphi|^2 + (2|\varphi|^2 + 4|\alpha|^2)|\varphi|^2 = 0 \ .$$

(4.9)

To exploit this equation, fix $R > 1$ and let $\chi_R$ denote the function of $z = x_1 + ix_2$ given by the rule $z \to \chi(\tfrac{\ln(|z|)}{\ln(R)})$. (Remember that $\chi$ is non-increasing and that it obeys $\chi(s) = 1$ for $s \le \tfrac{1}{4}$ and $\chi(s) = 0$ for $s \ge \tfrac{3}{4}$.) Multiply both sides of (4.9) by $\chi_R^2$ and integrate the result over any given constant $(t, x_3)$ slice of $(0,\infty) \times \mathbb{R}^2 \times \mathbb{R}$. Having done this, integrate by parts twice with each $x_1$ and $x_2$ derivative and use the fact that $|\varphi|$ is uniformly bounded to obtain the following inequalities (the first is for $t \le t_0$ and the other is for $t \ge t_0$):

- $-\tfrac{1}{2}\tfrac{d^2}{dt^2}\int_{\mathbb{R}^2}\chi_R^2|\varphi|^2 + 2\int_{\mathbb{R}^2}\chi_R^2|\mathfrak{a}|^2|\varphi|^2 \le c_0 \tfrac{\mathfrak{z}^2}{t_0^2}\tfrac{1}{\ln R}$ *for* $t \le t_0$.
- $-\tfrac{1}{2}\tfrac{d^2}{dt^2}\int_{\mathbb{R}^2}\chi_R^2|\varphi|^2 + 2\int_{\mathbb{R}^2}\chi_R^2|\mathfrak{a}|^2|\varphi|^2 \le c_0 \tfrac{\mathfrak{z}^2}{t^2}\tfrac{1}{\ln R}$ *for* $t \ge t_0$.

(4.10)

To elaborate: The bound on the right hand side of (4.10) is obtained by using integration by parts twice to bound the integral of the product of $\chi_R$ with the $-\tfrac{1}{2}(\tfrac{\partial^2}{\partial x_1^2} + \tfrac{\partial^2}{\partial x_2^2})|\varphi|^2$ term in (4.9). All of these derivatives on $|\varphi|^2$ are moved to $\chi_R$ which results in a function of $|z|$

whose norm is bounded by $c_0 (\frac{1}{\ln R})^2 \frac{1}{|z|^2}$ and whose support lies where $R^{1/4} \le |z| \le R^{3/4}$. The bound $|\varphi|^2 \le \frac{\mathfrak{z}^2}{t^2}$ is then used as well.

Fix some very large T (much greater than 1) and integrate both sides of (4.10) from a given value of $t \in (0, T)$ to T. Then, use the identity $\frac{\partial}{\partial t}|\varphi|^2 = 4\alpha|\varphi|^2$ in the resulting inequality to see that

$$\tfrac{1}{2} (\int_{\mathbb{R}^2} \chi_R^2 4\alpha |\varphi|^2 )|_t + \int_t \int_{\mathbb{R}^2} \chi_R^2 |\mathfrak{a}|^2 |\varphi|^2 \le c_0 (\tfrac{\mathfrak{z}^2}{t_0^2} + \tfrac{\mathfrak{z}^2}{t_0}) \tfrac{1}{\ln R} + c_0 \tfrac{\mathfrak{z} R^2}{T^3} .$$

(4.11)

Look at what this says when $t \le t_0$: If $t < t_0$, the left most term on the left hand side of (4.11) is positive (by assumption). This is true for any choice of R. Therefore, taking R very large and then T very much larger than R leads to the conclusion that $\varphi$ must be identically zero. This conclusion is nonsensical because it runs afoul of the assumptions of the proposition.

*Part 3*: Having fixed $\varepsilon > 0$, let $\Omega_\varepsilon$ denote the subset of $(0,\infty) \times \mathbb{R}^2 \times \mathbb{R}$ where $\alpha$ obeys $\alpha \ge \varepsilon$. This is a smooth submanifold with boundary for all but a measure zero set of values for $\varepsilon$. If $\varepsilon$ is not from this measure zero set, then the outward pointing normal derivative of $\alpha$ along $\partial \Omega_\varepsilon$ is negative. It is also the case that $t \ge t_0$ on $\Omega_\varepsilon$ (because of what is said in Part 2); and this implies that $\alpha$ is uniformly bounded on $\Omega_\varepsilon$ with an $\varepsilon$-independent upper bound (in particular, $\alpha \le \frac{\mathfrak{z}}{t_0}$ on $\Omega_\varepsilon$ for some fixed number $\mathfrak{z}$).

With the preceding understood, fix $R > 1$ and multiply both sides of (4.8) by $\chi_R \alpha$ with $\chi_R$ as before. Then, fix some very large T (much greater than $t_0$ and 1) and integrate the resulting identity over the $\Omega_\varepsilon$ part of $[t_0, T] \times \mathbb{R}^2 \times \{x_3 = 0\}$ to obtain the following:

$$\int_{((0,T]\times\mathbb{R}^2)\cap\Omega_\varepsilon} \chi_R 4\alpha^2 |\varphi|^2 \le c_0 (\tfrac{\mathfrak{z}^2}{t_0^2} + \tfrac{\mathfrak{z}^2}{t_0}) \tfrac{1}{\ln R} + c_0 \tfrac{\mathfrak{z} R^2}{T^3} .$$

(4.12)

With regards to the derivation: The product of $\alpha$ with $-(\frac{\partial^2}{\partial t^2} + \frac{\partial^2}{\partial x_1^2} + \frac{\partial^2}{\partial x_2^2})\alpha$ is written as

$$-\tfrac{1}{2} (\tfrac{\partial^2}{\partial t^2} + \tfrac{\partial^2}{\partial x_1^2} + \tfrac{\partial^2}{\partial x_2^2}) \alpha^2 + |\tfrac{\partial}{\partial t} \alpha|^2 + |d\alpha|^2 .$$

(4.13)

Integration by parts is employed twice to see that integral of the product of $\chi_R$ with (4.13) leads to a positive contribution to the left hand side of (4.48) but for an error of size at most $c_0 (\frac{\mathfrak{z}^2}{t_0^2} + \frac{\mathfrak{z}^2}{t_0}) \frac{1}{\ln R}$. As was the case with (4.11), this error arises when all of the derivatives from the $-\frac{1}{2}(\frac{\partial^2}{\partial x_1^2} + \frac{\partial^2}{\partial x_2^2})\alpha^2$ term in (4.13) land on $\chi_R$ to produce a term whose norm is bounded by $c_0(\frac{1}{\ln R})^2 \frac{1}{|z|^2} \frac{\mathfrak{z}^2}{t^2}$ and whose support lies where $R^{1/4} \leq |z| \leq R^{3/4}$.

Taking R very large and then take T very much larger than R in (4.12) leads to the conclusion that $\Omega_\varepsilon$ must be empty. This is to say that $\alpha < \varepsilon$ for any positive $\varepsilon$. Thus, $\alpha$ can't be positive; and since it can't take on a non-positive local minimum it follows that $\alpha$ must be strictly negative.

*Part 4*: This part of the proof with Parts 5-9 explain why $t|\alpha|$ has a positive lower bound on the whole of $(0,\infty) \times \mathbb{R}^2 \times \mathbb{R}$ when $t|\alpha|$ has a positive lower bound where t is less than some positive $t_0$. To do this, note again that there is no need to consider Nahm pole imposters because $t|\mathfrak{a}_3| = \frac{1}{2}$ for Nahm pole imposters. Supposing that $(A, \mathfrak{a})$ is not a Nahm pole imposter, then $\mathfrak{a}_3$ can be written as $\alpha\sigma$ with $\alpha$ a smooth function and with $\sigma$ being a unit length, A-covariantly constant section of $\text{ad}(P_0)$. The function $\alpha$ obeys (4.8). Let $q = t\alpha$. By virtue of (4.8), this function $q$ obeys

$$-(\frac{\partial^2}{\partial t^2} + \frac{\partial^2}{\partial x_1^2} + \frac{\partial^2}{\partial x_2^2})q + \frac{2}{t}\frac{\partial}{\partial t}q + (4|\varphi|^2 - \frac{2}{t^2})q = 0 \; .$$

(4.14)

Noting that $q < 0$ everywhere (because $\alpha$ is nowhere zero), the maximum principle forbids a local maximum for $q$ where $|\varphi|^2 \leq \frac{1}{2t^2}$.

*Part 5*: The preceding observation is exploited momentarily with the help of the following, second observation:

*There exists $\delta_* > 0$ such that $|\varphi|^2 \leq \frac{1}{4t^2}$ where $q > -\delta_*$.*

(4.15)

To prove the claim in (4.15), assume to the contrary that there is no such $\delta_*$ to derive nonsense. If there is no such $\delta_*$, then the pull-backs of $(A, \mathfrak{a})$ by translations along the $\mathbb{R}^2$ factor of $(0,\infty) \times \mathbb{R}^2 \times \mathbb{R}$ and coordinate rescalings of $(0,\infty) \times \mathbb{R}^2 \times \mathbb{R}$ can be used to construct a sequence of pairs $\{(A^n, \mathfrak{a}^n)\}_{n \in \mathbb{N}}$ obeying (1.4) and with the following additional properties: First, each integer n version of $\mathfrak{a}_n$ obeys $|\mathfrak{a}^n| \leq \frac{\mathfrak{z}}{t}$ with $\mathfrak{z}$ being independent of n. Second, each $\mathfrak{a}_n$ version of $\alpha$ is strictly negative. Third, the integer n version of the function $q$ at the point $(t=1, z=0, x_3 = 0)$ is greater than $-\frac{1}{n^2}$ whereas the $(A^n, \mathfrak{a}^n)$ version of the function $\varphi$ at this same point obeys $|\varphi|^2 \geq \frac{1}{4}$.

Lemma 2.3 can be used to extract a subsequence in $\mathbb{N}$ (to be denoted by $\Lambda$) such that the corresponding sequence $\{(A^n, \mathfrak{a}^n)\}_{n \in \Lambda}$ converges on compact subsets of $(0,\infty) \times \mathbb{R}^2 \times \mathbb{R}$ (after termwise applications of elements in $\text{Aut}(P)$) to a solution to (1.4) that is henceforth denoted by $(A^\infty, \mathfrak{a}^\infty)$. This solution has $|\mathfrak{a}^\infty| \le \frac{3}{t}$ also. The corresponding version of $\alpha$ for this pair (to be denoted by $\alpha^\infty$) is a non-negative function that obeys (4.8) with $\varphi$ replaced by $\varphi^\infty = \mathfrak{a}^\infty{}_1 - i\mathfrak{a}^\infty{}_2$. (It is non-negative because all $\alpha^n$ are positive.) But this function $\alpha^\infty$ vanishes at $(t=1, z=0, x_3=0)$. Thus it must vanish identically on $(0,\infty) \times \mathbb{R}^2 \times \mathbb{R}$ (see Parts 1-3 of this proof). Meanwhile, the corresponding version of the function $\varphi$ (denoted by $\varphi^\infty$) is not identically zero because it is not zero at $(t=1, z=0, x_3=0)$; its norm there is greater than $\frac{1}{4}$. But now note that these last two facts, $\alpha^\infty \equiv 0$ and $\varphi^\infty \ne 0$ are nonsense together because $\alpha^\infty \equiv 0$ implies that the time derivative of $|\varphi^\infty|$ is zero (the top bullet of (3.9)); and that is not compatible with the top the bound $|\mathfrak{a}^\infty| \le \frac{3}{t}$ unless $\varphi^\infty \equiv 0$.

*Part 6*: Define a function on $(0, \infty)$ to be denoted by $q_*$ by the rule whereby $q_*(t)$ is the supremum of the values of $q$ on $\{t\} \times \mathbb{R}^2 \times \mathbb{R}$. Parts 7-10 prove that $q_*$ is a Lipshitz function (which implies that it is $C^1$ on a full measure set in $(0,\infty)$). By assumption, there are positive numbers $t_0$ and $\varepsilon$ such that $q_* < -\varepsilon$ for $t < t_0$. Let $\delta$ denote the smaller of the numbers $\varepsilon$ and $\delta_*$ (with $\delta_*$ coming from (4.15)). As explained also in Parts 7-10, the function $q_*$ has the following property: Given $t > 0$ such that $q_*(t) \in (-\delta, 0)$ and such that $\frac{d}{dt} q_*$ is defined at $t$, and then given sufficiently small but positive $\Delta$ such that $\frac{d}{dt} q_*$ is also defined at $t + \Delta$, then

$$(\tfrac{d}{dt} q_*)_{t+\Delta} \ge \frac{(t+\Delta)^2}{t^2} (\tfrac{d}{dt} q_*)_t .$$

(4.16)

(Remember that $(\frac{d}{dt} q_*)_t$ is defined on a full measure subset of $(0,\infty)$.) As is explained directly, this inequality leads to nonsense which can be avoided only in the event that the requirement $q_* > -\delta$ is never satisfied. And if that requirement is never satisfied, then $|\mathfrak{a}_3|$ (which is $|\alpha|$ and thus $\frac{1}{t}|q|$) is no less than $\frac{\delta}{t}$ for all $t$.

To derive nonsense from (4.16), let $t_1$ denote the infimum of the times $t$ in $(t_0, \infty)$ where $q > -\delta$. If $q$ is ever less than $-\delta$, then there is some interval in $[t_1, \infty)$ where $q$'s values are between $-\delta$ and $-\frac{1}{2}\delta$ and $(\frac{d}{dt} q_*)$ is positive on a full measure set. Let $t_2$ denote a point in that full measure set and let $r > 0$ denote the value of $\frac{d}{dt} q_*$ at $t_2$. It then follows that $(\frac{d}{dt} q_*)_{t_2+\Delta} \ge r$ for any positive $\Delta$ where $(\frac{d}{dt} q_*)_{t+\Delta}$ exists. This implies that $q_*$ is greater than $-\delta + r(t - t_2)$ for all $t > t_2$.

This last conclusion implies that $q_*$ is zero at some time prior to the time $t_2 + r^{-1}\delta$. And that is nonsensical for the following reasons: This value 0 for $q_*$ at finite time can't be taken on at an actual point in $\mathbb{R}^2$ because that event would run afoul of the first assertion of the lemma (that $\alpha < 0$ everywhere). If there is a non-convergent sequence of points $\{z_n\}_{n\in\mathbb{N}}\subset\mathbb{R}^2$ such that $\lim_{n\to\infty}q(t_*,z_n)=0$ for some fixed positive time $t_*$, then $(A,\mathfrak{a})$ could be translated along the $\mathbb{R}^2$ factor by first $z_1$, then $z_2$, and so on to obtain a corresponding sequence $\{(A^n,\mathfrak{a}^n)\}_{n\in\mathbb{N}}$ with corresponding sequence $\{q_n\}_{n\in\mathbb{N}}$ obeying $\lim_{n\to\infty}q_n(t_*,z=0,x_3=0)=0$. The sequence $\{(A^n,\mathfrak{a}^n)\}_{n\in\mathbb{N}}$ would have a convergent subsequence (invoke Lemma 2.3) whose limit (denoted by $(A^\infty,\mathfrak{a}^\infty)$ would obey (1.4) with $|\mathfrak{a}^\infty|\leq\frac{3}{t}$. The limit would also obey $\alpha^\infty(t_*,0,0)=0$ and $\alpha^\infty\geq 0$; and it would also obey $t|\alpha^\infty|\geq\varepsilon$ for $t\leq t_0$ because the t coordinate is not rescaled to make the sequence $\{(A^n,\mathfrak{a}^n)\}_{n\in\mathbb{N}}$. If $\varphi^\infty$ is not identically zero, then this hypothetical $(A^\infty,\mathfrak{a}^\infty)$ runs afoul of the first assertion in Proposition 4.3. If $\varphi^\infty$ is identically zero, then the function $\alpha^\infty$ must obey the equation

$$-\left(\frac{\partial^2}{\partial t^2}+\frac{\partial^2}{\partial x_1^2}+\frac{\partial^2}{\partial x_2^2}\right)\alpha^\infty=0,$$

(4.17)

because each $\alpha^n$ obeys (4.8). But a solution to (4.17) which is non-negative can't have a zero unless it is identically zero, which is forbidden because $t|\alpha^\infty|\geq\varepsilon$ for $t\leq t_0$.

*Part 7*: This part of the proof and Parts 8-10 prove that $q_*$ is a Lipschitz function and that it obeys (4.16) when t and $t+\Delta$ are drawn from a full measure subset. To begin: The $C^\infty$ bounds from Lemma 2.2 can be used with Cerf's constructions [C] to uniformly approximate $q$ in the $C^\infty$ topology on any region of the form $[t_1,t_2]\times\mathbb{R}^2\times\mathbb{R}$ when $0<t_1<t_2$ by an $x_3$-independent, negative function that achieves its supremum for each $t\in[t_1,t_2]$ on $\{t\}\times\mathbb{R}^2\times\mathbb{R}$. Moreover, the supremum is achieved at a single point but for a finite set of times where it is achieved at two points. Whether there be one maximum or two at a given time t, the Hessian matrix $\mathcal{H}$ at any maxima of $q$ is non-degenerate along the $\mathbb{R}^2$ factor of $(0,\infty)\times\mathbb{R}^2\times\mathbb{R}$ for such an approximation. Moreover, in the event that there are two points where $q$ is equal to its maximum value, then the t-derivative of the value of the approximation to $q$ at one of these maxima is strictly greater than the t-derivative of the value the approximation to $q$ at the other.

*Part 8*: With the preceding as background, set $t_1=\frac{1}{2}t_0$, fix $t_2>t_0$ and then choose $\varepsilon>0$. Having done this, let $\mathfrak{q}$ denote an approximation to $q$ on $[t_1,t_2]\times\mathbb{R}^2\times\mathbb{R}$ of the sort described in Part 7 with the $C^3$ norm of $\mathfrak{q}-q$ being less than $\varepsilon$. Fix $t_*>0$ and let z denote a point in $\{t_*\}\times\mathbb{R}^2\times\mathbb{R}$ where $\mathfrak{q}$ takes on its maximum value. By construction, the matrix of second derivatives of $\mathfrak{q}$ at z in the directions tangent to the $\mathbb{R}^2$ factor of $[t_1,t_2]\times\mathbb{R}^2\times\mathbb{R}$ is non-degenerate. Under these circumstances, there will be a unique local maximum of $\mathfrak{q}$ near z on $\{t\}\times\mathbb{R}^2\times\{0\}$ for times near $t_*$ that moves smoothly as a function of t. Let $I \subset (0,\infty)$ denote an open interval where this occurs. The position of the local maximum is a smooth function of t for $t \in I$ denoted by z(t). It is obtained by solving the equation

$$\mathcal{H} \cdot \dot{z} + d(\tfrac{\partial}{\partial t} q)|_{(t,z(t),0)} = 0 \tag{4.18}$$

where $\dot{z}$ denotes the time derivative of z(t) and $\mathcal{H}$ denotes the Hessian at the point z(t) for the function on $\mathbb{R}^2$ that is defined by the rule $z \to q(t,z,x_3 = 0)$. Let $f$ denote the function that is given by the rule $t \to q(t,z(t),0)$. The first and second derivatives of the $f$ can be written as

- $\frac{d}{dt} f = (\frac{\partial}{\partial t} q)|_{(t,z(t),0)}$.
- $\frac{d^2}{dt^2} f = (\frac{\partial^2}{\partial t^2} q)|_{(t,z(t),0)} - \langle \dot{z}, \mathcal{H} \cdot \dot{z} \rangle$

$$\tag{4.19}$$

Given that $q$ obeys (4.14) it then follows from (4.19) that

$$- \tfrac{d^2}{dt^2} f + \tfrac{2}{t} \tfrac{d}{dt} f - (4|\varphi|^2|_{(t,z(t),0)} - \tfrac{2}{t^2})q - (\text{tr}(\mathcal{H}) + \langle \dot{z}, \mathcal{H} \cdot \dot{z} \rangle) = \mathfrak{e} \tag{4.20}$$

where $\mathfrak{e}$ obeys the bound $|\mathfrak{e}| \le c_0 \varepsilon \tfrac{1}{t}$. This leads to the inequality

$$- \tfrac{d^2}{dt^2} f + \tfrac{2}{t} \tfrac{d}{dt} f < 0 \tag{4.21}$$

when both $|\varphi|^2 \le \tfrac{1}{4t^2}$ and $\varepsilon \le c_0^{-1} \inf_{t \in [t_1, t_2]} |\alpha(t,z(t),0)|$. (The matrix $\mathcal{H}$ is negative definite because z(t) is a local maximum for the function q(t).) Assume henceforth that $\varepsilon$ is small enough so that (4.21) is obeyed. (The upper bound for a suitable $\varepsilon$ is positive. This is because $|\alpha(t,z(t),0)|$ is not less than the infimum of $|\alpha|$ on the whole of $[t_1, t_2] \times \mathbb{R}^2 \times \mathbb{R}$ which is bounded away from zero since it is the infimum of $t|q_*|$ on this interval.)

The inequality in (4.21) can be written equivalently as

$$\tfrac{d}{dt}(\tfrac{1}{t^2} \tfrac{d}{dt} f) > 0 \tag{4.22}$$

which is the form that is used below.

*Part 9*: Suppose that $t \in [t_1, t_2)$ and that $\Delta > 0$ is such that $t + \Delta < t_2$. Then integrating (4.22) from t to $t+\Delta$ leads the inequality

$$(\tfrac{d}{dt} f)(t+\Delta) \ge \tfrac{(t+\Delta)^2}{t^2} (\tfrac{d}{dt} f)(t) \tag{4.23}$$

Now suppose that $t_{\ddagger} \in (t_1, t_2)$ is such that there are two maxima for $q$ at $t_{\ddagger}$, one at the point $z_1(t_{\ddagger})$ and the other at the point $z_2(t_{\ddagger})$. These have corresponding function $f_1$ and $f_2$ which are equal at $t = t_{\ddagger}$. The Cerf construction of $q$ is such that the following is true: If $\mu > 0$ and sufficiently small (require among other things that $t_{\ddagger}-\mu > t_1$ and $t_{\ddagger}+\mu < t_2$), then the derivative of one of these functions (take it to be $f_1$) is greater than that of the other on the whole of $(t_{\ddagger}-\mu, t_{\ddagger}+\mu)$. Thus

- $f_1(t_{\ddagger}) = f_2(t_{\ddagger})$.
- $\frac{d}{dt} f_1 > \frac{d}{dt} f_2$ on $(t_{\ddagger}-\mu, t_{\ddagger}+\mu)$.

(4.24)

Now define a function $\mathfrak{f}$ on $(t_{\ddagger}-\mu, t_{\ddagger}+\mu)$ by the rule whereby $\mathfrak{f} = f_2$ for $t \leq t_{\ddagger}$ and $\mathfrak{f} = f_1$ for $t > t_{\ddagger}$. By virtue of (4.24), this function $\mathfrak{f}$ is given by $\mathfrak{f}(t) = \max(f_1(t), f_2(t))$ which is to say that $\mathfrak{f}(t)$ is $\sup_z \{q(t,z,0)\}$. (It is therefore $\varepsilon$ close to $\sup_x \{q(t,z,0)\}$.)

Note that $\mathfrak{f}$ is smooth on the complement of $t_{\ddagger}$ and it is Lipshitz with Lipshitz constant independent of the approximation $q$ (the latter by virtue of (4.19)). The function $\mathfrak{f}$ also obeys (4.23) so long as neither $t$ nor $t+\Delta$ equal $t_{\ddagger}$. Indeed, this follows directly from (4.23) when $t$ and $t+\Delta$ are on the same side of $t_{\ddagger}$, and it follows from (4.23) and the second bullet of (4.24) when they are on different sides of $t_{\ddagger}$.

*Part 10*: To complete the proof, consider now a sequence of $q$'s as just described that converge uniformly on $[t_1, t_2] \times \mathbb{R}^2 \times \mathbb{R}$ to $q$. The sequence is denoted by $\{q_n\}_{n \in \mathbb{N}}$ and nothing is lost by requiring that for each n, the $q_n$ version of $\varepsilon$ is less than $\frac{1}{n^2}$. There is a corresponding version of $\mathfrak{f}$ for each $q_n$ which is denoted by $\mathfrak{f}_n$. The sequence $\{\mathfrak{f}_n\}_{n \in \mathbb{N}}$ converges uniformly on $[t_1, t_2]$ to $q_*$. Since the sequence $\{\mathfrak{f}_n\}_{n \in \mathbb{N}}$ is uniformly Lipshitz (by virtue of (4.19)), the function $q_*$ is also Lipshitz.

Each $\mathfrak{f}_n$ is smooth on the complement of a finite set and it obeys (4.23) supposing that both $t$ and $t+\Delta$ avoid this same finite set. Call this set $U_n$ (it might depend on n). It then follow that $q_*$ is differentiable on the complement of $\cup_{n \in \mathbb{N}} U_n$ and its derivative there is the limit of the derivatives of the functions $q_n$. This being the case, $q_*$ also obeys (4.23) on the complement of $\cup_{n \in \mathbb{N}} U_n$; which is to say that it obeys (4.16) on the complement of $\cup_{n \in \mathbb{N}} U_n$. With that understood, the proof ends by noting that $\cup_{n \in \mathbb{N}} U_n$ is the countable union of finite sets, so it is countable and thus has measure zero in $[t_1, t_2]$.

**c) Rescalings and translations**

Proposition 4.3 implies in part that the following constraint set is equivalent to CONSTRAINT SET 1 in the sense that the set of solutions that are described by one is the same set of solutions that are described by the other.

CONSTRAINT SET 1´: *A pair $(A, \mathfrak{a})$ is described by this constraint set when*
- *The function $\mathfrak{t}|\mathfrak{a}|$ is uniformly bounded on $(0,\infty)\times\mathbb{R}^2\times\mathbb{R}$.*
- *There exists $\varepsilon > 0$ that $\mathfrak{t}|\mathfrak{a}_3| > \varepsilon$ on $(0,\infty)\times\mathbb{R}^2\times\mathbb{R}$.*
- *$\varphi$ is not identically zero.*
- *$\langle \mathfrak{a}_3 \varphi \rangle$ is identically zero.*

(4.25)

This differs from CONSTRAINT SET 1 in the second bullet. The important point with regards to CONSTRAINT SET 1´ is this: If $(A, \mathfrak{a})$ is described by CONSTRAINT SET 1´, and if $\mathfrak{z}$ and $\varepsilon$ are such that

$$\frac{\varepsilon}{\mathfrak{t}} \leq |\mathfrak{a}_3| \quad and \quad |\mathfrak{a}| \leq \frac{\mathfrak{z}}{\mathfrak{t}}$$

(4.26)

on the whole of $(0,\infty)\times\mathbb{R}^2\times\mathbb{R}$, then any solution to (1.4) that is obtained from $(A, \mathfrak{a})$ by the action via pull-back of any translation along the $\mathbb{R}^2$ factor of $(0,\infty)\times\mathbb{R}^2\times\mathbb{R}$ or any coordinate rescaling also obeys the bounds in (4.26) with the *same* $\varepsilon$ and $\mathfrak{z}$. This fact is exploited to prove an assertion to the effect that the action on the set $\mathcal{M}_f$ (via pull-backs) of the group of translations along the $\mathbb{R}^2$ factor and coordinate rescalings has closed orbits. (Remember that $\mathcal{M}_f$ consists of pairs $(A, \mathfrak{a})$ obeying (1.4), satisfying CONSTRAINT SET 1 and with the corresponding set Z being finite. By virtue of Proposition 4.3, elements in the set $\mathcal{M}_f$ are described by CONSTRAINT SET 1´ also.)

**Lemma 4.4**: *Fix a pair $(A, \mathfrak{a}) \in \mathcal{M}_f$ and fix a sequence $\{(\lambda_n, a_n)\}_{n\in\mathbb{N}}$ from $(0,\infty)\times\mathbb{R}^2$. Let $\{(A^n, \mathfrak{a}^n)\}_{n\in\mathbb{N}}$ denote the sequence in $\mathcal{M}_f$ whose k'th element is the pull-back of $(A, \mathfrak{a})$ by the diffeomorphism of $(0,\infty)\times\mathbb{R}^2\times\mathbb{R}$ that is defined by the rule whereby $(\mathfrak{t}, z, x_3)$ is sent to $(\lambda_n \mathfrak{t}, \lambda_n z + a_n, \lambda_n x_3)$. There is a subsequence in $\mathbb{N}$ (to be denoted by $\Lambda$), and a corresponding sequence of automorphisms of P (denoted by $\{g_n\}_{n\in\Lambda}$) such that $\{(g_n{}^*A^n, g_n{}^*\mathfrak{a}^n)\}_{n\in\Lambda}$ converges in the $C^\infty$ topology on compact subsets of $(0,\infty)\times\mathbb{R}^2\times\mathbb{R}$ to an element in $\mathcal{M}_f$. Moreover, if $(A, \mathfrak{a})$ is not a Nahm pole imposter, then the two $\Lambda$-indexed sequences whose n'th elements are the respective $(g_n{}^*A^n, g_n{}^*\mathfrak{a}^n)$ versions of $\sigma$ and $\alpha$ converge in the $C^\infty$ topology on compact subsets of $(0,\infty)\times\mathbb{R}^2\times\mathbb{R}$ to the versions of $\sigma$ and $\alpha$ that are defined by the limit in $\mathcal{M}_f$.*

*Proof of Lemma 4.4*: The proof has three parts. Part 1 considers the convergence assertions without regard to whether the limit is in $\mathcal{M}_f$. Parts 2 and 3 prove that the version of $\varphi$ for the limit solution has a finite number of zeros in each constant $(\mathfrak{t}, x_3)$ slice of $(0,\infty)\times\mathbb{R}^2\times\mathbb{R}$.

*Part 1*: Lemma 2.3 supplies the subsequence $\Lambda \subset \mathbb{N}$, the sequence $\{g_n\}_{n \in \Lambda}$ of automorphisms of $P_0$, and the limit pair $(A^\infty, \mathfrak{a}^\infty)$ for the sequence $\{(g_n*A^n, g_n*\mathfrak{a}^n)\}_{n \in \Lambda}$. According to Lemma 2.3, this limit pair obeys the first two bullets of CONSTRAINT SET 1´ and the fourth bullet.

With regards to the convergence of the $\sigma$ and $\alpha$ sequences when $(A, \mathfrak{a})$ is not a Nahm pole imposter: A key point is that $\sigma = -\frac{\mathfrak{a}_3}{|\mathfrak{a}_3|}$ is $\nabla_A$-covariantly constant. With this understood, let $\sigma^n$ and $\alpha^n$ for $n \in \Lambda$ denote the $(A^n, \mathfrak{a}^n)$ versions of $\sigma$ and $\alpha$. Since $(A, \mathfrak{a})$ is not a Nahm pole imposter, the section $\sigma^n$ is $\nabla_{A^n}$-covariantly constant. This implies that the sequence $\{g_n*\sigma^n\}_{n \in \Lambda}$ converges in the $C^\infty$ topolology on compact subsets of $(0, \infty) \times \mathbb{R}^2 \times \mathbb{R}$ to a unit normed, $\nabla_{A^\infty}$-covariantly constant section of $ad(P)$. (This is because each $n \in \Lambda$ version of $g_n*\sigma^n$ is covariantly constant for the connection $g_n*A^n$ and because $\{g_n*A^n\}_{n \in \Lambda}$ converges in the $C^\infty$ topology on compact subsets to $A^\infty$.) Denote this unit normed section by $\sigma^\infty$. Since each $n \in \Lambda$ version of $g_n*\mathfrak{a}^n_3$ can be written as $\alpha^n g_n*\sigma^n$, it follows that $\mathfrak{a}^\infty_3$ can be written as $\alpha^\infty \sigma^\infty$ with $\alpha^\infty$ being a smooth function on $(0, \infty) \times \mathbb{R}^2 \times \mathbb{R}$. It then follows from the manner of $C^\infty$ convergence of $\{g_n*\sigma^n\}_{n \in \Lambda}$ and $\{g_n*\mathfrak{a}^3_n\}_{n \in \Lambda}$ that the function sequence $\{\alpha^n\}_{n \in \Lambda}$ converges in the $C^\infty$ topology to $\alpha^\infty$.

*Part 2*: The zero locus of $\varphi^\infty$ is the concern in this part of the proof. In this regard, there are, a priori, five possibilities for the zero locus of $\varphi^\infty$ in any given $\{t\} \times \mathbb{R}^2 \times \{x_3\}$ slice: The first one is that $\varphi^\infty$ is identically zero; and the remaining four are as follow:

- *There are no zeros of $\varphi^\infty$.*
- *The zeros of $\varphi^\infty$ are in 1-1 correspondence with those of $\varphi$ with corresponding zeros having the same local degree.*
- *There is just one zero of $\varphi^\infty$ and it corresponds to just a single zero of $\varphi$, its degree being the same as that zero.*
- *There is just one zero of $\varphi^\infty$ and its local degree is the sum of the degrees of the zeros of $\varphi$.*

(4.27)

(If $\varphi^\infty$ is not identically zero, then $(A^\infty, \mathfrak{a}^\infty)$ is described by Proposition 3.1. In particular, if a given bullet of (4.27) holds in any one slice constant $(t, x_3)$ slice of $(0, \infty) \times \mathbb{R}^2 \times \mathbb{R}$, then that bullet must hold in all of them.) To see why (4.27) plus the $\varphi^\infty \equiv 0$ option acccounts for all possibilities, assume for the moment that $\varphi^\infty$ is not identically zero. It is then a $\nabla_{\hat{A}^{(\infty)}}$-holomorphic section of the corresponding version of $\mathcal{L}^+$ and the zeros in any $\{t\} \times \mathbb{R}^2 \times \{x_3\}$ slice are a discrete set. Fix $p \in \mathbb{R}^2$ so that $\varphi^\infty$ is zero at $(1, p, 0)$. The function $|\varphi^\infty|$ is then strictly positive on any given, sufficiently small but positive radius circle about $p$; and therefore $\varphi^\infty/|\varphi^\infty|$ has some strictly positive degree on such a circle (which is independent of the radius if the latter is small). This degree is defined as

follows: Fix a small radius disk about p whose closure is disjoint from all other zeros of $\varphi^\infty$. Let $\mathcal{L}^{(\infty)+}$ denote the $(A^\infty, \mathfrak{a}^\infty)$ version of $\mathcal{L}^+$. Fix an isomorphism from $\mathcal{L}^{(\infty)+}$ on this disk with the product $\mathbb{C}$ bundle. This isomorphism identifies $\varphi^\infty$ as a $\mathbb{C}$-valued function which has an expansion near the origin that has the form $\eta(z-p)^m + \mathcal{O}(|z-p|^{m+1})$ with $\eta \neq 0$ and with m being a positive integer. This integer n is the degree.

This paragraph is a digression about $\mathcal{L}^{(\infty)+}$. View this complex line bundle as a line subbundle in $ad(P)_\mathbb{C}$. Viewed as such, it is the pointwise limit of the sequence of complex line subbundles (parametrized by $\Lambda$) whose n'th member is the pull-back via $g_n$ of the $(A^n, \mathfrak{a}^n)$ version of $\mathcal{L}^+$. To see that this is so, let $\sigma^\infty$ denote the $(A^\infty, \mathfrak{a}^\infty)$ version of $\sigma$. Then $\mathcal{L}^{(\infty)+}$ is the +1 eigenspace of the endomorphism $[\frac{i}{2}\sigma^\infty, \cdot]$. Keep in mind also that $\nabla_{A^\infty} \sigma^\infty \equiv 0$. Meanwhile, if $n \in \Lambda$, then the $(A^n, \mathfrak{a}^n)$ version of $\mathcal{L}^+$ is the +1 eigenspace of the endomorphism $[\frac{i}{2}\sigma^n, \cdot]$ with $\sigma^n$ denoting here and in what follows the $(A^n, \mathfrak{a}^n)$ version of $\sigma$. The key point now (as proved in Part 1) is that the sequence $\{g_n{}^*\sigma^n\}_{n \in \Lambda}$ converge to $\sigma^\infty$ in the $C^\infty$-topology on compact subsets of $(0, \infty) \times \mathbb{R}^2 \times \mathbb{R}$. Since $\mathcal{L}^{(\infty)+}$ is the +1 eigenspace of $[\frac{i}{2}\sigma^\infty, \cdot]$ and each $n \in \Lambda$ version of $\mathcal{L}^+$ is the +1 eigenspace of $[\frac{i}{2}\sigma^n, \cdot]$, the automorphism $g_n h_n$ then identifies (in $\mathbb{C}$-linear fashion) the large $n \in \Lambda$ versions of $\mathcal{L}^+$ on the given open set with $\mathcal{L}^{\infty+}$ on that open set. Such an identification is implicitly used in what follows to view the sufficiently large n versions of $\varphi$ as being sections of $\mathcal{L}^{(\infty)+}$ on the given open set.

With the digression now over, let C denote a very small radius circle centered at p such that p is the only zero of $\varphi^\infty$ in the disk of radius 100 times C's radius centered at p. Given the $C^\infty$-nature of the convergence, the preceding observations about $\varphi^\infty$ on C implies that each sufficiently large n version of $\varphi^n$ must have positive norm on C, and it implies that $\varphi^n/|\varphi^n|$ must have the same positive degree around C as does $\varphi^\infty/|\varphi^\infty|$ when n is large (What is denoted here by $\varphi^n$ is the $(A^n, \mathfrak{a}^n)$ version of $\varphi$.) This implies in turn that $\varphi^n$ has one or more zeros inside C whose degrees sum up to the degree of $\varphi^\infty$ at p.

Now each $n \in \Lambda$ version of $\varphi^n$ is, after all, just a translation and rescaling of $\varphi$, so $|\varphi|$ is strictly positive on the circle that is obtained from C by undoing the translation and rescaling. And, the sum of the degrees of $\varphi$ inside that circle (call it $C_n$) must sum up to the degree of $\varphi^\infty$ at p. Granted this, then bullets 1-4 of (4.27) account for the possibilities for the collection of circles $\{C_n\}_{n \in \Lambda}$ (which depends in turn on the behavior of $\{\lambda_n\}_{n \in \Lambda}$ and $\{a_n/\lambda_n\}_{n \in \Lambda}$.) Keep in mind here that a point $p_* \in \mathbb{R}^2$ is a zero of $\varphi$ if and only if the point $\lambda_n^{-1}(p_* - a_n)$ is a zero of $\varphi^n$.

It remains now to rule out the possibility that $\varphi^\infty$ is identically zero. This is done in the next part of the proof.

*Part 3*: If $(A, \mathfrak{a})$ is a Nahm pole imposter, then it is Aut(P) equivalent to any rescaling or translation of itself. This implies that each $n \in \mathbb{N}$ version of $|\varphi^n|$ is identical to $|\varphi|$. And, that implies in turn that $\varphi^\infty \equiv 0$ if and only if $\varphi \equiv 0$ (which is to say that $(A, \mathfrak{a})$ is aut(P) equivalent to the $|w| = 1$ limit of the Nahm pole imposters).

With the preceding understood, assume that $(A, \mathfrak{a})$ is not a Nahm pole imposter. As explained in Part 1, the sequence $\{g_n * \sigma^n\}_{n \in \Lambda}$ converges on compact subsets of $(0, \infty) \times \mathbb{R}^2 \times \mathbb{R}$ to a unit length, $\nabla_{A^\infty}$-covariantly constant section of ad(P). This section is again denoted by $\sigma^\infty$; it is relevant here because $\mathfrak{a}_3^\infty$ can then be written as $\alpha^\infty \sigma^\infty$ with $\alpha^\infty$ being the limit of the sequence $\{\alpha^n\}_{n \in \Lambda}$. This function $\alpha^\infty$ obeys the $(A^\infty, \mathfrak{a}^\infty)$ version of (4.8) which says this

$$--(\tfrac{\partial^2}{\partial t^2} + \tfrac{\partial^2}{\partial x_1^2} + \tfrac{\partial^2}{\partial x_2^2}) \alpha^\infty = 0 \; .$$

(4.28)

if $\varphi^\infty$ is identically zero. Thus, $\alpha^\infty$ is harmonic.

In the case at hand, the bounds in (4.26) imply

$$\tfrac{\varepsilon}{t} \leq |\alpha^\infty| \leq \tfrac{\mathfrak{z}}{t} \; ,$$

(4.29)

with $\varepsilon$ and $\mathfrak{z}$ being positive numbers. (Note in this regard that $\alpha^\infty < 0$ because each $\alpha^n$ is strictly negative.) As explained in the next paragraph, the bounds in (4.29) are not possible if (4.28) holds.

Fix $T > 4$. Let $f_T$ denote the function of just the coordinate t given by the rule

$$t \to f_T(t) = \varepsilon \tfrac{t}{T} - \varepsilon (1 - \tfrac{1}{T})$$

(4.30)

This function is equal to $-\varepsilon$ at $t = 1$ and it is equal to $-\tfrac{\varepsilon}{T}$ at $t = T$. As a consequence, it is greater than or equal to $\alpha^\infty$ at both $t = 1$ and at $t = T$. Because both $f_T$ and $\alpha^\infty$ are harmonic, the function $f_T$ must be greater than or equal $\alpha$ on the interval $[t, T]$. Now, by virtue of (4.29), $\alpha^\infty \geq -\tfrac{\mathfrak{z}}{t}$ at all points. Therefore, the inequality

$$\varepsilon \tfrac{t}{T} - \varepsilon(1 - \tfrac{1}{T}) \geq -\tfrac{\mathfrak{z}}{t}$$

(4.31)

must hold for $t \in (1, T)$. In particular, if this is to hold at $t = \tfrac{1}{2} T$, then the bound

$$-\tfrac{1}{4} \varepsilon \geq -2 \tfrac{\mathfrak{z}}{t}$$

(4.32)

must hold; which is not possible if $T > \tfrac{8 \mathfrak{z}}{\varepsilon}$.

### d) Proof of Theorem 4.1

The proof of Theorem 4.1 has six parts. The assumption at the start and througout is that $(A, \mathfrak{a})$ and $(A´, \mathfrak{a}´)$ are two elements in the space $\mathcal{M}_f$ with corresponding $\varphi$ and $\varphi´$ having the same zero locus in any given $\{t\} \times \mathbb{R}^2 \times \{x_3\}$. In addition, the degree of vanishing of $\varphi$ at any given zero is the same as that of $\varphi´$.

*Part 1*: Write $\mathfrak{a}_3$ as $\alpha \sigma$ with $|\sigma| = 1$ and with $\alpha < 0$. Meanwhile, $\mathfrak{a}_3´$ can be written as $f \alpha \sigma´$ with $f$ denoting a strictly positive function on $(0, \infty) \times \mathbb{R}^2 \times \mathbb{R}$. (But it is independent of $x_3$ coordinate, as is $\alpha$.) Since neither $\alpha$ nor $f$ are zero, there is an automorphism of P on $(0, \infty) \times \mathbb{R}^2 \times \mathbb{R}$ that identifies $\sigma$ with $\sigma´$. After acting by this automorphism, write $\mathfrak{a}_3´$ as $f \alpha \sigma$. (Both $\sigma$ and $\sigma´$ can be identified with the $\nabla_{\theta_0}$-constant element $\sigma_3$ that appears in (1.6) if needs be.) There is still the freedom to act by automorphisms of P that fix $\sigma$. This can be used to write $\varphi´$ as $h \varphi$ with h being a real and strictly positive function. Such an automorphism exists because $\varphi$ and $\varphi´$ have to the same set of zeros and because $\varphi$'s degree at each of these zeros is the same as that $\varphi´$.

Meanwhile, the ad(P)-valued 1-form $A´ - A$ is proportional to $\sigma$ because neither $(A, \mathfrak{a})$ nor $(A´, \mathfrak{a}´)$ is a Nahm pole imposter. (This is where the latter assumption is used in the proof.) As a consequence, the connection $A´$ can be written as $A + (b_t dt + b) \sigma$ with $b_t$ denoting a real function on $(0, \infty) \times \mathbb{R}^2 \times \mathbb{R}$ and with b denoting an $\mathbb{R}$-valued 1-form on $(0, \infty) \times \mathbb{R}^2 \times \mathbb{R}$ that annilihates the tangent vectors to the $(0, \infty)$ factor.

*Part 2*: The $(A´, \mathfrak{a}´)$ version of the equation in in the top bullet of (3.9) when written using $b_t$, b, h and $f$ asserts this:

$$\tfrac{\partial}{\partial t} h - 2i b_t h + 2\alpha h = 2\alpha f h \ .$$

(4.33)

The only imaginary term is $b_t h$ and since h is not zero, $b_t$ must vanish identically. Meanwhile, the real part of (4.33) relates $f$ and h:

$$\tfrac{\partial}{\partial t} h = 2\alpha (f - 1) h \ .$$

(4.34)

Now write $b_1 + i b_2$ as $\mu$. Then the prime version of the second bullet in (3.9) asserts that

$$\tfrac{\partial}{\partial \bar{z}} h = i \mu h \ .$$

(4.35)

Since h is strictly positive, it has a natural logarithm everywhere. This logarithm is denoted by 2w. The identities in (4.36) and (4.37) can be written using w as:

$$\tfrac{\partial}{\partial t} w = \alpha(f - 1) \quad \text{and} \quad \tfrac{\partial}{\partial z} w = \tfrac{i}{2}\mu .$$

(4.36)

Finally, the primed version of (3.10) can be written (with the help of (4.36)) as:

$$-\left(\tfrac{\partial^2}{\partial t^2} + \tfrac{\partial^2}{\partial x_1^2} + \tfrac{\partial^2}{\partial x_2^2}\right) w = (1 - e^{4w}) |\varphi|^2$$

(4.37)

(See Section 7a for the more details about the derivation.)

The key observation with regards to (4.37) comes via the maximum principle:

*If* w *obeys* (4.37), *then it has neither positive local maximima nor negative local minima.*

(4.38)

Here is the argument for (4.38) with regards to local maxima: The integral of $|\varphi|^2$ about any ball in $(0,\infty) \times \mathbb{R}^2 \times \mathbb{R}$ is positive. This implies that the integral of the right hand side of (4.37) on any sufficiently small radius ball centered on any $w > 0$ point is negative. On the other hand, if the ball has sufficiently small radius and center on a local maximum of w, then the integral of the left hand side of (4.37) is a negative multiple of the radial derivative of the average of w over balls with the given center. This derivative is non-positive for small radius balls if w has a local maximum at the center; so the integral of the left hand side of (4.37) is non-negative for small radius balls centered on local maxima of w. (The argument is simpler if $\varphi$ is non-zero at a hypothetical positive local maximum. In this case, the right hand side of (4.37 is negative at the positive local maximum whereas the left hand side would be non-negative.)

But for notation (changing positive to negative) the same argument rules out the existence of negative local minimum.

*Part 3*: This part of the proof establishes Theorem 4.1 for those pairs in $\mathcal{M}_f$ where the corresponding version of $\varphi$ is non-vanishing. This is done by comparing any given such $(A', \mathfrak{a}')$ with the Nahm pole solution. Remember that the latter is the $m = 0$ model solution in Section 1c; it has $\mathfrak{a}_3 = -\tfrac{1}{2t}\sigma$ and $|\varphi| = \tfrac{1}{\sqrt{2}t}$. The pair $(A', \mathfrak{a}')$ can be any element in $\mathcal{M}_f$ whose corresponding $\varphi'$ has no zeros. Given what is said in Parts 1 and 2, it is sufficient to prove that the function w is identically zero. To this end, note that w in this case is half of the natural logarithm of $\sqrt{2}t|\varphi'|$; and (4.37) in this case is

$$-\left(\tfrac{\partial^2}{\partial t^2} + \tfrac{\partial^2}{\partial x_1^2} + \tfrac{\partial^2}{\partial x_2^2}\right) w = (1 - e^{4w}) \tfrac{1}{2t^2} .$$

(4.39)

The proof that $w = 0$ has five steps.

Step 1: This step and the next two steps prove the following assertion:

*If w is not identically zero, then* w *can not be a bounded function on* $(0, \infty) \times \mathbb{R}^2 \times \mathbb{R}$.

(4.40)

To see this, suppose for the sake of argument that w > 0 at one or more points, and that it is bounded. As explained directly, this assumption leads to nonsense, a violation of (4.38). To start the explanation, let $\mathfrak{w}$ denote the supremum of w. Fix $x_3 = 0$ and fix a sequence of points, $\{(t_n, z_n)\}_{n \in \mathbb{N}} \subset (0, \infty) \times \mathbb{R}^2$ such that the numerical sequence $\{w(t_n, z_n, 0)\}_{n \in \mathbb{N}}$ is increasing and converging as $n \to \infty$ to $\mathfrak{w}$. The sequence $\{(t_n, z_n)\}_{n \in \mathbb{N}}$ has no convergent subsequence because then w would take on a positive local maximum. But, the equation in (4.39) for an $x_3$-independent w is invariant with respect to translations along the $\mathbb{R}^2$ factor (for example $z \to z + a$) and with respect to rescalings $(t, z) \to (\lambda t, \lambda z, \lambda x_3)$ for $\lambda \in (0, \infty)$. Therefore, rather than considering w on a sequence of points, one can consider a sequence of solutions $\{w_n(t, z)\}_{n \in \mathbb{N}}$ to (4.39) such that each function $w_n$ is bounded above by $\mathfrak{w}$ and such that $\lim_{n \to \infty} w_n(t = 1, z = 0) = \mathfrak{w}$.

Step 2: Of course, for each index n, the function $w_n$ is half of the natural log of $\sqrt{2t}|\varphi^n|$ where $(A^n, \mathfrak{a}^n)$ is the rescaled and translated version of $(A', \mathfrak{a}')$. By virtue of Lemma 4.4 (using $(A', \mathfrak{a}')$ for $(A, \mathfrak{a})$), the sequence $\{(A^n, \mathfrak{a}^n)\}_{n \in \mathbb{N}}$ has a subsequence (labeled by a subsequence $\Lambda \subset \mathbb{N}$) that converges after termwise action of automorphisms of P to a pair $(A^\infty, \mathfrak{a}^\infty)$ from $\mathcal{M}_f$. Since $\varphi'$ has no zeros, the first bullet of (4.27) is the only viable option; thus $\varphi^\infty$ has no zeros.

The convergence just described implies that the sequence $\{w_n\}_{n \in \Lambda}$ converges in the $C^\infty$ topology to $\frac{1}{2} \ln(\sqrt{2t}|\varphi^\infty|)$ which is the version of the function w using the Nahm pole solution for $(A, \mathfrak{a})$ and using $(A^\infty, \mathfrak{a}^\infty)$ in lieu of $(A', \mathfrak{a}')$. This new version of w (call it $w_\infty$) takes on the value $\mathfrak{w}$ and this is its maximum. This contradicts (4.38).

Step 3: The same argument (except for sign changes) as the one just given proves that the solution to (4.39) coming from $(A', \mathfrak{a}')$ can not be bounded from below if it is negative at any point.

Step 4: The function w does indeed have an upper bound. This is because $|\varphi'|$ is $e^{2w}|\varphi|$ which is $e^{2w} \frac{1}{\sqrt{2t}}$; and so an unbounded w would run afoul of the condition in the top bullet of CONSTRAINT SET 1 (or CONSTRAINT SET 1´ in (4.25)). Thus, w must be negative if it is not identically zero.

Step 5: This step proves that w can't be negative. To this end, suppose for the sake of argument that it were (so as to derive nonsense). Choose a sequence from $(0, \infty) \times \mathbb{R}^2$ to be denoted by $\{(t_n, z_n)\}_{n \in \mathbb{N}}$ with the property that $\{w(t_n, z_n)\}_{n \in \mathbb{N}}$ is decreasing and unbounded. For each index n, translate z so that the origin in $\mathbb{R}^2$ is mapped to $z_n$, and then rescale coordinates $(t,z) \to (t_n t, t_n z, t_n x_3)$. Pull back $(A´, \mathfrak{a}´)$ by the concatenation of these two maps to get a new pair from the space $\mathcal{M}_f$. The index n version of this new pair is denoted by $(A^n, \mathfrak{a}^n)$.

Appeal to Lemma 4.4 to obtain a subsequence $\Lambda \subset \mathbb{N}$ such that $\{(A_n, \mathfrak{a}_n)\}_{n \in \Lambda}$ converges to a pair in $\mathcal{M}_f$ after term-wise applications of automorphisms of P. The limit is denoted by $(A^\infty, \mathfrak{a}^\infty)$. As before, the only viable option in (4.27) for $\varphi^\infty$ is that it be non-zero. But, this conclusion is nonsensical because it runs afoul of the fact that $|\varphi^\infty|$ at $(t = 1, z = 0, x_3 = 0)$ is the limit (for $n \in \Lambda$ and $n \to \infty$) of the values of $|\varphi^n|$ at this point, and that limit is zero since $\lim_{i \in \Lambda} w(t_n, z_n) \to -\infty$.

*Part 4*: The case just considered (where $\varphi$ has no zeros) has the following implications for the general case:

**Lemma 4.5**: *Assume that $(A, \mathfrak{a})$ is from $\mathcal{M}_f$ and not a Nahm pole imposter. Given $\varepsilon > 0$, there exists $\kappa > 1$ with the following significance: Fix $t \in (0, \infty)$. The components $\{\mathfrak{a}_i\}_{i=1,2,3}$ of $\mathfrak{a}$ obey $|1 - 2t|\mathfrak{a}_i|| < \varepsilon$ at all points in $\{t\} \times \mathbb{R}^2 \times \mathbb{R}$ where the distance to $\varphi^{-1}(0)$ is greater than $\kappa t$.*

*Proof of Lemma 4.5*: If there were no such $\kappa$, then suitably chosen sequences of translations along the $\mathbb{R}^2$ factor of $(0, \infty) \times \mathbb{R}^2 \times \mathbb{R}$ and rescalings of the coordinates could be used to construct a sequence of elements in $\mathcal{M}_f$ that converges to an element in $\mathcal{M}_f$ with two contradictory properties: The limit pair is not equivalent to the Nahm pole solution via an automorphism of P; and its version of $\varphi$ is non-vanishing. (The convergence follows from Lemma 4.4 and the non-vanishing of $\varphi^\infty$ follows from (4.27)).

*Part 5*: This part of the proof considers the case of Theorem 4.1 where the zero locus of the corresponding $\varphi$ for any pair is $(0, \infty) \times \{0\} \times \mathbb{R}$ and the vanishing degree there is some positive integer (to be denoted by $m$). Take $(A, \mathfrak{a})$ to be the integer $m$ version of the model solution $(A^{(m)}, \mathfrak{a}^{(m)})$. Suppose that $(A´, \mathfrak{a}´)$ is a second element in $\mathcal{M}_f$ with $\varphi´$ vanishing only on $(0, \infty) \times \{0\}$ and with degree $m$ also. The proof that $(A´, \mathfrak{a}´)$ is equivalent via an automorphism of P to $(A^{(m)}, \mathfrak{a}^{(m)})$ has three steps.

<u>Step 1</u>: Define w to be half of the natural logarithm of $|\varphi´|/|\varphi^{(m)}|$. Thus w obeys (4.37) with $\varphi = \varphi^{(m)}$. Suppose for the sake of argument that w is not identically zero. (There would be nothing to say if $w \equiv 0$.) If w is positive somewhere, let $\mathfrak{w}$ denote either its supremum on $(0, \infty) \times \mathbb{R}^2 \times \mathbb{R}$ or $+\infty$ if w has no upper bound on $(0, \infty)$. If w is non-positive, let $\mathfrak{w}$ denote either the infimum of w on $(0, \infty) \times \mathbb{R}^2 \times \mathbb{R}$ or $-\infty$ if w is not bounded from below. Because w lacks positive local maximum and negative local minimum (see (4.40)), there is a non-convergent sequence $\{(t_n, z_n)\}_{n \in \mathbb{N}} \subset (0, \infty) \times \mathbb{R}^2$ such that the sequence $\{w(t_n,z_n,0)\}_{i\in\mathbb{N}}$ is monotonically approaching $\mathfrak{w}$ (increasing or decreasing as the case may be) and such that it is either unbounded (if $\mathfrak{w}$ is $\pm\infty$) or converges to $\mathfrak{w}$ if $\mathfrak{w}$ is finite. By virtue of the top bullet of CONSTRAINT SET 1 (the top bullet of (4.25)) and Lemma 4.5, the sequence $\{|z_n|/t_n\}_{n\in\mathbb{N}}$ must be bounded. (If not, then $\{w(t_n,z_n,0)\}_{n\in\mathbb{N}}$ would have a subsequence with limit zero.) Therefore, no generality is lost by choosing the sequence $\{(t_n,z_n)\}_{n\in\mathbb{N}}$ so that $\{z_n/t_n\}_{n\in\mathbb{N}}$ converges in $\mathbb{R}^2$.

Step 2: Fix an index n and let $w_n(t,z,x_3)$ denote $w(t_n t, t_n z, t_n x_3)$. These rescaled versions of w are such that $\lim_{n\to\infty} w_n(1,z_*,0) = \mathfrak{w}$. Meanwhile, for each index n, let $(A^n, \mathfrak{a}^n)$ denote the pull-back of $(A´,\mathfrak{a}´)$ via this rescaling map. Note that $w_n$ is equal to half of the natural logarithm of $|\varphi^n|/|\varphi^{(m)}|$ because $(A^{(m)}, \mathfrak{a}^{(m)})$ is scale invariant (see the remarks subsequent to (1.9)).

According to Lemma 4.4, there is a subsequence $\Lambda \subset \mathbb{N}$ such that $\{(A^n,\mathfrak{a}^n)\}_{n\in\Lambda}$ converges (after termwise action by automorphisms of P) to a pair $(A^\infty,\varphi^\infty)$ from $\mathcal{M}_f$. The only viable option from the list in (4.27) is the case where $\varphi^\infty$ vanishes only along $(0,\infty)\times\{0\}\times\mathbb{R}$ with degree $m$. (This is equivalent to the assertion that $\lim_{n\to\infty} z_n/t_n = 0$.)

Step 3: It follows from what is said in Step 2 and from the manner of convergence ($C^\infty$ convergence on compact subsets) that $\{w_n\}_{n\in\Lambda}$ must converge to $\frac{1}{2}\ln(|\varphi^\infty|/|\varphi^{(m)}|)$ which is the $(A^\infty,\varphi^\infty)$ version of w. This implies in turn that $\mathfrak{w}$ must be finite; and then it must be zero because it is either the maximum or minimum of $\frac{1}{2}\ln(|\varphi^\infty|/|\varphi^{(m)}|)$ and the latter has neither (by virtue of (4.38)).

*Part 6*: This part prove the remaining cases of Theorem 4.1. These are the cases where $\varphi$ (and thus $\varphi´$) has two or more zeros in any given $\{t\}\times\mathbb{R}^2\times\{x_3\}$ slice of $(0,\infty)\times\mathbb{R}^2\times\mathbb{R}$. Assume as before that $w = \ln(|\varphi´|/|\varphi|)$ is not identically zero (if it is, then $(A´,\mathfrak{a}´)$ is equivalent via an automorphism of P to $(A,\mathfrak{a})$). The three steps that follow derive nonsense from this assumption.

Step 1: Define $\mathfrak{w}$ as in the previous part (if w is positive, it is either the supremum of w, or $+\infty$; and if w is non-positive, then it is either the infimum of w or $-\infty$). As done in Part 5, there is a non-convergent sequence $\{(t_n,z_n)\}_{n\in\mathbb{N}} \subset (0,\infty)\times\mathbb{R}^2$ such that $\{w(t_n,z_n,0)\}_{n\in\mathbb{N}}$ is monotonically approaching $\mathfrak{w}$ (increasing or decreasing as the case may be) and such that it is either unbounded (if $\mathfrak{w}$ is $\pm\infty$) or converges to $\mathfrak{w}$ if $\mathfrak{w}$ is finite.

No generality is lost by choosing the sequence $\{(t_n,z_n)\}_{n\in\mathbb{N}}$ so that either $\{t_n\}_{n\in\mathbb{N}}$ is decreasing with limit zero; or $\{t_n\}_{n\in\mathbb{N}}$ is increasing and unbounded. Indeed, if $\{t_n\}_{i\in\mathbb{N}}$ is bounded, then $\{z_n\}_{n\in\mathbb{N}}$ must diverge, in which case $\{\text{dist}(z_n,\varphi^{-1}(0))/t_n\}_{i\in\mathbb{N}}$ diverges, in which case $\mathfrak{w}$ is zero (see Lemma 4.5).

In fact, by virtue of Lemma 4.5, the sequence $\{\text{dist}(z_n, \varphi^{-1}(0))/t_n\}_{n\in\mathbb{N}}$ must be bounded. Therefore, no generality is lost by assuming that there exists $p_* \in \mathbb{R}^2$ such that $\varphi$ is zero on $(0,\infty)\times\{p_*\}\times\mathbb{R}$ and such that $\{(p_*-z_n)/t_n\}_{n\in\mathbb{N}}$ converges. If $\{t_n\}_{n\in\mathbb{N}}$ diverges, then this will be the case for any choice of $p_*$; but if $\{t_n\}_{n\in\mathbb{N}}$ converges to 0, then there is a unique such $p_*$.

Step 2: For each index n, define $w_n$ by the rule $w_n(t,z,x_3) = w(t_n t, t_n z + p_*, t_n x_3)$. Also, let $(A^n, \mathfrak{a}^n)$ denote the corresponding pull-back of $(A, \mathfrak{a})$; and let $(A'^n, \mathfrak{a}'^n)$ denote the corresponding pull-back of $(A', \mathfrak{a}')$. Lemma 4.4 can be invoked to obtain a subsequence in $\mathbb{N}$ (to be denoted by $\Lambda$) such that both $\{(A^n, \mathfrak{a}^n)\}_{n\in\Lambda}$ and $\{(A'^n, \mathfrak{a}'^n)\}_{n\in\Lambda}$ converge to elements in $\mathcal{M}_f$ after suitable termwise action on each by automorphisms of P. The respective limits are denoted by $(A^\infty, \mathfrak{a}^\infty)$ and $(A'^\infty, \mathfrak{a}'^\infty)$. It follows from (4.27) that both $\varphi^\infty$ and $\varphi'^\infty$ vanish at a single point and that their respective degrees of vanishing at that point are equal. (If $\{t_n\}_{n\in\Lambda}$ limits to zero, then the third bullet of (4.27) describes the situation; and if $\{t_n\}_{n\in\Lambda}$ diverges, then the fourth bullet of (4.27) does.)

Step 3: It follows from what was just said in Step 2 (and from the $C^\infty$ convergence on compact subsets) that the sequence $\{w_n\}_{n\in\Lambda}$ converges to $\tfrac{1}{2}\ln(|\varphi'^\infty|/|\varphi^\infty|)$. This is necessarily a smooth function on $(0,\infty)\times\mathbb{R}^2\times\mathbb{R}$ because $(A^\infty, \mathfrak{a}^\infty)$ and $(A'^\infty, \mathfrak{a}'^\infty)$ obey (1.4) and they are described by the top bullet in CONSTRAINT SET 1, and because $\varphi^\infty$ and $\varphi'^\infty$ have the same single zero on any $\{t\}\times\mathbb{R}^2\times(x_3)$ slice and vanish there with identical degrees. This convergence implies that $\mathfrak{w}$ is the maximum (if positive) or minimum (if negative) of $\tfrac{1}{2}\ln(|\varphi'^\infty|/|\varphi^\infty|)$.

Granted the preceding, then $\mathfrak{w}$ must be zero (which is the desired nonsense) because of the following observations: The function $\tfrac{1}{2}\ln(|\varphi'^\infty|/|\varphi^\infty|)$ (call it $w^\infty$) obeys the equation in (4.37) with $\varphi$ replaced by $\varphi^\infty$. Therefore, the conclusions of (4.38) hold; which is to say that $w^\infty$ *does not* have a local maximum, nor a local minimum (let alone global ones).

## 5. CONSTRAINT SET 2 solutions: The norm of $|\varphi|$

This section and the next establish various properties of CONSTRAINT SET 2 solutions to (1.4) which are used in the subsequent sections. The principle observation in this section is stated by the proposition that follows directly.

**Proposition 5.1**: *Suppose that* $(A, \mathfrak{a})$ *is a solution to (1.4) that is described by* CONSTRAINT SET 2. *Given* $\delta > 0$, *there exists* $\kappa_\delta > 1$ *with the following significance: The function* $|\varphi|$ *obeys* $|\varphi| > \tfrac{1}{\sqrt{2t}}(1-\delta)$ *where* $|z| > \kappa_\delta t$ *on* $(0,\infty)\times\mathbb{R}^2\times\mathbb{R}$.

To put this proposition in perspective: The top bullet of (2.1) makes this same assertion for the case when $(A, \mathfrak{a})$ is one of the $m \geq 0$ model solutions from Section 1c.

By way of an outline for what is to come in this section: Section 5a states a lemma which is essentially parenthetical. It is included because of the insight it gives with regards to the integrals of $|F_A|^2$ over other sorts of unbounded sets in $(0,\infty) \times \mathbb{R}^2 \times \mathbb{R}$ when $(A, \mathfrak{a})$ obeys the first and third bullets of CONSTRAINT SET 2. Sections 5b and 5c supply various preliminary observations about CONSTRAINT SET 2 solutions (and other solutions) to (1.4) that are used in Section 5d's proof of Proposition 5.1.

There is an extra Sections 5e which uses Proposition 5.1 to establish an a priori upper bound for the vanishing degree of $\varphi$ (from a CONSTRAINT SET 2 version of $(A, \mathfrak{a})$) at the origin in the constant $(t, x_3)$ slices of $(0,\infty) \times \mathbb{R}^2 \times \mathbb{R}$. See Lemma 5.6.

### a) Integrals of $|F_A|^2$

The lemma that follows momentarily assumes that $(A, \mathfrak{a})$ is a solution to (1.4) that obeys the first and third bullets of CONSTRAINT SET 2. This lemma concerns the integral of $|B_{A3}|^2 + |E_{A1}|^2 + |E_{A2}|^2$ on $t > 0$ versions of the domain $[t,\infty) \times \mathbb{R}^2 \times \mathbb{R}$ and the domains where $x$ is bounded away from zero. (Remember that $x$ is shorthand for $(t^2 + |z|^2)^{1/2}$.)

With regards to notation: Supposing that $(A, \mathfrak{a})$ denotes a solution to (1.4) that is described by the third bullet of CONSTRAINT SET 2, there is a function on $(0,\infty)$ that is defined by the rule

$$t \to f(R) = \int_{(0,\infty) \times (\mathbb{R}^2 - D_R) \times \{0\}} (|B_{A3}|^2 + |E_{A1}|^2 + |E_{A2}|^2)$$

(5.1)

with $D_R$ denoting here the $|z| < R$ disk in $\mathbb{R}^2$. The third bullet of CONSTRAINT SET 2 asserts that the integral depicted on the right hand side of (5.4) is bounded by an R-independent multiple of $\frac{1}{R}$. Granted the proceding, the rest of this section and the subsequent ones use $\Xi$ to denote a positive number which is chosen so that the bound $f(R) < \frac{\Xi}{R}$ holds for all R. When $\Xi$ appears, it necessarily depends a priori on a given solution to (1.4) that is described by the third bullet in CONSTRAINT SET 2. This dependence will often be assumed implicitly.

**Lemma 5.2**: *Given positive numbers $\mathfrak{z}$ and $\Xi$, there exists $\kappa > 1$ with the following significance: Suppose that $(A, \mathfrak{a})$ is a solution to (1.4) with $|\mathfrak{a}| \leq \frac{\mathfrak{z}}{t}$ and with the function $f$ from (5.1) obeying $f(R) \leq \frac{\Xi}{R}$ for all $R > 0$. Granted these assumptions:*
- *The integral of $|B_{A3}|^2 + |E_{A1}|^2 + |E_{A2}|^2$ over any $t > 0$ version of $[t,\infty) \times \mathbb{R}^2 \times \{0\}$ is finite; and this integral is bounded by $\frac{\kappa}{t}$.*

- *If* $x > 0$, *then the integral of* $|B_{A3}|^2 + |E_{A1}|^2 + |E_{A2}|^2$ *over the part of* $(0,\infty) \times \mathbb{R}^2 \times \{0\}$ *where* $x > \mathrm{x}$ *is finite; and this integral is bounded by* $\frac{K}{\mathrm{x}}$.

*Proof of Lemma 5.2*: The second bullet of the lemma follows from the first one and from the $f(R) < \frac{\Xi}{R}$ bound because the $x > \mathrm{x}$ part of $(0,\infty) \times \mathbb{R}^2 \times \{0\}$ is contained in the union of the part of $(0,\infty) \times \mathbb{R}^2 \times \{0\}$ where $t \geq \frac{1}{\sqrt{2}} \mathrm{x}$ and the part where $|z| \geq \frac{1}{\sqrt{2}} \mathrm{x}$.

The argument for the top bullet is this: Lemma 2.2 can be invoked because $|\mathfrak{a}| \leq \frac{\mathfrak{z}}{t}$ on $(0,\infty) \times \mathbb{R}^2 \times \mathbb{R}$. According to Lemma 2.2, the integrand $|B_{A3}|^2 + |E_{A1}|^2 + |E_{A2}|^2$ on any $s > 0$ version of $\{s\} \times \mathbb{R}^2 \times \{0\}$ is bounded by $c_\mathfrak{z} \frac{1}{s^4}$ with $c_\mathfrak{z}$ being independent of s. Since the area of the $|z| \leq t$ disk in $\mathbb{R}^2$ is $\pi t^2$, the preceding bound for the integrand leads to the integral bound:

$$\int_{[t,\infty) \times D_t \times \{0\}} (|B_{A3}|^2 + |E_{A1}|^2 + |E_{A2}|^2) \leq c_\mathfrak{z} \frac{1}{t}.$$

(5.2)

Meanwhile, the integral of $|B_{A3}|^2 + |E_{A1}|^2 + |E_{A2}|^2$ over the domain $[t,\infty) \times (\mathbb{R}^2 - D_{R=t}) \times \{0\}$ is no greater than $\frac{\Xi}{t}$.

### b) Preliminary lower bounds for $|\mathfrak{a}_3|$ and $|\varphi|$

Th upcoming lemma assert lower bounds for $|\mathfrak{a}_3|$ and $|\varphi|$ at various points on the constant t slices of $(0,\infty) \times \mathbb{R}^2 \times \mathbb{R}$. To set the background for the lemma: The top bullet of CONSTRAINT SET 2 says that $|\mathfrak{a}| \leq \frac{\mathfrak{z}}{t}$ on the whole of $(0,\infty) \times \mathbb{R}^2 \times \mathbb{R}$ with $\mathfrak{z}$ being constant. Meanwhile, the second bullet of CONSTRAINT SET 2 asserts that $|\mathfrak{a}| \geq \frac{\varepsilon}{t}$ where t is less than $t_0$, with $t_0$ and $\varepsilon$ being positive numbers. Then, the third bullet of CONSTRAINT SET 2 says that the function $f$ in (5.1) obeys $f(R) \leq \frac{\Xi}{R}$ with $\Xi$ being independent of R.

**Lemma 5.3**: *Given positive numbers $\mathfrak{z}$ and $\varepsilon$ and $\Xi$, there exists $\kappa > \sqrt{2}$ with the following significance*: *Suppose that* $(A, \mathfrak{a})$ *is a solution to (1.4) on* $(0,\infty) \times \mathbb{R}^2 \times \mathbb{R}$ *that obeys* CONSTRAINT SET 2 *with* $|\mathfrak{a}| \leq \frac{\mathfrak{z}}{t}$ *for all positive t and with* $|\mathfrak{a}| > \frac{\varepsilon}{t}$ *when* $t \leq t_0$ *with* $t_0 > 0$. *Finally, assume that (5.1)'s function f obeys* $f(R) < \frac{\Xi}{R}$ *for all* $R > 0$.
- *If* $t \leq t_0$, *then* $|\mathfrak{a}_3| > \frac{1}{\kappa t}$ *on the* $|z| > \kappa t$ *part of* $\{t\} \times \mathbb{R}^2 \times \{0\}$.
- *If* $t \leq t_0$, *then* $|\varphi| > \frac{1}{\kappa t}$ *on the* $|z| > \kappa t$ *part of* $\{t\} \times \mathbb{R}^2 \times \{0\}$.
- *If* $t > t_0$, *then there exists* $\kappa' \geq \kappa$ *(determined by* $(A,\mathfrak{a})$ *and t) such that the subset in* $\{t\} \times \mathbb{R}^2 \times \{0\}$ *where* $|\varphi| < \frac{1}{\kappa' t}$ *is contained in the* $|z| < \kappa t$ *disk*.

By way of an immediate corollary: Supposing that $(A, \mathfrak{a})$ obeys (1.4) and CONSTRAINT SET 2, then Lemma 5.3 implies that the the zero locus of $\varphi$ in any constant $(t, x_3)$ slice of $(0, \infty) \times \mathbb{R}^2 \times \mathbb{R}$ is either empty or it is the origin in $\mathbb{R}^2$. (Keep in mind that the zero set in any such slice is independent of the slice.)

*Proof of Lemma 5.3*: The proof has three parts. Parts 1 and 2 prove the top two bullets of the lemma and Part 3 proves the lemma's third bullet.

*Part 1*: Suppose to the contrary that either the first or the second bullet of the lemma were false so as to derive nonsense. In this event, there would exists a sequence $\{(t_n, z_n)\}_{n \in \mathbb{N}} \subset (0, t_0) \times \mathbb{R}^2$ such that the following conditions are obeyed for any given positive integer n:

- $|\mathfrak{a}_3| < \frac{1}{n t_n}$ at $(t_n, z_n, 0)$ *if the lemma's first bullet is false*.
- $|\varphi| < \frac{1}{n t_n}$ at $(t_n, z_n, 0)$ *if the lemma's second bullet is false*.
- $|z_n| > n t_n$.

(5.3)

Let $\{(A^n, \mathfrak{a}^n)\}_{n \in \mathbb{N}}$ denote the sequence of solutions to (1.4) whose n'th term is the pull-back of the given solution $(A, \mathfrak{a})$ by the coordinate rescaling and $\mathbb{R}^2$-factor translation diffeomorphism $(t, z, x_3) \to (t_n t, t_n z + z_n, x_3)$.

To see the significance of this sequence, fix $n \in \mathbb{N}$ for the moment. Let $\mathfrak{a}^n_3$ and $\varphi^n$ denote the respective $(A^n, \mathfrak{a}^n)$ versions of $\mathfrak{a}_3$ and $\varphi$. If the first bullet's assertion is false, then the norm of $\mathfrak{a}^n_3$ is less than $\frac{1}{n}$ at the point $(1, 0, 0)$ whereas the norm of $\varphi^n$ is greater than $\varepsilon$ at this point. If the second bullet's assertion is false, then the norm of $\varphi^n$ is less than $\frac{1}{n}$ at the point $(1, 0, 0)$ whereas the norm of $\mathfrak{a}^n_3$ at this point is greater than $\varepsilon$. In any event, the integral of $|F_A|^2$ on the domain $(0, \infty) \times D_{n/2} \times \{0\}$ is no greater than $\frac{2 \Xi}{n}$.

*Part 2*: Use the sequence $(A^n, \mathfrak{a}^n)$ as input for Lemma 2.3 and let $(A^\infty, \mathfrak{a}^\infty)$ denote one of its limits. By virtue of the manner of $C^\infty$ convergence that is dictated by Lemma 2.3, the solution $(A^\infty, \mathfrak{a}^\infty)$ has the following properties:

- $|\mathfrak{a}^\infty| \leq \frac{3}{t}$.
- $|\mathfrak{a}| \geq \frac{\varepsilon}{t}$ *where* $t < t_0$.
- *If the lemma's first bullet is false, then* $\mathfrak{a}^\infty_3 = 0$ *and* $|\varphi^\infty| \geq \varepsilon$ *at the point* $(1, 0, 0)$.
- *If the lemma's second bullet is false, then* $\varphi^\infty = 0$ *and* $|\mathfrak{a}^\infty_3| \geq \varepsilon$ *at the point* $(1, 0, 0)$.
- $F_{A^\infty} \equiv 0$.

(5.4)

To obtain nonsense from this when the lemma's first bullet is false, use the top bullet of (3.10) to see that $\frac{\partial}{\partial t}\alpha^\infty \geq |\varphi^\infty|^2$ on the whole of $(0,\infty)\times\mathbb{R}^2\times\mathbb{R}$. (This $\alpha^\infty$ is the $(A^\infty, \mathfrak{a}^\infty)$ version of the function $\alpha$.) If $\alpha^\infty$ is zero at the point $(1, 0, 0)$, then it must be positive for $t > 1$ on the $z = 0$ locus. What with the bullet of (3.9), this implies that $|\varphi^\infty|$ is *increasing* for $t > 1$ on the $z = 0$ locus since $|\varphi^\infty| > \varepsilon$ when $t = 1$ on the $z = 0$ locus. Thus, $|\varphi^\infty| > \varepsilon$ at all times $t > 1$ on the $z = 0$ locus. This last conclusion is the desired nonsense because it runs afoul of the top bullet in (5.4).

Now suppose that the lemma's second bullet is false. Then $\varphi^\infty = 0$ on the whole of the $z = 0$ locus by virtue of the top bullet in (1.4). To obtain nonsense from this fact, use the third bullet of (1.4) and the fifth bullet of (5.4) to see that $\frac{\partial}{\partial t}|\mathfrak{a}^\infty{}_3| = 0$ on the whole of the $z = 0$ locus. This last conclusion implies that $|\mathfrak{a}^\infty{}_3| > \varepsilon$ on the whole of the $z = 0$ locus which is the desired nonsense because it fouls the top bullet in (5.4).

*Part 3*: The third bullet of the lemma follows from the second bullet given that $\nabla_{\hat{A}t}\varphi = 2\alpha\varphi$ and given that $\alpha \geq -\frac{\mathfrak{z}}{t}$. To see why, let $\kappa_*$ denote the version of the number $\kappa$ from the lemma's second bullet. Now fix $t > t_0$. If $z$ has norm greater than $\kappa_* t$, then the equation $\nabla_{At}\varphi = 2\alpha\varphi$ and the lower bound $\alpha > -\frac{\mathfrak{z}}{t}$ imply that

$$|\varphi| \geq \left(\frac{t_0}{t}\right)^{2\mathfrak{z}-1}\frac{1}{\kappa_* t}$$

(5.5)

at the point $(t, z, x_3 = 0)$. Thus, if $\kappa'$ is set equal to $\kappa_*\left(\frac{t}{t_0}\right)^{2\mathfrak{z}-1}$, and if $|\varphi| \leq \frac{1}{\kappa' t}$ at a given point $(t, z, 0) \in (t_0, \infty)\times\mathbb{R}^2\times\{0\}$, then $|z|$ must be less than $\kappa t$.

**c) The $|z| \to \infty$ limit of $|\varphi|$**

The upcoming Lemma 5.4 makes a formal assertion to the effect that any given CONSTRAINT SET 2 solution to (1.4) looks more and more like the Nahm pole solution from Section 1c as $|z| \to \infty$ on $\{t\}\times\mathbb{R}^2\times\mathbb{R}$ for any fixed, positive time $t$. (The Nahm pole solution is the $m = 0$ model solution.) It has a corollary (which is Lemma 5.5) to the effect that $|\varphi|$ limits to $\frac{1}{\sqrt{2}t}$ as $|z| \to \infty$ on $\{t\}\times\mathbb{R}^2\times\mathbb{R}$

**Lemma 5.4**: *Let $(A, \mathfrak{a})$ denote a given CONSTRAINT SET 2 solution to (1.4). Fix a sequence $\{\lambda_n\}_{n\in\mathbb{N}} \subset (0,\infty)$ that is bounded from above; and fix a sequence $\{z_n\}_{n\in\mathbb{N}} \subset \mathbb{R}^2$ with no convergent subsequences. For any given $n \in \mathbb{N}$, let $(A^n, \mathfrak{a}^n)$ denote the solution to (1.4) that is obtained from $(A, \mathfrak{a})$ by pulling back the latter via the diffeomorphism $(t, z, x_3) \to (\lambda_n t, \lambda_n(z + z_n), x_3)$ of $(0,\infty)\times\mathbb{R}^2\times\mathbb{R}$. The sequence $\{(A^n, \mathfrak{a}^n)\}_{n\in\mathbb{N}}$ after possibly* *termwise applications of elements in* Aut(P) *converges in the* $C^\infty$*-topology on compact subsets of* $(0,\infty) \times \mathbb{R}^2 \times \mathbb{R}$ *to the Nahm pole solution from Section 1c.*

The proof is given momentarily. The next lemma is an immediate corollary:

**Lemma 5.5**: *Let* $(A, \mathfrak{a})$ *denote a given* CONSTRAINT SET 2 *solution to (1.4). Given positive numbers* t *and* $\delta$, *there exists* $r > 0$ *such that if* $s \in (0, t]$, *then*

- $\frac{(1-\delta)}{\sqrt{2}s} < |\varphi| < \frac{(1+\delta)}{\sqrt{2}s}$ ,
- $\frac{(1-\delta)}{2s} < |\mathfrak{a}_3| < \frac{(1+\delta)}{2s}$ ,
- $|\beta| < \frac{\delta}{s}$

*where* $|z| > rs$ *on the slice* $\{s\} \times \mathbb{R}^2 \times \mathbb{R}$.

Keep in mind that there is no claim at this point that $r$ is independent of t.

*Proof of Lemma 5.5*: Imagine if Lemma 5.5 were not true: There would exist positive numbers t and $\delta$, a bounded sequence $\{t_n\} \in (0, t]$ and a sequence of points $\{z_n\}_{n \in \mathbb{N}} \subset \mathbb{R}^2$ with the following two properties: The bounds from Lemma 5.5 are violated at the points $\{(t_n, z_n, 0)\}_{n \in \mathbb{N}}$; and the sequence $\{\frac{|z_n|}{t_n}\}_{k \in \mathbb{N}}$ diverges. If the sequences $\{\lambda_n = t_n\}_{n \in \mathbb{N}}$ and $\{z_n = \frac{z_n}{t_n}\}_{n \in \mathbb{N}}$ are used as input for Lemma 5.4, then no subsequence of the resulting $\{(A^n, \mathfrak{a}^n)\}_{n \in \mathbb{N}}$ can converge to the Nahm pole solution in the manner indicated (after any termwise action of elements in Aut(P)) because $|\varphi| \equiv \frac{1}{\sqrt{2}t}$ and $|\mathfrak{a}_3| \equiv \frac{1}{2t}$ and $\beta \equiv 0$ for the Nahm pole solution.

*Proof of Lemma 5.4*: By way of notation, positive numbers $\mathfrak{z}, \varepsilon, t_0$ and $\Xi$ are chosen so that $|\mathfrak{a}| \leq \frac{\mathfrak{z}}{t}$ for all $t > 0$, so that $|\mathfrak{a}| \geq \frac{\varepsilon}{t}$ where $t \leq t_0$ and so that $f(R) < \frac{\Xi}{R}$ for all $R > 0$. The proof also uses $\kappa_*$ to denote Lemma 5.3's version of $\kappa$. Keep in mind the assertions of Lemma 5.3 to the effect that the norms of both $\mathfrak{a}_3$ and $\varphi$ are greater than $\frac{1}{\kappa_* t}$ on the $t < t_0$ and $|z| > \kappa_* t$ part of $(0, \infty) \times \mathbb{R}^2 \times \mathbb{R}$.

Since no generality is lost by assuming that each $n \in \mathbb{N}$ version of $\lambda_n$ is less than 1, this upper bound is henceforth assumed.

To start the proof: The sequence $\{(A^n, \mathfrak{a}^n)\}_{n \in \mathbb{N}}$ has the following properties:

- $|\mathfrak{a}^n| \leq \frac{\mathfrak{z}}{t}$ *for all* $t > 0$.
- $|\mathfrak{a}^n_3| \geq \frac{1}{\kappa_* t}$ *where both* $t \leq \frac{1}{\lambda_n} t_0$ *and* $|z| > |z_n| - \kappa_* t$

- $\int_{(0,\infty)\times D_R\times\{0\}} |F_{A^k}|^2 < \frac{2\,\Xi}{|z_n|}$ *if* $R < \frac{1}{2}|z_n|$.

- $|\varphi^n| > \frac{1}{\kappa_* t}$ *where both* $t < \frac{1}{\lambda_n} t_0$ *and* $|z| > |z_n| - \kappa_* t$. *In general, there exists* $\delta(t) > 0$ *such that* $|\varphi^n| > \delta(t)$ *where* $|z| < |z_n| - \kappa_* t$.

(5.6)

The properties in the first and third bullets follow by rescaling and translating the bounds from CONSTRAINT SET 2. The properties in the second and fourth bullet follow from Lemma 5.3.

It follows from Lemma 2.3 that this sequence has a convergent subsequence (after possibly termwise action of elements from Aut(P)) that converges on compact subsets of $(0,\infty)\times\mathbb{R}^2\times\mathbb{R}$. The limit is denoted by $(A^\infty, \mathfrak{a}^\infty)$. It obeys the following:

- $|\mathfrak{a}^\infty| \leq \frac{3}{t}$ *for all* $t > 0$.
- $|\mathfrak{a}^\infty_3| \geq \frac{1}{\kappa_* t}$ *for all* $t \leq t_0$.
- $F_{A^\infty} \equiv 0$.
- $|\varphi^k| > \frac{1}{\kappa_* t}$ *where* $t < t_0$ *and* $|\varphi| > 0$ *in general*

(5.7)

These follow by taking limits of the corresponding bullets in (5.6). This solution is therefore described by CONSTRAINT SET 1 if $\langle \mathfrak{a}^\infty_3 \varphi^\infty \rangle \equiv 0$ (which is to say that $\beta^\infty \equiv 0$ since $\varphi^\infty$ is non-vanishing). But note that the top bullets of (3.11) and (3.13) imply that

$$\frac{\partial}{\partial t}\langle \varphi^\infty \beta^\infty \rangle \equiv 0,$$

(5.8)

which implies that $\langle \varphi^\infty \beta^\infty \rangle \equiv 0$ since $\langle \varphi^\infty \beta^\infty \rangle \to 0$ as $t \to 0$ (because of the top bullet in (5.7)). Thus, $\beta^\infty \equiv 0$ and the solution is indeed described by CONSTRAINT SET 1. According to Theorem 1, it $(A^\infty, \mathfrak{a}^\infty)$ must be the Nahm pole solution.

### d) Proof of Proposition 5.1

The proof has three parts.

*Part 1*: Fix a positive number to be denoted by $\delta$. For each positive time $t$, there is an associated number $r(t)$ that is defined as follows: The number $r(t)$ is the smallest choice for $r$ such that if $s \in (0, t]$, then $|\varphi| > \frac{1}{\sqrt{2}s}(1-\delta)$ at all points where $|z| > rs$ on the slice $\{s\}\times\mathbb{R}^2\times\mathbb{R}$. Lemma 5.5 guarantees that there is a finite $r(t)$ for each $t$. Because of the definition, there exists $(s,z) \in (0, t]\times\mathbb{R}^2$ with $|z| = r(t)s$ such that $|\varphi|(s,z,0) = \frac{1}{\sqrt{2}s}(1-\delta)$.

With the preceding understood, suppose for the sake of an argument that there is a sequence $\{t_n\}_{n\in\mathbb{N}}$ with the corresponding sequence $\{r(t_n)\}_{n\in\mathbb{N}}$ having no convergent subsequences. Nonsense will be generated shortly from the existence of such a sequence of times. To this end, choose (for each positive integer n) a point $(s_n, z_n) \in (0, t_k) \times \mathbb{R}^2$ with $|z_n| = r(t_n) s_n$ with $|\varphi|(s_n, z_n, 0) = \frac{1}{\sqrt{2} s_n}(1-\delta)$. By virtue of the definition of the function $r(\cdot)$, the sequence $\{s_n\}_{n\in\mathbb{N}}$ can't have convergent subsequences and neither can the sequence $\{z_n = \frac{z_n}{s_n}\}_{k\in\mathbb{N}}$. On the other hand, no generality is lost by assuming that the sequence $\{\frac{z_k}{|z_k|}\}_{n\in\mathbb{N}}$ of unit length vectors in $\mathbb{R}^2$ converges to the vector $(1,0) \in \mathbb{R}^2$. (This is because pull-backs of CONSTRAINT SET 2 solutions of (1.4) by constant angle rotations about the origin of the $\mathbb{R}^2$ factor of $(0,\infty)\times\mathbb{R}^2\times\mathbb{R}$ are also CONSTRAINT SET 2 solutions of (1.4).)

Use $(A^n, \mathfrak{a}^n)$ to denote the pull-back of $(A, \mathfrak{a})$ by the combination of coordinate rescaling and $\mathbb{R}^2$-factor translation of $(0, \infty)\times\mathbb{R}^2\times\mathbb{R}$ that is defined by the rule whereby $(t, z, x_3)$ is sent to $(s_n t, s_n z + z_n, s_n x_3)$. (This diffeomorphism is chosen so that it sends the point $(1, 0, 0)$ to $(s_n, z_n, 0)$.)

Defined in this way, the pair $(A^n, \mathfrak{a}^n)$ has the following properties:

- $|\mathfrak{a}^n| \leq \frac{3}{t}$ *for all* $t > 0$.
- $|\varphi^n| = \frac{(1-\delta)}{\sqrt{2}}$ *at the point* $(1, 0\ 0)$ *in* $(0, \infty)\times\mathbb{R}^2\times\mathbb{R}$.
- $|\varphi^n| \geq \frac{(1-\delta)}{\sqrt{2}t}$ *at any* $(t, z, 0) \in (0, 1]\times\mathbb{R}^2\times\mathbb{R}$ *if* $t \leq 1$ *and* $|z + \frac{z_n}{s_n}| > r(t_n) t$.
- $\int_{(0,\infty)\times D_R\times\{0\}} |F_{A^n}|^2 \leq \frac{2\Xi}{r(t_n)}$ *if* $R < \frac{1}{2} r(t_n)$.

(5.9)

The top bullet condition holds because $(A, \mathfrak{a})$ obeys the first bullet condition in CONSTRAINT SET 2. The fourth bullet condition in (5.9) results from writing the coordinate rescaling and $\mathbb{R}^2$ factor translation of the third bullet condition in CONSTRAINT SET 2 using the fact that $r(t_n)$ is $\frac{|z_n|}{s_n}$. The second and third bullet conditions in (5.16) hold by virtue of the definition of $(A^n, \mathfrak{a}^n)$.

*Part 2*: According to Lemma 2.3, the sequence $\{(A^n, \mathfrak{a}^n)\}_{n\in\mathbb{N}}$ has a subsequence that converges in the $C^\infty$ topology on compact subsets of $(0, \infty)\times\mathbb{R}^2\times\mathbb{R}$ (after possible termwise actions of elements in Aut(P)) to a solution to (1.4). The latter is denoted by $(A^\infty, \mathfrak{a}^\infty)$. By virtue of the manner of convergence and by virtue of the conditions in (5.9), this solution has the following properties:

- $|\mathfrak{a}^\infty| \leq \frac{3}{t}$ *for all* $t > 0$.

- $|\varphi^\infty| = \frac{(1-\delta)}{\sqrt{2t}}$ *at any point of the form* $(t, 0\ x_3)$ *in* $(0, 1) \times \mathbb{R}^2 \times \mathbb{R}$.
- $|\varphi^\infty| \geq \frac{(1-\delta)}{\sqrt{2t}}$ *at any* $(t, z = x_1 + ix_2, x_3) \in (0, 1] \times \mathbb{R}^2 \times \mathbb{R}$ *with* $t \leq 1$ *and* $x_1 < 0$.
- $F_{A^\infty} = 0$.

(5.10)

It is also the case that $|\varphi^\infty| > 0$ everywhere. This is because the zeros of any $n \in \mathbb{N}$ version of $\varphi^n$ on the $(t = 1, x_3 = 0)$ slice of $(0, \infty) \times \mathbb{R}^2 \times \mathbb{R}$ are at the point $-\frac{\zeta_n}{s_n}$ in $\mathbb{R}^2$ which has distance $r(t_n)$ from the origin.

As explained in Part 2, the pair $(A^\infty, \mathfrak{a}^\infty)$ must be a Nahm pole solution (this is suggested by the fact that $\varphi^\infty$ is nowhere zero and $A^\infty$ is flat). The conclusion that $(A^\infty, \mathfrak{a}^\infty)$ is a Nahm pole solution constitutes the nonsense that proves Proposition 5.1 because the norm of $\varphi$ for the Nahm pole solution slice is precisely $\frac{1}{\sqrt{2t}}$ which fouls the second bullet in (5.17).

*Part 3*: The proof that $(A^\infty, \mathfrak{a}^\infty)$ is a Nahm pole solution will invoke Theorem 1. Since $|\varphi^\infty| > 0$, it is necessary that $\langle \mathfrak{a}^\infty_3 \varphi^\infty \rangle \equiv 0$ (equivalently, that $\beta^\infty \equiv 0$) and that the function $t|\mathfrak{a}_3^\infty|$ is bounded from below by a positive number where t is small. The fact that $\beta^\infty \equiv 0$ follows because (5.8) also holds in this case and because $\langle \varphi^\infty \beta^\infty \rangle \to 0$ as $t \to 0$ (due to the top bullet in (5.10)). The proof that $t|\mathfrak{a}_3^\infty|$ is bounded from below by some positive number when t is given in the subsequent paragraphs.

Since $B_{A_3^\infty} \equiv 0$, the equation $(\nabla_{A_1^\infty} + i \nabla_{A_2^\infty}) \varphi^\infty$ from the second bullet of (1.4) implies (by differentiating it) the identity

$$\tfrac{1}{2} \left( \tfrac{\partial^2}{\partial x_1^2} + \tfrac{\partial^2}{\partial x_1^2} \right) |\varphi^\infty|^2 = |\nabla_{A_1^\infty} \varphi^\infty|^2 + |\nabla_{A_2^\infty} \varphi^\infty|^2 .$$

(5.11)

Fix $R > 1$ and let $\chi_R$ denote the function on $\mathbb{R}^2$ given by the rule $z \to \chi(\tfrac{\ln|z|}{\ln R} - 1)$ where $\chi$ is the standard bump function on $\mathbb{R}^2$. This function $\chi_R$ is equal to 1 where $|z| < R$ and equal to 0 where $|z| > R^2$. Multiply both sides of (5.11) by $\chi_R$, integrate the result over the any $t > 0$ version of $\{t\} \times \mathbb{R}^2 \times \{0\}$, then integrate twice by parts on the left hand side and use the fact that $|\varphi^\infty| \leq \mathfrak{z}$ on this slice to see that

$$\int_{\{t\} \times \mathbb{R}^2 \times \{0\}} \chi_R (|\nabla_{A_1^\infty} \varphi^\infty|^2 + |\nabla_{A_2^\infty} \varphi^\infty|^2) \leq c_0 \mathfrak{z}^2 \tfrac{1}{\ln R} .$$

(5.12)

Take the $R \to \infty$ limit of this last inequality to conclude that $\nabla_A^\perp \varphi^\infty \equiv 0$. This implies in particular that $|\varphi^\infty|$ is a function only of the coordinate t.

Since $|\varphi^\infty|$ is a function of only t, the third bullet of (5.10) implies that $|\varphi^\infty| \geq \frac{(1-\delta)}{\sqrt{2}t}$ on the whole of $\{t\} \times \mathbb{R}^2 \times \mathbb{R}$ if $t \leq 1$. Meanwhile, the third bullet of (1.4) says in this case that $\frac{\partial}{\partial t}\alpha^\infty = |\varphi^\infty|^2$ since $\beta^\infty \equiv 0$ and $B_{A_3^\infty} \equiv 0$. Thus

$$\frac{\partial}{\partial t}\alpha^\infty \geq \frac{(1-\delta)^2}{2t^2}$$

(5.13)

where $t \leq 1$. Integrate the latter inequality from any given $t < 1$ to 1 (and use the fact that the norm of $\alpha^\infty$ at $t = 1$ is bounded by $\mathfrak{z}$ because of the top bullet in (5.10)) to see that

$$\alpha^\infty|_t \leq -\frac{(1-\delta)^2}{2t} + \mathfrak{z}.$$

(5.14)

This inequality implies in particular that $|\alpha^\infty{}_3| > \frac{1}{4t}$ where $t < \frac{(1-\delta)^2}{4\mathfrak{z}}$.

### e) The vanishing degree of $\varphi$

The lemma that follows in this section asserts an upper bound for the degree of vanishing of $\varphi$ at the origin in any constant $(t, x_3)$ slice of $(0, \infty) \times \mathbb{R}^2 \times \mathbb{R}$ given certain assumptions that hold for CONSTRAINT SET 2 solutions by virtue of Proposition 5.1.

**Lemma 5.6**: *Given positive numbers $\mathfrak{z}, \delta, r$, there exist a positive integer to be denoted by $\hat{m}$ with the following significance: Suppose that $(A, \mathfrak{a})$ is a solution to (1.4) such that $|\mathfrak{a}| \leq \frac{\mathfrak{z}}{t}$ on the whole of $(0, \infty) \times \mathbb{R}^2 \times \mathbb{R}$, and such that $|\varphi| \geq \frac{\delta}{t}$ on the $|z| \geq rt$ part $(0, \infty) \times \mathbb{R}^2 \times \mathbb{R}$. Given these properties, then the degree of vanishing of $\varphi$ at the $z = 0$ point on any constant $(t, x_3)$ slice of $(0, \infty) \times \mathbb{R}^2 \times \mathbb{R}$ is at most $\hat{m}$.*

*Proof of Lemma 5.6*: The lemma is proved by assuming it to be false so as to generate nonsense. In particular, if the lemma is false for a given data set $(\mathfrak{z}, \delta, r)$, then there is a sequence $\{(A^n, \mathfrak{a}^n)\}_{n \in \mathbb{N}}$ of solutions to (1.4) that obey the assumptions of the lemma and with the following additional property: The $(A^n, \mathfrak{a}^n)$ version of $\varphi$ (which is denoted henceforth by $\varphi^n$) vanishes with degree greater than n at the $z = 0$ point in each constant $(t, x_3)$ slice of $(0, \infty) \times \mathbb{R}^2 \times \mathbb{R}$. A point to keep in mind is this: If the vanishing degree of $\varphi^n$ at the $z = 0$ point in each constant $(t, x_3)$ slice of $(0, \infty) \times \mathbb{R}^2 \times \mathbb{R}$ is $k$, then the $\nabla_A$-covariant derivatives of $\varphi^n$ to order at least $k-1$ must also vanish at these $z = 0$ points.

Lemma 2.3 can be invoked for this sequence (because $|\mathfrak{a}^n| \leq \frac{\mathfrak{z}}{t}$ for each n) to extract a subsequence $\Lambda \subset \mathbb{N}$ and a limit $(A^\infty, \mathfrak{a}^\infty)$ which obeys (1.4). It also obeys $|\mathfrak{a}^\infty| \leq \frac{\mathfrak{z}}{t}$ on all of $(0, \infty) \times \mathbb{R}^2 \times \mathbb{R}$ and it obeys $|\varphi^\infty| \geq \frac{\delta}{t}$ on the $|z| \geq rt$ part $(0, \infty) \times \mathbb{R}^2 \times \mathbb{R}$. In particular, $\varphi^\infty$ is not identically zero. In addition, the $C^\infty$ convergence of $\{(A^n, \mathfrak{a}^n)\}_{n \in \mathbb{N}}$ after termwise action of Aut(P) (which is dictated by Lemma 2.3) implies that the covariant derivatives of $\varphi^\infty$ to all orders must vanish at the $z = 0$ point in each constant $(t, x_3)$ slice of $(0, \infty) \times \mathbb{R}^2 \times \mathbb{R}$. But then $\varphi^\infty$ must be identically zero which is the desired nonsense.

## 6. Norms of $\beta$ and $\mathfrak{b}$

Supposing that $(A, \mathfrak{a})$ is a given solution to (1.4), this section is concerned with the pointwise norms of the corresponding versions of $\beta$ and $\mathfrak{b}$ and integrals of the squares of these norms on constant $(t, x_3)$ slices. Each subsection concerns some sort of bound on $\beta$ and $\mathfrak{b}$ (and in some cases $\hat{A}$ and $\alpha$ and $\varphi$) given certain assumptions about $(A, \mathfrak{a})$. One or more lemmas in each subsection (and in two cases a proposition) summarize what is done there. Each result in this section requires only that $(A, \mathfrak{a})$ obey a specified subset of the following conditions:

- $|\mathfrak{a}| \leq \frac{3}{t}$ on $(0, \infty) \times \mathbb{R}^2 \times \mathbb{R}$.
- *There are positive numbers $r$ and $\delta$ such that $|\varphi| \geq \frac{\delta}{t}$ where $|z| > rt$ on $(0, \infty) \times \mathbb{R}^2 \times \mathbb{R}$.*
- *The function $|F_A|^2$ is integrable on any $R > 0$ version of $(0, \infty) \times (\mathbb{R}^2 - D_R) \times \{0\}$, and its integral there is bounded by $\frac{\Xi}{R}$ with $\Xi$ being independent of $R$.*

(6.1)

CONSTRAINT SET 2 solutions to (1.4) are a priori described by the first and third bullets of (6.1); and Proposition 5.1 guarantees that they are also described by the second bullet with $\delta$ being any positive number less than $\frac{1}{\sqrt{2}}$ and with $r$ sufficiently large given $\delta$.

### a) Pointwise bounds for the norms of $\beta$ and $\mathfrak{b}$

The following lemma is the center-piece of this subsection.

**Lemma 6.1**: *Given positive numbers $\mathfrak{z}, \delta, r$, there exists a sequence $\kappa > 1$ and a sequence $\{\kappa_q\}_{q \in \mathbb{N}}$ with the following significance: Suppose that $(A, \mathfrak{a})$ is a solution to (1.4) that is described by the top two bullets in (6.1) with the given $(\mathfrak{z}, \delta, r)$. Then the following hold for any $t > 0$:*
- $|\beta| < \frac{\kappa}{t}$ *on the whole of $(0, \infty) \times \mathbb{R}^2 \times \mathbb{R}$.*
- *If $q \in \mathbb{N}$, then $|\mathfrak{b}| < \frac{\kappa_q}{t}$ on the $|z| > \frac{1}{q} t$ part of $(0, \infty) \times \mathbb{R}^2 \times \mathbb{R}$.*
- *The set $\{\kappa_q\}_{q \in \mathbb{N}}$ from the previous bullet is bounded from above (which is to say that the function $t|\mathfrak{b}|$ is bounded) when the following holds: Every coordinate rescaling* *limit version of φ from Section 2d has the same degree of vanishing at the z = 0 point in constant $(t, x_3)$ slices as the original $(A, \mathfrak{a})$ version of φ.*

***Proof of Lemma 6.1***: The bound for β follows from the top bullet in (6.1) because the norm of |β| is at most that of |$\mathfrak{a}_3$|. The proof for the bound for |$\mathfrak{b}$| has six parts. Parts 1-3 prove the assertions in the second bullet and Parts 4-6 prove that of the third bullet. With regards to notation used here and subsequently in this paper: What is denoted in what follows by $c_{\mathfrak{z}}$ denotes a number that depends only on the number $\mathfrak{z}$ that appears in the top bullet of (6.1). It is greater than 1 and it can be assumed to increase between appearances.

*Part 1*: On the one hand, $|\langle \varphi^* \mathfrak{b} \rangle| \leq |\nabla_A \varphi|$, and on the other, $|\langle \varphi^* \mathfrak{b} \rangle| = |\varphi||\mathfrak{b}|$ because $\mathfrak{b}$ is a section of $\mathcal{L}^+$. Thus,

$$|\mathfrak{b}| \leq \tfrac{1}{|\varphi|} |\nabla_A \varphi| .$$

(6.2)

Since $|\varphi| > \tfrac{\delta}{t}$ where $|z| > rt$ by virtue of the second bullet of (6.1) and since $|\nabla_A \varphi| < c_{\mathfrak{z}} \tfrac{1}{t^2}$ (see Lemma 2.2), this last bound implies that $|\mathfrak{b}| \leq \tfrac{c_{\mathfrak{z}}}{\delta} \tfrac{1}{t}$ where $|z| \geq rt$ on $(0, \infty) \times \mathbb{R}^2 \times \mathbb{R}$.

*Part 2*: Lemma 5.6 asserts that there is a positive integer (to be denoted by $\hat{m}$) that depends only on $\mathfrak{z}, \delta$ and $r$ from the top two bullets of (6.1), and with the following significance: If $(A, \mathfrak{a})$ is described by these same top two bullets for the given $\mathfrak{z}, \delta$ and $r$, then the degree of vanishing of φ at the z = 0 point on any constant $(t, x_3)$ slice of $(0, \infty) \times \mathbb{R}^2 \times \mathbb{R}$ is at most $\hat{m}$. With this understood, then Taylor's theorem with remainder implies this: There exists a function ε: $(0, \infty) \to (0, 1)$ such that the norm of φ obeys $|\varphi| \geq \varepsilon(t) \tfrac{|z|^{\hat{m}}}{t^{\hat{m}+1}}$ on the $|z| \leq rt$ part of $(0, \infty) \times \mathbb{R}^2 \times \mathbb{R}$. As explained momentarily, this function ε(·) has a strictly positive lower bound on any constant, positive μ version of the domain in $(0, \infty) \times \mathbb{R}^2 \times \mathbb{R}$ where $|z| \geq \mu t$. Supposing that to be true, then the desired bound for |$\mathfrak{b}$| follows again from (6.2) since $|\nabla_A \varphi| < c_{\mathfrak{z}} \tfrac{1}{t^2}$.

To prove that there is a positive, t-independent lower bound for ε(t) where $|z| \geq \mu t$, assume it to be false so as to derive nonsense. Supposing that there is no such positive lower bound, then there is a sequence $\{t_n, z_n\}_{n \in \mathbb{N}}$ in $(0, \infty) \times \mathbb{R}^2$ with the following properties: First, $\mu t_n \leq |z_n| \leq r t_n$. Second,

$$|\varphi|(t_n, z_n, x_3) \leq \tfrac{1}{n} \tfrac{|z|^{\hat{m}}}{t^{\hat{m}+1}} .$$

(6.3)

Part 3 explains how this sequence runs afoul of what is said in Lemmas 2.3 and 5.6.

*Part 3*: Supposing that n is a positive integer, let $(A^n, \mathfrak{a}^n)$ denote the pull-back of $(A, \mathfrak{a})$ via the coordinate rescaling map $(t, z, x_3) \to (t_n t, t_n z, t_n x_3)$. By virtue of this definition and by virtue of (6.3), the norm of the corresponding $\varphi^n$ obeys

$$|\varphi^n|(1, q_n, 0) \leq \tfrac{1}{n} |q_n|^{\hat{m}}$$

(6.4)

where $q_n$ denotes the point $\tfrac{z_n}{t_n}$. Keep in mind that $q_n$ obeys $\mu \leq |q_n| \leq r$.

Now, let $(A^\infty, \mathfrak{a}^\infty)$ denote one of Lemma 2.3's limits for the sequence $\{(A^n, \mathfrak{a}^n)\}_{n \in \mathbb{N}}$. The norm of the corresponding $\varphi^\infty$ obeys $|\varphi^\infty| \geq \varsigma |z|^{\hat{m}}$ where $|z| < r$ on the slice $\{1\} \times \mathbb{R}^2 \times \mathbb{R}$ with $\varsigma$ being positive. This holds for suitable $\varsigma$ because of Lemma 5.6. The same $\hat{m}$ appears because $\varphi^\infty$ is also described by the top two bullets of (6.1) with the same data $\mathfrak{z}, \delta$ and $r$. Because $\mu \leq |q|_n \leq r$, the bound for $|\varphi^\infty|$ is not compatible with (6.4) when n is large given the manner of convergence that is dictated by Lemma 2.3 (which is $C^\infty$ convergence on compact sets). This incompatibility is the desired nonsense.

*Part 4*: This part and Part 5 of the proof establishes the third bullet of Lemma 6.1. This task requires a depiction of the t-derivative of $\langle \varphi^* \mathfrak{b} \rangle$ which is this one:

$$\tfrac{\partial}{\partial t} \langle \varphi^* \mathfrak{b} \rangle = \langle \varphi^* (E_{A1}^+ + i E_{A2}^+) \rangle .$$

(6.5)

To explain: The equation in (6.5) is derived by first using the second and third bullets in (3.13) to obtain the equation

$$\nabla_{\hat{A}t} \mathfrak{b} + 2\alpha \mathfrak{b} = E_{A1}^+ + i E_{A2}^+ .$$

(6.6)

Then take the inner product of both sides with $\varphi^*$ and use (3.9)'s first bullet.

Integrate both sides of (6.5) on the domain $(t, \infty)$ and use the fundamental theorem of calculus (with the fact that $|\langle \varphi^* \mathfrak{b} \rangle| \leq c_{\mathfrak{z}} \tfrac{1}{t^2}$ so $|\langle \varphi^* \mathfrak{b} \rangle|$ limits to 0 as $t \to \infty$) to see that

$$|\mathfrak{b}|(t, z, \cdot) \leq c_{\mathfrak{z}} \tfrac{1}{|\varphi|_t} \int_t^\infty \tfrac{1}{s^2} |\varphi|_s \, ds$$

(6.7)

where $|\varphi|_s$ is shorthand for the value of $|\varphi|$ at $(s, z, 0)$ and where $|\varphi|_t$ has the analogous definition. (Note that the bound $|E_A| \leq c_{\mathfrak{z}} \tfrac{1}{t^2}$ is also used to obtain (6.7).)

*Part 5*: Let *m* denote the degree of vanishing of $\varphi$ on any given constant $(t, x_3)$ slice of $(0, \infty) \times \mathbb{R}^2 \times \mathbb{R}$. Because $|\varphi|$ obeys the bound $|\varphi|(s, z, x_3) \leq c_\delta \frac{|z|^m}{s^{m+1}}$ (because of Lemma 2.2's bounds), the inequality in (6.7) leads to the following one:

$$|\flat|(t, z, x_3) \leq c_\delta \frac{1}{m+2} \frac{1}{|\varphi|_t} \frac{|z|^m}{t^{m+2}} .$$

(6.8)

This last bound leads in turn to the desired bound $|\flat| \leq \kappa_* \frac{1}{t}$ with $\kappa_*$ being independent of the variables t and z (with $|z| \leq rt$), if there exists a positive $\varepsilon_*$ such that $|\varphi|$ obeys the bound $|\varphi|(t, z, x_3) \geq \varepsilon_* \frac{|z|^m}{t^{m+1}}$ at all points $(t, z)$ with $|z| \leq rt$.

To see about this $\varepsilon_*$, let $\varsigma(t)$ denote the value of $|\nabla_A^{\otimes m}\varphi|$ at the point $(t, z=0, x_3=0)$. If there is a positive lower bound for the function $\varsigma(\cdot)$ on $(0, \infty)$, then there is a positive $\varepsilon_*$ with the desired properties. That this is so follows using Taylor's theorem with remainder (and using Lemma 2.2 to bound order $(m+1)$, $\nabla_A$-covariant derivatives of $\varphi$).

*Part 6*: Suppose for the sake of argument that there is either no postive lower bound for $\mu(\cdot)$ for a given $(A, \mathfrak{a})$ obeying the third bullet's assumptions. This assumption is used in what follows to derive nonsense. In this event, there is a sequence $\{t_n\}_{n \in \mathbb{N}}$ such that $|\nabla_A^{\otimes m}\varphi|(t_n, 0, 0) \leq \frac{1}{n}$. Use the sequence $\{t_n\}_{n \in \mathbb{N}}$ to define a corresponding sequence $\{(A^n, \mathfrak{a}^n)\}_{n \in \mathbb{N}}$ whose n'th term is the pull-back of $(A, \mathfrak{a})$ by the coordinate rescaling map $(t, z, x_3) \to (t_n t, t_n z, t_n x_3)$. By construction, the $(A^n, \mathfrak{a}^n)$ version of $\varphi$ (denoted by $\varphi^n$) obeys

$$|\nabla_A^{\otimes m}\varphi^n|(1, 0, 0) \leq \frac{1}{n} .$$

(6.9)

The sequence $\{(A^n, \mathfrak{a}^n)\}_{n \in \mathbb{N}}$ can now be used as input to Lemma 2.3. Let $(A^\infty, \mathfrak{a}^\infty)$ denote a limit as described by that lemma. Because of (6.9) and because of the $C^\infty$ convergence dictated by Lemma 2.3, the corresponding version of $\varphi$ must have vanishing degree greater than *m* at the $z = 0$ point on the $(t = 1, x_3 = 0)$ slice of $(0, \infty) \times \mathbb{R}^2 \times \mathbb{R}$.

The preceding conclusion is the desired nonsense because it contradicts the explicit assumptions of Lemma 6.1's third bullet.

**b) Bounds for covariant derivatives of $\beta$ and $\flat$**

The next lemma is can be used in conjunction with Lemma 6.1 to obtain bounds on derivatives of $\alpha$ and $\nabla_{\hat{A}}$-covariant derivatives of $\varphi, \beta$ and $\flat$, and bounds for $\langle \sigma F_{\hat{A}} \rangle$ and its derivatives. By way of an explanation, bounds for $\nabla_A$-covariant derivatives of $\mathfrak{a}$ and $F_A$ are given in Lemma 2.2, but these don't a priori lead to bounds for $\nabla_{\hat{A}}$-derivatives.

**Lemma 6.2**: *Given $\mathfrak{z} > 0$, there exists $\kappa > 1$ with the following significance: Suppose that $(A, \mathfrak{a})$ is a solution to (1.4) such that $|\mathfrak{a}| \leq \frac{\mathfrak{z}}{t}$ on $(0, \infty) \times \mathbb{R}^2 \times \mathbb{R}$ (the top bullet of (6.1) holds). Then both $|\nabla_{\hat{A}t}\beta|$ and $|\frac{\partial}{\partial t}\alpha|$ are bounded by $\frac{\kappa}{t^2}$. The following is also true: Fix $t \in (0, 1)$ and $\gamma \in (0, \frac{1}{2})$ and let B denote a ball in $(0, \infty) \times \mathbb{R}^2 \times \{0\}$ of radius $\gamma t$ centered at a point in $\{t\} \times \mathbb{R}^2 \times \{0\}$. Fix $\mathcal{Z} \geq 0$ so that $|\mathfrak{b}| \leq \frac{\mathcal{Z}}{t}$ on B. For any given nonnegative integer n, there exists $\kappa_n \geq 1$ that depends only on $\gamma, \mathfrak{z}$ and $\mathcal{Z}$ such that the bounds below hold on the concentric ball with radius $\frac{1}{2}\gamma t$.*

- $|\nabla^{\otimes n} F_{\hat{A}}| \leq \kappa_n \frac{1}{t^{n+2}}$,
- $|\nabla^{\otimes n} \alpha| \leq \kappa_n \frac{1}{t^{n+1}}$,
- $|(\nabla_{\hat{A}})^{\otimes n} \varphi| \leq \kappa_n \frac{1}{t^{n+1}}$
- $|(\nabla_{\hat{A}})^{\otimes n} \mathfrak{b}| + |(\nabla_{\hat{A}})^{\otimes n} \beta| \leq \kappa_n \frac{1}{t^{n+2}} \left( \int_{B \times \{0\}} (|\beta|^2 + |\mathfrak{b}|^2) \right)^{1/2}$ *(which is at most $4\kappa_n (\mathcal{Z} + \mathfrak{z}) \frac{1}{t^{n+1}}$).*

*Proof of Lemma 6.2*: The assertions about the t-derivative of $\alpha$ and $\nabla_{\hat{A}t}\beta$ follow from Lemmas 2.2 and the bound $\mathfrak{a}| \leq \frac{\mathfrak{z}}{t}$ because (3.7) and the top bullet of (3.8) imply that

$$\langle \sigma \nabla_{At} \mathfrak{a}_3 \rangle = \frac{\partial}{\partial t} \alpha - 4|\beta|^2 .$$

(6.10)

The remaining bounds can be obtained from the initial pointwise bounds on $|\varphi|$, $|\alpha|, |\beta|$ (via the top bullet of (6.1)) and $|\mathfrak{b}|$ (via the assumption), and from the bounds from Lemma 2.2 by applying elliptic bootstrapping arguments to the equations in (3.9)–(3.13). The paragraphs that follow describe how this works.

The identity in the top bullet of (3.13) when written using the connection A instead of $\hat{A}$ has the schematic form $-2i \overline{\partial}_A \mathfrak{b} = \frac{1}{2}(B_{A3} + [\frac{1}{2}\sigma, B_{A3}]) - 4|\mathfrak{b}|^2 \sigma$ with $\overline{\partial}_A$ denoting the covariant $\overline{\partial}$-operator that is defined by A. Since A is smooth and $B_{A3}$ is smooth with a priori bounds from Lemma 2.2, relatively standard elliptic boot-strapping arguments produce bounds for the norms of the $\nabla_A^\perp$-covariant derivatives of $\mathfrak{b}$ to any given order in any disk in any constant $(t, x_3)$ slice of $(0, \infty) \times \mathbb{R}^2 \times \mathbb{R}$ where $|\mathfrak{b}|$ is a priori bounded on a concentric, larger radius disk. (Keep in mind when doing the bootstrapping that the $\nabla_A^\perp$ covariant derivatives of $\sigma$ can be written as linear combinations of $\mathfrak{b}$ and $\mathfrak{b}^*$.) Section 6d has an example (Lemma 6.7) of what can be said using this $\overline{\partial}_A$-depiction of the top bullet of (3.11) when there isn't an a priori point-wise initial bound for $|\mathfrak{b}|$.

Supposing that $\mathfrak{b}$ and its $\nabla_A^\perp$-covariant derivatives have a priori bounds on a disk in a constant $(t, x_3)$ slice $(0, \infty) \times \mathbb{R}^2 \times \mathbb{R}$, then the lower two bullets of (3.13) can be used to obtain analogous bounds for $\beta$ and its $\nabla_A^\perp$-covariant derivatives on any smaller radius concentric disk. This is because these two bullets when written using A instead of Â say that $2i\partial_A\beta = \frac{1}{2}((E_{A1}-iE_{A2}) + \frac{1}{2}[\sigma, (E_{A1}-iE_{A2})]) - 4\langle\mathflat^*\beta\rangle\sigma$. Note also: Since $\nabla_A^\perp$ and $\nabla_{\hat{A}}^\perp$ differ by a zero order endomorphism that is linear in $\mathflat$, the Lemma 2.2 bounds for the $\nabla_A^\perp$-covariant derivatives of $\varphi$ lead immediately to a prior bounds for the corresponding $\nabla_{\hat{A}}^\perp$-covariant derivatives. They also lead (via the top bullet in (3.12)) for bounds on the derivative to any given order of $\langle\sigma B_{\hat{A}3}\rangle$.

The $\nabla_{\hat{A}t}$ derivatives of $\beta$ and $\mathflat$ (and those of $\alpha$) can be bounded now using the equations in (3.11) and the top bullet in (3.10) with the initial bound for $\partial_t\alpha$ from the first part of Lemma 6.2. Meanwhile, a priori bounds for $\nabla^\perp$ derivative of $\alpha$ can be obtained from the lower bullet in (3.12) and the lower bullet in (3.10). Mixed $\nabla_{\hat{A}t}$ and $\nabla_{\hat{A}}^\perp$ covariant derivatives can be likewise be obtained from the equations in (3.9)–(3.13).

### c) Constant $(t, x_3)$ integrals of $|\beta|^2$ and $|\mathflat|^2$

The following lemma is the center-piece of this subsection.

**Lemma 6.3**: *Given positive numbers $\mathfrak{z}$, $\delta$, $r$ and $\Xi$, there exists $\kappa > 1$ with the following significance: Suppose that $(A, \mathfrak{a})$ is a solution to (1.4) that is described by the three bullets in (6.1) with the given values of $\mathfrak{z}$, $\delta$, $r$ and $\Xi$. Then following hold for any $t > 0$:*

- $|\beta|^2 + |\mathflat|^2$ *is integrable on* $\{t\}\times\mathbb{R}^2\times\{0\}$ *and the resulting integral is at most $\kappa$;*
- *If $R \geq t$, then the integral of $|\beta|^2 + |\mathflat|^2$ on $\{t\}\times(\mathbb{R}^2-D_R)\times\{0\}$ is at most $\kappa\frac{t}{R}$.*
- *If $R \geq t$, then $|\beta| + |\mathflat| \leq \kappa(\frac{t}{R})^{1/4}\frac{1}{t}$*

*Proof of Lemma 6.3*: The proof has three parts. Part 1 proves the second bullet of the lemma (it is needed to prove the first bullet); Part 2 proves the first bullet of the lemma; and Part 3 proves the last bullet.

*Part 1*: The proof of the second bullet invokes (6.5) and an analogous equation with $\beta$ replacing $\mathflat$:

$$\tfrac{\partial}{\partial t}\langle\varphi^*\beta\rangle = \langle\varphi^*B_{A3}^+\rangle.$$

(6.11)

This $\beta$ analog of (6.5) comes from the top bullets in (3.11) and (3.13) by taking the inner product of both sides with $\varphi^*$.

Supposing that $S > R$, let $D_{R,S}$ denote the annulus in $\mathbb{R}^2$ where $R \leq |z| \leq S$. Since $|\langle\varphi^*\beta\rangle| = |\varphi||\beta|$ and $|\langle\phi^*\mathflat\rangle| = |\phi||\mathflat|$, the identities in (6.5) and (6.11) lead directly to this:

$$\frac{d}{dt}\Big(\int_{\{t\}\times D_{R,S}\times\{0\}} |\varphi|^2(|\beta|^2+|\mathfrak{b}|^2)\Big)^{1/2} \le \Big(\int_{\{t\}\times D_{R,S}\times\{0\}} |\varphi|^2 |F_A|^2\Big)^{1/2}$$

(6.12)

The right hand side of this last inequality is integrable on $(t,\infty)$ since its integrand is at most $c_{\mathfrak{z}} \frac{S}{t^3}$. More to the point, by virtue of the first and third bullets in (6.1), the $(t,\infty)$ integral of the right hand side of (6.1) is at most $\mathfrak{z} \Xi \frac{1}{\sqrt{tR}}$ because $D_{R,S}$ is contained in the $|z|\ge R$ part of $\mathbb{R}^2$. (Note in particular that this bound is independent of S). Meanwhile, the integral of $|\varphi|^2(|\beta|^2+|\mathfrak{b}|^2)$ over $D_{R,S}$ is at most $c_{\mathfrak{z}} \frac{S^2}{t^4}$ because of the a priori bounds from the top bullet of (3.1) (to bound $|\varphi||\beta|$), and because $|\varphi||\mathfrak{b}|\le |\nabla_A \varphi|$ which Lemma 2.2 bounds by $c_{\mathfrak{z}} \frac{1}{t^2}$. The important take-away is that the integral on the left hand side has limit 0 as $t\to\infty$. Granted this, then (6.12) when integrated leads to the bound

$$\int_{\{t\}\times D_{R,S}\times\{0\}} |\varphi|^2(|\beta|^2+|\mathfrak{b}|^2) \le \mathfrak{z}^2 \Xi^2 \frac{1}{tR}.$$

(6.13)

This last inequality then leads to the following one if $R > rt$:

$$\int_{\{t\}\times D_{R,S}\times\{0\}} (|\beta|^2+|\mathfrak{b}|^2) \le \frac{1}{\delta^2}\mathfrak{z}^2 \Xi^2 \frac{t}{R}.$$

(6.14)

because $|\varphi| \ge \frac{\delta}{t}$ where $|z|\ge rt$ on $\{t\}\times\mathbb{R}^2\times\mathbb{R}$ (see the second bullet in (6.1)).

Since the right hand side of (6.14) is independent of S, the $S\to\infty$ limit of the left hand side of (6.14) can be taken to obtain the desired a priori bound for the $\mathbb{R}^2 - D_R$ integral of $|\beta|^2+|\mathfrak{b}|^2$ when $R\ge rt$. The latter bound leads to the bound in the second bullet of Lemma 6.3 because the area of the annulus in $\{t\}\times\mathbb{R}^2\times\{0\}$ where $t\le |z|\le rt$ is bounded by $\pi r^2 t^2$ and because both $|\beta|^2$ and $|\mathfrak{b}|^2$ are bounded on this annulus by a t-independent multiple of $\frac{1}{t^2}$. (The bound for $\beta$ follows from the top bullet in (6.1) and the bound for $\mathfrak{b}$ follows from Lemma 6.1.)

*Part 2*: The assertion in the first bullet follows from the assertion in the second bullet if there is a t-independent bound for the integrals of $|\beta|^2$ and $|\mathfrak{b}|^2$ on the $|z|\le t$ disk in $\{t\}\times\mathbb{R}^2\times\{0\}$. Such a bound holds for the $|\beta|^2$ integral because the area of the disk is $\pi t^2$ and because $|\beta|^2$ is bounded $\frac{\mathfrak{z}^2}{t^2}$. The argument for the $|\mathfrak{b}|^2$ integral is far more devious since $|\mathfrak{b}|^2$ has no analogous pointwise upper bound. This argument has four steps.

Step 1: Let $\chi_t$ denote the function on $\mathbb{R}^2$ given by the rule $z \to \chi_*(z) = \chi(2r - \frac{|z|}{t})$. This function is equal to 1 where $|z| \geq 2rt$ and it is equal to zero where $|z| \leq rt$. Define a new connection on P (to be denoted by $\mathcal{A}$) as follows:

$$\mathcal{A} = \hat{A} - \tfrac{i}{4} \chi_* \tfrac{1}{|\varphi|^2} (\langle \varphi^* \nabla_{\hat{A}} \varphi \rangle - \langle \varphi \nabla_{\hat{A}} \varphi^* \rangle) \sigma \ .$$

(6.15)

The curvature 2-form for $\mathcal{A}$ is

$$F_\mathcal{A} = (1 - \chi_*) F_{\hat{A}} - \tfrac{i}{4} (\tfrac{\partial}{\partial t} \chi_* dt + \tfrac{\partial}{\partial x_1} \chi_* dx_1 + \tfrac{\partial}{\partial x_2} \chi_* dx_2) \wedge \tfrac{1}{|\varphi|^2} (\langle \varphi^* \nabla_{\hat{A}} \varphi \rangle - \langle \varphi \nabla_{\hat{A}} \varphi^* \rangle) \sigma \ .$$

(6.16)

The connection $\mathcal{A}$ is smooth and equal to $\hat{A}$ where $|z| \leq rt$ (because $\chi_*$ is zero there). By way of comparison, it has zero curvature and $\nabla_\mathcal{A} \varphi \equiv 0$ where $|z| \geq 2rt$ (because $\chi_* \equiv 1$ there). This is by design.

A crucial point is that the $F_\mathcal{A} \equiv 0$ property where $|z| \geq 2rt$ implies that the integral of $\langle \sigma F_\mathcal{A} \rangle$ on the constant $(t, x_3)$ slices of $(0, \infty) \times \mathbb{R}^2 \times \mathbb{R}$ is independent of t and $x_3$. (The $\nabla_\mathcal{A} \varphi = 0$ property where $|z| \geq 2rt$ implies that this integral is equal to $2\pi$ times the degree of vanishing of $\varphi$ at the $z = 0$ point of these constant $(t, x_3)$ slices.)

Step 2: As noted in Step 1, the curvature 2-form of $\mathcal{A}$ pulls back to any given constant $(t, x_3)$ slice of $(0, \infty) \times \mathbb{R}^2 \times \{0\}$ so as to be zero where $|z| \geq 2rt$. It is also equal to $B_{\hat{A}3} dx_1 \wedge dx_2$ where $|z| \leq rt$. In between (on the $rt \leq |z| \leq 2rt$ annulus) it obeys the bound $|F_\mathcal{A}| \leq \mu \tfrac{1}{t^2}$ with $\mu$ being independent of t and bounded a priori given $\mathfrak{z}, \delta$ and $r$. Such a bound follows from the formula in (6.16) and the formula in (3.12) given the upper bounds for $|\hat{b}|$ on the $rt \leq |z| \leq 2rt$ annulus (from Lemma 6.1 for $|\hat{b}|$) and the bounds from Lemma 2.2 for $|\nabla_\mathcal{A} \varphi|$ and $|B_{\hat{A}3}|$; and from the lower bound $|\varphi| \geq \tfrac{\delta}{t}$ where $|z| > rt$ from the second bullet in (6.1).

Step 3: Since $|F_\mathcal{A}| \leq \mu \tfrac{1}{t^2}$ where $rt \leq |z| \leq 2rt$, the contribution to the integral of $\langle \sigma F_\mathcal{A} \rangle$ over a given constant $(t, x_3)$ slice of $(0, \infty) \times \mathbb{R}^2 \times \mathbb{R}$ from its $|z| \geq rt$ part is bounded by $\mu r^2$ which is independent of t. Since the total integral over the slice is proportional to the degree of vanishing of $\varphi$ at the $z = 0$ point (which is independent of t), it follows as a consequence that the integral of $\langle \sigma B_{\hat{A}3} \rangle$ over the $|z| \leq rt$ slice is independent of t and bounded a priori given $\mathfrak{z}, r, \delta, \Xi$ and the degree of vanishing of $\varphi$ at $z = 0$. (Keep in mind that the pull-back of $F_\mathcal{A}$ to the $|z| \leq rt$ part of the slice is $B_{\hat{A}3} dx_1 dx_2$.) With the preceding understood, then a t-independent upper bound for the integral of $|\hat{b}|^2$ over the $|z| \leq rt$ part of the slice (determined by $\mathfrak{z}, r, \delta$ and $\Xi$ and the degree of vanishing of $\varphi$ at the $z = 0$ point) follows from the identity in the top bullet of (3.12) since $|B_{A3}| \leq c_\mathfrak{z} \frac{1}{t^2}$ (from Lemma 2.2) and since the area of the $|z| \leq rt$ disk in the slice is $\pi r^2 t^2$.

Step 4: Granted what is said in Step 3, all that remains is to note that the degree of vanishing of $\varphi$ at the $z = 0$ point in any constant $(t, x_3)$ slice of $(0, \infty) \times \mathbb{R}^2 \times \mathbb{R}$ has an upper bound determined a priori by the numbers $\mathfrak{z}, \delta, r$. (This is what is said by Lemma 5.6).

*Part 3*: This part of the proof gives the argument for the third bullet of the lemma. To this end, fix $t > 0$ and $\varepsilon \in (0, 1)$ and suppose that $p \in \mathbb{R}^2$ is a point with norm greater than R and such that $|\beta| + |\mathfrak{b}|$ at $(t, p, 0)$ is greater than $\frac{\varepsilon}{t}$. Because $R > t$, the norm of $|\mathfrak{b}|$ on the disk of radius $\frac{1}{2} t$ in $\{t\} \times \mathbb{R}^2 \times \{0\}$ centered at $(t, p, 0)$ is bounded by $\frac{\mathcal{Z}}{t}$ with $\mathcal{Z}$ being independent of t and p and R. This bound leads via Lemma 6.2 to a $\frac{\mathcal{Z}_1}{t^2}$ bound for the norm of $\nabla^\perp(|\beta| + |\mathfrak{b}|)$ on the concentric, radius $\frac{1}{4} t$ disk with $\mathcal{Z}_1$ being independent of t, p and R. As a consequence $|\beta| + |\mathfrak{b}|$ is greater than $\frac{\varepsilon}{2t}$ in the concentric disk of radius $\frac{\varepsilon}{4\mathcal{Z}_2} t$ in $\{t\} \times \mathbb{R}^2 \times \{0\}$ with $\mathcal{Z}_2$ denoting the maximum of $\mathcal{Z}_1$ and 1. These bounds lead in turn to a $\frac{\varepsilon^4}{64(\mathcal{Z}')^2}$ lower bound for the integral of $|\beta|^2 + |\mathfrak{b}|^2$ over this same disk. That lower bound runs a foul of the bound in the second bullet of the lemma if $\varepsilon > \mathcal{Z}_3 (\frac{t}{R})^{1/4}$ with $\mathcal{Z}_3$ being independent of t and p and R.

### d) The $t \to 0$ and $t \to \infty$ limits of constant $(t, x_3)$ integrals of $|\beta|^2 + |\mathfrak{b}|^2$

The central lemma for this subsection concerns the behavior of the constant $(t, x_3)$ integrals of $|\beta|^2 + |\mathfrak{b}|^2$ in the event that the function $t|\alpha|$ is bounded from below by a positive number where t is either very small or very large.

**Lemma 6.4**: *Suppose that $(A, \mathfrak{a})$ is a solution to (1.4) that is described by (6.1).*
- *Fix $\delta > 0$ and use it to define a function on $(0, \infty)$ by the rule*

$$t \to t^{-4\delta} \int_{\{t\} \times \mathbb{R}^2 \times \{0\}} (|\beta|^2 + |\mathfrak{b}|^2) \ .$$

*This function is bounded on $(0, t_0]$ if the $(A, \mathfrak{a})$ version of $\alpha$ is less than $-\frac{\delta}{t}$ on $(0, t_0]$.*
- *Fix $\delta > 0$ and use it to define a function on $(0, \infty)$ by the rule*

$$t \to t^{4\delta} \int_{\{t\} \times \mathbb{R}^2 \times \{0\}} (|\beta|^2 + |\mathfrak{b}|^2) \ .$$

*This function is bounded on $(t_0, \infty)$ if the $(A, \mathfrak{a})$ version of $\alpha$ is less than $-\frac{\delta}{t}$ on $[t_0, \infty)$.*

The proof of this lemma exploits an equation that can be derived from the equations in the two bullets of (3.11):

$$\tfrac{\partial}{\partial t}(|\beta|^2 - |\mathfrak{b}|^2) + i(\tfrac{\partial}{\partial x_1} - i\tfrac{\partial}{\partial x_2})\langle\beta^*\mathfrak{b}\rangle - i(\tfrac{\partial}{\partial x_1} + i\tfrac{\partial}{\partial x_2})\langle\mathfrak{b}^*\beta\rangle = -4\alpha(|\beta|^2 + |\mathfrak{b}|^2).$$
(6.17)

This equation is obtained by first taking the inner product of both sides of the top bullet's equation in (3.12) with $\beta^*$; and taking the inner product of both sides of the second bullet's equation in (3.12) with $\mathfrak{b}^*$. Having done that, subtract the first equation from the second and take the real part of what results.

*Proof of Lemma 6.4*: The proof has four parts.

*Part 1*: Fix $t \in (0,\infty)$ and then fix $R \gg t$. Let $\chi_R$ denote the function on the constant t slice of $(0,\infty)\times\mathbb{R}^2\times\mathbb{R}$ that is given by the rule $\chi_R(z) = \chi(\tfrac{|z|}{R} - 1)$. This function is equal to 1 where $|z| < R$ and it is equal to zero where $|z| > 2R$. Multiply both sides of (6.17) by $\chi_R$ and then integrated over the slice of $\{t\}\times\mathbb{R}^2\times\{0\}$. Integrate by parts move the $x_1$ and $x_2$ derivatives from $\langle\beta^*\mathfrak{b}\rangle$ and $\langle\mathfrak{b}^*\beta\rangle$ to the function $\chi_R$. Since the norms of the derivatives of $\chi_R$ are bounded by $c_0\tfrac{1}{R}$, this justifies the assertion that the norms of the $\langle\beta^*\mathfrak{b}\rangle$ and $\langle\mathfrak{b}^*\beta\rangle$ integrals are bounded by a (t, R) independent multiple of $\tfrac{1}{R}$. (Remember that $|\beta|^2$ and $|\mathfrak{b}|^2$ have integrals on $\{t\}\times\mathbb{R}^2\times\{0\}$ with a t-independent upper bound.) Therefore, taking $R \to \infty$ on both sides of the family of finite R identities leads to this:

$$\tfrac{d}{dt}\int_{\{t\}\times\mathbb{R}^2\times\{0\}}(|\beta|^2 - |\mathfrak{b}|^2) = -4\int_{\{t\}\times\mathbb{R}^2\times\{0\}}\alpha(|\beta|^2 + |\mathfrak{b}|^2).$$
(6.18)

Now suppose that $\delta$ is positive and that the function $\alpha$ obeys the bound $\alpha < -\tfrac{\delta}{t}$ where t lies in some given interval in $(0, \infty)$. In this event, (6.18) leads to this

$$\tfrac{d}{dt}\int_{\{t\}\times\mathbb{R}^2\times\{0\}}(|\beta|^2 - |\mathfrak{b}|^2) \geq \tfrac{4\delta}{t}\int_{\{t\}\times\mathbb{R}^2\times\{0\}}(|\beta|^2 + |\mathfrak{b}|^2)$$
(6.19)

on the given interval in $(0, \infty)$.

*Part 2*: This part and Part 3 prove the top bullet's assertion. To start, note that the assumptions leading to (6.19) hold on the interval $(0, t_0]$ with $t_0 > 0$. Thus, (6.19) is valid on this interval. Let $f(t)$ denote the $\{t\}\times\mathbb{R}^2\times\{0\}$ integral of $|\beta|^2 - |\mathfrak{b}|^2$. As explained directly, this function can not be negative at any $t_1 \in (0, t_0]$. Indeed, if it is negative at $t=t_1$, then (6.19) implies that it is negative on the whole of the interval $(0, t_1]$ and that

$$\tfrac{d}{dt} f \geq -\tfrac{4\delta}{t} f \tag{6.20}$$

on $(0, t_1]$. This implies in turn that $f$ obeys $f(t) \leq (\tfrac{t_1}{t})^{4\delta} f(t_1)$ which leads to nonsense since $|f|$ has a t-independent upper bound.

Given that $f \geq 0$ on $(0, t_0]$, then (6.19) leads to this instead:

$$\tfrac{d}{dt} f \geq \tfrac{4\delta}{t} f \quad on \ (0, t_0] \ . \tag{6.21}$$

And, (6.21) implies that $f(t) \leq (\tfrac{t}{t_0})^{4\delta} f(t_0)$ for $t \leq t_0$. Granted this, then integrating (6.19) leads in turn to the assertion that

$$\int_{(0,t]\times \mathbb{R}^2 \times \{0\}} \tfrac{1}{s}(|\beta|^2 + |\mathfrak{b}|^2) \leq \tfrac{1}{4\delta}(\tfrac{t}{t_0})^{4\delta} f(t_0) \ . \tag{6.22}$$

The latter bound implies the following: Fix $t \in (0, t_0]$ and $\mathcal{Z} \geq 2$. Then the subset in $[\tfrac{1}{2} t, t]$ where the $\{s\} \times \mathbb{R}^2 \times \{0\}$ integral of $|\beta|^2 + |\mathfrak{b}|^2$ is greater than $\mathcal{Z} \tfrac{1}{4\delta} (\tfrac{t}{t_0})^{4\delta} f(t_0)$ has measure at most $\tfrac{t}{\mathcal{Z}}$.

*Part 3*: The concluding observation in Part 2 is almost but not quite the assertion made by top bullet of Lemma 6.4. The full assertion is proved with the help of a Bochner-Weitzenboch formula for $\beta$ that is obtained by differentiating the two bullets in (3.11). This formula is as follows:

$$\nabla_{\hat{A}t}^2 \beta + \nabla_{\hat{A}1}^2 \beta + \nabla_{\hat{A}2}^2 \beta = (-2\langle \sigma B_{\hat{A}3} \rangle - 2 \tfrac{\partial}{\partial t} \alpha + 4\alpha^2) \beta - 4i(\tfrac{\partial}{\partial x_1} \alpha - i \tfrac{\partial}{\partial x_2} \alpha) \mathfrak{b} \ . \tag{6.23}$$

Take the inner product of both sides of this equation with $\beta^*$ (and use the complex conjugate equation also) to obtain this:

$$-\tfrac{1}{2}(\tfrac{\partial^2}{\partial t^2} + \tfrac{\partial^2}{\partial x_1^2} + \tfrac{\partial^2}{\partial x_2^2})|\beta|^2 + |\nabla_{\hat{A}t}\beta|^2 + |\nabla_A^\perp \beta|^2 = (-2\langle \sigma B_{\hat{A}3} \rangle - 2\partial_t \alpha + 4\alpha^2)|\beta|^2$$
$$- 2i(\tfrac{\partial}{\partial x_1} \alpha - i \tfrac{\partial}{\partial x_2} \alpha)\langle \beta^* \mathfrak{b} \rangle + 2i(\tfrac{\partial}{\partial x_1} \alpha + i \tfrac{\partial}{\partial x_2} \alpha)\langle \mathfrak{b}^* \beta \rangle \ . \tag{6.24}$$

With (6.24) in hand, use the identities in (3.12) to see that (6.24) leads to the inequality

$$-\tfrac{1}{2}\left(\tfrac{\partial^2}{\partial t^2}+\tfrac{\partial^2}{\partial x_1^2}+\tfrac{\partial^2}{\partial x_2^2}\right)|\beta|^2+|\nabla_{\hat{A}t}\beta|^2+|\nabla_A^{\perp}\beta|^2 \leq \tfrac{\varsigma}{t^2}(|\beta|^2+|\flat|^2)$$

(6.25)

with $\varsigma$ being independent of t. The number $\varsigma$ comes from Lemma 2.2 and from using the top bullet in Lemma 6.4 to bound the norm of $\tfrac{\partial}{\partial t}\alpha$ and from using the top bullet in (6.1) to bound $|\beta|$ and $|\alpha|$. A crucial fact: The terms on right hand side of (6.24) are at most quadratic with regards to $\flat$ when $B_{\hat{A}3}$ and the derivatives of $\alpha$ are written in terms of the curvature of A and quadratic terms in $\varphi$, $\beta$ and $\flat$. Some of these quadratic appearances of $|\flat|$ on the right hand side of (6.24) are multipled by $|\beta|^2$, but this is of little concern because $|\beta|$ has an a priori bound proportional to a t-independent multiple of $\tfrac{1}{t}$. This is important because there is no t-independent, pointwise upper bound for $|\flat|$; there is only a t-independent bound for the $\{t\}\times\mathbb{R}^2\times\{0\}$ integral of $|\flat|^2$. (A cubic or quartic term in $|\flat|$ on the right hand side of (6.24) would be evil if it appeared with the wrong sign. The upcoming (6.30) depicts a Bochner-Weitzenboch formula for $\flat$ which is analogous to (6.23) and which leads in turn to a version of (6.24) where $|\flat|^4$ appears with an evil sign.)

Integrate both sides of (6.25) over a given constant t version of $\{t\}\times\mathbb{R}^2\times\{0\}$. Multiple instances of integration by parts on the left hand side (which can be justified using the cut-off functions $\{\chi_R\}_{R>t}$ from Step 1) can be used to see the following:

$$-\tfrac{1}{2}\tfrac{\partial^2}{\partial t^2}\int_{\{t\}\times\mathbb{R}^2\times\{0\}}|\beta|^2 + \int_{\{t\}\times\mathbb{R}^2\times\{0\}}(|\nabla_{\hat{A}t}\beta|^2+|\nabla_{\hat{A}}^{\perp}\beta|^2) \leq \tfrac{\varsigma}{t^2}\int_{\{t\}\times\mathbb{R}^2\times\{0\}}(|\beta|^2+|\flat|^2).$$

(6.26)

To procede from here, fix $t\in(0,t_0]$ and then use the function $\chi$ to construct a bump function on $(0,\infty)$ which is non-negative, equal to 1 on $[\tfrac{1}{2}t,\tfrac{3}{2}t]$ and equal to 0 on the domains $(0,\tfrac{1}{4}t]$ and $[2t,\infty)$. This function can and should be constructed so that the norm of its second derivative is bounded by $c_0\tfrac{1}{t^2}$. Multiply both sides of (6.26) by this function, then integrate the resulting inequality over $(0,\infty)$. What with (6.22), two instances of integration by parts on the left hand side leads to an inequality of the form

$$\int_{[\tfrac{1}{2}t,\tfrac{3}{2}t]\times\mathbb{R}^2\times\{0\}}(|\nabla_{\hat{A}t}\beta|^2+|\nabla_{\hat{A}}^{\perp}\beta|^2) \leq \mathfrak{z}\,t^{4\delta-1}$$

(6.27)

with $\mathfrak{z}$ being independent of t. The bound here for the t-derivative of $\beta$ implies (using the fundamental theorem of calculus) that

$$\left(\int_{\{s\}\times\mathbb{R}^2\times\{0\}}|\beta|^2\right)^{1/2} - \left(\int_{\{s'\}\times\mathbb{R}^2\times\{0\}}|\beta|^2\right)^{1/2} \leq \sqrt{\mathfrak{z}}\,t^{2\delta}$$

(6.28)

for any pair (s, s´) from $[\frac{1}{2}t, \frac{3}{2}t]$. The latter bound with what is said at the end of Part 2 gives the lemma's bound for the $\{t\}\times\mathbb{R}^2\times\{0\}$ integral of $|\beta|^2$ when $t \in (0, t_0]$. And, such a bound for the integral of $|\beta|^2$ leads to the same bound for that of $|\mathfrak{b}|^2$ because the former integral is greater than the latter on $(0, t_0]$.

*Part 4*: This part proves the assertion made by the second bullet of Lemma 6.4. The argument here is analogous up to a point. First, arguments from Part 2 can be redone when (6.19) holds on $[t_0, \infty)$ and they prove that the function $f$ must be *negative* on that interval rather than positive (which is to say that the $\{t\}\times\mathbb{R}^2\times\mathbb{R}$ integral of $|\mathfrak{b}|^2$ must be greater than that of $|\beta|^2$ when $t \in [t_0, \infty)$.)

Granted the preceding changes, then arguments that differ only cosmetically from those used subsequent to (6.20) in Part 2 lead to the bound $|f(t)| \leq (\frac{t_0}{t})^{4\delta}|f|(t_0)$ when $t \geq t_0$ and to the following analog of (6.22):

$$\int_{[t,\infty)\times\mathbb{R}^2\times\{0\}} \tfrac{1}{s}(|\beta|^2 + |\mathfrak{b}|^2) \leq \tfrac{1}{4\delta}(\tfrac{t_0}{t})^{4\delta}|f|(t_0)$$

(6.29)

which holds for $t \geq t_0$. The latter bound implies in turn that the $\{t\}\times\mathbb{R}^2\times\{0\}$ integral of $|\beta|^2 + |\mathfrak{b}|^2$ must be less than $\frac{1000}{\delta}(\frac{t_0}{t})^{4\delta}|f|(t_0)$ on more than a fraction greater than $\frac{3}{4}$ of any interval of the form $[t, 2t)$ in $[t_0, \infty)$.

The argument used in Part 3 can now be repeated (using $\beta$ still) to obtain a $\mathfrak{z}t^{-4\delta}$ bound for the $\{t\}\times\mathbb{R}^2\times\{0\}$ integral of $|\beta|^2$ for any $t \in [t_0, \infty)$ with $\mathfrak{z}$ being independent of the choice of t. This bound for $|\beta|^2$ leads to a $2\mathfrak{z}\, t^{-4\delta}$ bound for the $\{t\} \times \mathbb{R}^2\times \{0\}$ integral of $|\mathfrak{b}|^2$ when the latter integral is less than twice that of $|\beta|^2$. When the $\{t\} \times \mathbb{R}^2\times \{0\}$ integral of $|\mathfrak{b}|^2$ is greater than twice that of $|\beta|^2$, then it is less than twice the value of $f(t)$ which is less than $(\frac{t_0}{t})^{4\delta}|f|(t_0)$, which is also a t-independent multiple of $(\frac{1}{t})^{4\delta}$.

### e) Pointwise bounds for $|\beta|^2 + |\mathfrak{b}|^2$ when its $\{t\}\times\mathbb{R}^2\times\{0\}$ integral is small

The central lemma in this subsection says in effect that $t|\beta|$ and $t|\mathfrak{b}|$ are both poinwise small at times t when the $\{t\}\times\mathbb{R}^2\times\{0\}$ integral of $|\beta|^2+|\mathfrak{b}|^2$ is small.

**Lemma 6.5**: *Given $\mathfrak{z} > 0$, there exists $\kappa > 1$ with the following significance: Suppose that $(A, \mathfrak{a})$ is a solution to (1.4) such that $|\mathfrak{a}| \leq \frac{3}{t}$ on $(0,\infty)\times\mathbb{R}^2\times\mathbb{R}$. Given $t > 0$, if the corresponding $\beta$ and $\mathfrak{b}$ obey*

$$\int_{\{s\}\times\mathbb{R}^2\times\{0\}} (|\beta|^2 + |\mathfrak{b}|^2) \leq \tfrac{1}{\kappa}$$

*for all* $s \in [\tfrac{1}{4}t, 2t]$, *then*

$$|\beta| + |\mathfrak{b}| \leq \tfrac{\kappa}{t} \big( \sup_{s\in[\tfrac{1}{2}t,\tfrac{3}{2}t]} \int_{\{s\}\times\mathbb{R}^2\times\{0\}} (|\beta|^2 + |\mathfrak{b}|^2) \big)^{1/2}$$

*on the whole of* $[\tfrac{3}{4}t, \tfrac{5}{4}t] \times \mathbb{R}^2 \times \mathbb{R}$.

Note that Lemmas 6.4 and 6.5 lead immediately to very small pointwise bounds for $t|\beta|$ and $t|\mathfrak{b}|$ where either $t \leq (0, t_0]$ or $t \geq [t_0, 0)$ if $t\alpha$ is bounded away from zero on the given interval. Here is a precise statement:

**Corollary 6.6**: *Suppose that* $(A, \mathfrak{a})$ *is a solution to (1.4) that is described by (6.1). If* $\delta$ *is positive and if the corresponding version of* $\alpha$ *is less than* $-\tfrac{\delta}{t}$ *on an interval of the form* $(0, t_0]$, *then* $|\beta| + |\mathfrak{b}| \leq \kappa_*(\tfrac{1}{t})^{1-2\delta}$ *on this interval with* $\kappa_*$ *being independent of* $t$. *By the same token, if* $\alpha < -\tfrac{\delta}{t}$ *on an interval of the form* $[t_0, \infty)$, *then* $|\beta| + |\mathfrak{b}| \leq \kappa_*(\tfrac{1}{t})^{1+2\delta}$ *on this same interval with* $\kappa_*$ *likewise independent of* $t$.

Lemma 6.5 is proved with the help of another lemma which says (roughly) that the sup norm of $t|\mathfrak{b}|$ can be bounded uniformly on any disk in $\{t\}\times\mathbb{R}^2\times\mathbb{R}$ where the integral of $|\mathfrak{b}|^2 + |\beta|^2$ is small.

**Lemma 6.7**: *Given* $\mathfrak{z} > 0$, *there exists* $\kappa > 1$ *with the following significance: Suppose that* $(A, \mathfrak{a})$ *is a solution to (1.4) such that* $|\mathfrak{a}| \leq \tfrac{\mathfrak{z}}{t}$ *on* $(0,\infty)\times\mathbb{R}^2\times\mathbb{R}$. *Fix* $t \in (0, \infty)$ *and fix a disk* $D \subset \mathbb{R}^2$ *of radius* $\gamma t$ *with* $\gamma \leq 1$. *Let* $D' \subset D$ *denote the concentric disk with radius* $\tfrac{1}{2}\gamma t$. *Fix* $\delta \in (0, 1)$ *such that* $|\mathfrak{b}|$ *at some point in* $D'$ *is greater than* $\delta$ *times its supremum on* $D$. *If the integral of* $|\mathfrak{b}|^2$ *over* $\{t\}\times D\times\{0\}$ *is less than* $\tfrac{1}{\kappa}\delta^2$, *then* $|\mathfrak{b}| \leq \tfrac{\kappa}{\delta\gamma t}$ *in the disk* $D$.

Lemma 6.5 is proved first assuming Lemma 6.7.

***Proof of Lemma 6.5***: If the $\{s\} \times \mathbb{R}^2 \times \{0\}$ integrals of $|\mathfrak{b}|^2$ for $s \in [\tfrac{1}{4}t, 2t]$ are bounded by a suitably small positive number supplied by Lemma 6.7, then Lemma 6.7 leads to an a priori bound $\tfrac{\mathcal{Z}}{t}$ bound for $|\mathfrak{b}|$ on these $\{s\} \times \mathbb{R}^2 \times \{0\}$ domains with $\mathcal{Z}$ independent of $t$. (Lemma 6.1 bounds $|\mathfrak{b}|$ by a t-independent multiple of $\tfrac{1}{t}$ where $|z| \geq rs$ on $\{s\}\times\mathbb{R}^2\times\{0\}$ when $s \in [\tfrac{1}{4}t, 2t]$.) With that understood, take D for Lemma 6.7 to be a radius $\tfrac{1}{100}t$ disk centered at a point in $\{s\}\times\mathbb{R}^2\times\{0\}$ where $|\mathfrak{b}|$ is greater than $\tfrac{1}{2}$ its supremum on $\{s\}\times\mathbb{R}^2\times\{0\}$.) This $\tfrac{\mathcal{Z}}{t}$ bound for $|\mathfrak{b}|$ is the input from Lemma 6.7 to Lemma 6.5's proof.

The remainder of the proof is given in four parts. By way of notation, this proof uses $\varsigma$ to denote a number that is independent of t and greater than 1. It's value can be assumed to increase between successive appearances. (It can depend on (A, $\mathfrak{a}$) however.)

*Part 1*: The arguments in Part 3 of the proof of Lemma 6.4 can be repeated starting from the Bochner-Weitzenboch formula in (6.23) to obtain the differential inequality in (6.25). (The derivation of (6.25) does not use the $\frac{z}{t}$ bound for $|\mathfrak{b}|$.) There is also a Bochner-Weitzenboch formula from (3.11) for $\mathfrak{b}$ which is this one:

$$\nabla_{\hat{A}t}^2 \mathfrak{b} + \nabla_{\hat{A}1}^2 \mathfrak{b} + \nabla_{\hat{A}2}^2 \mathfrak{b} = (2\langle \sigma B_{\hat{A}3}\rangle + 2\tfrac{\partial}{\partial t}\alpha + 4\alpha^2)\mathfrak{b} + 4i(\tfrac{\partial}{\partial x_1}\alpha + i\tfrac{\partial}{\partial x_2}\alpha)\beta \;.$$

(6.30)

Taking the inner product of both sides of this with $\mathfrak{b}^*$ leads to a corresponding analog of (6.24). The latter then leads to the following analog of (6.25):

$$-\tfrac{1}{2}(\tfrac{\partial^2}{\partial t^2} + \tfrac{\partial^2}{\partial x_1^2} + \tfrac{\partial^2}{\partial x_2^2})|\mathfrak{b}|^2 + |\nabla_{\hat{A}t}\mathfrak{b}|^2 + |\nabla_A^\perp \mathfrak{b}|^2 \leq \tfrac{\varsigma}{t^2}(|\beta|^2 + |\mathfrak{b}|^2)$$

(6.31)

where $\varsigma$ is again independent of t. The derivation of (6.31) <u>does</u> use the $\frac{z}{t}$ bound for $|\mathfrak{b}|$; it is used to bound an evil signed multiple of $|\mathfrak{b}|^4$ by $(\frac{z}{t})^2|\mathfrak{b}|^2$. (This evil $|\mathfrak{b}|^4$ would appear on the right hand side of (6.31) with a positive sign. It comes from using (3.12) to write $\langle \sigma B_{\hat{A}3}\rangle$ in (6.31) as $\langle \sigma B_{A3}\rangle - 4|\mathfrak{b}|^2$.)

A first important take-away from (6.25) and (6.30) is this: The function $|\beta| + |\mathfrak{b}|$ obeys the differential inequality

$$-\tfrac{1}{2}(\tfrac{\partial^2}{\partial t^2} + \tfrac{\partial^2}{\partial x_1^2} + \tfrac{\partial^2}{\partial x_2^2})(|\beta| + |\mathfrak{b}|) \leq \tfrac{\varsigma}{t^2}(|\beta| + |\mathfrak{b}|)$$

(6.32)

on $[\tfrac{1}{4}t, 2t] \times \mathbb{R}^2 \times \{0\}$. This inequality will be used momentarily with the Green's function for the Laplace operator to obtain the desired pointwise bounds for $|\beta|$ and $|\mathfrak{b}|$.

*Part 2*: To exploit (6.32), fix a point p in $[\tfrac{3}{4}t, \tfrac{5}{4}t] \times \mathbb{R}^2 \times \{0\}$ and let B denote the radius $\tfrac{1}{100}t$ ball centered at p. Introduce by way of notation $G_p$ to denote the Dirichelet Green's function for the Laplacian on B with pole at the point p. This function is positive on B−{p} where it obeys the bound $G_p(\cdot) \leq \tfrac{1}{4\pi}\tfrac{1}{|(\cdot) - p|}$.

Another function is needed, which is a bump function constructed using $\chi$ which must be non-negative, compactly supported in B and equal to 1 in the ball concentric to B with radius $\tfrac{1}{200}t$. This function can and should be constructed so that the norms of its second derivatives are bounded by $c_0 \tfrac{1}{t^2}$.

Multiply both sides of (6.32) by the product of this bump function and $G_p$. Then integrate by parts on the left hand side to see that

$$|\beta|(p) + |\mathfrak{b}|(p) \leq c_0 \frac{1}{t^3} \int_{B\times\{0\}} (|\beta| + |\mathfrak{b}|) + \frac{\varsigma}{t^2} \int_B \frac{1}{|(\cdot) - p|}(|\beta| + |\mathfrak{b}|) .$$

(6.33)

The left most integral on the right hand side of this inequality is bounded by

$$c_0 \frac{1}{t^{3/2}} \Big( \int_{B\times\{0\}} (|\beta|^2 + |\mathfrak{b}|^2) \Big)^{1/2} .$$

(6.34)

The next steps explain why the right most integral in (6.33) has a similar bound (with $c_0$ replaced by $c_0\varsigma$). Lemma 6.5 follows directly if the right hand side of (6.33) is no greater than $\varsigma$ times larger what is written in (6.34).

*Part 3*: This part of the proof sets the stage for bounding the right most integral on the right hand side of (6.33) by $\varsigma$ times what is written in (6.34). To start, repeat the arguments in Part 3 of the proof of Lemma 6.4 from (6.25) through (6.27) to obtain the following analog of (6.27):

$$\int_{[\frac{1}{2}t,\frac{3}{2}t]\times\mathbb{R}^2\times\{0\}} (|\nabla_{\hat{A}t}\beta|^2 + |\nabla_{\hat{A}}^{\perp}\beta|^2) \leq \frac{\varsigma}{t^2} \int_{[\frac{1}{4}t,2t]\times\mathbb{R}^2\times\{0\}} (|\beta|^2 + |\mathfrak{b}|^2)$$

(6.35)

The bump function built from $\chi$ that is used to derive (6.35) can be used again with (6.31) to obtain the following $\mathfrak{b}$ version of (6.35):

$$\int_{[\frac{1}{2}t,\frac{3}{2}t]\times\mathbb{R}^2\times\{0\}} (|\nabla_{\hat{A}t}\mathfrak{b}|^2 + |\nabla_{\hat{A}}^{\perp}\mathfrak{b}|^2) \leq \frac{\varsigma}{t^2} \int_{[\frac{1}{4}t,2t]\times\mathbb{R}^2\times\{0\}} (|\beta|^2 + |\mathfrak{b}|^2) .$$

(6.36)

Note in particular that the left hand side of (6.35) bounds the integral of $|\nabla|\beta||^2$ over the indicated domain; and the left hand side of (6.36) bounds that of $|\nabla|\mathfrak{b}||^2$.

*Part 4*: There is what is called *Hardy's* inequality for functions on a ball in $\mathbb{R}^3$ whose squares are integrable and whose derivatives have integrable squares. This inequality asserts the following: Supposing that $r > 0$, let B denote a ball of radius r in $\mathbb{R}^3$ centered at a given point p. Suppose that $f$ is a function on B with $f^2$ and $|\nabla f|^2$ being integrable on $B_r$. Then

$$\int_B \frac{1}{|(\cdot) - p|^2} f^2 \leq c_0 \int_B (|\nabla f|^2 + \tfrac{1}{r^2} f^2) \ .$$

(6.37)

This inequality will be applied momentarily using $f = |\beta| + |\flat|$ with B denoting a given radius $\frac{1}{100}$ t ball in $[\tfrac{3}{4}t, \tfrac{5}{4}t] \times \mathbb{R}^2 \times \{0\}$. The right hand side of (6.37) for these versions of $f$ and B obeys

$$\int_B \frac{1}{|(\cdot) - p|^2} (|\beta| + |\flat|)^2 \leq \frac{\varsigma}{t^2} \int_{[\tfrac{1}{4}t, 2t] \times \mathbb{R}^2 \times \{0\}} (|\beta|^2 + |\flat|^2)$$

(6.38)

This is because of (6.35) and (6.36) and what is said at the end of Part 3.

The bound in (6.38) leads directly to the desired bound (by $\varsigma$ times what is depicted in (6.31)) for the right most integral on the hand side of (6.34).

*Proof of Lemma 6.7*: Write $\hat{A}$ as $A - \flat$ to write the top bullet in (3.13) as

$$-i(\nabla_{A1} - i\nabla_{A2})\flat = B_{A3}{}^+ + 4|\flat|^2 \sigma \ .$$

(6.39)

To exploit this equation, fix a positive number to be denoted by $\delta$ and suppose for the moment that there is a point in $D'$ which is such that $|\flat|$ at $(t, q, 0)$ is greater than $\delta$ times it supremum in $\{t\} \times D \times \{0\}$. Choose such a point in $D'$ and denote it by q. Then use parallel transport by A along the rays in $\{t\} \times \mathbb{R}^2 \times \{0\}$ from the point $(t, q, 0)$ to define an isomorphism between the bundle ad(P) on $\{t\} \times D \times \{0\}$ and the product lie algebra bundle. Use this isomorphism to write A as $\theta_0 + \mathfrak{k}$ with $\mathfrak{k}$ being a Lie algebra valued 1-form on $\{t\} \times D \times \{0\}$. The isomorphism then writes (6.39) as

$$-i(\tfrac{\partial}{\partial x_1} - i\tfrac{\partial}{\partial x_2})\flat = B_{A3}{}^+ + i[\mathfrak{k}_1 - i\mathfrak{k}_2, \flat] + 4|\flat|^2 \sigma \ .$$

(6.40)

To continue, let $\chi_q$ denote the function on D given by $\chi(\frac{100|z-q|}{\gamma t} - 1)$ which is a function with compact support in D that is equal to 1 on the disk in D with radius $\frac{1}{100}\gamma t$ centered at the point q. Multiply both sides of the equation in (6.41) by the cut-off Green's function $\frac{i}{2\pi(z-q)} \chi_q$ and integrate the result over D. An integration by parts then leads to an inequality that bounds $|\flat|$ at $(t,q,0)$ by $c_0$ times a sum of four integrals:

$$\int_{\{t\}\times D\times\{0\}} \tfrac{1}{|z-q|} |B_{A3}| + \int_{\{t\}\times D\times\{0\}} \tfrac{1}{|z-q|} |\mathfrak{k}| |\flat| + \int_{\{t\}\times D\times\{0\}} \tfrac{1}{|z-q|} |\nabla\chi_q| |\flat| + \int_{\{t\}\times D\times\{0\}} \tfrac{1}{|z-q|} |\flat|^2 \ .$$

(6.41)

The task now is to obtain useful bounds for these integrals. The paragraphs that follow do this moving from left most to right most.

The integral with $B_{A3}$ is bounded by $c_0 c_ɜ \gamma \frac{1}{t}$ because $|B_{A3}|$ is bounded by $c_ɜ \frac{1}{t^2}$ (see Lemma 2.2) and because the integral of the function $\frac{1}{|z-q|}$ over D is bounded by $c_0 \gamma t$.

To bound the integral with $|ƻ|$, note first that the 1-form $ƻ$ annihilates the tangent vectors to the rays in D from q. Because of this, it can be written at any given $z \in D$ as an integral along the ray from q to z using $B_{A3}$. This integral can be used to bound the norm of $ƻ$ by $|ƻ| \leq c_ɜ |z-q| \frac{1}{t^2}$. Therefore, the integral over $\{t\} \times D \times \{0\}$ of $\frac{1}{|z-q|} |ƻ||Ƅ|$ is at most

$$c_0 c_ɜ \frac{1}{t^2} \int_{\{t\} \times D \times \{0\}} |Ƅ|$$

(6.42)

which is, in turn, at most $c_0 c_ɜ \frac{1}{t} ( \int_{\{t\} \times D \times \{0\}} |Ƅ|^2 )^{1/2}$.

The integral with $|\nabla \chi_q|$ is bounded by $\frac{1}{\gamma}$ times what is written in (6.42) (with a perhaps larger $c_0$) because $|\nabla \chi_q| \leq c_0 \frac{1}{\gamma t}$ and the distance to q from where it is not zero is greater than $c_0^{-1} \gamma t$.

To obtain a bound for the left most term in (6.41), introduce by way of notation ¥ to denote the supremum of the function $|Ƅ|$ on the disk $\{t\} \times D \times \{0\}$. There is nothing to prove if $¥ < \frac{200}{\gamma t}$, so suppose that this is not the case. Fix for the moment $\varepsilon \in (0, 1)$ and let $r = \frac{\varepsilon}{¥}$. (Thus, $r < \frac{\varepsilon}{200} \gamma t$.) Break the integral in (6.41) into the part where $|z-q| \leq r$ and the part where $|z-q| \geq r$. The former (the part where $|z-q| \leq r$) is no greater than $2\pi ¥^2 r$ which is equal to $2\pi \varepsilon ¥$ by virtue of the definition of r. The part where $|z-q| \geq r$ is no greater than $\frac{1}{r}$ times the integral of $|Ƅ|^2$ over $\{t\} \times D \times \{0\}$ which is

$$\frac{1}{\varepsilon} ¥ \int_{\{t\} \times D \times \{0\}} |Ƅ|^2 .$$

(6.43)

Since $|Ƅ|(t,q,0) > \delta ¥$ (by assumption), the preceding four paragraphs lead to this:

$$\delta ¥ \leq c_0 c_ɜ \frac{1}{\gamma t} ( \int_{\{t\} \times D \times \{0\}} |Ƅ|^2 )^{1/2} + c_0 \varepsilon (1 + \frac{1}{\varepsilon^2} \int_{\{t\} \times D \times \{0\}} |Ƅ|^2 ) ¥ .$$

(6.44)

This is an inequality with ¥ appearing on both sides; it appears on the right with a factor δ, and it appears on the left with a factor that depends on the choice of ε. Note in particular that if $\varepsilon < c_0^{-1} \delta$ and if $\int_{\{t\} \times D \times \{0\}} |Ƅ|^2 \leq \varepsilon^2$, then the factor that multiplies ¥ on the right hand side of (6.44) will be less than $\frac{1}{2}\delta$. And, if this is the case, then (6.44) bounds the supremum of $|\mathfrak{b}|$ on $\{t\}\times D\times\{0\}$ by $c_0 c_{\mathfrak{z}} \frac{1}{\delta\gamma t}$.

### f) When the very small t or large t integrals of $|\beta|^2 + |\mathfrak{b}|^2$ on $\{t\}\times \mathbb{R}^2 \times \{0\}$ are small.

The central result in this section is a proposition that says in effect that if either the very small t or very large t integrals of $|\beta|^2 + |\mathfrak{b}|^2$ on $\{t\}\times\mathbb{R}^2\times\{0\}$ are small (less than a fixed positive bound), then these integrals limit to zero almost like $t^2$ as $t\to 0$, or almost like $\frac{1}{t^2}$ as $t\to\infty$ as the case may be.

**Proposition 6.8**: *Given positive numbers $\mathfrak{z}$, $\delta$ and $r$, there exists $\kappa > 1$ with the following significance: Let $(A,\mathfrak{a})$ denote a solution to (1.4) that obeys the first two bullets of (6.1) with the given $\mathfrak{z}$, $\delta$, and $r$. If $\int_{\{t\}\times\mathbb{R}^2\times\{0\}} (|\beta|^2 + |\mathfrak{b}|^2) \leq \frac{1}{\kappa}$ for all $t < t_0$ for some $t_0 > 0$, then*

$$\lim_{t\to 0} t^{-2+\varepsilon} \int_{\{t\}\times\mathbb{R}^2\times\{0\}} (|\beta|^2 + |\mathfrak{b}|^2) = 0 \quad \text{and} \quad \lim_{t\to 0} t^{\varepsilon}(|\beta|^2 + |\mathfrak{b}|^2) = 0$$

*for any $\varepsilon > 0$. By the same token, if $\int_{\{t\}\times\mathbb{R}^2\times\{0\}} (|\beta|^2 + |\mathfrak{b}|^2) \leq \frac{1}{\kappa}$ for all $t \in [t_0, \infty)$, then*

$$\lim_{t\to\infty} t^{2-\varepsilon} \int_{\{t\}\times\mathbb{R}^2\times\{0\}} (|\beta|^2 + |\mathfrak{b}|^2) = 0 \quad \text{and} \quad \lim_{t\to\infty} t^{4-\varepsilon}(|\beta|^2 + |\mathfrak{b}|^2) = 0$$

*for any $\varepsilon > 0$.*

The proof is given directly. Note that some of the lemmas along the way to the proof are used later also.

***Proof of Proposition 6.8***: The idea behind the proof is this: If the integral of $|\beta|^2 + |\mathfrak{b}|^2$ on a given constant $(t, x_3)$ slice of $(0,\infty)\times\mathbb{R}^2\times\mathbb{R}$ is sufficiently small (meaning less than some positive, t-independent upper bound), then $t|\beta|$ and $t|\mathfrak{b}|$ are correspondingly small (see Lemma 6.5). The plan is to prove that if $t|\beta|$ and $t|\mathfrak{b}|$ are suitably small where $t \leq t_0$ (or $t \geq t_0$), then $(A,\mathfrak{a})$ will be close to one of the model solutions from Section 1c where $t < t_1$ with $t_1 \ll t_0$ (or $t > t_1$ with $t_1 \gg t_0$). And, if $(A,\mathfrak{a})$ is suitably close to one of the model solutions for very small t (or very large t), then $t\alpha$ will be uniformly bounded away from zero for these values of t. This in turn implies (via Lemma 6.4) that the integrals of $|\beta|^2 + |\mathfrak{b}|^2$ on the constant $t \ll t_1$ and constant $x_3$ slices of $(0,\infty)\times\mathbb{R}^2\times\mathbb{R}$ limit to zero as $t\to 0$ (or as $t\to\infty$) which implies that $t|\beta|$ and $t|\mathfrak{b}|$ also limit to zero (by Corollary 6.6). The latter fact will be used to prove that $(A,\mathfrak{a})$ actually converges to one of the model solution as $t\to 0$ (or $t\to\infty$). And, if $(A,\mathfrak{a})$ does indeed converge as $t\to 0$ (or t → ∞) to one of the model solutions, then the t → 0 (or t → ∞) limit tα is no smaller than $-\frac{1}{2}$ because no model solution version of tα is greater than this. Supposing that tα does converge everywhere to something no less than $-\frac{1}{2}$, then the assertion of the lemma follows from Lemma 6.4 and Corollary 6.6.

The details of the proof are presented directly in six parts.

*Part 1*: This first part of the proof constitutes a digression of sorts to state and then prove a lemma about sequences of solutions to (1.4). This upcoming lemma assumes somewhat more than Lemma 2.3 to draw a correspondingly stronger conclusion. The lemma introduces terminology: A pair (A, 𝔞) is said to be *decomposable* when it can be written in terms of a data set $(\hat{A}, \flat, \alpha, \beta, \varphi, \sigma)$ with $\mathfrak{a}_3$ having a smooth decomposition as depicted in (3.7) and with φ being a section of σ's version of the line bundle $\mathcal{L}^+$. In addition, A can be written as $A = \hat{A} + \flat$ so that $\nabla_{\hat{A}} \sigma = 0$ and so that $\flat$ is described by (3.8).

**Lemma 6.9**: *Let $\{(A^n, \mathfrak{a}^n)\}_{n \in \mathbb{N}}$ denote a sequence of decomposable solutions to (1.4) with the following property: There is a positive numbers $\mathfrak{z}$ and a sequence of positive numbers $\{\mathcal{Z}_k\}_{k \in \mathbb{N}}$ such that if $n \in \mathbb{N}$, then*

- $|\mathfrak{a}^n| \leq \frac{\mathfrak{z}}{t}$ *on* $(0, \infty) \times \mathbb{R}^2 \times \mathbb{R}$.
- *Supposing that $n \geq k$, then the $(A^n, \mathfrak{a}^n)$ version of $\flat$ obeys $|\flat| \leq \frac{\mathcal{Z}_k}{t}$ where $t \in [\frac{1}{k}, k]$ and $|z| \leq k$ on $(0, \infty) \times \mathbb{R}^2 \times \mathbb{R}$.*

*Granted these assumptions, there exists a decomposable solution to (1.4) (to be denoted by $(A^\infty, \mathfrak{a}^\infty)$), a subsequence of $\mathbb{N}$ (denoted by $\Lambda$), and a corresponding sequence of automorphisms of P (to be denoted by $\{g_n\}_{n \in \Lambda}$); and these are such that the $\Lambda$-indexed sequence whose n'th term is the $(A^n, \mathfrak{a}^n)$ version of the set $(\hat{A}, \flat, \alpha, \beta, \sigma, \varphi)$ converges in the $C^\infty$-topology on compact subsets of $(0, \infty) \times \mathbb{R}^2 \times \mathbb{R}$ to the $(A^\infty, \mathfrak{a}^\infty)$ version of $(\hat{A}, \flat, \alpha, \beta, \sigma, \varphi)$. Moreover, if all of the $n \in \mathbb{N}$ version of φ obey the second bullet of (6.1) with the same δ and r, then so does the $(A^\infty, \mathfrak{a}^\infty)$ version of φ. In this event, all of the $n \in \Lambda$ versions of φ and the $(A^\infty, \mathfrak{a}^\infty)$ version will vanish only at the z = 0 point in each constant $(t, x_3)$ slice of $(0, \infty) \times \mathbb{R}^2 \times \mathbb{R}$; and the vanishing degrees of the $n \in \Lambda$ versions at this point will have an $n \to \infty$ limit which is the vanishing degree at that point of the $(A^\infty, \mathfrak{a}^\infty)$ version of φ.*

*Proof of Lemma 6.9*: The existence of a limit and the convergence in the indicated manner follows from the a priori bounds that are given by Lemma 6.2. (These bounds hold because of the n-independent bounds on $|\mathfrak{a}^n|$ and on the $(A^n, \mathfrak{a}^n)$ version of $|\flat|$.) With regards to convergence of the sequence of $(A^n, \mathfrak{a}^n)$ versions of σ, keep in mind that the $(A^n, \mathfrak{a}^n)$ version of $\mathfrak{b}$ (which is determined by the corresponding version of $\mathfrak{b}$) is given in terms of the corresponding version of $\sigma$ by the formula $\mathfrak{b} = -\frac{1}{4}[\sigma, \nabla_A \sigma]$. Thus, convergence of the sequence of $\mathfrak{b}$'s implies convergence of the covariant derivatives of the sequence of $\sigma$'s; and likewise for sequences of higher order covariant derivatives.

The assertion that the $(A^\infty, \mathfrak{a}^\infty)$ version of $\varphi$ obeys the second bullet of (6.1) if the $n \in \mathbb{N}$ versions of $\varphi$ do follows from what is said by Lemma 2.3. Assuming that this bullet is obeyed by the various $n \in \mathbb{N}$ versions of $\varphi$, then those version of $\varphi$ (and the $(A^\infty, \mathfrak{a}^\infty)$ version) can only vanish at the origin in each constant $(t, x_3)$ slice of $(0, \infty) \times \mathbb{R}^2 \times \mathbb{R}$. (A zero at any other point would have to foul the second bullet of (6.1). See Proposition 3.1)

Consider now the assertion about the degree of vanishing of the $(A^\infty, \mathfrak{a}^\infty)$ version of $\varphi$ at the $z = 0$ point in each constant $(t, x_3)$ slice of $(0, \infty) \times \mathbb{R}^2 \times \mathbb{R}$. In this regard: Lemma 5.6 supplies an upper bound to the degrees of vanishing of the $n \in \mathbb{N}$ versions of $\varphi$. Therefore, $\Lambda$ can always be *chosen* so that these degrees are independent of n. But, as explained in the four steps that follow, the $C^\infty$-convergence dictated by the first parts of the lemma imply that the $n \in \Lambda$ versions of the vanishing degree must be independent of n when n is sufficiently large, and that this asymptotic value is the vanishing degree of of the $(A^\infty, \mathfrak{a}^\infty)$ version of $\varphi$.

Step 1: Suppose that $\sigma$ and $\sigma'$ are two unit length sections of ad(P) on a given subset in $(0, \infty) \times \mathbb{R}^2 \times \mathbb{R}$. They determine corresponding versions of the line bundle $\mathcal{L}^+$ (the $\sigma'$ version is denoted by $\mathcal{L}^{+'}$) which are complex lines subbundles of $\mathrm{ad}(P)_\mathbb{C}$. If $\langle \sigma \sigma' \rangle \neq -1$, then the orthogonal projection in $\mathrm{ad}(P)_\mathbb{C}$ from the line $\mathcal{L}^{+'}$ to the line $\mathcal{L}^+$ defines a $\mathbb{C}$-linear isomorphism between them. This happens in particular when if $\sigma$ and $\sigma'$ are nearly the same. Supposing that $\langle \sigma \sigma' \rangle \neq -1$, then the orthogonal projection from $\mathcal{L}^{+'}$ to $\mathcal{L}^+$ is the restriction to $\mathcal{L}^{+'}$ of the endomorphism of $\mathrm{ad}(P)_\mathbb{C}$ that can be defined by the rule $\eta \to \lambda \langle \lambda^* \eta \rangle$ with $\lambda$ being any unit norm section of $\mathcal{L}^+$.

Step 2: Now suppose that the subset where $\sigma$ and $\sigma'$ are defined is a disk centered at the $z = 0$ point in some constant $(t, x_3)$ slice of $(0, \infty) \times \mathbb{R}^2 \times \mathbb{R}$. If $\varphi'$ is a section of $\mathcal{L}^{+'}$ that vanishes only at the $z = 0$ point in this disk, then $\lambda \langle \lambda^* \varphi' \rangle$ will be a section of $\mathcal{L}^+$ that likewise vanishes only at the $z = 0$ point (assuming that $\langle \sigma \sigma' \rangle > -1$ on the whole disk). Moreover, the degree of vanishing of $\lambda \langle \lambda^* \varphi' \rangle$ at this point (as a section of $\mathcal{L}^+$) is the same as that of $\varphi'$ (as a section of $\mathcal{L}^{+'}$).

Step 3: Let $\varphi$ denote a section of $\mathcal{L}^+$ on the disk in question and assume that the only zero of $\varphi$ on this disk is also at the $z = 0$ point. Then $\langle \varphi^* \varphi' \rangle$ is a map from the disk to $\mathbb{C}$ that vanishes only at the origin. Thus, it maps the complement of the $z = 0$ point in the disk to $\mathbb{C}-\{0\}$. As such, its winding number is equal to the degree of vanishing of $\varphi´$ at the origin minus that of $\varphi$. Keep in mind in this regard that $\langle \varphi^*\varphi´\rangle$ can be written as $\langle \varphi^*\lambda\rangle\langle \lambda^*\varphi´\rangle$ which implies that $\langle \varphi^*\varphi´\rangle$ does indeed compute the difference between the degrees of vanishing of the two sections ($\varphi$ and $\lambda\langle\lambda^*\varphi´\rangle$) of $\mathcal{L}^+$.

In particular, if C is a constant $|z| > 0$ circle in the disk and if $|\varphi´ - \varphi| < \frac{1}{100}|\varphi|$ on C, then $\langle \varphi^*\varphi´\rangle$ on C will have strictly positive real part. And, if this is the case, then the vanishing degrees of $\varphi$ and $\varphi´$ at the $z = 0$ point will be the same since the map $\langle \varphi^*\varphi´\rangle$ from C to $\mathbb{C}-\{0\}$ will have vanishing winding number.

Step 4: To bring to bear what is said by Steps 1-3, let $\{g_n\}_{n\in\Lambda}$ denote the automorphisms of P that are supplied with $(A^\infty, \mathfrak{a}^\infty)$ by Lemma 6.9. Let $\{\sigma^n\}_{n\in\Lambda}$ denote the $g_n$ pullback of the $(A^n, \mathfrak{a}^n)$ version of $\sigma$. Likewise, let $\varphi^n$ denote the $g_n$ pullback of the $(A^n, \mathfrak{a}^n)$ version of $\varphi$. Meanwhile, let $\sigma^\infty$ denote the $(A^\infty, \mathfrak{a}^\infty)$ version of $\sigma$ and let $\varphi^\infty$ denote the corresponding version of $\varphi$. Since the sequence $\{\sigma^n\}_{n\in\Lambda}$ converges pointwise to $\sigma$ on the radius $2r$ disk about the $z = 1$ point in the $(t = 1, x_3 = 0)$ slice of $(0,\infty)\times\mathbb{R}^2\times\mathbb{R}$, the map $\langle(\varphi^\infty)^*\varphi^n\rangle$ from the $|z| = r$ circle in this disk to $\mathbb{C}-\{0\}$ computes the difference between the vanishing degrees at the $z = 0$ point of $\varphi^n$ and $\varphi^\infty$. By virtue of the fact that $|\varphi^n|$ and $|\varphi^\infty|$ are bounded away from zero on this circle, and by virtue of the fact that $\{\varphi^n\}_{n\in\Lambda}$ converges to $\varphi^\infty$ on this circle, it follows from what is said in Step 2 that $\langle(\varphi^\infty)^*\varphi^n\rangle$ has strictly positive real part on this circle when n is large. As a consequence $\varphi^n$ and $\varphi^\infty$ have the same vanishing degree at the $z = 0$ point in the disk when n is large.

*Part 2*: This part of the proof states and then proves a lemma that is for the most part just a corollary to Lemma 6.9. (This lemma is used in subsequent sections also.)

**Lemma 6.10**: *Suppose that* $(A, \mathfrak{a})$ *is a solution to (1.4) that obeys the three bullets in (6.1). Assume in addition that* $|\mathfrak{b}| \leq \frac{\mathcal{Z}}{t}$ *where either* $t \leq 1$ *or* $t \geq 1$ *on* $(0,\infty)\times\mathbb{R}^2\times\mathbb{R}$ *with* $\mathcal{Z}$ *being independent of* t. *Let m denote the degree of vanishing of* $\varphi$ *at the z = 0 point in each constant* $(t, x_3)$ *slice of* $(0,\infty)\times\mathbb{R}^2\times\mathbb{R}$. *There exists* $\varepsilon_* > 0$ *such that* $|\varphi| \geq \varepsilon_* \frac{|z|^m}{t^{m+1}}$ *on the part of* $(0,\infty)\times\mathbb{R}^2\times\mathbb{R}$ *where* $|z| \leq rt$ *and either* $t \leq 1$ *or* $t \geq 1$ *as the case may be.*

*Proof of Lemma 6.10*: The proof has two steps.

Step 1: This step explains why the asserted lower bound follows if there exists $\mu > 0$ such that $t^{m+1}|(\nabla_A^{\otimes m})\varphi| \geq \mu$ at the $z = 0$ point in each constant $(t, x_3)$ slice of $(0,\infty)\times\mathbb{R}^2\times\mathbb{R}$ where either $t \leq 1$ or $t \geq 1$. To start, note that in a $(t, x_3)$ slice of

$(0,\infty) \times \mathbb{R}^2 \times \mathbb{R}$ where this bound holds at $z = 0$, Taylor's theorem with remainder (with the bounds from Lemma 2.2) can be used to see that $|\varphi|$ on that slice will be greater than $\frac{1}{2}\mu \frac{|z|^m}{t^{m+1}}$ where $|z| \leq c_{\mathfrak{z}}^{-1}\mu t$. Meanwhile, $|\varphi|$ is greater than $\frac{\delta}{t}$ where $|z| \geq rt$ by assumption; and $|\varphi|$ must be greater than a fixed multiple of $\frac{1}{t}$ (which is determined by $\mathfrak{z}$, $\delta$, $r$ and $\mu$) where $|z|$ is between $c_{\mathfrak{z}}^{-1}\mu t$ and $rt$ (supposing that $c_{\mathfrak{z}}^{-1}\mu \leq r$). The existence of such a lower bound follows from the top bullet in (3.9), the bound $|\mathfrak{a}_{\mathfrak{z}}| \leq \frac{\mathfrak{z}}{s}$ for $s \leq t$ and the various $s \leq t$ versions of the second bullet of (6.1). Indeed, if $r$ is given and less than $r$, then $|\varphi|$ must be greater than or equal to $\frac{\delta}{s}$ where $s \leq \frac{r}{r}t$ and $|z| = rt$ (which is $rs$) so as not to foul the second bullet of (6.1) at time $s$. Granted that lower bound, then the top bullet of (3.9) and the top bullet of (6.1) lead to the lower bound $|\varphi| \geq (\frac{r}{r})^{\mathfrak{z}-1}\frac{\delta}{t}$ at time $t$ and where $|z| = rt$.

Step 2: With Step 1's lesson understood, suppose that there is no positive lower bound to the norm of $t^{m+1}|(\nabla_A^{\otimes m})\varphi|$ at the $z = 0$ locus where $t \leq 1$ (or where $t \geq 1$). The plan is to derive nonsense from this assumption. In particular, if there is no such positive lower bound, then there is a sequence $\{t_n\}_{n \in \mathbb{N}}$ (all less than 1 or all greater than 1) where the corresponding versions of $t^{m+1}|(\nabla_A^{\otimes m})\varphi|$ at the $z = 0$ point limit to zero as $n \to \infty$. Pulling back $(A, \mathfrak{a})$ by the coordinate rescaling diffeomorphisms $(t, z, x_3) \to (t_n t, t_n z, t_n x_3)$ constructs a sequence $\{(A^n, \mathfrak{a}^n)\}_{n \in \mathbb{N}}$ that can be used for the input sequence in Lemma 6.9. This sequence would converge in the manner dictated by Lemma 6.9 to a solution to (1.4) whose corresponding $\varphi$ would vanish to degree greater than $m$ at the $z = 0$ point in the $(t=1, x_3 = 0)$ slice. This is because it and its covariant derivatives to order $m$ or less all vanish at this $z = 0$ point. The greater than $m$ degree of vanishing is the desired nonsense because it flat-out contradicts what is said at the end of Lemma 6.9.

*Part 3*: This part of the proof contains a lemma that can be viewed in part as an application of Lemma 6.10 (so it is also something of a corollary to Lemma 6.9). This lemma concerns solutions with $\beta \equiv 0$ and $\mathfrak{b} \equiv 0$ (so it harks back to Section 4).

**Lemma 6.11**: *Any solution to (1.4) that obeys the conditions in the first two bullets of (6.1) with vanishing $\beta$ and $\mathfrak{b}$ is one of the model solutions from Section 1c.*

*Proof of Lemma 6.11*: Let $(A, \mathfrak{a})$ denote a solution to (1.4) of requisite sort. If the corresponding $\varphi$ vanishes on the constant $(t, x_3)$ slices of $(0, \infty) \times \mathbb{R}^2 \times \mathbb{R}$, then this can happen only at the $z = 0$ point (because of the second bulet in (6.1)). In any event, let $m$ denote the degree of vanishing of $\varphi$ at these $z = 0$ points.

Write $|\varphi|$ as $e^{2w}|\varphi^{(m)}|$ where $\varphi^{(m)}$ is the version of $\varphi$ from the model solution in Section 1c whose vanishing degree at the $z = 0$ point in any constant $(t, x_3)$ slice of $(0, \infty) \times \mathbb{R}^2 \times \mathbb{R}$ is the same as that of $\varphi$. There are two important facts to note about w: First, it obeys the following version of (4.37):

$$-(\tfrac{\partial^2}{\partial t^2} + \tfrac{\partial^2}{\partial x_1^2} + \tfrac{\partial^2}{\partial x_2^2})w + (e^{4w} - 1)|\varphi^{(m)}|^2 = 0 \ .$$

(6.47)

Second, w is bounded from above and from below. This is proved momentarily. The next paragaph uses these facts to finish the proof of the lemma.

As just noted, w obeys (6.47) which is an equation that has a corresponding maximum principle. In particular, this maximum principle implies that w has neither positive local maxima nor negative local minima. Granted that |w| is a priori bounded and <u>not</u> identically zero, then a sequence of solutions to (1.4) can be constructed as pullbacks of combinations of translations of the $\mathbb{R}^2$ factor of $(0, \infty) \times \mathbb{R}^2 \times \mathbb{R}$ and coordinate rescalings to obtain (via Lemma 6.9) a new version of w that obeys (6.47) (using either $|\varphi^{(m)}|$ or using $\tfrac{1}{\sqrt{2}t}$ in lieu of $|\varphi^{(m)}|$) whose value at the point $(t = 1, z = 0)$ is a local maximum or minimum; either the supremum or infimum of the original w. (The translations and coordinate rescalings that define the input sequence for Lemma 6.9 are chosen so that the values of corresponding sequence of w's at $(t = 1, z = 0)$ approaches either the supremum or infimum of the original w.) Thus, the new version of w is is non-zero, it obeys a version of (6.47) and it takes a local positive maximum or negative minimum where $t = 1$ and $z = 0$. This last conclusion is the desired nonsense. It can be avoided only if w is identically zero which implies that $(A, \mathfrak{a})$ is the integer $m$ model solution from Section 1c. (Supposing that $w \equiv 0$, then $|\varphi| = |\varphi^{(m)}|$ which implies that $\alpha$ is model version $\alpha^{(m)}$ because of (3.9a); and it implies that $\langle \sigma E_A \rangle$ and $\langle \sigma B_A \rangle$ are the same as their model counterparts because of (3.10). These facts imply that $(A, \mathfrak{a})$ is $\mathrm{Aut}(P)$ equivalent to the integer $m$ version of what is depicted in (1.7) and (1.8).)

To prove that |w| is bounded, note first that w is bounded from above and below where $|z| > rt$ by virtue of (6.46). To see that w is bounded from above where $|z| \le rt$, invoke Taylor's theorem with remainder on any given constant $(t, x_3)$ slice with the fact that $\varphi$ vanishes to degree $m$ at the $z = 0$ point (and Lemma 2.2's bound $|(\nabla_A^{\otimes m})\varphi| \le c_\mathfrak{z} \tfrac{1}{t^{m+1}}$) to conclude that $|\varphi| \le c_\mathfrak{z} \tfrac{|z|^m}{t^{m+1}}$ on the slice. Meanwhile, $|\varphi^{(m)}|$ is observedly greater than $c_0^{-1} \tfrac{|z|^m}{t^{m+1}}$ on this slice (look at the second bullet in (1.7).) The assertion that w is bounded from below where $|z| \le rt$ follows from Lemma 6.10's lower bound $|\varphi| \ge \varepsilon \tfrac{|z|^m}{t^{m+1}}$ and the observed upper bound $|\varphi^{(m)}| \le c_0(m+1) \tfrac{|z|^m}{t^{m+1}}$).

*Part 4*: The lemma that follows is the central result for this part of the proof.

**Lemma 6.12**: *Given positive numbers $\mathfrak{z}$, $\delta$, $r$ and $\Xi$, and also a number $\mu \in (0,1)$, there exists $r_\mu > 1$ with the following significance: Let $(A, \mathfrak{a})$ denote a solution to (1.4) that obeys the three bullets of (6.1) with the given $\mathfrak{z}$, $\delta$, $r$ and $\Xi$. The $(A, \mathfrak{a})$ version of the function $\alpha$ obeys $\alpha < -\frac{(1-\mu)}{2t}$ on the $|z| \geq r_\mu t$ part of $(0, \infty) \times \mathbb{R}^2 \times \mathbb{R}$.*

*Proof of Lemma 6.12*: The plan is to assume that there exists $\mu$ with no corresponding $r_\mu$ so as to derive a nonsensical conclusion. If there is no $r_\mu$ for a given choice of $\mu$, then there exist the following two sequences: First, a sequence $\{(A_\ddagger^n, \mathfrak{a}_\ddagger^n)\}_{n \in \mathbb{N}}$ of solutions to (1.4) that are described by the three bullets of (6.1) with the given $\mathfrak{z}$, $\delta$, $r$ and $\Xi$. Second, a sequence $\{(t_n, z_n)\}_{n \in \mathbb{N}} \subset (0, \infty) \times \mathbb{R}^2$ with any given $n \in \mathbb{N}$ member obeying $|z_n| \geq n t_n$. These two sequences are such that any given $n \in \mathbb{N}$ version of $\alpha$ (denoted by $\alpha_\ddagger^n$) obeys $\alpha_\ddagger^n \geq -\frac{(1-\mu)}{2t_n}$ at the point $(t_n, z_n, 0)$. The two steps that follow derive nonsense from this.

    Step 1: Given $n \in \mathbb{N}$, define a new solution to (1.4) to be denoted by $(A^n, \mathfrak{a}^n)$ to be the pull-back of $(A_\ddagger^n, \mathfrak{a}_\ddagger^n)$ by the translation/coordinate rescaling diffeomorphism of $(0, \infty) \times \mathbb{R}^2 \times \mathbb{R}$ that is defined by the rule $(t, z, x_3) \to (t_n t, t_n z + z_n, t_n x_3)$. The corresponding version of $\alpha$ for this new pair (denoted by $\alpha^n$) obeys $\alpha^n \geq -\frac{(1-\mu)}{2}$ at the point $(t = 1, z = 0, x_3 = 0)$. Meanwhile, the corresponding version of $\varphi$ (denoted by $\varphi^n$) obeys $|\varphi^n| \geq \frac{\delta}{t}$ where $|z| < \frac{1}{100} n$ and $t \leq \frac{1}{2r} n$ due to the second bullet in (6.1). In addition, the versions of $\beta$ and $\mathfrak{b}$ for this new pair obey $|\beta^n| + |\mathfrak{b}^n| \leq \mathcal{Z}(\frac{1}{n})^{1/4} \frac{1}{t}$ on this domain with $\mathcal{Z}$ being independent of $t$, $z$ and $n$. (See the third bullet of Lemma 6.3.) Finally, note that the integral of $|F_{A^n}|^2$ over the $|z| \leq \frac{1}{100} n$ part of $(0, \infty) \times \mathbb{R}^2 \times \mathbb{R}$ is bounded by $\frac{\Xi}{n}$ by virtue of the third bullet of (6.1).

    Step 2: Lemma 6.9 can be invoked using the input sequence $\{(A_n, \alpha_n)\}_{n \in \mathbb{N}}$. Let $(A^\infty, \mathfrak{a}^\infty)$ denote a resulting limit pair. This limit pair has $|\varphi^\infty| \geq \frac{\delta}{t}$ on the whole of $(0, \infty) \times \mathbb{R}^2 \times \mathbb{R}$ whereas $\beta^\infty \equiv 0$ and $\mathfrak{b} \equiv 0$ on the whole of $(0, \infty) \times \mathbb{R}^2 \times \mathbb{R}$. The curvature of the connection $A^\infty$ must also vanish identically. The final crucial point is that $\alpha^\infty$ at the $(t = 1, z = 0)$ is no smaller than $-\frac{(1-\mu)}{2}$. (These conclusions all follow from the manner of convergence that is dictated by Lemma 6.8.)

    Granted the preceding, it then follows from Lemma 6.11 that this pair $(A^\infty, \mathfrak{a}^\infty)$ must be the Nahm pole solution. That conclusion is the desired nonsense because the version of $\alpha$ for the Nahm pole solution is equal to $-\frac{1}{2}$ where $t = 1$ and $z = 0$.

*Part 5*: Lemma 6.12 describes the function $\alpha$ where $\frac{|z|}{t}$ is large. The next lemma is a step towards describing $\alpha$ where $\frac{|z|}{t}$ is not especially large but t is very small or large.

**Lemma 6.13**: *Given positive numbers $\mathfrak{z}$, $\delta$, $r$ and $\Xi$, and also a number $\mu \in (0,1)$, there exists $\kappa_\mu > 1$ with the following significance: Let $(A, \mathfrak{a})$ denote a solution to (1.4) that obeys the three bullets of (6.1) with the given $\mathfrak{z}$, $\delta$, $r$ and $\Xi$. Suppose in addition that*

$$\int_{\{t\}\times\mathbb{R}^2\times\{0\}} (|\beta|^2 + |\mathfrak{b}|^2) \leq \tfrac{1}{\kappa_\mu}$$

*when t is positive and sufficiently small (or sufficiently large). Then there exists postive $t_0$ such that $\alpha < -\frac{(1-\mu)}{2t}$ where $t < t_0$ (or where $t > t_0$ as the case may be).*

*Proof of Lemma 6.13*: The small t assertion is proved in the three steps that follow. The proof of the large t assertion is not given because it differs only cosmetically from what is said below for the case of small t.

<u>Step 1</u>: Assume to the contrary that here exists $\mu \in (0,1)$ without a corresponding $\kappa_\mu$. There is in this case three sequences, the first being a sequence $\{(A_\ddagger^n, \mathfrak{a}_\ddagger^n)\}_{n\in\mathbb{N}}$ of solutions to (1.4) that all obey (6.1) with the given values of $\mathfrak{z}$, $\delta$, $r$ and $\Xi$. The second is a sequence $\{(t_{0n}, t_{1n}, z_n)\}_{n\in\mathbb{N}}$ with the n'th member being a pair in $(0,1)$ and a point in $\mathbb{R}^2$. Moreover, the n'th member $(t_{0n}, t_{1n}, z_n)$ is such that $t_{0n} < \frac{1}{n}$ and $t_{1n} < \frac{1}{n} t_{0n}$ and $|z_n| \leq r_{\mu/2} t_{1n}$ with $r_{\mu/2}$ from Lemma 6.12. These sequences have the following additional properties: Supposing that $n \in \mathbb{N}$, then the integral of the $(A_\ddagger^n, \mathfrak{a}_\ddagger^n)$ version of $|\beta|^2 + |\mathfrak{b}|^2$ on the constant $(t = t_n, x_3 = 0)$ slice of $(0,\infty)\times\mathbb{R}^2\times\mathbb{R}$ is smaller than $\frac{1}{n}$. Meanwhile, the $(A_\ddagger^n, \mathfrak{a}_\ddagger^n)$ version of $\alpha$ is greater than $-\frac{(1-\mu)}{2t_{1n}}$ where $t = t_{1n}$ and $z = z_n$.

Keep in mind for what follows that the properties listed above with Lemma 6.5 imply this: There is a positive number (to be denoted here by K) which depends only on $\mathfrak{z}$, $\delta$, $r$ and $\Xi$ such that if $n \in \mathbb{N}$ and $n \geq K$, then the $(A_\ddagger^n, \mathfrak{a}_\ddagger^n)$ versions of $|\beta|$ and $|\mathfrak{b}|$ are bounded pointwise by $\frac{K}{nt}$ on the $t \leq t_{0n}$ part of $(0,\infty)\times\mathbb{R}^2\times\mathbb{R}$.

<u>Step 2</u>: Given $n \in \mathbb{N}$, define $(A^n, \mathfrak{a}^n)$ to be the pull-back of $(A_\ddagger^n, \mathfrak{a}_\ddagger^n)$ by the coordinate rescaling diffeomorphism given by the rule $(t, z, x_3) \to (t_{1n}t, t_{1n}z, t_{1n}x_3)$. Use the corresponding sequence $\{(A^n, \mathfrak{a}^n)\}_{n\in\mathbb{N}}$ as input for Lemma 6.9 which then outputs a limit pair $(A^\infty, \mathfrak{a}^\infty)$. Given the manner of $C^\infty$ convergence dictated by Lemma 6.9 and given the properties listed in Step 1, this limit pair $(A^\infty, \mathfrak{a}^\infty)$ must obey the following on the whole of $(0,\infty)\times\mathbb{R}^2\times\mathbb{R}$:

- $|\mathfrak{a}^\infty| \leq \frac{\mathfrak{z}}{t}$ *on the whole of* $(0,\infty) \times \mathbb{R}^2 \times \mathbb{R}$.
- $|\varphi^\infty| \geq \frac{\delta}{t}$ *where* $|z| \geq rt$ *on* $(0,\infty) \times \mathbb{R}^2 \times \mathbb{R}$.
- $\alpha^\infty \leq -\frac{(1-\frac{1}{2}\mu)}{2t}$ *where* $|z| \geq r_{1/2}t$ *on* $(0,\infty) \times \mathbb{R}^2 \times \mathbb{R}$.
- $\beta^\infty \equiv 0$ *and* $\mathfrak{b}^\infty \equiv 0$
- $F_{A^\infty} \equiv 0$.
- $\alpha^\infty \geq -\frac{(1-\mu)}{2}$ *at a* $t = 1$ *point*.

(6.46)

Step 3: The first and second bullets and fourth bullets of (6.46) imply via Lemma 6.11 that $(A^\infty, \mathfrak{a}^\infty)$ is one of the model solutions from Section 1c; and then the vanishing of all of the curvature (the fifth bullet in (6.46)) implies that it must be the $m = 0$ solution, a Nahm pole. Whether the $m = 0$ model or not, it then follows that $\alpha^\infty$ is nowhere greater than $-\frac{1}{2t}$ (see the second bullet of (2.1)). This $\alpha^\infty \leq -\frac{1}{2t}$ bound is the desired nonsense because it explicitly contradicts the fourth bullet in (6.46).

*Part 6*: Lemma 6.13 supplies a positive $\kappa_*$ that depends only on $\mathfrak{z}, \delta, r$ and $\Xi$ such that the following is true: If $(A, \mathfrak{a})$ is a solution to (1.4) that obeys the three bullets in (6.1) with the given $\mathfrak{z}, \delta, r$ and $\Xi$, and if the $\{t\} \times \mathbb{R}^2 \times \{0\}$ integrals of $|\beta|^2 + |\mathfrak{b}|^2$ are less than $\frac{1}{\kappa_*}$ for all t sufficiently small (or for all t sufficiently large) then the corresponding version of the function $\alpha$ will be less than $-\frac{1}{4t}$ where t is sufficiently small (or sufficiently large as the case may be). In this event, Lemma 6.4 can be brought to bear to see that the $t \to 0$ limit of the $\{t\} \times \mathbb{R}^2 \times \{0\}$ integrals of $|\beta|^2 + |\mathfrak{b}|^2$ is zero (or the $t \to \infty$ limit of these integrals). This implies via Lemma 6.13 that the function $\alpha$ is less than any given $\mu > 0$ version of $-\frac{(1-\mu)}{2t}$ on the whole of $\{t\} \times \mathbb{R}^2 \times \mathbb{R}$ when t is sufficiently small (or sufficiently large). Granted this last fact, then the assertion in Proposition 6.8 follows directly from Lemma 6.4 and Corollary 6.6.

**g) A special case of Theorem 2**

The following proposition is a special case of what is asserted by Theorem 2. It says in effect that a solution to (1.4) obeying the constraints in (6.1) is necessarily one of the model solutions if both the very small t and very large t integrals of $|\beta|^2 + |\mathfrak{b}|^2$ on constant $(t, x_3)$ slices of $(0, \infty) \times \mathbb{R}^2 \times \mathbb{R}$ are small.

**Proposition 6.14**: *Given positive numbers $\mathfrak{z}$, $\delta$, $r$ and $\Xi$, there exists $\kappa > 1$ with the following significance: Let $(A, \mathfrak{a})$ denote a solution to (1.4) that obeys the three bullets of (6.1) with the given $\mathfrak{z}$, $\delta$, $r$ and $\Xi$. Suppose in addition that*

$$\int_{\{t\} \times \mathbb{R}^2 \times \{0\}} (|\beta|^2 + |\mathfrak{b}|^2) \le \tfrac{1}{\kappa}$$

*when $t$ is positive and sufficiently small <u>and</u> when $t$ is sufficiently large. Then $(A, \mathfrak{a})$ is one of the model solutions from Section 1c.*

***Proof of Proposition 6.14***: The number $\kappa$ is chosen so that Proposition 6.8's a priori bounds can be used where $t$ is both very small and very large. The proof that $(A, \mathfrak{a})$ is a model solution when Proposition 6.8 is in play has five parts. With regards to notation: What is denoted by $\mathcal{Z}_\varepsilon$ in the proof is a number that is greater than 1 which depends on a given choice of another positive number $\varepsilon$. This $\mathcal{Z}_\varepsilon$ can depend on $(A, \mathfrak{a})$ also. The important point is that it has no dependence on any given point $(t, z, x_3)$ in $(0, \infty) \times \mathbb{R}^2 \times \mathbb{R}$. Its value can be assumed to increase between successive appearances.

*Part 1*: Take the inner product of the equations in (3.11) with $\varphi^*$ and then use (3.9) to write the resulting equations on the $z \ne 0$ part of $(0, \infty) \times \mathbb{R}^2 \times \mathbb{R}$ as follows:

- $\frac{\partial}{\partial t} \langle \varphi^* \beta \rangle = -i 2 \frac{\partial}{\partial z} \langle \varphi^* \mathfrak{b} \rangle$ .
- $\frac{\partial}{\partial t} ( \frac{1}{|\varphi|^2} \langle \varphi^* \mathfrak{b} \rangle ) = -i 2 \frac{\partial}{\partial \bar{z}} ( \frac{1}{|\varphi|^2} \langle \varphi^* \mathfrak{b} \rangle )$ .

(6.47)

Here, $\frac{\partial}{\partial z}$ is the standard shorthand for $\frac{1}{2} ( \frac{\partial}{\partial x_1} - i \frac{\partial}{\partial x_2} )$; and likewise $\frac{\partial}{\partial \bar{z}}$ is shorthand for $\frac{1}{2} ( \frac{\partial}{\partial x_1} + i \frac{\partial}{\partial x_2} )$. These equations will be exploited by introducing a $\mathbb{C}$-valued function to be denoted by $\nu$ that is defined at any given $(t, z, x_3)$ by the rule

$$\nu = - \int_t^\infty \langle \varphi^* \mathfrak{b} \rangle .$$

(6.48)

This integral converges for any $(t, z, x_3)$ because $|\varphi| \le \frac{\mathfrak{z}}{t}$ (due to the top bullet in (6.1)) and because there exists a number $\mathcal{Z}$ which is independent of $t$ and such that $|\mathfrak{b}| \le \frac{\mathcal{Z}}{t}$ (see Lemmas 6.1 and Proposition 6.8). In fact Lemma 6.1 and Proposition 6.8 imply the following: Given $\varepsilon > 0$, there exists $\mathcal{Z}_\varepsilon$ (which is independendent of $t$, $z$ and $x_3$) such that

$$|\nu| \le \mathcal{Z}_\varepsilon \tfrac{1}{t^\varepsilon} \text{ where } t \le 1 \quad \text{and} \quad |\nu| \le \mathcal{Z}_\varepsilon \tfrac{1}{t^{2-\varepsilon}} \text{ where } t \ge 1 .$$

(6.49)

(A very small value for ε less than $\frac{1}{100}$ should be fixed and used throughout.) Other crucial properties of $v$ are established in the subsequent parts of the proof. What follows directly explains how $v$ will be used to prove Proposition 6.14.

The function $v$ is introduce by virtue of it obeying the following identities:

- $\frac{\partial}{\partial t} v = \langle \varphi^* \mathfrak{b} \rangle$ .
- $-i 2 \frac{\partial}{\partial \bar{z}} v = \langle \varphi^* \beta \rangle$ .

(6.50)

The top identity follows from the definition of $v$ in (6.48). As explained in Part 2 of this proof, the lower identity in (6.50) follows from the top bullet equation in (6.47) since that and the top bullet idenity in (6.50) say that the sum of $\langle \varphi^* \beta \rangle$ and $2i \frac{\partial}{\partial \bar{z}} v$ is independent of t for any fixed z and $x_3$. (Part 2 also proves that $v$ is differentiable.)

The identities in (6.50) can be used to write the lower bullet in (6.47) as a second order equation for $v$:

$$\frac{\partial}{\partial t} ( \frac{1}{|\varphi|^2} \frac{\partial}{\partial t} v) + 4 \frac{\partial}{\partial \bar{z}} ( \frac{1}{|\varphi|^2} \frac{\partial}{\partial z} v) = 0.$$

(6.51)

To expoit this identity, mutliply both sides by $\bar{v}$ and use the identities in (6.50) to write the result as follows:

$$|\beta|^2 + |\mathfrak{b}|^2 = \frac{\partial}{\partial t} ( \frac{1}{|\varphi|^2} \bar{v} \langle \varphi^* \mathfrak{b} \rangle) + 2i \frac{\partial}{\partial \bar{z}} ( \frac{1}{|\varphi|^2} \bar{v} \langle \varphi^* \beta \rangle) .$$

(6.52)

(The derivation also uses the fact that $\mathcal{L}^+$ is a complex *line* bundle to write a given section (call it η) on the $\varphi \neq 0$ part of $(0, \infty) \times \mathbb{R}^2 \times \mathbb{R}$ as $\frac{1}{|\varphi|^2} \langle \varphi^* \eta \rangle \varphi$.)

Fix T >> 1 and R >> T; and then imagine integrating both sides of (6.52) over the part of $(0, \infty) \times \mathbb{R}^2 \times \{0\}$ where $\frac{1}{T} \leq t \leq T$ and where $0 < |z| \leq R$. If there are no issues with the z = 0 locus with regards to the $\frac{1}{|\varphi|}$ factor on the right hand side of (6.52), then an integration by parts on the right hand side identifies the integral of $|\beta|^2 + |\mathfrak{b}|^2$ on the $\frac{1}{T} \leq t \leq T$ and $|z| \leq R$ part of $(0, \infty) \times \mathbb{R}^2 \times \{0\}$ with an integral on this domain's boundary.

Subsequent parts of the proof explain why these boundary integrals limit to zero as T and R get ever larger subject to R being significantly larger than T (and why the z=0 locus has no integration by parts contribution). Supposing that all of this happens, and given what the left hand side integral is, this limit of zero for the right hand side integral proves that β and $\mathfrak{b}$ have to vanish identically. And, supposing that both $\beta \equiv 0$ and $\mathfrak{b} \equiv 0$, then Lemma 6.11 asserts that (A, a) must be a model solution from Section 1c.

*Part 2*: This part of the proof first establishes that $v$ is differentiable and it then justifies the equality in the second bullet of (6.50). The proof that $v$ is differentiable starts with an appeal to Lemma 6.2 on the domain $(0,\infty) \times \mathbb{R}^2 \times \mathbb{R}$ which can be done because $t|\mathfrak{b}|$ has an apriori upper bound on this domain. This appeal results in the following bounds on $(0,\infty) \times \mathbb{R}^2 \times \mathbb{R}$ (with K being constant):

- $|\nabla_{\hat{A}}\varphi| \leq \frac{K}{t^2}$ and $|\nabla_{\hat{A}}\nabla_{\hat{A}}\varphi| \leq \frac{K}{t^3}$
- $|\nabla_{\hat{A}}\mathfrak{b}| \leq \frac{K}{t^2}$ and $|\nabla_{\hat{A}}\nabla_{\hat{A}}\varphi| \leq \frac{K}{t^3}$

(6.53)

These bound lead in turn to bounds for the derivatives of $\langle \varphi^* \mathfrak{b} \rangle$ (with L being constant):

- $|\nabla \langle \varphi^* \mathfrak{b} \rangle| \leq \frac{L}{t^3}$
- $|\nabla \nabla \langle \varphi^* \mathfrak{b} \rangle| \leq \frac{L}{t^4}$

(6.54)

To construct the derivative of $v$ along the $\mathbb{R}^2$-factor of $(0,\infty) \times \mathbb{R}^2 \times \mathbb{R}$, fix points z and u in $\mathbb{R}^2$ and then use (6.54) and Taylor's theorem with remainder to conclude that

$$v|_{(t,z+u)} = v|_{(t,z)} + \sum_{i=1,2} \int_t^\infty u^i \nabla_i \langle \varphi^* \mathfrak{b} \rangle + \mathfrak{E}$$

(6.55)

where the norm of $\mathfrak{E}$ is bounded by $|u|^2 \frac{L}{t^3}$. Note in particular that the explicit linear term in u that appears in (6.55) defines a linear functional on $T^*\mathbb{R}^2$ whose norm is bounded by $\frac{L}{t^2}$ (this is courtesy of (6.54).) What is written in (6.55) demonstrates two facts: First, that $v$ is differentiable along the $\mathbb{R}^2$ factor directions; and second that its differential in these directions at any given $(t, z, x_3)$ is bounded by $\frac{L}{t^2}$. Similar arguments using the bounds from Lemma 6.2 establish that $v$ is in fact $C^\infty$. These arguments are omitted.

With regards to the lower identity in (6.50): As noted in Part 1, when $\langle \varphi^* \mathfrak{b} \rangle$ is written using the top bullet in (6.50), then the top bullet of (6.47) makes the assertion that $\langle \varphi^* \mathfrak{b} \rangle + i2 \frac{\partial}{\partial z} v$ is independent of t for any fixed choice of $z \in \mathbb{R}^2$ and $x_3 \in \mathbb{R}$. With this understood, then the second bullet of (6.50) follows if this this function of t is zero; and that will be the case if the $t \to \infty$ limit of $\frac{\partial}{\partial z} v$ is zero since the $t \to \infty$ limit of $\langle \varphi^* \mathfrak{b} \rangle$ is zero (due to the top bullet in (6.1) and Proposition 6.8.) Meanwhile, the $t \to \infty$ limit of $\frac{\partial}{\partial z} v$ is indeed zero because (as noted in the previous paragraph), the norm of the derivative of $v$ on any given constant t slice of $(0,\infty) \times \mathbb{R}^2 \times \mathbb{R}$ is bounded by $\frac{L}{t^2}$.

*Part 3*: The hypothetical integration of (6.52) on $(0,\infty)\times(\mathbb{R}^2-0)\times\mathbb{R}$ and then integration by parts on the left hand side must deal with the $z = 0$ locus in $(0,\infty)\times\mathbb{R}^2\times\mathbb{R}$ where the factor $\frac{1}{|\varphi|^2}$ diverges when $\varphi$ has a zero. Since $\frac{1}{|\varphi|}|\langle\varphi^*\beta\rangle| = |\beta|$, any contribution from the $z = 0$ locus when integrating by parts will be due to a divergence there of $\frac{1}{|\varphi|}\nu$. This being the only issue, then there is no $z = 0$ contribution to an integration by parts because $\frac{1}{|\varphi|}\nu$ is, in fact, bounded on each constant t slice of $(0,\infty)\times\mathbb{R}^2\times\mathbb{R}$:

$$\tfrac{1}{|\varphi|}|\nu| \leq \mathcal{Z}_\varepsilon t^{1-\varepsilon} \text{ where } t \leq 1 \quad \text{and} \quad \tfrac{1}{|\varphi|}|\nu| \leq \mathcal{Z}_\varepsilon \tfrac{1}{t^{1-\varepsilon}} \text{ where } t \geq 1.$$

(6.56)

The proof that (6.56) holds follows directly.

The bound in (6.56) holds where $|z| \geq rt$ because of (6.49) and because $|\varphi|$ there is no smaller than $\frac{\delta}{t}$. With this understood, consider now where $|z| \leq rt$. To obtain the desired bound here, introduce by way of notation $m$ to denote the degree of vanishing of $\varphi$ at the $z = 0$ point on the constant $(t,x_3)$ slices of $(0,\infty)\times\mathbb{R}^2\times\mathbb{R}$. To prove the bound, fix $s \geq t$ for the moment and invoke Taylor's theorem with remainder plus Lemma 2.2 to obtain the bound $|\varphi| \leq c_{\mathfrak{z}} \frac{|z|^m}{s^{m+1}}$ at any point $(s,z,x_3)$. Meanwhile, Lemma 6.10's lower bound leads to the upper bound $\frac{1}{|\varphi|} \leq \varepsilon_*^{-1} \frac{t^{m+1}}{|z|^m}$ where $|z| \leq rt$. Using these bounds with (6.48) and with Proposition 6.8's bound for $|\mathfrak{b}|$ where $t \geq t_0$ and the bound $|\mathfrak{b}| \leq \frac{z}{t}$ otherwise leads to the following $t \geq 1$ and $|z| \leq rt$ bound:

$$\tfrac{1}{|\varphi|}|\nu| \leq \varepsilon_*^{-1} c_{\mathfrak{z}} \mathcal{Z}_\varepsilon t^{m+1} \int_t^\infty \tfrac{1}{s^{m+3-\varepsilon}}\, ds$$

(6.57)

This bound leads directly to the $t \geq 1$ bound in (6.56).

In the case where $t \leq 1$, use Proposition 6.8's bound where $t \leq t_0$ and the bound $|\mathfrak{b}| \leq \frac{z}{t}$ otherwise (and the afore-mentioned bound $\frac{1}{|\varphi|} \leq \varepsilon_*^{-1} \frac{t^{m+1}}{|z|^m}$ at $(t,z,x_3)$ and the bound $|\varphi| \leq c_{\mathfrak{z}} \frac{|z|^m}{s^{m+1}}$ at $(s, z, x_3)$ for $s \geq t$) o see that

$$\tfrac{1}{|\varphi|}|\nu| \leq \varepsilon_*^{-1} c_{\mathfrak{z}} \mathcal{Z}_\varepsilon t^{m+1} \int_t^1 \tfrac{1}{s^{m+1+\varepsilon}}\, ds + \tfrac{1}{|\varphi|}(|\nu|\big|_{s=1})$$

(6.58)

This bound leads directly to the $t \leq 1$ bound in (6.56).

*Part 4*: Looking ahead to the boundary terms in the planned integration by parts, this part of the proof derives the following bound for the integral of $\frac{1}{|\varphi|^2}|v|^2$ on any given constant $(t, x_3)$ slice of $(0, \infty) \times \mathbb{R}^2 \times \mathbb{R}$:

$$\int_{\{t\} \times \mathbb{R}^2 \times \{0\}} \frac{1}{|\varphi|^2}|v|^2 \leq \mathcal{Z}_\varepsilon t^{2-\varepsilon} \text{ where } t \leq 1 \quad \text{and} \quad \int_{\{t\} \times \mathbb{R}^2 \times \{0\}} \frac{1}{|\varphi|^2}|v|^2 \leq \mathcal{Z}_\varepsilon t^\varepsilon \text{ where } t \geq 1.$$

(6.59)

To prove these inequalities, use the bounds in (6.56) to see that the contribution to the respective integrals from the $|z| \leq rt$ part of $\{t\} \times \mathbb{R}^2 \times \{0\}$ are no greater than what is depicted on the right hand sides of the inequalities in (6.59). This is because the area of the $|z| \leq rt$ part of the integration domain is bounded by $\pi r^2 t^2$.

To consider the $|z| \geq rt$ part of the integration domain, use the fact that $|\varphi| \geq \frac{\delta}{t}$ there to obtain the preliminary bound:

$$\int_{\{t\} \times (\mathbb{R}^2 - D_{rt}) \times \{0\}} \frac{1}{|\varphi|^2}|v|^2 \leq \frac{t^2}{\delta^2} \int_{\{t\} \times \mathbb{R}^2 \times \{0\}} |v|^2 .$$

(6.60)

Meanwhile, from (6.48):

$$\int_{\{t\} \times \mathbb{R}^2 \times \{0\}} |v|^2 \leq c_\partial \left( \int_t^\infty \frac{1}{s} \left( \int_{\{s\} \times \mathbb{R}^2 \times \{0\}} |b|^2 \right)^{1/2} ds \right)^2$$

(6.61)

Invoke Proposition 6.1's bound for the case when $t \geq 1$ to see that the square of the $(t, \infty)$ integral that appears on the right hand side of (6.61) is bounded by $\mathcal{Z}_\varepsilon \frac{1}{t^{2-\varepsilon}}$ when $t \geq 1$. Using this bound with (6.60) leads to a $\mathcal{Z}_\varepsilon t^\varepsilon$ bound for the $|z| \geq rt$ contribution to the $t \geq 1$ version of the integral in (6.59).

Proposition 6.1's bound for the case $t \leq 1$ leads to a $\mathcal{Z}_\varepsilon t^{-\varepsilon}$ bound for the right hand side of (6.61). And, using that bound in (6.60) then leads to a $\mathcal{Z}_\varepsilon t^{2-\varepsilon}$ bound for the $|z| \geq rt$ contribution to any given $t \leq 1$ version of the integral in (6.59)

*Part 5*: This last part of the proof implements the integration by parts strategy that is outlined in Part 1. To this end, fix some large T (much greater than 1) and then fix $R \gg T$. Let $\chi_R$ denote here the function on $\mathbb{R}^2$ that is given by the rule $z \to \chi(\frac{|z|}{R} - 1)$. This function is equal to 1 where $|z| \leq R$ and it is equal to zero where $|z| > 2R$. Multiply both sides of (6.52) by $\chi_R$ and then integrate the resulting identity over the $\frac{1}{T} \leq t \leq T$ part of $(0, \infty) \times \mathbb{R}^2 \times \{0\}$. Then integrate by parts to obtain an inequality that has the form

$$\int_{[\frac{1}{T},T]\times\mathbb{R}^2\times\{0\}} (|\beta|^2 + |\mathfrak{b}|^2) \leq \int_{\{T\}\times D_{2R}\times\{0\}} \frac{1}{|\varphi|} |\nu||\mathfrak{b}| + \int_{\{\frac{1}{T}\}\times D_{2R}\times\{0\}} \frac{1}{|\varphi|} |\nu||\mathfrak{b}|$$

$$+ \frac{1}{R} \int_{[\frac{1}{T},T]\times(D_{2R}-D_R)\times\{0\}} \frac{1}{|\varphi|} |\nu||\beta| \ .$$

(6.62)

The integrals that appear on the right hand side of (6.62) are bounded in turn.

With regards to the right most term: Ignoring the factor of $\frac{1}{R}$, Proposition 6.8 and (6.59) imply that there is an R and T independent upper bound for the contribution to the integral giving the right most term in (6.62) from the $t \in [1, T]$ part of the integration domain. They also imply that there is an R and T independent bound for the contribution from the $t \in [\frac{1}{T}, 1]$ part of the integration domain. Putting back the $\frac{1}{R}$ factor, it follows therefore that the right most term on the right hand side of (6.62) is no greater than $\mathcal{Z}_* \frac{1}{R}$ with $\mathcal{Z}_*$ being independent of t and R.

With regards to the other terms on the right hand side of (6.62): The small t bounds from Proposition 6.8 and the small t bound in (6.59) lead to a $\mathcal{Z}_* \frac{1}{T^{2-\varepsilon}}$ bound for middle term on the right hand side of (6.62). Here, $\mathcal{Z}_*$ is again independent of T and R. By the same token, the large t bounds from Proposition 6.8 and the large t bound in (6.59) lead to a $\mathcal{Z}_* \frac{1}{T^{1-\varepsilon}}$ bound for the left most term on the right hand side of (6.62).

It follows from what was just said that the $R \to \infty$ and then $T \to \infty$ limit of the right hand side of (6.62) is equal to zero. Therefore, $\beta$ and $\mathfrak{b}$ must be everywhere zero.

### h) The $\frac{|z|}{t} \to \infty$ limit

The central proposition in this subsection is hinted at by Lemma 6.12. By way of notation, the proposition uses $(A^{(0)}, \mathfrak{a}^{(0)})$ to denote the Nahm pole solution from Section 1c (it is the $m = 0$ version of what is depicted in (1.7).) The proposition also refers to an automorphism of $P/\{\pm 1\}$. The latter is a principal SO(3) bundle. As noted just prior to the statement of Theorem 2, automorphisms of $P/\{\pm 1\}$ also act on the space of connections on P and the space of sections of ad(P).

**Proposition 6.15**: *Let $(A, \mathfrak{a})$ denote a solution to (1.4) that is described by (6.1). There is an automorphism of $P/\{\pm 1\}$ on the $|z| > t$ part of $(0, \infty) \times \mathbb{R}^2 \times \mathbb{R}$ to be denoted by $g_>$ with the following significance: If k is a non-negative integer, then*

$$\lim_{|z|/t \to \infty} (|(\nabla_{A^{(0)}})^{\otimes k} (g_>^* A - A^{(0)})| + |(\nabla_{A^{(0)}})^{\otimes k} (g_>^* \mathfrak{a} - \mathfrak{a}^{(0)})|) = 0$$

*with the limit being uniform with respect to the value of $x = (|z|^2 + t^2)^{1/2}$.*

Note that the third bullet of Theorem 2 follows from Proposition 5.1 and Proposition 6.15.

*Proof of Proposition 6.15*: The argument for Lemma 6.12 can be repeated with only cosmetic changes to see that

$$\lim_{|z|/t \to \infty} |t\alpha - \tfrac{1}{2}| = 0 \quad and \quad \lim_{|z|/t \to \infty} |t|\varphi| - \tfrac{1}{\sqrt{2}}| = 0$$

(6.63)

with both limits being uniform with respect to $x$. Meanwhile, $\lim_{|z|/t \to \infty} (|\beta| + |\mathfrak{b}|) = 0$ (uniformly with respect to $x$) because of Lemma 6.3. Keeping this in mind, fix an isomorphism between P and the product principal SU(2) bundle over $(0,\infty) \times \mathbb{R}^2 \times \mathbb{R}$. Let $\theta_0$ denote the product connection and let $\{\sigma_1, \sigma_2, \sigma_3\}$ a corresponding $\theta_0$-covariantly constant orthonormal frame that obeys the commutation conditions in (1.6). Since $\langle \varphi^2 \rangle = 0$ and $\langle \varphi \sigma \rangle = 0$, there is an automorphism of $P/\{\pm 1\}$ on $\pi_2((0,\infty) \times (\mathbb{R}^2 - 0) \times \mathbb{R}) = 0$ that writes $\mathfrak{a}_3$ with $\alpha \sigma_3 + \beta$ and it writes $\varphi$ as $|\varphi|(\sigma_1 - i\sigma_2)$ where $|z| > rt$. Use $g_>$ to denote this automorphism. Compare (6.63) with the $m = 0$ version of (1.7) to see that

$$\lim_{|z|/t \to \infty} (t|g_>{}^*A - A^{(0)}| + t|g_>{}^*\mathfrak{a} - \mathfrak{a}^{(0)}|) = 0.$$

(6.64)

This is the $k = 0$ statement of Proposition 6.15. Granted this $k = 0$ statement, then standard elliptic regularity techniques used in conjunction with the third bullet of (6.1) can be used to prove $k \geq 1$ assertions.

## 7. The natural logarithm of $|\varphi|$

As in the previous section, $(A, \mathfrak{a})$ denotes a solution to (1.4) that is described by (6.1). This section returns to the milieu of Section 5 to refine some of what is said there about the behavior of the norm of $|\varphi|$. Of particular interest here is a function to be denoted by $w$ which is defined on the complement of the $z = 0$ locus in $(0,\infty) \times \mathbb{R}^2 \times \mathbb{R}$ by writing $|\varphi|^2$ as $\frac{1}{2t^2} e^{4w}$. This section derives bounds for the integrals of $w$ and $w^2$ on the constant $(t, x_3)$ slices of $(0,\infty) \times \mathbb{R}^2 \times \mathbb{R}$ and for the integrals of $|dw|^2$ on the $|z| \geq rt$ part of these slices.

### a) The equation for the natural logarithm

The desired bounds $w$ come about by exploiting a second order equation for $w$ on the complement of the $z = 0$ locus that follows from (3.9) and (3.10):

$$-\left(\tfrac{\partial^2}{\partial t^2} + \tfrac{\partial^2}{\partial x_1^2} + \tfrac{\partial^2}{\partial x_2^2}\right)w + \tfrac{1}{2t^2}(e^{4w} - 1) + 4(|\beta|^2 + |\mathfrak{b}|^2) = 0.$$

(7.1)

(The equation in (4.37) is the $\beta \equiv 0$, $\mathfrak{b} \equiv 0$ version of (7.1).) With regards to the $z = 0$ locus exclusion: If $\varphi$ is not everywhere non-zero, then $w$ has a logarithmic singularity along the $z = 0$ locus; it is not differentiable on this locus.

With regards to the derivation of (7.1): The first step is to rewrite the equations in (3.9) using $w$. In this regard: The first bullet of (3.9) relates $\alpha$ to $w$:

$$\alpha = -\tfrac{1}{2t} + \tfrac{\partial}{\partial t} w.$$
(7.2)

Meanwhile, the second bullet in (3.9) is equivalent to the following identities

$$\nabla_{\hat{A}1} \tfrac{\varphi}{|\varphi|} = 2 \tfrac{\partial w}{\partial x_2} \tfrac{\varphi}{|\varphi|} \quad and \quad \nabla_{\hat{A}2} \tfrac{\varphi}{|\varphi|} = -2 \tfrac{\partial w}{\partial x_1} \tfrac{\varphi}{|\varphi|}$$
(7.3)

where $\varphi \neq 0$. (Here is how to prove (7.3): Fix a point $(t, z, x_3)$ where $z \neq 0$ and then fix a product structure for $\mathcal{L}^+$ on a small radius disk in the slice $\{t\} \times \mathbb{R}^2 \times \{x_3\}$ centered at $z$ with $\varphi$ nowhere zero on the disk. There is an isomorphism between $\mathcal{L}^+$ and the product $\mathbb{C}$ bundle over the disk that identifies $\varphi$ with $e^{2w}$. This is to say that $\tfrac{\varphi}{|\varphi|}$ is identified with the constant function 1. The connection $\hat{A}$ when viewed as a connection on $\mathcal{L}^+$ is identified with a connection on the product $\mathbb{C}$ bundle that has the form $\theta_0 - i2\hat{a}$ where $\theta_0$ denotes the product connection and $\hat{a}$ denotes an $\mathbb{R}$-valued 1-form on the disk. Having written $\hat{A}$ and $\varphi$ in this way, then the respective real and imaginary parts of the second bullet in (3.9) are equivalent to the respective right hand and left hand equations in (7.3).

To continue with the derivation of (7.1): Differentiate (7.2) with respect to $t$ to write the $t$-derivative of $\alpha$ as $\tfrac{1}{2t^2} + \tfrac{\partial^2}{\partial t^2} w$ and then use the top bullet in (3.10) to identify that with $\langle \sigma B_{\hat{A}3} \rangle + |\varphi|^2 + 4(|\beta|^2 + |\mathfrak{b}|^2)$. Thus, $\langle \sigma B_{\hat{A}3} \rangle$ is obtained from $\tfrac{1}{2t^2} + \tfrac{\partial^2}{\partial t^2} w$ by subtracting $|\varphi|^2 + 4(|\beta|^2 + |\mathfrak{b}|^2)$. Meanwhile, taking the $\nabla_{\hat{A}2}$-covariant derivative of the left hand equation in (7.3) and subtracting the $\nabla_{\hat{A}1}$-covariant derivative of the right hand equation in (7.3) writes $\langle \sigma B_{\hat{A}3} \rangle$ in terms of $w$ (where $z \neq 0$):

$$\langle \sigma B_{\hat{A}3} \rangle = -(\tfrac{\partial^2}{\partial x_1^2} + \tfrac{\partial^2}{\partial x_2^2}) w .$$
(7.4)

The preceding two depictions of $\langle \sigma B_{\hat{A}3} \rangle$ lead directly to (7.1).

**b) Pointwise behavior of $w$**

With regards to a lower bound for $w$: There isn't one if $\varphi$ has non-zero vanishing degree at the $z = 0$ point on the constant $(t, x_3)$ slices of $(0, \infty) \times \mathbb{R}^2 \times \mathbb{R}$. Indeed, $w$ near the $z = 0$ point has the form

$$w = \tfrac{m}{2} \ln(\tfrac{|z|}{t}) + \mathcal{X}$$

(7.5)

with *m* being the vanishing degree and with $\mathcal{X}$ being a smooth function. On the other hand, the second bullet in (6.1) leads to an a priori lower bound for *w* where $|z| \geq rt$ on each constant $(t, x_3)$ slice of $(0, \infty) \times \mathbb{R}^2 \times \mathbb{R}$. The next two lemmas give more bounds for *w*.

**Lemma 7.1**: *Given positive numbers $\mathfrak{z}, \delta, r, \Xi$ and then $\varepsilon$, there exists $r_\varepsilon \geq 1$ with the following significance: Suppose that $(A, \mathfrak{a})$ is a solution to (1.4) that is described by (6.1) with the given parameters. Then $(A, \mathfrak{a})$'s version of w obeys $|w| \leq \varepsilon$ where $|z| \geq r_\varepsilon t$. In general, $|w| \leq r_1(1 + |\ln \tfrac{|z|}{t}|)$ where $|z| > 0$.*

With regards to an upper bound: The first and second bullets in (6.1) plus Lemma 2.2 leads to an a priori upper bound for *w* where $\varphi \neq 0$ (which is the complement of the $z = 0$ locus). Here is a far stronger upper bound:

**Lemma 7.2**: *Let $(A, \mathfrak{a})$ denote a solution to (1.4) that is described by (6.1). Let m denote the vanishing degree of $\varphi$ at the $z = 0$ point on the constant $(t, x_3)$ slices of $(0, \infty) \times \mathbb{R}^2 \times \mathbb{R}$. Introduce by way of notation $w^{(m)}$ to denote $(A^{(m)}, \mathfrak{a}^{(m)})$ version of the function w (with $(A^{(m)}, \mathfrak{a}^{(m)})$ denoting the integer m model solution from Section 1c. Then $w \leq w^{(m)}$ with equality if and only if $(A, \mathfrak{a})$ is Aut(P)-equivalent to $(A^{(m)}, \mathfrak{a}^{(m)})$.*

The function $w^{(m)}$ is given by this formula:

$$w^{(m)} = \tfrac{1}{2} \ln(\tfrac{(m+1)\sinh(\Theta)}{\sinh((m+1)\Theta)})$$

(7.6)

with $\Theta$ as in (1.5). Note in particular that $w^{(m)} < 0$ and that $w^{(m)}$ has limit zero as $\tfrac{|z|}{t} \to \infty$.

*Proof of Lemma 7.1*: The proof has two parts. The first part proves the general $|z| > 0$ bound and the second proves the assertion concerning the domain where $|w| < \varepsilon$.

*Part 1*: To prove the general upper bound $|w| \leq r_1(1 + |\ln \tfrac{|z|}{t}|)$, note first that $|w|$ is bounded where $|z| \geq rt$ by the maximum of the norms of the natural logarithms of $\tfrac{\delta}{2}$ and $\tfrac{3}{2}$. This follows from the top two bullets of (6.1). With this bound understood, fix *z* with $0 < |z| < rt$. Since $|w|$ at time $s = \tfrac{1}{r}|z|$ is bounded by $\max(|\ln \tfrac{\delta}{2}|, |\ln \tfrac{3}{2}|)$ and since $|\tfrac{\partial}{\partial t} w| \leq \tfrac{3}{t}$ in any event (because $\alpha = -\frac{1}{2t} + \frac{\partial}{\partial t} w$), it follows by integrating from s to t (at fixed z) that $|w(z,t) - w(s,t)| \leq 3 \ln(\frac{t}{s})$. The bound in the lemma follows directly from this.

*Part 2*: The argument for the existence of the desired $r_\varepsilon$ is similar to the argument in the proof of Lemma 6.12. To start, assume to the contrary that there exists positive $\varepsilon$ with no corresponding $r_\varepsilon$. No generality is lost by assuming that $\varepsilon < \frac{1}{10,000}$. In this event, there is a sequence $\{(A_\ddagger^n, \mathfrak{a}_\ddagger^n)\}_{n \in \mathbb{N}}$ of solutions to (1.4) and an associated sequence $\{(t_n, z_n)\}_{n \in \mathbb{N}} \subset (0, \infty) \times \mathbb{R}^2$ with the following properties: First, each $n \in \mathbb{N}$ version of $(A_\ddagger^n, \mathfrak{a}_\ddagger^n)$ satisfies (6.1) with the $\mathfrak{z}$, $\delta$, $r$ and $\Xi$ being independent of n. Second, each $n \in \mathbb{N}$ version of $(t_n, z_n)$ are such that $|z_n| > n t_n$. Finally, for each $n \in \mathbb{N}$, the norm of the $(A_\ddagger^n, \mathfrak{a}_\ddagger^n)$ version of $\varphi$ at the point $(t_n, z_n, 0)$ differs from $\frac{1}{\sqrt{2} t_n}$ by more than $\frac{1}{100} \varepsilon$.

For each $n \in \mathbb{N}$, let $(A^n, \mathfrak{a}^n)$ denote the pull-back of $(A_\ddagger^n, \mathfrak{a}_\ddagger^n)$ by the coordinate rescaling and $\mathbb{R}^2$ factor translation map $(t, z, x_3) \to (t_n t, t_n z + z_n, t_n x_3)$. A subsequence of the sequence $\{(A^n, \mathfrak{a}^n)\}_{n \in \mathbb{N}}$ can be used as input for Lemma 6.9 if it is relabeled and numbered consecutively from 1. (The second bullet assumption in Lemma 6.9 is obeyed for a suitably chosen and then consecutively renumbered subsequence by virtue of the third bullet in Lemma 6.3.) Lemma 6.9 supplies a limit pair $(A^\infty, \mathfrak{a}^\infty)$ which obeys (1.4) and has vanishing $\beta$ and $\mathfrak{b}$. It also has $|\varphi^\infty| \geq \frac{\delta}{t}$ at every point. Granted these conditions, Lemma 6.11 asserts that $(A^\infty, \mathfrak{a}^\infty)$ is the Nahm pole solution (which is the $m = 0$ model solution from Section 1c). But this is nonsensical because the norm of the Nahm pole version of $\varphi$ at $(t = 1, z = 0)$ is $\frac{1}{\sqrt{2}}$ whereas $|\varphi^\infty|$ must differ from $\frac{1}{\sqrt{2}}$ by more than $\frac{1}{100} \varepsilon$ at this point. (This is because of the $C^\infty$ convergence dictated by Lemma 6.9 and because $(A^n, \mathfrak{a}^n)$ versions of $|\varphi|$ differ from $\frac{1}{\sqrt{2}}$ by more than $\frac{1}{100} \varepsilon$ at that same point.)

***Proof of Lemma 7.2***: Let u denote $w - w^{(m)}$. This is a smooth function on $(0, \infty) \times \mathbb{R}^2 \times \mathbb{R}$ because (7.5) with an appropriate version of $\mathcal{X}$ is obeyed by $w^{(m)}$. Subtracting the $(A^{(m)}, \mathfrak{a}^{(m)})$ version of (7.1) from $(A, \mathfrak{a})$'s version leads to the following differential equation:

$$-\left(\frac{\partial^2}{\partial t^2} + \frac{\partial^2}{\partial x_1^2} + \frac{\partial^2}{\partial x_2^2}\right) u + (e^{4u} - 1)|\varphi^{(m)}|^2 + 4(|\beta|^2 + |\mathfrak{b}|^2) = 0.$$

(7.7)

This equation also has an associated maximum principle; and it implies that u can not have a positive local maximum. As a consequence, each component in $(0, \infty) \times \mathbb{R}^2 \times \{0\}$ of the $u > 0$ set is non-compact. Moreover, Lemma 7.1 implies that the coordinate t has either no upper bound or no positive lower bound on any component of the $u > 0$ set.

With the preceding in mind, suppose that u > 0 at some point in $(0,\infty)\times\mathbb{R}^2\times\{0\}$. The five steps that follow derive nonsense from this assumption. (The proof of a lemma in Step 3 comes after the last step.)

Step 1: Supposing that u > 0 at some point in $(0,\infty)\times\mathbb{R}^2\times\{0\}$, then there is a sequence $\{(t_n, z_n)\}_{n\in\mathbb{N}}$ in $(0,\infty)\times\mathbb{R}^2$ with the following two properties: First, the sequence $\{u(t_n, z_n, 0)\}_{n\in\mathbb{N}}$ is a sequence of positive numbers that either diverges (if u has no positive upper bound) or limits as $n \to \infty$ to the least upper bound for the values of u on $(0,\infty)\times\mathbb{R}^2\times\mathbb{R}$. Second, the sequence the $\{\frac{z_n}{t_n}\}_{n\in\mathbb{N}}$ is a bounded sequence in $\mathbb{R}^2$ (so as not to foul Lemma 7.1).

For each $n \in \mathbb{N}$, let $(A^n, \mathfrak{a}^n)$ denote the pull-back of $(A, \mathfrak{a})$ by the coordinate rescaling diffeomorphism $(t, z, x_3) \to (t_n t, t_n z, t_n x_3)$. Use the sequence $\{(A^n, \mathfrak{a}^n)\}_{n\in\mathbb{N}}$ as the input sequence for Lemma 2.3 and let $(A^\infty, \mathfrak{a}^\infty)$ denote a limit solution of (1.4) that is supplied by that lemma. What is denoted by $\Lambda$ in what follows is the associated subsequence of $\mathbb{N}$ that is supplied by Lemma 2.3. The pair $(A^\infty, \mathfrak{a}^\infty)$ also obeys the conditions in (6.1) with the same values for $\mathfrak{z}, \delta, r$ and $\Xi$. Because of this, the corresponding version of $\varphi$ (which is henceforth denoted by $\varphi^\infty$) vanishes only at the z = 0 point in each constant $(t, x_3)$ slice of $(0,\infty)\times\mathbb{R}^2\times\mathbb{R}$ if it vanishes at all. Because $(A^\infty, \mathfrak{a}^\infty)$ is described by (6.1), the degree of vanishing of $\varphi^\infty$ at z = 0 in these slices is no smaller than $m$. Denote the degree of vanishing of $\varphi^\infty$ on the z = 0 locus by $k$.

There is one other point to make; this is with regards to the sequence $\{\frac{z_n}{t_n}\}_{n\in\Lambda}$. Because this is a bounded sequence in $\mathbb{R}^2$, no generality is lost by assuming henceforth that this sequence converges. Its limit is denoted by q.

Step 2: Supposing that n is a positive integer, let $w^n$ denote the corresponding version of $w$ and let $u^n$ denote the corresponding version of u. There is also an $(A^\infty, \mathfrak{a}^\infty)$ version of $w$ which will be denoted by $w^\infty$. The manner of convergence dictated by Lemma 2.3 is such that $\{w^n\}_{n\in\Lambda}$ converges to $w^\infty$ on compact subsets of $(0,\infty)\times(\mathbb{R}^2-0)\times\mathbb{R}$.

There is also an $(A^\infty, \mathfrak{a}^\infty)$ version of u which is denoted by $u^\infty$. It is defined to be $w^\infty - w^{(k)}$. The sequence $\{u^n\}_{n\in\Lambda}$ converges in the $C^\infty$ topology on compact subsets of $(0,\infty)\times(\mathbb{R}^2-0)\times\mathbb{R}$ to the function $w^{(k)} - w^{(m)} + u^\infty$. (This is because $\{w^n\}_{n\in\Lambda}$ converges on compact subsets of $(0,\infty)\times(\mathbb{R}^2-0)\times\mathbb{R}$ to $w^\infty$.)

Step 3: The following lemma will be used to deal with the $k = m$ case:

**Lemma 7.3**: *Let* $\{(A^n, \mathfrak{a}^n)\}_{n \in \mathbb{N}}$ *denote a sequence of solutions to (1.4) that all obey the same* $\mathfrak{z}, \delta, r$ *and* $\Xi$ *version of (6.1). Use this sequence as input to Lemma 2.3; and let* $(A^\infty, \mathfrak{a}^\infty)$ *denote a corresponding limit that is supplied by Lemma 2.3. Assume that all of the* $\{(A^n, \mathfrak{a}^n)\}_{n \in \mathbb{N}}$ *versions of* $\varphi$ *and also the* $(A^\infty, \mathfrak{a}^\infty)$ *version have the same degree of vanishing along the* $z = 0$ *locus in* $(0, \infty) \times \mathbb{R}^2 \times \mathbb{R}$. *Let* $\Lambda$ *denote the associated subsequence from Lemma 2.3. The sequence indexed by* $\Lambda$ *whose n'th element is the* $(A^n, \mathfrak{a}^n)$ *version of the function* $\mathrm{u}$ *has a subsequence that converges in the* $C^\infty$ *topology on compact subsets of* $(0, \infty) \times \mathbb{R}^2 \times \mathbb{R}$ *to the* $(A^\infty, \mathfrak{a}^\infty)$ *version of* $\mathrm{u}$.

This lemma is proved momentarily. The next paragraph uses it to finish the proof of Lemma 7.2 in the case when $k = m$.

Supposing that $k = m$, let $\{u^n\}_{n \in \Lambda'}$ denote the subsequence from Lemma 7.3 with $\Lambda'$ being the relevant subsequence of $\Lambda$. Lemma 7.3 asserts that

$$u^\infty(1, q, 0) = \lim_{n \in \Lambda'} u^n(1, q_n, 0) \ , \tag{7.8}$$

which is the supremum of the values of $u$ on $(0, \infty) \times \mathbb{R}^2 \times \mathbb{R}$. This also the supremum of the values of $u^\infty$ on $(0, \infty) \times \mathbb{R}^2 \times \mathbb{R}$ (this also follows from Lemma 7.3). Since $u^\infty$ obeys its own version of (7.7), it follows by the maximum principle that $u^\infty$ is constant. This implies that it is identically zero because $u^\infty \to 0$ as $|z| \to \infty$ at fixed t. (Lemma 7.1 applies to $u^\infty$ also.) Therefore, the supremum of $u$ on $(0, \infty) \times \mathbb{R}^2 \times \mathbb{R}$ must be zero which is the desired nonsense since it was assumed a priori that $u$ was positive at some point.

Step 4: This step assumes that $k > m$. Some preliminary observations: First, since $u^\infty$ on the $t = 1$ slice of $(0, \infty) \times \mathbb{R}^2 \times \mathbb{R}$ is bounded in any event, and since $w^{(k)} - w^{(m)}$ is negative for $k > m$ and it is unbounded from below as $|z| \to 0$, the point $q$ can't be at the origin. This implies that $\{u^n + w^{(m)} - w^{(k)}\}_{n \in \Lambda}$ converges in the $C^\infty$ topology to $u^\infty(1, q, 0)$ on some fixed (n-independent) ball in $(0, \infty) \times \mathbb{R}^2 \times \{0\}$ centered at $(1, q, 0)$. The latter observation implies in turn that $u^\infty(1, q, 0)$ is strictly greater than the supremum of the values of $u$ on $(0, \infty) \times \mathbb{R}^2 \times \mathbb{R}$. Thus, the supremum of $u^\infty$ on $(0, \infty) \times \mathbb{R}^2 \times \mathbb{R}$ is also positive.

Step 5: Given Step 4's observations, one can replace the original $(A, \mathfrak{a})$ by $(A^\infty, \mathfrak{a}^\infty)$ in Steps 1-4 to obtain a second limit $(A^{\infty 2}, \mathfrak{a}^{\infty 2})$. As in the original case, nonsense ensues if the vanishing degree of the corresponding version of $\varphi$ on the $z = 0$ locus is equal to $k$. If it is greater than $k$ (call it $k_2$), then one can repeat Steps 1-4 once again, and perhaps yet again, and so on. But note that this process of repetition has to end after a finite number of steps because the value of $k$ at each step increases whereas Lemma 5.6 puts an uppper bound on the size of any given version of $k$. Therefore, at some point, some version of *k* is equal to a previous version and, at that point, Step 3's nonsense ensues.

*Proof of Lemma 7.3*:  The proof is a riff on the proof of the third bullet in Lemma 6.1. To start:  Renumber the sequence $\{(A^n, \mathfrak{a}^n)\}_{n \in \Lambda}$ consecutively from 1.  It is proved momentarily that a subsequence of this sequence satisfies the assumptions of Lemma 6.9 after renumbering the subsequence consecutively from 1.  Assume that this is the case, renumber consecutively from 1 and then let $(A^\infty, \mathfrak{a}^\infty)$ denote a subsequence that is supplied by Lemma 6.9.  Let $\Lambda \subset \mathbb{N}$ denote the associated subsequence of integers (this is not the original version of $\Lambda$).  Fix $n \in \Lambda$, and let $\sigma^n$ and $\hat{A}^n$ denote the $(A^n, \mathfrak{a}^n)$ version of the data $\sigma, \hat{A}$.  Likewise, let $\sigma^\infty$ and $\hat{A}^\infty$ denote the $(A^\infty, \mathfrak{a}^\infty)$ version.  Lemma 6.9 implies that the sequence $\{\langle \sigma^n B_{\hat{A}^n 3}\rangle\}_{n \in \Lambda}$ converges to $\langle \sigma^\infty B_{\hat{A}^\infty 3}\rangle$ in the $C^\infty$ topology on compact subsets of $(0, \infty) \times \mathbb{R}^2 \times \mathbb{R}$.

With the preceding in mind, fix $n \in \Lambda$ and subtract the $(A^n, \mathfrak{a}^n)$ version of (7.4) from the $(A^\infty, \mathfrak{a}^\infty)$ version to obtain this:

$$-\left(\frac{\partial^2}{\partial x_1^2} + \frac{\partial^2}{\partial x_2^2}\right)(u^\infty - u^n) = \langle \sigma^\infty B_{\hat{A}^\infty 3}\rangle - \langle \sigma^n B_{\hat{A}^n 3}\rangle \ .$$

(7.9)

This identity holds on any given constant $(t, x_3)$ slice of $(0, \infty) \times \mathbb{R}^2 \times \mathbb{R}$.  Fix such a slice and let $D$ denote the disk of radius $3rt$ about the $z = 0$ point in that slice.  Having chosen a point $p$ in the concentric radius $2rt$ disk, let $\mathfrak{G}_p$ denote the Dirichelet Green's function for the operator $-\left(\frac{\partial^2}{\partial x_1^2} + \frac{\partial^2}{\partial x_2^2}\right)$ on $D$ with pole at $p$.  Multiply both sides of (7.9) by $\mathfrak{G}_p$ and then integrate the resulting identity over $D$.  Integration by parts then writes $(u^\infty - u^n)|_{(t,p,0)}$ as the sum of two integrals.  The first is an integral on the $|z| = 2rt$ circle of the product of the radial derivative of $\mathfrak{G}_p$ and $(u^\infty - u^n)$.  The second is an integral over $D$ of the product of $\mathfrak{G}_p$ and the expression on the right hand side of (7.9).  Since both integrals converge to zero uniformly with respect to the point $p$ and the chosen time $t$ (Lemma 6.9 here), the sequence $\{u^n\}_{n \in \Lambda}$ converges uniformly to $u^\infty$ in a neighorhood in $(0, \infty) \times \mathbb{R}^2 \times \{0\}$ of the concentric, radius $rt$ disk in $\{t\} \times \mathbb{R}^2 \times \{0\}$.

The $\partial D$ integral and the integral of $\mathfrak{G}_p$ times the right hand side of (7.9) can be differentiated with respect to the chosen time $t$ and point $p$ to obtain integrals for the derivatives of the difference $(u^\infty - u^n)$.  Lemma 6.9 guarantees that these also converge uniformly to zero.

The three steps that follow prove that Lemma 6.9's assumptions are met by Lemma 7.3's sequence $\{(A^n, \mathfrak{a}^n)\}_{n \in \mathbb{N}}$.

Step 1: At issue is whether the renumbered sequence $\{(A^n, \mathfrak{a}^n)\}_{n \in \Lambda}$ meets the requirements of Lemma 6.9's second bullet. Key points to keep in mind are these: First, $\{(A^n, \mathfrak{a}^n)\}_{n \in \Lambda}$ converges in the $C^\infty$-topology on compact subsets of $(0, \infty) \times \mathbb{R}^2 \times \mathbb{R}$ to $(A^\infty, \mathfrak{a}^\infty)$ after termwise application of elements in Aut(P). Second, the vanishing degrees of all $n \in \Lambda$ versions of $\varphi^n$ along the $z = 0$ locus are the same and this one vanishing degree is that of $\varphi^\infty$ along the $z = 0$ locus.

Step 2: Fix $T > 1$ to see about a bound for the various $(A^n, \mathfrak{a}^n)$ versions of $t|\mathfrak{b}|$ on the domain where $\frac{1}{T} \leq t \leq T$ and $|z| \leq T$. Keep in mind that Lemma 6.1 asserts in part the following: Given $\mu > 0$, there is an n-independent upper bound for each $n \in \Lambda$ version of $t|\mathfrak{b}|$ on the $|z| \geq \mu t$ part of $(0, \infty) \times \mathbb{R}^2 \times \mathbb{R}$.

Let $\varphi^\infty$ denote the $(A^\infty, \varphi^\infty)$ version of $\varphi$. Let $m$ denote the degree of vanishing of $\varphi^\infty$ at the $z = 0$ point of each constant $(t, x_3)$ slice of $(0, \infty) \times \mathbb{R}^2 \times \mathbb{R}$. Then there exists a positive $\varepsilon_T$ whose significance is this: If $t \in [\frac{1}{T}, T]$, then

$$|(\nabla_{A^\infty})^{\otimes m} \varphi^\infty|_{(t, z=0)} \geq \varepsilon_T \frac{1}{t^{m+1}} \; .$$

(7.10)

The first role of $(A^\infty, \mathfrak{a}^\infty)$ in the proof is here, to deliver this last bound.

Step 2: If $n \in \Lambda$, let $\varphi^n$ denote the $(A^n, \mathfrak{a}^n)$ version of $\varphi$. If n is sufficiently large (given T), then the manner of convergence dictated by Lemma 2.3 implies that

$$|(\nabla_{A^n})^{\otimes m} \varphi^n|_{(t, z=0)} \geq \tfrac{1}{2} \varepsilon_T \frac{1}{t^{m+1}} \; .$$

(7.11)

when $t \in [\frac{1}{T}, T]$. (This is the second role for $(A^\infty, \mathfrak{a}^\infty)$.)

The bound in (7.10) with the $c_3 \frac{1}{t^{m+2}}$ bound from Lemma 2.2 for the order $(m+1)$ $\nabla_{A^n}$-covariant derivative of $\varphi^n$ implies (via Taylor's theorem with remainder) that there exist positive, n-independent numbers $\delta_T$ and $\mu_T$ such that

$$|\varphi^n| \geq \delta_T \frac{|z|^k}{t^{k+1}}$$

(7.12)

for $t \in [\frac{1}{T}, T]$ and $|z| \leq \mu_T t$.

Step 3: The plan is to use the $(A^n, \mathfrak{a}^n)$ versions of (6.7) and (6.8) to bound the corresponding version of $|\mathfrak{b}|$. To do this, first assume that $n \in \Lambda$ is large so that (7.11) holds. Then, use the fact that $\varphi^n$ also has vanishing degree $m$ on the $z = 0$ locus to conclude that function $|\varphi|_s$ in (6.7) (which is $|\varphi^n|$ at times $s > t$) obeys $|\varphi|_s \leq c_3 \frac{|z|^k}{s^{k+1}}$ where $|z| \le t$ (which is also where $|z| \le s$ because $s \ge t$.) Use of this last upper bound and (7.11)'s lower bound in (6.7) leads directly to the desired n-independent upper bound for the norm of the $(A^n, \mathfrak{a}^n)$ version of $\flat$ on the part of $(0, \infty) \times \mathbb{R}^2 \times \mathbb{R}$ where $t \in [\frac{1}{T}, T]$ and $|z| \le \mu_T t$.

### c) Bounds for integrals of $w^2$

Suppose again that $(A, \mathfrak{a})$ is a solution to (1.4) that is described by the three bullets in (6.1). The central lemma for this subsection states an a priori bounds for the integral of $w^2$ on the $|z| > rt$ part of any constant $(t, x_3)$ slice of $(0, \infty) \times \mathbb{R}^2 \times \mathbb{R}$. The notation is such that $D_{rt}$ denotes the $|z| < rt$ disk in $\mathbb{R}^2$. The notation also has $\nabla w$ denoting the full gradient of $w$ on $(0, \infty) \times \mathbb{R}^2 \times \mathbb{R}$ (it has the t-derivative and the $x_1$ and $x_2$ derivatives).

**Lemma 7.4**: *Given positive numbers $\mathfrak{z}, \delta, r$ and $\Xi$, there exists $\kappa > 1$ with the following significance: Let $(A, \mathfrak{a})$ denote a solution to (1.4) that is described by (6.1) with the given data $\mathfrak{z}, \delta, r$ and $\Xi$. If $t \in (0, \infty)$, then*

$$\int_{\{t\} \times (\mathbb{R}^2 - D_{rt}) \times \{0\}} w^2 \le \kappa t^2 \quad \text{and} \quad \int_{[\frac{1}{2}t, 2t] \times (\mathbb{R}^2 - D_{rt}) \times \{0\}} |\nabla w|^2 \le \kappa t$$

*Proof of Lemma 7.4*: The proof has eight steps.

<u>Step 1</u>: Let $h$ denote the function $1 - e^{4w}$. An important point to keep in mind is that $h$ is nonnegative (Lemma 7.2) and that there exists $r_* \ge 2r$ (which depends only on $\mathfrak{z}$, $\delta$, $r$ and $t$) such that $h \ge \frac{1}{2}|w|$ on the $|z| \ge r_* t$ part of $(0, \infty) \times \mathbb{R}^2 \times \mathbb{R}$ (Lemma 7.1). Note also that $h \le |w|$ in any event; and thus $h$ limits uniformly to zero as $\frac{|z|}{t} \to \infty$ because $w$ does (see Lemma 7.1 again). In particular, the number $r_*$ can be chosen (to conserve notation) so that $|h| \le \frac{1}{100}$ on the $|z| \ge r_* t$ part of $(0, \infty) \times \mathbb{R}^2 \times \mathbb{R}$.

<u>Step 2</u>: Multiply both sides of (7.1) by $e^{4w}$ to see that $h$ obeys the following:

$$-\tfrac{1}{4}(\tfrac{\partial^2}{\partial t^2} + \tfrac{\partial^2}{\partial x_1^2} + \tfrac{\partial^2}{\partial x_2^2})h + \tfrac{1}{2t^2} h(1-h) - 4e^{4w}(|\nabla w|^2 + |\beta|^2 + |\flat|^2) = 0.$$

(7.12)

Multiply the latter equation again by $h$ and rearrange appropriately to obtain this

$$-\tfrac{1}{8}(\tfrac{\partial^2}{\partial t^2} + \tfrac{\partial^2}{\partial x_1^2} + \tfrac{\partial^2}{\partial x_2^2})h^2 + \tfrac{1}{2t^2} h^2(1-h) + 8e^{4w}(e^{4w} - \tfrac{1}{2})|\nabla w|^2 = he^{4w}(|\beta|^2 + |\flat|^2).$$

(7.13)

This equation will be used momentarily.

Step 3: Fix $R > 0$ and let $¥_R$ denote the function on $\mathbb{R}^2$ that is given by the rule whereby $¥_R = 1$ where $|z| \leq R$ and $¥ = e^{1 - |z|/R}$ where $|z| \geq R$. Both $¥_R$ and $|d¥_R|$ have finite integrals on the constant $(t, x_3)$ slices of $(0, \infty) \times \mathbb{R}^2 \times \mathbb{R}$. Note in this regard that $d¥_R$ vanishes where $|z| < R$ and it is equal to $-\frac{1}{R} ¥_R d|z|$ where $|z| > R$.

Fix $t \in (0, \infty)$ and supposing that $R > r_* t$, multiply both sides of (7.13) by $¥_R$ and then integrate the result over the constant $(t, x_3 = 0)$ slice of $(0, \infty)$. There is no worry in this regard about integrating the $e^{4w}(e^{4w} - \frac{1}{2})|dw|^2$ term near the $z = 0$ locus; it is bounded by $c_{\mathfrak{z}} \frac{1}{t^2}$ because

$$e^{4w}|\nabla w|^2 = 2|\nabla(t|\varphi|)|^2$$

(7.14)

which is bounded in turn by $c_{\mathfrak{z}} \frac{1}{t^2}$ (use the top bullet of (6.1) and Lemma 2.2). Keep in mind also that $e^{4w} \leq 1$ because $w \leq w^{(m)}$ and $w^{(m)} \leq 0$.

To put the results into manageable form, introduce by way of notation $f_R$ to denote the function of $t \in (0, \infty)$ given by the integral of $¥_R \hbar^2$ on the slice $\{t\} \times \mathbb{R}^2 \times \{0\}$. Meanwhile, let $g_R$ denote the integral of $¥_R e^{4w}|\nabla w|^2$ on that same slice. As is shown in the next step, doing this integration and then a single integration by parts along the $\mathbb{R}^2$ factor to put one factor of $\frac{\partial}{\partial x_1}$ and one of $\frac{\partial}{\partial x_2}$ on $¥_R$ leads to this inequality:

$$-\tfrac{1}{8} \tfrac{d^2}{dt^2} f_R + \tfrac{1}{2t^2}(1 - \tfrac{1}{100}) f_R - c_0 \tfrac{1}{R^2} f_R + 2g_R \leq \varsigma_1 + 4 \int_{\{t\} \times (\mathbb{R} - D_{r_* t}) \times \{0\}} (|\beta|^2 + |\mathfrak{b}|^2)$$

(7.15)

where $\varsigma_1$ is a number that depends only on the data $\mathfrak{z}, \delta, r$ and $\Xi$. (In any event, it is independent of $R$ and $t$).

Let $c_*$ denote the version of the number $c_0$ that appears on the left hand side of (7.15). The inequality in (7.15) where $t^2 \leq \frac{1}{1000 c_*} R^2$ implies this

$$-\tfrac{d^2}{dt^2} f_R + \tfrac{3}{t^2} f_R + 2g_R \leq \varsigma_2$$

(7.16)

where $\varsigma_2$ is a number that depends only on $\mathfrak{z}, \delta, r$ and $\Xi$. (The number $\varsigma_2$ is the sum of $8\varsigma_1$ and 32 times the upper bound from Lemma 6.3 for the integral of $|\beta|^2 + |\mathfrak{b}|^2$ on the constant $(t, x_3)$ slices of $(0, \infty) \times \mathbb{R}^2 \times \mathbb{R}$.)

Step 4: This step says more about $\varsigma_1$ and the derivation of (7.15). First of all, the term $\frac{1}{2t^2} h^2(1-h)$ from (7.13) is in any event nonnegative and larger than $\frac{1}{2t^2}(1-\frac{1}{100}) h^2$ where $|z| \geq r_* t$. This understood, add the integral of $\frac{1}{2t^2} \yen_R (1 - \frac{1}{100}) h^2$ to both sides of (7.13) to see that the integral of the term $\frac{1}{2t^2} h^2(1-h)$ on the left hand side of (7.12) accounts for the factor of $\frac{1}{2t^2}(1-\frac{1}{100}) f_R$ on the left hand side of (7.15) and a factor on the right that is no greater than the integral of $\frac{1}{2t^2} h^2$ over the $|z| \leq r_* t$ part of the slice $\{t\} \times \mathbb{R}^2 \times \{0\}$. That integral contributes to $\varsigma_1$ because it is bounded by $c_0 r_*^2$.

Consider next the contribution to the integral of the left hand side of (7.13) from the term $-\frac{1}{2} \yen_R (\frac{\partial^2}{\partial x_1^2} + \frac{\partial^2}{\partial x_2^2}) h^2$. One integration by parts for each of the respect $x_1$ and $x_2$ partial derivatives can be used to bound the contribution of that term from below by the integral of $-c_0 |\nabla^\perp \yen_R| |h| |\nabla^\perp h|$. The latter integral is in turn bounded from below by the sum of $-c_0 \frac{1}{R^2} f$ and $-\frac{1}{1000} g_R$. This accounts for the $-c_0 \frac{1}{R^2} f$ term on the left hand side of (7.15). The $-\frac{1}{1000} g_R$ term is discussed further in the next paragraph.

The term $8e^{4w}(e^{4w} - \frac{1}{2})|\nabla w|^2$ in (7.13) is greater than $3 \yen_R e^{4w} |\nabla w|^2$ where $|z| \geq r_* t$. This understood, add the integral of $3 \yen_R e^{4w} |\nabla w|^2$ to both sides of the integral of (7.13) to see that the integal of $\yen_R 8 e^{4w}(e^{4w} - \frac{1}{2})|\nabla w|^2$ and $-\frac{1}{1000} g_R$ (from the preceding paragraph) account for the factor of $2 g_R$ on the left hand side of (7.14) plus a contribution to $\varsigma_1$ which is no greater than the integral of $c_0 e^{4w} |\nabla w|^2$ on the $|z| \leq r_* t$ part of the slice. The fact that this right hand side integral is bounded in turn by $c_3$ follows from (7.14).

Step 5: Fix $\mu > 0$ for the moment. The operator that appears on the left hand side of (7.16) is the $\mu = 3$ version of the operator $-\frac{d^2}{dt^2} + \frac{\mu}{t^2}$. This operator for any given positive $\mu$ has the Green's function that is depicted below in (7.17). The notation has $G_s(t)$ denoting the version of the Green's function with second derivative singularity at the point $s \in (0, \infty)$. Also, $\gamma$ is short hand for $\frac{1}{2}(1 + (1 + 4\mu)^{1/2})$.

- $G_s(t) = \frac{1}{2\gamma - 1} t^\gamma s^{-\gamma+1}$ for $t \leq s$.
- $G_s(t) = \frac{1}{2\gamma - 1} s^\gamma t^{-\gamma+1}$ for $s \leq t$.

(7.17)

The function $G_s(t)$ for the rest of this proof is the $\mu = 3$ version.

Step 6: Use the function $\chi$ to construct a smooth, nonincreasing function on $(0, \infty)$ that is equal to 1 where $t^2 \leq \frac{1}{2000 c_*} R^2$ and equal to 0 where $t^2 \geq \frac{1}{1000 c_*} R^2$. This function can and should be constructed so that the norm of its second derivative is bounded by $c_0 c_*^2 \frac{1}{R^2}$. Denote the function by $\varpi_R$. (The essential point is that (7.16) holds on $\varpi$'s support.) Having fixed $\varepsilon \ll R$, let $\varpi_\varepsilon$ denote the $R = \varepsilon$ version of $\omega_R$.

Step 7: Fix a point t with $\varepsilon^2 \le t^2 \le \frac{1}{3000 c_*} R^2$ which is a point where $\varpi_R = 1$ and $\varpi_\varepsilon = 0$. Multiply both sides of (7.16) by $(1-\varpi_\varepsilon)\varpi_R G_t(\cdot)$ and integrate the result over the domain $(0, \infty)$. (The integrands on both sides have compact support because of the factor $(1-\varpi_\varepsilon)\varpi_R$). Two instances of integration then leads to the following inequality:

$$f_R(t) + t \int_{t/2}^{2t} g_R \le c_0 \varsigma_2 t^2 + c_0 (\tfrac{t}{R})^\gamma \tfrac{1}{R} \int_{R/c_0}^{c_0 R} f_R + c_0 t^{-\gamma+1} \varepsilon^{\gamma-2} \int_{\varepsilon/c_0}^{c_0 \varepsilon} f_R \ .$$

(7.18)

The $c_0 \varsigma_2 t^2$ term on the right appears with no R or $\varepsilon$ dependent prefactor because the function $G_t(\cdot)$ is integrable on $(0, \infty)$ and thus so is $G_t(\cdot) \varsigma_2$. (The fact that $G_t(\cdot)$ is integrable is due to the fact that the relevant $\gamma$ for use in (7.17) is greater than 2.)

The first observation with regards to this inequality is that the rightmost term on the right hand of (7.18) has limit 0 as $\varepsilon \to 0$. To see this, note first that $f_R$ at any $t \in (0, R]$ is bounded by $\pi R^2$. Thus, the term with $\varepsilon$ is bounded a priori by $c_0 t^{-\gamma+1} \varepsilon^{\gamma-1} R^2$ which has limit zero as $\varepsilon \to 0$ because $\gamma > 1$.

The second observation is with regards to the middle term in (7.18). This term is bounded by $c_0 t^2$ because: The function $f_R$ is bounded by $\pi R^2$ so its integral is bounded by $c_0 R^3$. Hence, the term itself is bounded by $c_0 t^2 (\tfrac{t}{R})^{\gamma-2}$ which is bounded by $c_0 t^2$ since $\gamma > 2$ and $t < R$.

Step 8: What was just said in Step 7 about (7.18) implies two bounds: The first is that $f_R(t)$ for $t \le c_0^{-1} R$ is no greater than $c_0(1+\varsigma_2) t^2$. The second bound is that the $[\tfrac{1}{2} t, 2t]$ integral of $g_R$ is at most $c_0(1+\varsigma_2) t$. Since both these bounds are independent of R, the $R \to \infty$ limit can be taken in both cases for fixed t (invoke the dominated convergence theorem). These respective limits are denoted by $f_\infty$ and $g_\infty$.

The fact that $f_\infty(t) \le c_0(1+\varsigma_2) t^2$ implies Lemma 7.4's claim about the integral of $w^2$ because $h \ge \tfrac{1}{2} |w|$ where $|z| \ge r_* t$ and because $|w| \le c_0$ where $r_* t \ge |z| \ge rt$. For the same reasons, the fact that the integral of $g_\infty(\cdot)$ on $[\tfrac{1}{2} t, 2t]$ is bounded by $c_0(1+\varsigma_2) t$ implies Lemma 7.4's claim about the integral of $|\nabla w|^2$

### d) Upper bounds for the integral of w

As before, $(A, \mathfrak{a})$ denotes a solution to (1.4) that is described by (6.1). Fix $t > 0$ and then fix $R \gg t$. The central observation in this subsection supplies a *lower* bound for the integral of the $(A, \mathfrak{a})$ version of $|w|$ over $|z| < R$ part of $\{t\} \times \mathbb{R}^2 \times \{0\}$. A reference point to keep in mind in this regard is that of the model solution version of $w^{(m)}$ with $m$ being the degree of vanishing of the $(A,\mathfrak{a})$ version of $\varphi$ along the $z = 0$ locus. In particular, the $w^{(m)}$ versions of these integrals are unbounded as a function of R:

$$\int_{\{t\}\times D_R \times \{0\}} w^{(m)} = -\frac{\pi m(m+2)}{6} t^2 \ln(\tfrac{R}{t}) + \mathcal{O}(t^2) .$$

(7.19)

(This follows form the formula in (7.6).) Since $w \leq w^{(m)}$, the integral of $w$ over the domain $\{t\}\times D_R \times \{0\}$ has to be more negative than what is depicted in (7.19) unless the given pair $(A, \mathfrak{a})$ is Aut(P) equivalent to $(A^{(m)}, \mathfrak{a}^{(m)})$.

With the preceding understood, write $w$ as $w^{(m)} + u$. The proposition that follows talks about the integral of u. To set the notation: Having fixed a postive number R, the proposition uses $\chi_R$ to denote the function on $\mathbb{R}^2$ that is defined by the rule $z \to \chi(\tfrac{|z|}{R} - 1)$. This function is equal to 1 where $|z| \leq R$ and zero where $|z| \geq 2R$.

**Proposition 7.5**: *Given positive numbers $\mathfrak{z}, \delta, r$ and $\Xi$, there exists $\kappa > 1$ with the following significance: Let $(A,\mathfrak{a})$ denote a solution to (1.4) that is described by (6.1) with the given data $\mathfrak{z}, \delta, r$ and $\Xi$. Fix $t \in (0, \infty)$ and $R \geq rt$. Then*

- $\int_{\{t\}\times\mathbb{R}^2\times\{0\}} \chi_R u \leq -\tfrac{4}{3} t^2 \int_{[t,R]\times\mathbb{R}^2\times\{0\}} \tfrac{1}{s}(|\beta|^2 + |\mathfrak{b}|^2) + \kappa t^2 .$

- $\int_{\{t\}\times\mathbb{R}^2\times\{0\}} \chi_R \tfrac{\partial}{\partial t} u \leq -\tfrac{8}{3} t \int_{[t,R]\times\mathbb{R}^2\times\{0\}} \tfrac{1}{s}(|\beta|^2 + |\mathfrak{b}|^2) + \kappa t .$

*Proof of Proposition 7.5*: The proof has nine steps.

<u>Step 1</u>: The starting point is the equation in (7.7) for u. Write $e^{4u} - 1$ as $4u + \mathfrak{e}_1(u)$ and note two facts about $\mathfrak{e}_1(u)$: It is non-negative but less than $8u^2$ (keep in mind that u is negative). Thus, the term $(e^{4u} - 1)|\varphi^{(m)}|^2$ that appears in (7.7) can be written as

$$(e^{4u} - 1)|\varphi^{(m)}|^2 = 4u|\varphi^{(m)}|^2 + \mathfrak{e}_1(u)|\varphi^{(m)}|^2 .$$

(7.20)

Meanwhile, $|\varphi^{(m)}|^2$ can be written $\tfrac{1}{2t^2} + (|\varphi^{(m)}|^2 - \tfrac{1}{2t^2})$. Use this decomposition to write

$$(e^{4u} - 1)|\varphi^{(m)}|^2 = \tfrac{2}{t^2} u + \mathfrak{k}(u)$$

(7.21)

with $\mathfrak{k}(u)$ denoting here and subsequently $\mathfrak{e}_1(u)|\varphi^{(m)}|^2 + 4u(|\varphi^{(m)}|^2 - \frac{1}{2t^2})$. Note in particular that $\mathfrak{k}(u)$ is non-negative (remember that $\mathfrak{e}_1(u)$ is non-negative whereas $u$ and $(|\varphi^{(m)}|^2 - \frac{1}{2t^2})$ are non-positive). The equation (7.7) can now be written schematically as:

$$-\left(\frac{\partial^2}{\partial t^2} + \frac{\partial^2}{\partial x_1^2} + \frac{\partial^2}{\partial x_2^2}\right)u + \frac{2}{t^2}u = -4(|\beta|^2 + |\mathfrak{b}|^2) - \mathfrak{k}(u).$$

(7.22)

This identity will be used momentarily.

Step 2: To make something of (7.22), fix $t \in (0, \infty)$ for the moment; and supposing that $R > rt$, multiply both sides of (7.20) by the function $\chi_R$ and then integrate both sides over $\{t\} \times \mathbb{R}^2 \times \{0\}$. To put the result in manageable form, let $f_R$ now denote the integral of $\chi_R u$ over this domain. Two instances of integration by parts leads to this inequality:

$$-\frac{d^2}{dt^2} f_R + \frac{2}{t^2} f_R = -4 \int_{\{t\} \times \mathbb{R}^2 \times \{0\}} \chi_R(|\beta|^2 + |\mathfrak{b}|^2) - \int_{\{t\} \times \mathbb{R}^2 \times \{0\}} \chi_R \mathfrak{k}(u) + \frac{1}{R^2} \int_{\{t\} \times (D_{2R} - D_R) \times \{0\}} (\chi^{\Delta})_R u ,$$

(7.23)

where $(\chi^{\Delta})_R$ is short hand for the function on $\mathbb{R}^2$ that is obtained by evaluating the function $(\frac{d^2}{d\rho^2}\chi + \frac{1}{\rho}\frac{d}{d\rho})\chi$ on $\mathbb{R}$ at $\rho = \frac{|z|}{R} - 1$. This last equation will be written as

$$-\frac{d^2}{dt^2} f_R + \frac{2}{t^2} f_R = -4 \int_{\{t\} \times \mathbb{R}^2 \times \{0\}} (|\beta|^2 + |\mathfrak{b}|^2) + \mathfrak{E}$$

(7.24)

with $\mathfrak{E}$ introduced here to denote

$$\mathfrak{E} = 4 \int_{\{t\} \times \mathbb{R}^2 \times \{0\}} (1 - \chi_R)(|\beta|^2 + |\mathfrak{b}|^2) - \int_{\{t\} \times \mathbb{R}^2 \times \{0\}} \chi_R \mathfrak{k}(u) + \frac{1}{R^2} \int_{\{t\} \times (D_{2R} - D_R) \times \{0\}} (\chi^{\Delta})_R u .$$

(7.25)

Bounds for this $\mathfrak{E}$ term are given in the next two steps

Step 3: Consider first the left most integral on the right hand side of (7.25): Lemma 6.3 can be brought to bear to see that this term (it is non-negative) is no greater than $\varsigma_1 \frac{t}{R}$ with $\varsigma_1$ being a number that depends only on the data set $\mathfrak{z}, \delta, r$ and $\Xi$. A bound of the same form, $\varsigma_1 \frac{t}{R}$, also holds for the norm of the right most integral in (7.25). Indeed, a bound of this sort with a different $\varsigma_1$ (but depending only on $\mathfrak{z}, \delta, r$ and $\Xi$) follows from Lemma 7.4 using the Cauchy-Schwarz inequality.

With regards to the middle integral in (7.25): The most important point to keep in mind is that $\mathfrak{k}(u)$ is non-negative. Therefore, this term makes a *non-positve* contribution to the right hand side of (7.25). In any event, the integral of $\chi_R \mathfrak{k}(u)$ that is depicted in this term is a priori bounded by a number (denoted by $\varsigma_2$) that depends only on $\mathfrak{z}, \delta, r$ and $\Xi$. The next step derives such a bound.

Step 4: To derive a purely $\mathfrak{z}, \delta, r$ and $\Xi$ bound for the $\chi_R \mathfrak{k}(u)$ integral in (7.25), break the integral into the part where $|z| \geq rt$ and that where $|z| \leq rt$. Bounds for these two parts of the integral are given directly.

The integrand $\chi_R \mathfrak{k}(u)$ where $|z| \geq rt$ is no larger than $c_0 \frac{1}{t^2}(u^2 + u w^{(m)})$. Meanwhile, the integral of the latter function on the $|z| \geq rt$ part of $\{t\} \times \mathbb{R}^2 \times \{0\}$ is bounded by a number depending only on $\mathfrak{z}, \delta, r$ and t. Indeed, this follows from Lemma 5.6 (to bound m) and Lemma 7.4 and the formula in (7.6) for $w^{(m)}$.

To see about the $|z| \leq rt$ part of the integral of $\chi_R \mathfrak{k}(u)$, note first that $\chi_R \equiv 1$ here, so the question is about the integral of $\mathfrak{k}(u)$. Go back to the definition of $\mathfrak{k}(u)$ in Step 1 and write it $(e^{4u} - 1)|\varphi^{(m)}|^2 - \frac{2}{t^2} u$. The absolute value of the integral of $(e^{4u} - 1)|\varphi^{(m)}|^2$ over the domain in question is no greater than $\frac{\pi}{2} r^2$ because $(e^{4u} - 1)|\varphi^{(m)}|^2$ is no greater than $\frac{1}{2t^2}$ and the area of the domain is $\pi r^2 t^2$. Granted the preceding, then the claim about the $\chi_R \mathfrak{k}(u)$ integral in (7.25) being bounded by a purely $\mathfrak{z}, \delta, r$ and $\Xi$ dependent number follows if the integral of $|u|$ where $|z| \leq rt$ on $\{t\} \times \mathbb{R}^2 \times \{0\}$ is likewise bounded by a purely $\mathfrak{z}, \delta, r$ and $\Xi$ dependent multiple of $t^2$. The next lemma makes a formal assertion to this effect.

**Lemma 7.6**: *Given positive numbers $\mathfrak{z}, \delta, r$ and $\Xi$, there exists $\kappa > 1$ with the following significance: Let $(A, \mathfrak{a})$ denote a solution to (1.4) that is described by (6.1) with the given data $\mathfrak{z}, \delta, r$ and $\Xi$. Fix $t \in (0, \infty)$. Then the integral of $|u|$ on the $|z| \leq rt$ part of $\{t\} \times \mathbb{R}^2 \times \{0\}$ is bounded by $\kappa t^2$.*

***Proof of Lemma 7.6***: To start the proof, subtract the $(A^{(m)}, \mathfrak{a}^{(m)})$ version of (7.4) from the $(A, \mathfrak{a})$ version to obtain the following equation for u:

$$-(\tfrac{\partial^2}{\partial x_1^2} + \tfrac{\partial^2}{\partial x_2^2}) u = \langle \sigma B_{\hat{A}3} \rangle - \langle \sigma^{(m)} B_{A^{(m)}3} \rangle .$$

(7.26)

Then use the top bullet in (3.12) to write the latter equation as

$$-(\tfrac{\partial^2}{\partial x_1^2} + \tfrac{\partial^2}{\partial x_2^2}) u = \langle \sigma B_{A3} \rangle + 4|\mathfrak{b}|^2 - \langle \sigma^{(m)} B_{A^{(m)}3} \rangle .$$

(7.27)

To exploit (7.27), introduce $\chi_{rt}$ to denote the $R = rt$ version of $\chi_R$ and multiply both sides of (7.27) by $\chi_{rt}|z|^2$; then integrate the result over $\{t\} \times \mathbb{R}^2 \times \{0\}$. (The integrands for both sides have compact support where $|z| \leq 2rt$ because of $\chi_{rt}$). Integrating by parts twice on the left hand side identifies the integral of $\chi_{rt} u$ with an integral of the function $\chi_{rt}|z|^2(\langle\sigma B_{A3}\rangle + 4|\mathfrak{b}|^2 - \langle\sigma^{(m)} B_{A^{(m)}3}\rangle)$ plus integrals that are supported in the annular region where $rt < |z| < 2rt$. Each of these integrals are bounded by $\varsigma_3 t^2$ with $\varsigma_3$ depending only on $\mathfrak{z}, \delta, r$ and $t$. This is because $|z| \leq 2rt$ on the integration domain, and because of the bounds from Lemma 2.2 for $|B_{A3}|$ and those from Lemma 6.3 for the integral of $|\mathfrak{b}|^2$; and because $|u|$ is bounded a priori (given $\mathfrak{z}, \delta, r$ and $\Xi$) where $|z| \geq rt$ (use the second bullet of (6.1)).

Step 5: This step introduces three functions of the coordinate t for $(0, \infty)$ that are used to exploit the identity in (7.24). To define the first of these, fix $R > 0$. The first function is denoted by $\varpi_R$; it is given by the rule $t \to \chi(\frac{100 t}{R} - 1)$. This function equals 1 where $t \leq \frac{1}{100} R$ and it is equal to zero where $t \geq \frac{1}{50} R$. (The definition is slighly different than that given in Step 6 of the proof of Lemma 7.4. The latter version could be used here just as well.) The second function is defined using a positive number $\varepsilon$; it is denote by $\varpi_\varepsilon$ and it is defined to be $R = \varepsilon$ version of $\varpi_R$.

The third function is defined using a given choice of $s \in (0, \infty)$. This function is the $\mu = 2$ (hence $\gamma = 2$) version of the function $G_s(t)$ from (7.17). It is the Green's function with pole at s for the operator $-\frac{d^2}{dt^2} + \frac{2}{t^2}$.

Step 6: Multiply both sides of (7.24) by $(1 - \varpi_\varepsilon) \varpi_R G_s(t)$ and then integrate both sides of the resulting identity over $(0, \infty)$ (the integrands have compact support). Having done this, then integrate by parts to obtain a formula for $f_R(s)$:

$$f_R(t) = \mathcal{A}_\varepsilon + \mathcal{B}_R + C - \tfrac{3}{4} t^2 \int_{[t,R] \times \mathbb{R}^2 \times \{0\}} \tfrac{1}{s}(|\beta|^2 + |\mathfrak{b}|^2)$$

(7.28)

where the functions $\mathcal{A}_\varepsilon, \mathcal{B}_R$ and $C$ are as given below.

- $\mathcal{A}_\varepsilon = \tfrac{4}{3t} \int_{\tfrac{1}{100}\varepsilon}^{\varepsilon} s^2(\tfrac{d^2}{ds^2}\varpi_\varepsilon + \tfrac{2}{s}\tfrac{d}{ds}\varpi_\varepsilon) f_R + \tfrac{4}{3t} \int_0^\varepsilon \varpi_\varepsilon s^2 \mathfrak{E} + \tfrac{4}{3t} \int_{(0,\varepsilon] \times \mathbb{R}^2 \times \{0\}} \varpi_\varepsilon s^2(|\beta|^2 + |\mathfrak{b}|^2)$.

- $\mathcal{B}_R = \tfrac{3}{4} t^2 \int_{\tfrac{1}{100}R}^{R} \tfrac{1}{s}(\tfrac{d^2}{ds^2}\varpi_R - \tfrac{2}{s}\tfrac{d}{ds}\varpi_R) f_R + \tfrac{3}{4} s^2 \int_{[\tfrac{1}{100}R, R] \times \mathbb{R}^2 \times \{0\}} \tfrac{1}{s}(1 - \varpi_R)(|\beta|^2 + |\mathfrak{b}|^2)$.

- $\mathcal{C} = \frac{4}{3t}\int_0^t s^2\mathfrak{E} \;+\; \frac{3}{4}t^2\int_t^R \frac{1}{s}\varpi_R\mathfrak{E} \;-\; \frac{4}{3t}\int_{(0,t]\times\mathbb{R}^2\times\{0\}} s^2(|\beta|^2+|\mathfrak{b}|^2)$

(7.29)

Note that the integrand of the right most integral on the right hand side of (7.28) lacks a factor of $(1-\varpi_\varepsilon)$. This is also the case for the left most integral in (7.29)'s definition of $\mathcal{C}$. The integrals without this factor converge as the lower endpoint of the domain limits to 0 because the constant $(t,x_3)$ integrals of $(|\beta|^2+|\mathfrak{b}|^2)$ have a t-independent upper bound; and likewise so does $\mathfrak{E}$ on the domain $(0,R]$. The absence of the $(1-\varpi_\varepsilon)$ factor in the integrands is accounted for by the two right most integrals in the formula for $\mathcal{A}_\varepsilon$ in (7.29). By the same token, the integrand of the right most integral on the right hand side of (7.28) lacks a factor of $\varpi_R$. This lack is accounted for by the right most integral in the formula for $\mathcal{B}_R$ in (7.29).

*Step 7*: This step considers in turn what are denoted by $\mathcal{A}_\varepsilon$ and $\mathcal{B}_R$ in (7.28). The important point with regards to $\mathcal{A}_\varepsilon$ is that the $\varepsilon \to 0$ limit of $\mathcal{A}_\varepsilon$ is zero (for fixed s and R). This is because the left most integral that appear in the definition is $\mathcal{O}(\varepsilon)$ as $\varepsilon\to 0$ and the others are $\mathcal{O}(\varepsilon^3)$ as $\varepsilon \to 0$. Indeed, the left most integral is $\mathcal{O}(\varepsilon)$ because $f_R(t)$ is at most $\pi R^2$. This is because the absolute value of u is bounded by $c_0$ where $|z|\geq rt$ so the absolute value of the integral of $\chi_R u$ on the $rt < |z| \leq 2R$ part of the domain is at most $c_0\pi R^2$. Meanwhile, the integral of $|u|$ where $|z|\leq rt$ is bounded (courtesy of Lemma 7.6) by a t-independent multiple of $t^2$. The right most integrals in the definition of $\mathcal{A}_\varepsilon$ are $\mathcal{O}(\varepsilon^3)$ because $\mathfrak{E}$ is bounded where $t\leq R$ and so is any given $t\in(0,\infty)$ integral of $|\beta|^2+|\mathfrak{b}|^2$ on $\{t\}\times\mathbb{R}^2\times\{0\}$.

The same arguments as those just given (but for notation) prove that what is denoted by $\mathcal{B}_R$ in (7.29) can be written as

$$\mathcal{B}_R = t^2 \mathcal{Z}_R$$

(7.30)

where $\mathcal{Z}_R$ is number that is bounded by a purely $\mathfrak{z}, \delta, r$ and $\Xi$ dependent constant. Note in particular that this bound is independent of R. Note also that $\mathcal{Z}_R$ is independent of t.

*Step 8*: This step concerns the term that is denoted by $\mathcal{C}$ in (7.28). The first point to note concerns the left most and right most integrals in (7.29)'s definition. These are:

$$\int_0^t s^2\mathfrak{E} \quad and \quad \int_{(0,t]\times\mathbb{R}^2\times\{0\}} s^2(|\beta|^2+|\mathfrak{b}|^2) \,.$$

(7.31)

Both integrals are bounded by a purely $\mathfrak{z}, \delta, r$ and $\Xi$ dependent multiple of $t^3$ because $\mathfrak{E}$ is bounded by a purely $\mathfrak{z}, \delta, r$ and $\Xi$ dependent number where $s \leq R$, and because any constant $(s, x_3)$ integral of $|\beta|^2 + |\mathfrak{b}|^2$ is likewise bounded.

With regards to the middle integral in the definition of $\mathcal{C}$:

$$\int_t^R \tfrac{1}{s} \varpi_R \mathfrak{E}$$

(7.32)

This is a sum of three terms that correspond to the three integrals on the right hand side of (7.25). As noted in Step 3, the left most and right most integrals in (7.25) contribute respective factors to $\mathfrak{E}$ with absolute values bounded where $s \leq R$ by $\varsigma_1 \tfrac{s}{R}$. As in Step 3, the number $\varsigma_1$ depends only on $\mathfrak{z}, \delta, r$ and $\Xi$. The absolute values of their respective contributions to (7.32) are therefore bounded by $\varsigma_1$. Meanwhile, the contribution of the middle integral in (7.25) to (7.32) is non-positive. This integral's contribution to $\mathcal{C}$ is denoted in the next step by $-\mathcal{C}_*$

Step 9: It follows from what is said in Steps 7 and 8 about $\lim_{\varepsilon \to 0} \mathcal{A}_\varepsilon$ being zero and about $\mathcal{B}_R$ and about $\mathcal{C}$ that $f_R(s)$ can be written as

$$f_R(t) = -\tfrac{3}{4} t^2 \int_{[t,R] \times \mathbb{R}^2 \times \{0\}} \tfrac{1}{s}(|\beta|^2 + |\mathfrak{b}|^2) + t^2(\mathcal{Q}(t) - \mathcal{C}_*(t))$$

(7.33)

where $\mathcal{Q}$ and $\mathcal{C}_*$ have the following properties: First, the norms of $\mathcal{Q}$ and $t$ times its derivative are bounded by a purely $\mathfrak{z}, \delta, r$ and $\Xi$ dependent number. Second, $\mathcal{C}_*$ is non-negative and $t$ times the norm of its derivative is likewise bounded by a purely $\mathfrak{z}, \delta, r$ and $\Xi$ dependent number. This depiction of $f_R(t)$ (and its t-derivative) leads directly to the assertion of Proposition 7.5. (The t-derivative of the right hand side of (7.33) computes the integral in the proposition's second bullet's integral. There are no issues with regard to interchanging differentiation with respect to t and integration on the $(t, 0)$ slice because the integration domain on the slice is bounded.)

**e) Lower bounds for the integral of $w$**

The central result of this subsection is Lemma 7.7 which asserts a lower bound for the constant $(t, x_3)$ integrals of the function $\chi_R w$. This lower bound is a counterpoint to Proposition 7.5's upper bound.

**Lemma 7.7**: *Given positive numbers $\mathfrak{z}, \delta, r$ and $\Xi$, there exists $\kappa > 1$ with the following significance: Let $(A, \mathfrak{a})$ denote a solution to (1.4) that is described by (6.1) with the given data $\mathfrak{z}, \delta, r$ and $\Xi$. Fix $t \in (0, \infty)$ and $R \geq rt$. Then*

$$\int_{\{t\} \times \mathbb{R}^2 \times \{0\}} \chi_R u \geq -\kappa t^2 (1 + |\ln \tfrac{R}{t}|) \quad \text{and} \quad \int_{\{t\} \times \mathbb{R}^2 \times \{0\}} \chi_R \tfrac{\partial}{\partial t} u \geq -\kappa t (1 + |\ln \tfrac{R}{t}|)$$

*Proof of Lemma 7.7*: By way of notation for this proof and subsequent ones: What is denoted below by $\varsigma$ is a positive number that depends only on the given data $\mathfrak{z}, \delta, r$ and $\Xi$ from (6.1). It will always be greater than 1 and it can be assumed to increase between successive appearance.

As for the proof of the lemma: The assertion follows from (7.33) given a suitable upper bound for what is denoted by $\mathcal{C}_*(t)$ in that equation. (Remember that the constant $(s, x_3)$ integrals of $|\beta|^2 + |\mathfrak{b}|^2$ that appear on the right hand side of (7.33) are uniformly bounded; and thus the corresponding integral of these (with a factor $\tfrac{1}{s}$) for s between t and R is bounded by $\varsigma \ln(\tfrac{R}{t})$.)

With regards to $\mathcal{C}_*(t)$: Track back through the notation to see that it is bounded by

$$\varsigma \int_t^R \tfrac{1}{s^3} \Big( \int_{\{s\} \times \mathbb{R}^2 \times \{0\}} (|u|^2 + |u|(1 - 2s^2|\varphi^{(m)}|^2)) \Big) \; .$$

(7.34)

The key observation is that the integral of $|u|^2$ over the any given slice of $\{s\} \times \mathbb{R}^2 \times \{0\}$ is bounded by $\varsigma s^2$ and likewise that of $|u|(1 - 2s^2|\varphi^{(m)}|^2)$. These bounds follow from the bounds in Lemma 7.4 for $w$ and from the pointwise bound where $|z| > 0$ from Lemma 7.1 (use (7.6) also). (Lemma 7.4 is used to bound the relevant integrals on the $|z| > rs$ part of the constant $(s, x_3 = 0)$ slice of $(0, \infty) \times \mathbb{R}^2 \times \mathbb{R}$ whereas the pointwise bound from Lemma 7.1 is used to bound the relevant integrals on the $|z| < rs$ part of the slice.)

## 8. The $t \to \infty$ and $t \to 0$ limits of $(A, \mathfrak{a})$

The two propositions that follow effectively restate the first two bullets of Theorem 2. The first proposition concerns the $t \to \infty$ limit of a solution to (1.4) that obeys the constraints in (6.1), and the second concerns the $t \to 0$ limit.

**Proposition 8.1**: *Let $(A, \mathfrak{a})$ denote a solution to (1.4) that is described by (6.1). Let m denote the $z = 0$ locus vanishing degree of the corresponding version of $\varphi$. There is an automorphism of P on the $t > 1$ part of $(0, \infty) \times \mathbb{R}^2 \times \mathbb{R}$ to be denoted by $g_\infty$ such that*

- $\lim_{t\to\infty} (t|g_\infty * A - A^{(m)}| + t|g_\infty * \mathfrak{a} - \mathfrak{a}^{(m)}| + |\ln(\frac{|\varphi|}{|\varphi^{(m)}|})|) = 0$,

- $\lim_{t\to\infty} t^{k+1}(|(\nabla_{A^{(m)}})^{\otimes k}(g_\infty * A - A^{(m)})| + |(\nabla_{A^{(m)}})^{\otimes k}(g_\infty * \mathfrak{a} - \mathfrak{a}^{(m)})|) = 0$ *for any* $k \in \mathbb{N}$,

*( These limits are uniform with respect to the coordinate z of the $\mathbb{R}^2$ factor.*

The proof of Proposition 8.1 is in Section 8e. Sections 8a-8d supply tools that are used in its proof (and in the proof of the next proposition).

**Proposition 8.2**: *Let* $(A, \mathfrak{a})$ *denote a solution to (1.4) that is described by (6.1). Let* m *denote the* $z = 0$ *locus vanishing degree of the corresponding version of* $\varphi$. *There is a non-negative integer* p *and automorphism of* P *on the* $t > 1$ *part of* $(0, \infty) \times \mathbb{R}^2 \times \mathbb{R}$ *to be denoted by* $g_0$ *such that for any non-negative integer* k,

$$\lim_{t\to\infty} t^{k+1}(|(\nabla_{A^{(m+2p)}})^{\otimes k}(g_0 * A - A^{(m+2p)})| + |(\nabla_{A^{(m+2p)}})^{\otimes k}(g_\infty * \mathfrak{a} - \mathfrak{a}^{(m+2p)})|) = 0$$

*with the limit being uniform with respect to the coordinate z of the $\mathbb{R}^2$ factor.*

The proof of Proposition 8.2 is in Section 8g. Section 8f supplies a tool that is used in the proof.

**a) A divergence identity**

The divergence identity given momentarily in (8.2) plays a central role in the later parts of this section. To set the stage, let $(A, \mathfrak{a})$ denote a solution (1.4) whose version of $\varphi$ vanishes (if at all) only at the $z = 0$ point in each constant $(t, x_3)$ slice of $(0, \infty) \times \mathbb{R}^2$. Let $m$ denote the degree of vanishing of $\varphi$ at these $z = 0$ points. Write $|\varphi|$ as $\frac{1}{\sqrt{2t}} e^{2w}$ as before to define $w$ and then write $w$ as

$$w = \mathfrak{w} + \frac{m}{2} \ln \frac{|z|}{t} .$$

(8.1)

This function $\mathfrak{w}$ is a smooth function on the whole of $(0, \infty) \times \mathbb{R}^2 \times \mathbb{R}$.

Here is the promised identity:

$$\frac{\partial}{\partial t}(-\frac{1}{2}(\frac{1}{t}\frac{\partial}{\partial t}w + \frac{1}{t^2}w) + \frac{1}{2}(\frac{\partial}{\partial t}w)^2 - \frac{1}{2}|\nabla^\perp \mathfrak{w}|^2 - \frac{1}{8t^2}(e^{4w} - 1 - 4w) + |\beta|^2 - |\mathfrak{b}|^2) =$$
$$\frac{m+1}{2t}(\frac{\partial^2}{\partial x_1^2} + \frac{\partial^2}{\partial x_2^2})\mathfrak{w} - \frac{\partial}{\partial x_1}(\frac{\partial}{\partial t}\mathfrak{w} \frac{\partial}{\partial x_1}\mathfrak{w}) - \frac{\partial}{\partial x_2}(\frac{\partial}{\partial t}\mathfrak{w} \frac{\partial}{\partial x_2}\mathfrak{w})$$
$$- i(\frac{\partial}{\partial x_1} - i\frac{\partial}{\partial x_2})\langle \beta^* \mathfrak{b}\rangle + i(\frac{\partial}{\partial x_1} + i\frac{\partial}{\partial x_2})\langle \mathfrak{b}^* \beta\rangle .$$

(8.2)

An important point: This identity holds on the whole of $(0,\infty)\times\mathbb{R}^2\times\mathbb{R}$; the $z = 0$ locus is not excluded. (Remember that $\mathfrak{w}$ is smooth and so is the t-derivative of $w$.) This identity is called a divergence identity because a straight forward rearrangement of its terms depicts the function $\frac{m+1}{2t}\langle\sigma B_{\hat{A}3}\rangle$ as the divergence of a vector field on $(0,\infty)\times\mathbb{R}^2\times\mathbb{R}$.

The derivation of (8.2) is given below in two steps.

Step 1: Things to keep in mind: First, $\alpha$ is $-\frac{1}{2t} + \frac{\partial}{\partial t}w$ which is also $-\frac{m+1}{2t} + \frac{\partial}{\partial t}\mathfrak{w}$. Second, the operator $\frac{\partial^2}{\partial x_1^2} + \frac{\partial^2}{\partial x_2^2}$ annihilates the function $\ln|z|$ (and also t). As a consequence, $(\frac{\partial^2}{\partial x_1^2} + \frac{\partial^2}{\partial x_2^2})w = (\frac{\partial^2}{\partial x_1^2} + \frac{\partial^2}{\partial x_2^2})\mathfrak{w}$ on the complement of the $z = 0$ locus (which is where (7.1) holds.)

With the preceding understood, multiply both sides of (7.1) by $\alpha$ and notice that the result of doing this is the sum of four terms:

- $\frac{1}{2t}\frac{\partial^2}{\partial t^2}w - \frac{\partial}{\partial t}w\frac{\partial^2}{\partial t^2}w$
- $-\frac{1}{4t^3}(e^{4w} - 1) + \frac{1}{2t^2}\frac{\partial}{\partial t}w(e^{4w} - 1)$
- $\frac{m+1}{2t}(\frac{\partial^2}{\partial x_1^2} + \frac{\partial^2}{\partial x_2^2})\mathfrak{w} - \frac{\partial}{\partial t}\mathfrak{w}(\frac{\partial^2}{\partial x_1^2} + \frac{\partial^2}{\partial x_2^2})\mathfrak{w}$ .
- $4\alpha(|\beta|^2 + |\mathfrak{b}|^2)$.

(8.3)

The appearance of $\frac{m+1}{2t}$ and $\frac{\partial}{\partial t}\mathfrak{w}$ in the middle bullet (and not $-\frac{1}{2t}$ and $\frac{\partial}{\partial t}w$) are accounted for by writing of $\alpha$ as $-\frac{m+1}{2t} + \frac{\partial}{\partial t}\mathfrak{w}$ when multiplying the $(\frac{\partial^2}{\partial x_1^2} + \frac{\partial^2}{\partial x_2^2})w$ term that appears in (7.1) while keeping in mind that this term is $(\frac{\partial^2}{\partial x_1^2} + \frac{\partial^2}{\partial x_2^2})\mathfrak{w}$.

Step 2: The expression in the top bullet in (8.3) can be written as

$$\frac{\partial}{\partial t}(\frac{1}{2t}\frac{\partial}{\partial t}w) + \frac{1}{2t^2}\frac{\partial}{\partial t}w - \frac{1}{2}\frac{\partial}{\partial t}(\frac{\partial}{\partial t}w)^2 \ ;$$

(8.4)

and the second bullet can be written as

$$\frac{\partial}{\partial t}(\frac{1}{8t^2}(e^{4w} - 1)) - \frac{1}{2t^2}\frac{\partial}{\partial t}w \ .$$

(8.5)

Notice in particular that $\frac{1}{2t^2}\frac{\partial}{\partial t}w$ appears in both (8.4) and (8.5), but with + appearing before it in (8.4) and - appearing before it in (8.5). When added together, these two give

$$\frac{\partial}{\partial t}(\frac{1}{2t}(\frac{\partial}{\partial t}w + \frac{1}{t}w) + \frac{1}{8t^2}(e^{4w} - 1 - 4w) - \frac{1}{2}(\frac{\partial}{\partial t}w)^2) \ ,$$

(8.6)

which accounts for part of (8.1).

Step 3: The left most term in the third bullet of (8.3) (the term with $\frac{m+1}{2t}$) is left as is. Meanwhile, the right most term in the third bullet of (8.3) can be written after some rearrangements as

$$\tfrac{1}{2} \tfrac{\partial}{\partial t} ((\tfrac{\partial}{\partial x_1} \mathfrak{w})^2 + (\tfrac{\partial}{\partial x_2} \mathfrak{w})^2) - \tfrac{\partial}{\partial x_1}(\tfrac{\partial}{\partial t} \mathfrak{w} \tfrac{\partial}{\partial x_1} \mathfrak{w}) - \tfrac{\partial}{\partial x_2}(\tfrac{\partial}{\partial t} \mathfrak{w} \tfrac{\partial}{\partial x_2} \mathfrak{w}) \,.$$

(8.7)

This rewriting of the third bullet in (8.3) accounts for more of (8.1).

Step 4: The remaining parts of (8.1) are accounted for by writing $4\alpha(|\beta|^2+|\mathfrak{b}|^2)$ using the identity from (6.17).

**b) Integrating the right hand side of the divergence identity**

Supposing that R is a postive number, reintroduce the function $\chi_R$ from Proposition 7.5. (It is the function on $\mathbb{R}^2$ that is given by the rule $z \to \chi(\frac{|z|}{R} - 1)$. The lemma that follows concerns the integral of the product of $\chi_R$ with the right hand side of (8.2) over positive t versions of $\{t\} \times \mathbb{R}^2 \times \{0\}$. These integrals define a function of t which is denoted henceforth by $\mathcal{I}_R$.

**Lemma 8.3**: *Given positive numbers $\mathfrak{z}, \delta, r$ and $\Xi$, there exists $\kappa > 1$ with the following significance: Let $(A, \mathfrak{a})$ denote a solution to (1.4) that is described by (6.1) with the given data $\mathfrak{z}, \delta, r$ and $\Xi$. Fix $R > 0$. The function $\mathcal{I}_R(\cdot)$ can be written where t obeys $R \geq 2rt$ as*

$$\mathcal{I}_R(t) = -\tfrac{\pi m}{2} \tfrac{1}{t} + \mathfrak{C}_R$$

*with $\mathfrak{C}_R$ being a function of t that obeys $\int_t^{2t} |\mathfrak{C}_R| \leq \kappa \tfrac{t}{R}$.*

*Proof of Lemma 8.3*: The integral of $\chi_R$ times the right hand side of (8.2) has three sorts of terms. First, there are the terms that have $x_1$ or $x_2$ derivatives acting on either $\langle \beta^* \mathfrak{b} \rangle$ or $\langle \mathfrak{b}^* \beta \rangle$, and then there are the terms with a respective $x_1$ or $x_2$ derivatives acting on the product of $\tfrac{\partial}{\partial t} \mathfrak{w}$ with the $x_1$ or $x_2$ derivative of $\mathfrak{w}$. Finally, there is the term with second derivatives of $x_1$ and $x_2$ acting on $\mathfrak{w}$ (the term with $(\tfrac{\partial^2}{\partial x_1^2} + \tfrac{\partial^2}{\partial x_2^2})\mathfrak{w}$). The three parts that follow treat (in order) these three sorts of terms. Adding their contributions gives the claim that is made by Lemma 8.3.

By way of notation, these parts of the proof use $\varsigma$ to denote a number that is greater than 1 and that depends only on $\mathfrak{z}, \delta, r$ and $\Xi$. Its value can be assumed to increase between successive appearances.

*Part 1*: The integral of $\chi_R$ times an $x_1$ or $x_2$ derivative of either $\langle \beta^* \mathfrak{b} \rangle$ or $\langle \mathfrak{b}^* \beta \rangle$ can be written using integration by parts as the integral of either $\langle \beta^* \mathfrak{b} \rangle$ or $\langle \mathfrak{b}^* \beta \rangle$ times the corresponding $x_1$ or $x_2$ derivative of $\chi_R$. Since the norm of the derivative of $\chi_R$ is at most $c_0 \frac{1}{R}$ and since its support is where $|z| \geq R$, Lemma 6.3 can be brought to bear to see that the norms of these $\langle \beta^* \mathfrak{b} \rangle$ and $\langle \mathfrak{b}^* \beta \rangle$ integrals are at most $\varsigma \frac{t}{R^2}$. To summarize: The $\langle \beta^* \mathfrak{b} \rangle$ and $\langle \mathfrak{b}^* \beta \rangle$ terms contribute to the function $\mathfrak{C}_R$ in Lemma 8.3's depiction of $\mathcal{I}_R$.

*Part 2*: Integrate by parts to write the contribution to $\mathcal{I}_R$ from the terms that have products of $\frac{\partial}{\partial t} \mathfrak{w}$ with $x_1$ and $x_2$ derivatives of $\mathfrak{w}$ as

$$\int_{\{t\} \times \mathbb{R}^2 \times \{0\}} \langle \nabla^\perp \chi_R, \nabla^\perp \mathfrak{w} \rangle \tfrac{\partial}{\partial t} \mathfrak{w} \ .$$

(8.8)

To see about the contribution of this term, return to the definition of $\mathfrak{w}$ in (8.1) to write it as $w - \frac{m}{2} \ln \frac{|z|}{t}$. Writing it this way, $\nabla^\perp \mathfrak{w}$ becomes $\nabla^\perp w - \frac{m}{2} \frac{1}{|z|} d|z|$ and $\frac{\partial}{\partial t} \mathfrak{w}$ becomes $\frac{\partial}{\partial t} w + \frac{m}{2t}$. The integral in (8.8) then becomes a corresponding sum of four terms:

- $\frac{1}{R} \int_{\{t\} \times \mathbb{R}^2 \times \{0\}} \chi'_R \frac{\partial}{\partial |z|} w \frac{\partial}{\partial t} w$ .

- $\frac{m}{2Rt} \int_{\{t\} \times \mathbb{R}^2 \times \{0\}} \chi'_R \frac{\partial}{\partial |z|} w$ .

- $-\frac{m}{2R} \int_{\{t\} \times \mathbb{R}^2 \times \{0\}} \frac{1}{|z|} \chi'_R \frac{\partial}{\partial t} w$ .

- $\frac{\pi m^2}{2t}$ .

(8.9)

Here, $\chi'_R$ is the function on $\mathbb{R}^2$ that is obtained by evaluating the derivative of $\chi$ at $\frac{|z|}{R} - 1$. The paragraphs that follow discuss in turn these four contributions to (8.8). By way of a summary: All but the fourth bullet in (8.9) contribute to the $\mathfrak{C}_R$ term in Lemma 8.3's expression for $\mathcal{I}_R$.

With regards to the top bullet's integral in (8.9). The norm of the integral over $[t, 2t]$ of the function of t that is depicted by the top bullet is bounded by $\varsigma \frac{t}{R}$. This is a consequence of Lemma 7.4. (One could work harder to see that the integral itself is bounded by $\varsigma \frac{1}{R}$.)

With regards to the second bullet in (8.9): First integrate by parts to move the derivative off of $w$ to write that integral as:

$$-\frac{m}{2Rt} \int_{\{t\}\times\mathbb{R}^2\times\{0\}} (\frac{1}{R}\chi''_R + \frac{1}{|z|}\chi'_R)w \ .$$

(8.10)

Here, $\chi''_R$ denotes the function on $\mathbb{R}^2$ that is obtained by evaluating $\chi$'s second derivative at $\frac{|z|}{R} - 1$. The Cauchy-Schwarz inequality and Lemma 7.4's bound for the integral of $|w|^2$ can be used to bound what is written in (8.10) by $\varsigma \frac{1}{R}$.

With regards to the third bullet in (8.10): The Cauchy-Schwarz inequality and Lemma 7.4's bound lead directly to a $\varsigma \frac{t}{R}$ bound for the [t, 2t] integral of the function of t that is depicted in (8.10)'s third bullet.

With regards to the fourth bullet in (8.10): This is a direct calculation (the value is positive because $\chi$ is non-increasing).

*Part 3*: This last part of the proof of Lemma 8.3 discusses the contribution to $\mathcal{I}_R(t)$ from the $\frac{m+1}{2t}(\frac{\partial^2}{\partial x_1^2} + \frac{\partial^2}{\partial x_2^2})\mathfrak{w}$ term that appears on the right hand side of (8.2). To see about the integral of $\chi_R$ times this, integrate by parts to write it as

$$-\frac{m+1}{2Rt} \int_{\{t\}\times\mathbb{R}^2\times\{0\}} \chi'_R \frac{\partial}{\partial |z|}\mathfrak{w} \ .$$

(8.11)

Then write the derivative of $\mathfrak{w}$ as $\frac{\partial}{\partial |z|}w - \frac{m}{2}\frac{1}{|z|}d|z|$ as was done previously. Doing this writes (8.11) as a sum of an integral whose value is $-\frac{\pi m(m+1)}{2t}$ and one that is the same as that in the second bullet of (8.10) but for the replacement of $m$ with $-(m+1)$ out front. The latter integral contributes to the function $\mathfrak{C}_R$. Meanwhile, the $-\frac{\pi m(m+1)}{2t}$ contribution here and the $\frac{\pi m^2}{2t}$ contribution from the fourth bullet in (8.10) add to give the explicit term, $-\frac{\pi m}{2}\frac{1}{t}$, in Lemma 8.3's expression for $\mathcal{I}_R$.

c) **Integrating the left hand side of the divergence identity**

The integral on $\{t\}\times\mathbb{R}^2\times\{0\}$ of $\chi_R$ times the expression on the left hand side of (8.2) is written here as $\frac{d}{dt}\mathcal{J}_R$ with $\mathcal{J}_R$ being the function of t given by the rule

$$\mathcal{J}_R(t) = -\tfrac{1}{2} \int_{\{t\}\times\mathbb{R}^2\times\{0\}} \chi_R (\tfrac{1}{t}\tfrac{\partial}{\partial t} w + \tfrac{1}{t^2} w) - \tfrac{1}{2} \int_{\{t\}\times\mathbb{R}^2\times\{0\}} \chi_R |\nabla^\perp \mathfrak{w}|^2$$
$$+ \int_{\{t\}\times\mathbb{R}^2\times\{0\}} \chi_R (\tfrac{1}{2}(\tfrac{\partial}{\partial t} w)^2 + \tfrac{1}{8t^2}(1+4w-e^{4w}) + |\beta|^2 - |\mathfrak{b}|^2)$$

(8.12)

The lemma that follows momentarily describes this function. By way of notation, the lemma uses $D_{rt}$ (for $t > 0$) to denote the $|z| \le rt$ disk in $\mathbb{R}^2$.

**Lemma 8.4**: *Given positive numbers $\mathfrak{z}, \delta, r$ and $\Xi$, there exists $\kappa > 1$ with the following significance: Let $(A, \mathfrak{a})$ denote a solution to (1.4) that is described by (6.1) with the given data $\mathfrak{z}, \delta, r$ and $\Xi$. Fix $R > 0$. The function $\mathcal{J}_R(\cdot)$ can be written where $t$ obeys $R \ge 2rt$ as*

$$\mathcal{J}_R(t) = \tfrac{\pi m}{2} \ln(\tfrac{R}{t}) + 2 \int_{[t,R]\times\mathbb{R}^2\times\{0\}} \tfrac{1}{s}(|\beta|^2 + |\mathfrak{b}|^2) - \tfrac{1}{2} \int_{\{t\}\times D_{rt}\times\{0\}} |\nabla^\perp \mathfrak{w}|^2 + \mathfrak{Y}_R$$

*with $\mathfrak{Y}_R$ being a function of $t$ that obeys $-\kappa t \le \int_t^{2t} \mathfrak{Y}_R \le \kappa t(1 + |\ln \tfrac{R}{t}|)$.*

*Proof of Lemma 8.4*: The three parts of this proof explain how the three integrals in (8.12) contribute to the various terms in Lemma 8.4's formula.

*Part 1*: The left most term in (8.12)'s definition of $\mathcal{J}_R$ is this:

$$-\tfrac{1}{2} \int_{\{t\}\times\mathbb{R}^2\times\{0\}} \chi_R (\tfrac{1}{t}\tfrac{\partial}{\partial t} w + \tfrac{1}{t^2} w) \ .$$

(8.13)

Proposition 7.5 and (7.19) and Lemma 7.7 can be invoked to write the latter as follows:

$$\tfrac{\pi m(m+2)}{4} \ln(\tfrac{R}{t}) + 2 \int_{[t,R]\times\mathbb{R}^2\times\{0\}} \tfrac{1}{s}(|\beta|^2 + |\mathfrak{b}|^2) + \mathfrak{q}_{R1}$$

(8.14)

with $\mathfrak{q}_{R1}$ being a function of $t$ that is bounded from below by $-\varsigma$ and bounded from above (via Lemma 7.7) by $\varsigma(1+|\ln\tfrac{R}{t}|)$. Adding the $\mathfrak{q}_{R1}$ term in (8.14) to the lemma's function $\mathfrak{Y}_R$ is consistent with the lemma's requirement for $\mathfrak{Y}_R$ due to these norm bound.

*Part 2*: This part discusses the middle term on the right hand side of (8.12):

$$-\tfrac{1}{2} \int_{\{t\}\times\mathbb{R}^2\times\{0\}} \chi_R |\nabla^\perp \mathfrak{w}|^2 \ .$$

(8.15)

This expression is analyzed by writing $\nabla^\perp \mathfrak{w}$ as $\nabla^\perp w - \frac{m}{2} \frac{1}{|z|} d|z|$ where $|z| \geq 2rt$ which writes $|\nabla^\perp \mathfrak{w}|$ where $|z| \geq 2rt$ as follows:

$$|\nabla^\perp \mathfrak{w}|^2 = |\nabla^\perp w|^2 - \frac{m}{|z|} \frac{\partial}{\partial |z|} w + \frac{m^2}{4|z|^2}$$

(8.16)

This last decomposition then writes (8.15) as a sum of four integrals that correspond to the product of $\chi_R$ with the $|z| \geq rt$ integrals of the three terms in (8.16) plus this:

$$- \tfrac{1}{2} \int_{\{t\} \times D_{2rt} \times \{0\}} |\nabla^\perp \mathfrak{w}|^2 \ .$$

(8.17)

The latter appears explicitly in the lemma's formula for $\mathcal{J}_R$.

With regards to the integrals of the terms on the right hand side of in (8.16): The integral of the function $-\tfrac{1}{2} \chi_R |\nabla^\perp w|^2$ on the $|z| \geq 2rt$ part of $\{t\} \times \mathbb{R}^2 \times \{0\}$ contributes to $\mathfrak{Y}_R$. Lemma 7.4 guarantees that such an assignment is consistent with Lemma 8.2's requirement for $\mathfrak{Y}_\mathfrak{R}$. Meanwhile, Lemma 7.4 with the middle bullet in (6.1) guarantee that the assignment of the integral of $\tfrac{1}{2} \chi_R \frac{m}{|z|} \frac{\partial}{\partial |z|} w$ over the $|z| \geq 2rt$ part of $\{t\} \times \mathbb{R}^2 \times \{0\}$ can be likewise be assigned to the $\mathfrak{Y}_R$. (This is not directly obvious. To see why this assigment is consistent, first integrate by parts to write the integral in question as the sum of the integral of the function $-\tfrac{1}{2} \frac{m}{|z|} w$ over the $|z| = rt$ circle in $\{t\} \times \mathbb{R}^2 \times \{0\}$ and the integral over this circle of $\tfrac{1}{2} - (\frac{\partial}{\partial |z|} \chi_R) \frac{m}{|z|} w$. The norm of former is bounded by $\varsigma$ because of the second bullet in (6.1); and the norm of the latter is bounded by $\varsigma \frac{t}{R}$ because of what is said by Lemma 7.4 (use the Cauchy-Schwarz inequality to see this.)

Finally, there is the integral of $-\chi_R \frac{m^2}{8|z|^2}$ over the $|z| \geq rt$ part of $\{t\} \times \mathbb{R}^2 \times \{0\}$. This integral (up to adding a term bounded by $\varsigma$) is $-\frac{\pi m^2}{4} \ln(\frac{R}{t})$. Note that this term plus the left most term in (8.14) account for the explicit $\frac{\pi m}{2} \ln(\frac{R}{t})$ term in Lemma 8.4's depiction of the function $\mathcal{J}_R$.

*Part* 3: The whole of the right most integral in (8.12) can be safely assigned to the $\mathfrak{Y}_R$ term in Lemma 8.4's description. Indeed, this is so for the integral of $(\frac{\partial}{\partial t} w)$ by virtue of Lemma 7.4 and (to handle the $|z| \geq rt$ part of the integral) and by virtue of the fact that $|\frac{\partial}{\partial t} w| \leq c_3 \frac{1}{t}$ in any event (because $\alpha = -\frac{1}{\sqrt{2}t} + \frac{\partial}{\partial t} w$). Meanwhile, the integral of $\chi_R \frac{1}{8t^2} (1 + 4w - e^{4w})$ over $\{t\} \times \mathbb{R}^2 \times \{0\}$ is bounded by $c_0$ times that of $\frac{1}{8t^2} w^2$ which is in turn bounded by $\varsigma$ courtesy of Lemma 7.4 (for the $|z| \geq rt$ part of the integral) and courtesy of the assertion in Lemma 7.1 to the effect that $|w| \leq \varsigma (1 + |\ln(\frac{|z|}{t})|)$ where $|z| \leq rt$ for the (for the $|z| \le rt$ part of the integral). Then, last of all, there is the integral of $\chi_R(|\beta|^2 - |\mathfrak{b}|^2)$ over $\{t\} \times \mathbb{R}^2 \times \{0\}$ which is bounded by $\varsigma$ courtesy of Lemma 6.3.

### d) The $\{t\} \times \mathbb{R}^2 \times \{0\}$ integral of $|\beta|^2 + |\mathfrak{b}|^2$ as $t \to \infty$

This section states and then proves an assertion to the effect that the $\{t\} \times \mathbb{R}^2 \times \{0\}$ integrals of $|\beta|^2 + |\mathfrak{b}|^2$ limit to zero as $t \to \infty$. This fact is asserted formally by the proposition that follows

**Proposition 8.5**: *Suppose that* $(A, \mathfrak{a})$ *denotes a solution to (1.4) that is described by (6.1). Given* $\varepsilon > 0$, *there exists* $T_\varepsilon$ *such that* $\int_{\{t\} \times \mathbb{R}^2 \times \{0\}} (|\beta|^2 + |\mathfrak{b}|^2) \le \varepsilon$ *if* $t > T_\varepsilon$.

*Proof of Proposition 8.5*: The proof has five parts.

*Part 1*: Fix $R > 0$ and then fix $t > 0$ but so that $R > rt$. Multiply both sides of the divergence identity in (8.2) by $\chi_R$ and then integrate the result over $\{t\} \times \mathbb{R}^2 \times \{0\}$ to obtain the following identity:

$$\tfrac{d}{dt} \mathcal{J}_R = \mathcal{I}_R(t) \ .$$

(8.18)

Fix $t > 0$ but so that $R > 4rt$ and then fix $t_* \in [t, 2t]$. Meanwhile, fix $R_* \in [\tfrac{1}{2}R, R]$. Having done this, integrate both sides of (8.18) between $t_*$ and $R_*$ and use the fundamental theorem of calculus to write:

$$-\mathcal{J}_R(t_*) + \mathcal{J}_R(R_*) = \int_{t_*}^{R_*} \mathcal{I}_R \ .$$

(8.19)

The plan from here is to use Lemmas 8.3 and 8.4 to write the respective right hand and left hand sides of (8.19). In doing this, $\varsigma$ will again denote a number that is greater than 1 which depends only on $\mathfrak{z}, \delta, r$ and $\Xi$. Its value can increase between successive appearances. But $\varsigma$ is specifically independent of $t_*$, $R_*$ and $R$.

*Part 2*: Invoke Lemma 8.3 to write the right hand side of (8.19) as

$$\int_{t_*}^{R_*} \mathcal{I}_R = -\tfrac{\pi m}{2} \ln(\tfrac{R}{t}) + \mathfrak{X}_R$$

(8.20)

with $\mathfrak{X}_R$ obeying $|\mathfrak{X}_R| \leq \varsigma$. Then use Lemma 8.4 to write the left hand side of (8.19) as

$$-\mathcal{J}_R(t_*) + \mathcal{J}_R(R_*) = -\tfrac{\pi m}{2} \ln(\tfrac{R}{t}) - 2 \int_{[t_*,R_*]\times\mathbb{R}^2\times\{0\}} \tfrac{1}{s}(|\beta|^2 + |\mathfrak{b}|^2) - \tfrac{1}{2} \int_{\{R_*\}\times D_{rR_*}\times\{0\}} |\nabla^\perp \mathfrak{w}|^2$$
$$+ \tfrac{1}{2} \int_{\{t_*\}\times D_{rt_*}\times\{0\}} |\nabla^\perp \mathfrak{w}|^2 + \mathfrak{T}_R$$

(8.21)

with $\mathfrak{T}_R$ obeying $\mathfrak{T}_R \leq \varsigma$ if $t_*$ and $R_*$ are chosen from respective subsets in $[t, 2t]$ and $[\tfrac{1}{2}R, R]$ with measures $\tfrac{1}{2}t$ and $\tfrac{1}{4}R$ as the case may be.

*Part 3*: The explicit $-\tfrac{\pi m}{2} \ln(\tfrac{R}{t})$ terms in (8.20) and (8.21) cancel each other; and because of this, (8.19) leads to a curious inequality which is this:

$$2 \int_{[t_*,R_*]\times\mathbb{R}^2\times\{0\}} \tfrac{1}{s}(|\beta|^2 + |\mathfrak{b}|^2) + \tfrac{1}{2} \int_{\{R_*\}\times D_{rR_*}\times\{0\}} |\nabla^\perp \mathfrak{w}|^2 \leq \tfrac{1}{2} \int_{\{t_*\}\times D_{rt_*}\times\{0\}} |\nabla^\perp \mathfrak{w}|^2 + \varsigma$$

(8.22)

Here is what is curious about (8.22): This equation gives an R-independent upper bound for the integrals that appear on its left hand side since the integral on (8.22)'s right hand side is independent of R.

*Part 4:* As explained in Part 5, the curious feature of the inequality in (8.22) leads to a proof of Proposition 8.5. To set the stage for Part 5, let $\mathfrak{f}(\cdot)$ denote the function on $(0,\infty)$ whose value at any given t is the $\{t\}\times\mathbb{R}^2\times\{0\}$ integral of $|\beta|^2 + |\mathfrak{b}|^2$. An upper bound for the norm of the derivative of $\mathfrak{f}$,

$$|\tfrac{d}{dt}\mathfrak{f}| \leq \varsigma \tfrac{1}{t}\mathfrak{f} + \varsigma \mathfrak{f}^{1/2} \Big( \int_{\{t\}\times\mathbb{R}^2\times\{0\}} |F_A|^2 \Big)^{1/2},$$

(8.23)

can be had using the top bullets in (3.11) and (3.13) to write $\nabla_{\hat{A}t}\beta$ as $-2\alpha\beta + B_{A3}{}^+$; and using (6.6) to write $\nabla_{\hat{A}t}\mathfrak{b}$. Dividing both sides by $\mathfrak{f}^{1/2}$ and integrating from some time $t_0$ to a later value $t_1$ gives the following (by virtue of Lemma 5.3):

$$\mathfrak{f}^{1/2}(t_1) \geq \mathfrak{f}^{1/2}(t_0) - \varsigma \big(\ln \tfrac{t_1}{t_0} + (\tfrac{t_1}{t_0} - 1)^{1/2}\big).$$

(8.24)

Therefore, if $\varepsilon \in (0,1)$ and if $\mathfrak{f}(t_0) > \varepsilon$, then $\mathfrak{f}(t_1) > \tfrac{1}{2}\varepsilon$ if $t_1 \leq (1 + \tfrac{\varepsilon}{\varsigma})t_0$.

*Part 5*: Suppose that N is a given integer and that there are N times $\{t_0, t_1, \ldots, t_N\}$ with the following properties: First, $t_0 > t$, and then $t_i > 2t_{i-1}$ for each $i \in \{1, \ldots, N\}$. Second, the function $\mathfrak{f}$ at each $t_i$ is greater than $\varepsilon$. Then, by virtue of (8.24), one has

$$\int_{[t, 2t_N] \times \mathbb{R}^2 \times \{0\}} \tfrac{1}{s}(|\beta|^2 + |\mathfrak{b}|^2) \geq N \tfrac{\varepsilon^2}{\varsigma} \,.$$

(8.25)

Suppose for the sake of argument that N is large enough so that $N \tfrac{\varepsilon^2}{\varsigma}$ is greater than the number that is depicted on the right hand side of (8.22). Taking R on the left hand side of (8.22) to be greater than $2t_N$ then leads to a nonsensical inequality. (Keep in mind that the right hand side of (8.22) does not depend on R; it depends only on t.) Therefore, N is a priori bounded given t (which determines the value of the right hand side of (8.22)) and $\varepsilon$. Such an a priori bound on N leads directly proposition's assertion.

### e) Proof of Proposition 8.1

The proof has three steps.

Step 1: Proposition 8.5 is a given. Granted that, then Lemma 6.13 and Corollary 6.6 (or Proposition 6.8) imply that $\lim_{t \to \infty} t^{2-\varepsilon}(|\beta| + |\mathfrak{b}|) = 0$ for any $\varepsilon > 0$ with the limit being uniform with respect to the z coordinate.

Keeping this last observation in mind, let $\{\lambda_n\}_{n \in \mathbb{N}}$ denote any given increasing, unbounded sequence of positive numbers. For each $n \in \mathbb{N}$, let $(A^n, \mathfrak{a}^n)$ denote the pull-back of $(A, \mathfrak{a})$ by the coordinate rescaling diffeomorphism $(t, z, x_3) \to (\lambda_n t, \lambda_n z, \lambda_n x_3)$. Lemma 2.3 asserts that the the sequence $\{(A^n, \mathfrak{a}^n)\}_{n \in \mathbb{N}}$ has a subsequence that converges (after possibly termwise application of automorphisms of P) to a solution to (1.4) that also obeys the constraints in (6.1). This limit is denoted by $(A^\infty, \mathfrak{a}^\infty)$. Moreover, Lemma 6.9 can be brought to bear (because $t|\mathfrak{b}|$ limits to zero as $t \to \infty$) to see that the corresponding sequence $\{(\hat{A}^n, \mathfrak{b}^n, \alpha^n, \beta^n, \sigma^n, \varphi^n)\}_{n \in \mathbb{N}}$ converges in the $C^\infty$ topology on compact subsets of $(0, \infty) \times \mathbb{R}^2 \times \mathbb{R}$ (after possible termwise application of automorphisms of P) to the $(A^\infty, \mathfrak{a}^\infty)$ version of $(\hat{A}, \mathfrak{b}, \alpha, \beta, \sigma, \varphi)$. It then follows that the $(A^\infty, \mathfrak{a}^\infty)$ versions of $\beta$ and $\mathfrak{b}$ are zero; and so $(A^\infty, \mathfrak{a}^\infty)$ must be the integer *m* model solution $(A^{(m)}, \mathfrak{a}^{(m)})$ from Section 1c. (See Lemma 6.11 in this regard.)

Step 2: The limit in Step 1 is independent of the Lemma 2.2 subsequence and independent of the choice of scaling sequence $\{\lambda_n\}_{n \in \mathbb{N}}$ (as long it is increasing and unbounded.) As explained directly, this fact can be used to construct an automorphism of P on $(1, \infty) \times \mathbb{R}^2 \times \mathbb{R}$ to be denoted by $g_\infty$ such that

$$\lim_{t\to\infty} t(|g_\infty^*A - A^{(m)}| + |g_\infty^*\mathfrak{a} - \mathfrak{a}^{(m)}|) = 0$$

(8.26)

with the convergence being uniform with respect to the z coordinate of $\mathbb{R}^2$. There is also $C^k$ uniform convergence as $t \to \infty$ for $k \geq 1$ in the following sense:

$$\lim_{t\to\infty} t^{k+1}(|(\nabla_{A^{(m)}})^{\otimes k}(g_\infty^*A - A^{(m)})| + |(\nabla_{A^{(m)}})^{\otimes k}(g_\infty^*\mathfrak{a} - \mathfrak{a}^{(m)})|) = 0$$

(8.27)

To prove the preceding convergence assertions, invoke Lemma 2.3 with the sequence $\{\lambda_n = n\}_{n\in\mathbb{N}}$ as input. But for one small point, the assertion follows from what is said by Lemma 2.3 since the $\{t_n = \frac{1}{n}\}_{n\in\mathbb{N}}$ sequence of coordinate rescaling pull-backs of $(A, \mathfrak{a})$ converges to $(A^{(m)}, \mathfrak{a}^{(m)})$ after suitable termwise application of elements in Aut(P). There is no need to pass to a subsequence. The one small point is this: A literal reading of Lemma 2.3 yields the asserted bounds with a different automorphisms for each $n \in \mathbb{N}$ version of the domain where $t \in [\frac{1}{4n}, \frac{4}{n}]$. Assuming this to be the case, then the desired automorphism $g_\infty$ can be constructed from these $n \in \mathbb{N}$ labeled elements in Aut(P) by gluing successive ones together on the overlapping regions of the $\{t \in [\frac{1}{4n}, \frac{4}{n}]\}$ cover of the full domain $[1, \infty) \times \mathbb{R}^2 \times \{0\}$. This gluing construction is virtually identical to the construction that is used in the proof of Proposition 10.2 in [T] (which is modeled on a construction by Karen Uhlenbeck in [U].)

With regards to the z-coordinate uniformity of convergence: These instances of convergence, (8.26) and the various $k \in \mathbb{N}$ versions of (8.27), are uniform with respect to the z coordinate on $\mathbb{R}^2$ because Proposition 6.15's convergence is uniform with respect to the coordinate *x*.

Step 3: Given (8.26) and (8.27), then the only thing left to prove with regards to Proposition 8.1 is the assertion to the effect that the $t \to \infty$ limit of the natural logarithm of the ratio of $|\varphi|$ to $|\varphi^{(m)}|$ is zero. The latter assertion is an instance of Lemma 7.3 because $\varphi$ and $\varphi^{(m)}$ have the same $z = 0$ locus vanishing degree.

f) **The $\{t\} \times \mathbb{R}^2 \times \{0\}$ integral of $|\beta|^2 + |\mathfrak{b}|^2$ as $t \to 0$**

This subsection sets the stage for the proof of Proposition 8.2 (and the ultimate proof of Theorem 2) by stating and then proving three lemmas. The following (immediate) corollary to Propositions 8.5 and 6.14 should be kept in mind by way of motivation for the subsequent lemmas.

**Corollary 8.6**: *Given positive numbers $\mathfrak{z}$, $\delta$, $r$ and $\Xi$, there exists $\kappa > 1$ with the following significance: Let $(A, \mathfrak{a})$ denote a solution to (1.4) that obeys the three bullets of (6.1) with the given $\mathfrak{z}$, $\delta$, $r$ and $\Xi$. Suppose in addition that*

$$\int_{\{t\}\times\mathbb{R}^2\times\{0\}} (|\beta|^2 + |\mathfrak{b}|^2) \leq \tfrac{1}{\kappa}$$

*when $t$ is positive and sufficiently small. Then $(A, \mathfrak{a})$ is one of the model solutions from Section 1c.*

The first lemma in this subsection stems from the curious inequality in (8.22).

**Lemma 8.7**: *Let $(A, \mathfrak{a})$ denote a given solution to (1.4) that is described by (6.1). Then*

$$\limsup_{t\to 0} \int_{\{t\}\times\mathbb{R}^2\times\{0\}} (|\beta|^2 + |\mathfrak{b}|^2) = 0$$

*if there exists a number $\mathcal{E}$ and a decreasing sequence of times $\{t_n\}_{n\in\mathbb{N}} \subset (0, 1]$ with limit zero such that $\int_{\{t_n\}\times D_{r\, t_n}\times\{0\}} |\nabla^\perp \mathfrak{w}|^2 \leq \mathcal{E}$ for all $n \in \mathbb{N}$.*

*Proof of Lemma 8.7*: If the assumptions hold, then the argument in Part 5 of the proof of Proposition 8.5 can be applied with only cosmetic changes (taking $t_*$ to be successive versions of $t_n$ with R fixed and independent of n). Indeed, if N is an integer obeying $N\frac{\varepsilon^2}{\varsigma} \geq \mathcal{E}+\varsigma$, then (8.22) will be violated if there are N times $\{t_n\}_{n\in\{1,\ldots,N\}} \subset (0, \tfrac{1}{2} R)$ obeying $t_n < \tfrac{1}{2} t_{n-1}$ and $\mathfrak{f}(t_n) > \varepsilon$.

To set the stage for the second lemma, introduce by way of notation $g$ to denote the function on $(0, \infty)$ defined by the rule

$$t \to g(t) = \int_{\{t\}\times D_{r\, t}\times\{0\}} |\nabla^\perp \mathfrak{w}|^2$$

(8.28)

Nothing has been said at this point to exclude the event that $g$ has infinite lim-inf as $t \to 0$.

**Lemma 8.8**: *Let $(A, \mathfrak{a})$ denote a solution to (1.4) that is described by (6.1). Suppose that $\{t_n\}_{n\in\mathbb{N}} \subset (0, 1]$ is a decreasing sequence of times with limit zero. The values of the function $g$ on a subsequence of $\{t_n\}_{n\in\mathbb{N}}$ are bounded if and only if the following occurs: Let $\{(A^n, \mathfrak{a}^n)\}_{n\in\mathbb{N}}$ denote the sequence of solutions to (1.4) whose n'th term is the pull-back* *of* $(A, \mathfrak{a})$ *by the coordinate rescaling diffeomorphism* $(t, z, x_3) \to (t_n t, t_n z, t_n x_3)$. *Give the sequence* $\{(A^n, \mathfrak{a}^n)\}_{n \in \mathbb{N}}$ *to Lemma 2.3 and it outputs at least one limit solution with a corresponding version of* $\varphi$ *whose* $z = 0$ *locus vanishing degree is the same as the* $z = 0$ *locus vanishing degree of the* $(A, \mathfrak{a})$ *version of* $\varphi$.

**Proof of Lemma 8.8**: The proof of has three parts.

*Part 1*: Some background for what is to come: Having fixed $t \in (0, \infty)$, multiply both sides of (7.4) at time t by the function $\chi_{rt} \mathfrak{w}$ and integrate the result over $\{t\} \times \mathbb{R}^2 \times \{0\}$. (The function $\chi_{rt}$ equals 1 where $|z| \leq rt$ and it equals zero where $|z| \geq 2rt$.) An instance of integration by parts (and use of Lemmas 2.2 and 7.4) leads to this:

$$\int_{\{t\} \times D_{rt} \times \{0\}} |\nabla^{\perp} \mathfrak{w}|^2 \leq \varsigma + 4 \int_{\{t\} \times D_{rt} \times \{0\}} |\mathfrak{w}||\mathfrak{b}|^2 \ .$$

(8.29)

Observedly, if a given number $\mathcal{Z}$ is such that $\mathfrak{w}$ at time t obeys a bound $|\mathfrak{w}| \leq \mathcal{Z}$, then the right hand side of (8.29) will be less than $\varsigma(1 + \mathcal{Z})$; and thus so will the left hand side. By the same token, if a given number $\mathcal{Z}$ is such that the norm of $\mathfrak{b}$ at time t obeys a bound $|\mathfrak{b}| \leq \frac{\mathcal{Z}}{t}$, then the right hand side of (8.29) will again be less than $\varsigma(1 + \mathcal{Z})$. This is because the integral of $|\mathfrak{w}|$ on the $|z| \leq rt$ disk in $\{t\} \times \mathbb{R}^2 \times \{0\}$ is bounded by $\varsigma$ (see Lemma 7.6). Also, keep in mind that there is an a priori bound of this sort for $|\mathfrak{b}|$ if the integral of $|\mathfrak{b}|^2$ on $\{t\} \times \mathbb{R}^2 \times \{0\}$ is less than $\frac{1}{\varsigma}$ (see Lemmas 6.1 and 6.7).

*Part 2*: Returning now to the context of Lemma 8.8, let $(A^{\infty}, \mathfrak{a}^{\infty})$ denote a limit sequence that is supplied by Lemma 2.3 from the sequence $\{(A^n, \mathfrak{a}^n)\}_{n \in \mathbb{N}}$. Suppose in this part of the proof that the $(A^{\infty}, \mathfrak{a}^{\infty})$ and $(A, \mathfrak{a})$ versions of $\varphi$ have the same degree of vanishing on the $z = 0$ locus. Denote this degree by *m*. Let $\Lambda \subset \mathbb{N}$ denote the subsequence that is associated by Lemma 2.3 to the limit $(A^{\infty}, \mathfrak{a}^{\infty})$. Under the stated circumstances, the sequence that is indexed by $\Lambda$ whose n'th term is the norm of the *m*'th order covariant derivative of the $(A^n, \mathfrak{a}^n)$ version of $\varphi$ at $(t = 1, z = 0)$ has a positive lower bound. (This is because if the manner of $C^{\infty}$ convergence in Lemma 2.3 and because there is a positive lower bound for the $(A^{\infty}, \mathfrak{a}^{\infty})$ version of the norm for the analogous *m*'th order covariant derivative). The existence of such a positive lower bound with (8.1) implies via Taylor's theorem with remainder (and Lemma 2.2) that the sequence indexed by $\Lambda$ whose n'th element is the supremum norm on $\{1\} \times \mathbb{R}^2 \times \{0\}$ of the $(A^n, \mathfrak{a}^n)$ version of 𝔴 is bounded. This last observation leads via (8.29) to an an upper bound for the values of the function $g$ on the subsequence (see what is said in Part 1).

*Part 3*: Now suppose that the z = 0 locus vanishing degree of the $(A^\infty, \mathfrak{a}^\infty)$ version of φ is strictly greater than $m$. Denote this greater vanishing degree by $k$. If z ≠ 0, then the n ∈ Λ versions of the function $w$ will converge at (t = 1, z) to the $(A^\infty, \mathfrak{a}^\infty)$ versions (because $|\varphi| = \frac{1}{\sqrt{2t}} e^{2w}$ and the sequence $\{|\varphi^n|\}_{n \in \Lambda}$ converges to $|\varphi^\infty|$). This convergence of the $w$'s at t = 1 for any non-zero z can be viewed from the perspective of the respective $(A^\infty, \mathfrak{a}^\infty)$ and $\{(A^n, \mathfrak{a}^n)\}_{n \in \Lambda}$ versions of (8.1) at t = 1. There is some tension with regard to this perspective (because $k \neq m$) which resolves as follows: Given any z ≠ 0 with norm less than $\frac{1}{\varsigma}$, all sufficiently large n ∈ Λ versions of 𝔴 at (t = 1, z) obey $\mathfrak{w} \leq \frac{k-m}{4} \ln|z|$. Noting this, an instance of the fundamental theorem of calculus can be invoked (with the Cauchy-Schwarz inequality) to prove the following: If r ∈ (0, 1), and if some n ∈ Λ version of 𝔴 is less than $\frac{k-m}{4} \ln r$ at all points where t = 1 and |z| = r, then

$$c_{0-} |\ln r| \leq \int_{\{1\} \times D_1 \times \{0\}} |\nabla^\perp \mathfrak{w}|^2 .$$

(8.30)

(The norm of the integral of $\frac{\partial}{\partial |z|} |\mathfrak{w}|$ from r to $\frac{1}{\varsigma}$ along any ray from the origin in $\{1\} \times D_1 \times \{0\}$ is greater $\frac{k-m}{4} |\ln r|$.)

### g) Proof of Proposition 8.2

The proof of the proposition has four parts.

*Part 1*: A function on (0, 1] is called *purely unbounded* if its t → 0 lim-inf is infinite. If $(A, \mathfrak{a})$ is a solution to (1.4) that is described by (6.1), then the corresponding version of the function $g$ as defined in (8.28) is either purely unbounded or not. If it is not purely unbounded, then Lemma 8.7 in conjunction with Propositions 8.4 and 6.14 can be brought to bear to see that $(A, \mathfrak{a})$ is a model solution from Section 1c. The conclusions of Proposition 8.2 hold by default in this case. With this understood, assume henceforth that $(A, \mathfrak{a})$'s version of the function $g$ is purely unbounded.

A pair $(A^\ddagger, \mathfrak{a}^\ddagger)$ of connection on P and 1-form valued section of ad(P) is said to be a *t → 0 rescaling limit* of $(A, \mathfrak{a})$ when it arises as a Lemma 2.3 limit solution from for an input sequence that is defined as follows: Let $\{t_n\}_{n \in \mathbb{N}}$ denote a decreasing sequence in (0, 1] with limit zero. The n'th element of the desired input sequence is the pull-back of $(A, \mathfrak{a})$ via the coordinate rescaling diffeomorphism $(t, z, x_3) \to (t_n t, t_n z, t_n x_3)$. Let $\mathfrak{S}$ denote the set of all t → 0 rescaling limits of $(A, \mathfrak{a})$.

There are three features of $\mathfrak{S}$ that should be kept in mind: First, every pair from $\mathfrak{S}$ is described by (6.1) using the same data $\mathfrak{z}$, $\delta$, $r$ and $\Xi$ as used for $(A,\mathfrak{a})$. The second feature is this: If $(A^{\ddagger},\mathfrak{a}^{\ddagger})$ is from $\mathfrak{S}$, then the $z = 0$ locus vanishing degree of the corresponding version of $\varphi$ is strictly greater then the $z = 0$ locus vanishing degree of the $(A,\mathfrak{a})$ version of $\varphi$ (see Lemma 8.8). As for the third feature: If $(A^{\ddagger},\mathfrak{a}^{\ddagger})$ is from $\mathfrak{S}$, then every $t \to 0$ rescaling limit of $(A^{\ddagger},\mathfrak{a}^{\ddagger})$ is also in $\mathfrak{S}$. (This follows from the definition of the term *$t \to 0$ rescaling limit* and from the manner of convergence to the limit that is dictated by Lemma 2.3.) By virtue of these features (and by virtue of Lemma 5.6), the set $\mathfrak{S}$ can be said to have a *maximal* element. This an element whose version of $\varphi$ has the maximal $z = 0$ vanishing degree of the versions of $\varphi$ from the elements in $\mathfrak{S}$.

*Part 2*: This part of the proof states and then proves a proposition that characterizes maximal elements in the set $\mathfrak{S}$.

**Lemma 8.9**: *Let $(A,\mathfrak{a})$ denote a solution to (1.4) that is described by (6.1) whose version of the function $\mathfrak{g}$ is purely unbounded. Any maximal element in the corresponding version of the set $\mathfrak{S}$ is necessarily a model solution from Section 1c.*

Supposing that this lemma is correct, write the relevant model solution as $(A^{(k)},\mathfrak{a}^{(k)})$. Keep in mind that the integer $k$ is greater than the $z = 0$ locus degree of vanishing of the $(A, \mathfrak{a})$ version of $\varphi$ (which is henceforth denoted by $m$).

*Proof of Lemma 8.9*: Let $(A^{\ddagger},\mathfrak{a}^{\ddagger})$ denote a maximal element in $\mathfrak{S}$. Because $(A^{\ddagger},\mathfrak{a}^{\ddagger})$ is maximal, the $z = 0$ vanishing degree of its version of $\varphi$ can not be less than the $z = 0$ vanishing degree of the version of $\varphi$ that comes from any $t \to 0$ rescaling limit of $(A^{\ddagger},\mathfrak{a}^{\ddagger})$. Therefore, by virtue of Lemma 8.8, the $(A^{\ddagger},\mathfrak{a}^{\ddagger})$ version of the function $\mathfrak{g}$ is not purely unbounded. Therefore, Lemma 8.7 is in play (keep in mind that $(A^{\ddagger},\mathfrak{a}^{\ddagger})$ is described by (6.1)); and Lemma 8.7 plus Propositions 8.5 put Proposition 6.14 in play. Proposition 6.14 says that $(A^{\ddagger},\mathfrak{a}^{\ddagger})$ is a model solution from Section 1c.

*Part 3*: This part proves that the set $\mathfrak{S}$ has just the one element $(A^{(k)},\mathfrak{a}^{(k)})$. To this end, consider the sequence $\{t_n = \frac{1}{n}\}_{n \in \mathbb{N}}$. There is a subsequence $\Lambda \subset \mathbb{N}$ such that the corresponding pull-back sequence $\{(A^n,\mathfrak{a}^n)\}_{n \in \Lambda}$ converges in the $C^{\infty}$ topology on compact subsets of $(0,\infty) \times \mathbb{R}^2 \times \mathbb{R}$ to $(A^{(k)},\mathfrak{a}^{(k)})$ after termwise application of elements in $\mathrm{Aut}(P)$. Thus, given $\varepsilon > 0$, there exists an integer $N_\varepsilon$ such that if $n \in \Lambda$ and $n > N_\varepsilon$, then

$$|g_n{}^*A^n - A^{(k)}| + |g_n{}^*\mathfrak{a}^n - \mathfrak{a}^{(k)}| < \varepsilon$$

(8.31)

where $t \in (\varepsilon, \frac{1}{\varepsilon})$ and $|z| < \frac{1}{\varepsilon}$. Here, $g_n$ is a suitable automorphism of P.

Now suppose that there is element in $\mathfrak{S}$ that is not the maximal element $(A^{(k)}, \mathfrak{a}^{(k)})$ to derive nonsense. Proposition 8.1 applies to this element because it is described by (6.1). In particiular, it follows from Proposition 8.1 that the $t \to \infty$ limit of the corresponding version of $t|\mathfrak{a}_3|$ on the $z = 0$ locus can be written as $\frac{k'+1}{2}$ with $k'$ being an integer that is strictly less than $k$.

The preceding has the following implication: For each $n \in \Lambda$, there is a smallest integer $n' > n$ such that

$$|t|\mathfrak{a}_3^{n'}| - \frac{k+1}{2}| \geq \frac{1}{32}$$
(8.32)

at some point on the $z = 0$ locus where $t \in [1, 2]$. Denote this smallest $n'$ by $n'(n)$. An important note in this regard: The ratio $\frac{n}{n'(n)}$ must limit to zero as $n$ gets ever larger to be consistent with (8.31). Use the corresponding sequence $\{(A^{n'(n)}, \mathfrak{a}^{n'(n)})\}_{n \in \Lambda}$ as input for Lemma 2.3 and let $(A^\ddagger, \mathfrak{a}^\ddagger)$ denote a limit solution to (1.4) from Lemma 2.3. This solution obeys

- $|t|\mathfrak{a}_3^\ddagger| - \frac{k+1}{2}| \geq \frac{1}{32}$ *at some $z = 0$ point where $t \in [1, 2]$*,
- $|t|\mathfrak{a}_3^\ddagger| - \frac{k+1}{2}| \leq \frac{1}{32}$ *on the $t > 2$ part of the $z = 0$ locus*.

(8.33)

(The top bullet follows from (8.32) and the second bullet is due to the minimality of each $n'(n)$ and due to the fact that $\lim_{n \in \Lambda} \frac{n}{n'(n)}$ is zero.

Note in particular that the top bullet of (8.33) implies that $(A^\ddagger, \mathfrak{a}^\ddagger)$ can't be Aut(P) equivalent to $(A^{(k)}, \mathfrak{a}^{(k)})$. As a consequence, the $z = 0$ locus degree of vanishing of $(A^\ddagger, \mathfrak{a}^\ddagger)$ is some integer $k' < k$. Meanwhile, Proposition 8.1 applies to $(A^\ddagger, \mathfrak{a}^\ddagger)$ because $(A^\ddagger, \mathfrak{a}^\ddagger)$ is described by (6.1). Therefore, its $t \to \infty$ limit is necessarily a model solution from Section 1c whose version of $\varphi$ has the same $z = 0$ locus degree of vanishing as $\varphi^\ddagger$ (which is $k'$). In particular, Proposition 8.1 implies that $\lim_{t \to \infty} t|\mathfrak{a}_3^\ddagger| = \frac{k'+1}{2}$ which is the desired nonsense: It fouls the second bullet of (8.33) because $k' \neq k$.

*Part 4*: Except for one subtle point, the argument for Proposition 8.2 when $\mathfrak{S}$ contains only the one element $(A^{(k)}, \mathfrak{a}^{(k)})$ differs only cosmetically from what is said in Step 2 of Section 8e. The subtle point is the implication in Proposition 8.2 to the effect that $k - m$ is necessarily an *even* integer.

To prove that $k - m$ is even, fix some very small, positive $\varepsilon$ and then a time $t$ where the following holds on $\{t\} \times \mathbb{R}^2 \times \{0\}$ for $\ell = 0$ and for $\ell = 1$:

$$t^{k+1}|(\nabla_{A^{(k)}})^{\otimes \ell}(g_0^*A - A^{(k)})| + |(\nabla_{A^{(k)}})^{\otimes \ell}(g_\infty^*\mathfrak{a} - \mathfrak{a}^{(k)})| < \varepsilon$$

(8.34)

Identify an orthonormal frame for ad(P) at the $z = 0$ point in $\{t\}\times\mathbb{R}^2\times\{0\}$ with the basis vectors $\{\sigma_1, \sigma_2, \sigma_3\}$ that appear in (1.6). Parallel transporting this frame along the rays from the origin in $\{t\}\times\mathbb{R}^2\times\{0\}$ using the connection $A^{(k)}$ defines an orthonormal basis for ad(P) on the whole of $\{t\}\times\mathbb{R}^2\times\{0\}$. Denote this basis by $\{\sigma_1, \sigma_2, \sigma_3\}$ also. This basis will be used to view any given section of ad(P) over $\{t\}\times\mathbb{R}^2\times\{0\}$ as a map from $\{t\}\times\mathbb{R}^2\times\{0\}$ to $\mathbb{R}^3$. Just to be sure: If $\eta$ is a section of ad(P), then the $\mathbb{R}^3$ components of the corresponding map are $(\langle\sigma_1\eta\rangle, \langle\sigma_2\eta\rangle, \langle\sigma_3\eta\rangle)$.

There are two maps to $\mathbb{R}^3$ of interest. The first is given by $g_0^*\mathfrak{a}_3$ (which is a section of ad(P) and thus a map to $\mathbb{R}^3$.) This map is nowhere zero if $\varepsilon$ is small because $t|\mathfrak{a}_3^{(k)}| \geq \frac{1}{2}$. The second map is $g_0^*\sigma$; which has unit length so it too is nowhwere zero. The two maps $\frac{1}{|\mathfrak{a}_3|}g_0^*\mathfrak{a}_3$ and $g_0^*\sigma$ are nearly equal on the $|z| = \frac{1}{\sqrt{\varepsilon}}t$ circle in $\{t\}\times\mathbb{R}^2\times\{0\}$ if $\varepsilon$ is sufficiently small because $t|\beta| \leq c_0\varepsilon^{1/8}$ where $|z| = \frac{1}{\sqrt{\varepsilon}}t$. (See Lemma 6.3.) A map that agrees with $g_0^*\sigma$ where $|z| = \frac{1}{\sqrt{\varepsilon}}t$ can be made from the modification

$$\mathfrak{a}_3' = \mathfrak{a}_3 - (1-\chi(2\sqrt{\varepsilon}\tfrac{|z|}{t} - 1))\beta.$$

(8.35)

The resulting map is $\frac{1}{|\mathfrak{a}_3'|}g_0^*\mathfrak{a}_3'$. It is denoted by $\tau$. The map $\tau$ is likewise nowhere zero if $\varepsilon$ is small and it is equal to $g_0^*\sigma$ where $|z| = \frac{1}{\sqrt{\varepsilon}}t$. (Keep in mind that the composition of $(1-\chi)$ with $2\sqrt{\varepsilon}\frac{|z|}{t} - 1$ is equal to 0 where $|z| < \frac{1}{2\sqrt{\varepsilon}}t$ and it is equal to 1 where $|z| > \frac{1}{\sqrt{\varepsilon}}t$.)

Let $D^\tau$ denote the $|z| \leq \frac{1}{\sqrt{\varepsilon}}t$ disk in $\mathbb{R}^2$ and let $D^\sigma$ denote a second copy. View $S^2$ via stereographic projection as $\mathbb{C}\cup\{\infty\}$. The disk $D^\tau$ is viewed as the 'northern' hemisphere of $S^2$ via the map from $D^\tau$ to $\mathbb{C}\cup\{\infty\}$ that sends $z$ to $\frac{\sqrt{\varepsilon}z}{t}$. Meanwhile, $D^\sigma$ is viewed as the 'southern' hemisphere of $S^2$ via the map from $D^\sigma$ to $\mathbb{C}\cup\{\infty\}$ that sends $z$ to $\frac{tz}{\sqrt{\varepsilon}|z|^2}$. A map (to be denoted by $y$) from this same $S^2$ to the unit radius sphere in $\mathbb{R}^3$ is defined by the rule whereby $y$ on $D^\sigma$ is $g_0^*\sigma$ and $y$ on $D^\tau$ is $\tau$.

The section $\varphi^{(k)} - \tau\langle\tau\varphi^{(k)}\rangle$ is orthogonal to $\tau$ so it defines a section of $y^*TS^2$ on the northern hemisphere (the $D^\tau$ part of $S^2$). Meanwhile, the section $g_0^*\varphi$ is orthogonal to $g_0^*\sigma$ so it defines a section of $y^*TS^2$ on the southern hemisphere (the $D^\sigma$ part of $S^2$). Meanwhile, The former section, $\varphi^{(k)} - \tau\langle\tau\varphi^{(k)}\rangle$, vanishes only at the north pole and it vanishes there with degree $k$ if $\varepsilon$ is small. The latter section, $g_0^*\varphi$, vanishes only at the south pole and with vanishing degree is $-m$.

Near the equator of $S^2$ (where $|z| = \frac{1}{\sqrt{\varepsilon}} t$ on $D^\tau$ and $D^\sigma$), the both $\varphi^{(k)} - \tau\langle\tau\varphi^{(k)}\rangle$ and $g_0^*\varphi$ are uniformly large (their respective norms differ by at most $c_0 \frac{\varepsilon}{t}$ from $\frac{1}{\sqrt{2}t}$). Moreover, the two sections almost agree. This is because both have nearly the same norm on the $|z| = \frac{1}{\sqrt{\varepsilon}} t$ circle in $\{t\} \times \mathbb{R}^2 \times \mathbb{R}$ and because both are nearly $\nabla_{A^{(k)}}$- covariantly constant along this circle. (This last claim follows because

$$t^2 |\nabla_A^\perp \varphi| \le c_\jmath \varepsilon \;\; and \;\; t^2 |\nabla_{A^{(k)}}^\perp \varphi^{(k)}| \le c_\jmath \varepsilon \;\; and \;\; t|g_0^* A - A^{(k)}| \le \varepsilon$$

(8.36)

on the $|z| = \frac{1}{\sqrt{\varepsilon}} t$ circle in $\{t\} \times \mathbb{R}^2 \times \{0\}$ and because the radius of this circle is much less than $\frac{1}{\varepsilon} t$ when $\varepsilon$ is small.) Because $\varphi^{(k)} - \tau\langle\tau\varphi^{(k)}\rangle$ and $g_0^*\varphi$ are very close on this circle (assuming $\varepsilon$ is small), the function $\chi$ can be used (as was done in the construction of $\tau$) to obtain a section of $y^*TS^2$ on the whole of $S^2$ that agrees with $\varphi^{(k)} - \tau\langle\tau\varphi^{(k)}\rangle$ where $|z| < \frac{1}{2\sqrt{\varepsilon}} t$ on $D^\tau$, agrees with $g_0^*\varphi$ where $|z| < \frac{1}{2\sqrt{\varepsilon}} t$ on $D^\sigma$, and vanishes only at the poles. Denote this section by $\psi$.

The Euler class of the bundle $y^*TS^2$ equal to the sum of the vanishing degrees of $\psi$ at its zeros, thus $k - m$. The Euler class is also equal to the product of the degree of the map $y$ with the Euler class of $TS^2$. Since the latter is equal to 2, the integer $k - m$ is even.

## 9. Proof of Theorem 2

By virtue of Proposition 8.1, 8.2 and 6.15, any solution to (1.4) that is described by (6.1) is also described by the three bullets of Theorem 2. By virtue of Proposition 5.1, any CONSTRAINT SET 2 solution to (1.4) is described by (6.1). Thus, all CONSTRAINT SET 2 solutions to (1.4) obey the three bullets of Theorem 2.

Theorem 2 makes one additional assertion which is implied by the following proposition (via Proposition 5.1).

**Proposition 9.1**: *A solution to (1.4) that is described by (6.1) is a model solution from Section 1c if and only if there exists $t_0 > 0$ such that the corresponding version of $\langle \mathfrak{a}_3[\mathfrak{a}_1, \mathfrak{a}_2] \rangle$ is non-negative on the $t < t_0$ part of $(0, \infty) \times \mathbb{R}^2 \times \mathbb{R}$.*

This proposition is proved in Section 9b. Section 9a makes some preliminary observations that are used in the proof.

### a) Implications of t → 0 convergence

Let $(A, \mathfrak{a})$ denote here a solution to (1.4) that obeys the constraints in (6.1). Use $m$ to denote the $z = 0$ locus vanishing degree of $\varphi$. Proposition 8.2 is in play for $(A, \mathfrak{a})$ because (6.1) is obeyed. Proposition 8.2 gives an integer $k \geq m$ and a continuous, increasing function $\underline{\varepsilon}: (0, 1] \to (0, \infty)$ with limit zero as $t \to 0$, these such that

- $\left| |\mathfrak{a}_3| - \frac{1}{2t} \frac{(k+1)\sinh(\Theta)}{\sinh((k+1)\Theta)} \frac{\cosh((k+1)\Theta)}{\cosh(\Theta)} \right| \leq \frac{1}{t} \underline{\varepsilon}(t)$,
- $\left| |\varphi| - \frac{1}{2t} \frac{(k+1)\sinh(\Theta)}{\sinh((k+1)\Theta)} \right| \leq \frac{1}{t} \underline{\varepsilon}(t)$.
- $|\langle \mathfrak{a}_3 \varphi \rangle| \leq \frac{1}{t^2} \underline{\varepsilon}(t)$,
- $|\langle \mathfrak{a}_3 \nabla_A \varphi \rangle| \leq \frac{1}{t^3} \underline{\varepsilon}(t)$,

(9.1)

These inequalities follow from Proposition 8.2 because they hold with $\underline{\varepsilon} \equiv 0$ for the model solution $(A^{(k)}, \mathfrak{a}^{(k)})$ (see (1.7)). The preceding bounds have the following implication:

**Lemma 9.2**: *Given $\varepsilon > 0$, there exists $t_\varepsilon > 0$ such that if $t < t_\varepsilon$ and $|z| > \varepsilon t$ then*

- $\left| \alpha + \frac{1}{2t} \frac{(k+1)\sinh(\Theta)}{\sinh((k+1)\Theta)} \frac{\cosh((k+1)\Theta)}{\cosh(\Theta)} \right| \leq \frac{1}{t} \varepsilon$.
- $|\beta| \leq \frac{1}{t} \varepsilon$.
- $|\mathfrak{b}| \leq \frac{1}{t} \varepsilon$.

*Proof of Lemma 9.2*: The second bullet of the lemma follows from the second and third bullet in (9.1) because $|\langle \mathfrak{a}_3 \varphi \rangle| = |\beta| |\varphi|$. To prove the lemma's third bullet, write $\langle \mathfrak{a}_3 \nabla_A \varphi \rangle$ as the sum of $2i\alpha \langle \mathfrak{b} \varphi \rangle + \langle \beta \nabla_{\hat{A}}^\perp \varphi \rangle$. (Remember that $\mathfrak{b} = \frac{1}{4}[\sigma, \nabla_A \sigma]$.) Then, use the fact that $|\nabla_{\hat{A}}^\perp \varphi| \leq \varsigma \frac{1}{t^3}$ (which is a consequence of Lemma 2.2) and the fact that $t|\beta|$ is small where $t$ is small to see that $\langle \mathfrak{a}_3 \nabla_A \varphi \rangle$ will differ from $2i\alpha \langle \mathfrak{b} \varphi \rangle$ by at most a very small multiple of $\frac{1}{t^3}$ when $t$ is small (given a positive, $t$-independent lower bound for $|z|$.) Use this last fact with the top two bullets of (9.1) to see that $t|\mathfrak{b}|$ must be small when $t$ is small given a positive, $t$-independent upper bound for $|z|$. (Remember that $\mathfrak{b} = \frac{1}{2}(\mathfrak{b}_1 + i\mathfrak{b}_2)$.)

The top bullet of the lemma follows from the top bullet of (9.2) because if $|\beta|$ is small, then $\alpha$ is either very near one of $\pm \frac{1}{2t} \frac{(k+1)\sinh(\Theta)}{\sinh((k+1)\Theta)} \frac{\cosh((k+1)\Theta)}{\cosh(\Theta)}$; and it can't be the + sign version for $t$ very small (given a positive, $t$-independent upper bound for $|z|$) without fouling the second bullet in (9.1) after integrating the top bullet in (3.9).

The following consequence of (9.1) is needed for the proof of Proposition 9.1: Given $\varepsilon > 0$, there exists $t_\varepsilon > 0$ such that assertions below hold where $t < t_\varepsilon$ and $|z| < t$.

$$\nabla_{At} \mathfrak{a}_3 = -\frac{1}{t} \mathfrak{a}_3 + \mathfrak{e}_1 \frac{1}{t^2} \quad \text{with} \quad |\mathfrak{e}_1| \leq \varsigma \frac{|z|^2}{t^2} + \varepsilon.$$

(9.2)

This follows from (9.2) because it holds $(A^{(k)}, \mathfrak{a}^{(k)})$ with $\underline{\varepsilon} \equiv 0$. (See the formulas in (1.7).)

Ramifications of the identity in (9.2) are explored by writing $\mathfrak{a}_3$ and $A_t$ as done in Section 3: If $\mathfrak{a}_3$ is written as $\alpha\sigma + \beta + \beta^*$ and if $A_t$ is written as $\hat{A}_t - i(\beta - \beta^*)$, then the respective $\sigma$ and $\mathcal{L}^+$ parts of (9.2)'s identity can be written in turn as

- $\frac{\partial}{\partial t}\alpha - 4|\beta|^2 = -\frac{1}{t}\alpha + \mathfrak{e}_2 \frac{1}{t^2}$  with $|\mathfrak{e}_2| \leq \varsigma \frac{|z|^2}{t^2} + c_0 \varepsilon$.
- $\nabla_{\hat{A}t}\beta + 2\alpha\beta = -\frac{1}{t}\beta + \mathfrak{e}_3 \frac{1}{t^2}$  with $|\mathfrak{e}_3| \leq \varsigma \frac{|z|^2}{t^2} + c_0 \varepsilon$.

(9.3)

The top bullet in (9.2) plays the starring role in the proof of Proposition 9.1.

**b) Proof of Proposition 9.1**

The proof that $\langle \mathfrak{a}_3, [\mathfrak{a}_1, \mathfrak{a}_2]\rangle \geq 0$ for a model solution, say $(A^{(m)}, \mathfrak{a}^{(m)})$, is straightforward: Use (1.7) to see that this function is

$$\frac{1}{2t^3}\left(\frac{(m+1)\sinh(\Theta)}{\sinh((m+1)\Theta)}\right)^2 \frac{\cosh((m+1)\Theta)}{\cosh(\Theta)}$$

(9.4)

which is manifestly non-negative.

The converse assertion is this: If $(A, \mathfrak{a})$ is not a model solution from Section 1c, then the function $\langle \mathfrak{a}_1, [\mathfrak{a}_2, \mathfrak{a}_3]\rangle$ is somewhere negative on each sufficiently small, positive t version of $\{t\} \times \mathbb{R}^2 \times \mathbb{R}$. The proof of this converse assertion has four parts.

*Part 1*: If there is a positive time $t_0$ such that $\alpha < -\frac{1}{4t}$ where $t < t_0$, then Lemma 6.5 is in play: Lemma 6.5 says that the various $\{t\}\times\mathbb{R}^2\times\{0\}$ integrals of $|\beta|^2 + |\mathfrak{b}|^2$ limit to zero as t does. Granted that, then Corollary 8.6 is in play; and Corollary 8.6 says that $(A, \mathfrak{a})$ is one of the model solutions from Section 1c. Thus, if $(A, \mathfrak{a})$ is not a model solution from Section 1c, there is a decreasing sequence of positive times $\{t_n\}_{n\in\mathbb{N}}$ with limit 0 and a corresponding sequence of points $\{z_n\}_{n\in\mathbb{N}} \subset \mathbb{R}^2$ such that $\alpha(t_n, z_n, 0) \geq -\frac{1}{4t}$.

*Part 2*: The very small t part of the locus in $(0,\infty) \times \mathbb{R}^2 \times \mathbb{R}$ where $\alpha \geq -\frac{1}{4t}$ has to lie where $\frac{|z|}{t}$ is very small. That this must be so follows from the top bullet of Lemma 9.2 because Lemma 9.2's top bullet implies the following: Given $\varepsilon > 0$, there exists $t_\varepsilon$ such that $\alpha < -\frac{1}{2t}(1-\varepsilon)$ where $t < t_\varepsilon$ and $|z| > \varepsilon t$ on $(0, \infty)\times\mathbb{R}^2\times\mathbb{R}$.

The preceding has the following implication: If $\{(t_n, z_n)\}_{n\in\mathbb{N}} \in (0,\infty)\times\mathbb{R}^2$ is such that $\{t_n\}_{n\in\mathbb{N}}$ is decreasing to zero and such that $\alpha(t_n, z_n, 0) \geq -\frac{1}{4t_n}$, then $\lim_{n\to\infty} \frac{|z_n|}{t_n} = 0$.

*Part 3*: The top bullet in (9.3) can be rewritten using the top bullet in (9.1) to say the following: If $t < t_\varepsilon$ and $|z| \leq t$, then

$$\tfrac{\partial}{\partial t}\alpha + 2\alpha(\alpha+\tfrac{1}{2t}) - \tfrac{(k+1)^2}{2t^2} = \mathfrak{e}_4 \tfrac{1}{t^2} \quad \text{with} \quad |\mathfrak{e}_4| \leq \varsigma(\varepsilon + \tfrac{|z|^2}{t^2}).$$

(9.5)

This equation and the top bullet in Lemma 9.2 have the following immediate implication: There exists $t_\diamond > 0$ such that

$$\tfrac{\partial}{\partial t}\alpha \geq \tfrac{(k+1-\tfrac{1}{1000})^2}{2t^2}$$

(9.6)

where $-\tfrac{1}{4t} \leq \alpha \leq 0$ in the $t < t_\diamond$ part of $(0,\infty) \times \mathbb{R}^2 \times \mathbb{R}$ (assuming that it is not empty). This has the following consequence: If $t$ is very small (less than $t_\diamond$) and if $\alpha$ is greater than $-\tfrac{1}{4t}$ at a point $(t,z,0)$, then $\alpha(\cdot,z,0)$ is greater than 0 on the whole of the interval $[3t, t_\diamond]$. This is to say that if $\alpha$ is greater than $-\tfrac{1}{4t}$ at $(t,z,0)$, then it very rapidly (within time $2t$) becomes positive and stays positive until some $\mathcal{O}(1)$ time.

*Part 4*: Suppose that $\{(t_n, z_n)\}_{n \in \mathbb{N}} \in (0,\infty) \times \mathbb{R}^2$ is such that $\{t_n\}_{n \in \mathbb{N}}$ is decreasing to zero and such that $\alpha(t_n, z_n, 0) \geq -\tfrac{1}{4t_n}$. Then it follows from what is said in Part 3 that $\alpha(t, z=0, x_3=0) > 0$ on the whole interval $(0, t_\diamond)$. As a consequence of this (and because $\alpha$ is continuous): Given $t \in (0, t_\diamond]$, there exists a corresponding positive number $\rho_\diamond$ such that $\alpha > 0$ on the $|z| < \rho_\diamond t$ part of $\{t\} \times \mathbb{R}^2 \times \mathbb{R}$.

The positivity of $\alpha$ on the $|z| < \rho_\diamond t$ on $\{t\} \times \mathbb{R}^2 \times \mathbb{R}$ implies that $\langle \mathfrak{a}_3[\mathfrak{a}_1, \mathfrak{a}_2] \rangle$ is negative on the $0 < |z| < \rho_\diamond t$ part of $\{t\} \times \mathbb{R}^2 \times \mathbb{R}$ because $\langle \mathfrak{a}_3[\mathfrak{a}_1, \mathfrak{a}_2] \rangle = -2\alpha |\varphi|^2$.

**Appendix. The small t behavior of Kapustin-Witten instantons**

The assumption in what follows is that $(A, \mathfrak{a})$ is a solution to (1.4) that is described by (6.1) but not a model solution from Section 1c. This assumption is not stated explicitly in some of the subsequent lemmas. The integer *m* in what follows is the $z = 0$ vanishing degree of the corresponding version of $\varphi$. The integer *k* in what follows is the integer that appears in the $(A, \mathfrak{a})$ version of Proposition 8.2. Note in particular that *k* is greater than *m*. An upper bound for *k* can be given knowing only the numbers $\mathfrak{z}, \delta, r$ and $\Xi$ that appear in the $(A, \mathfrak{a})$ version of (6.1). See Lemma 5.6 in this regard.

**a) Small t behavior of $\alpha$**

The upcoming Lemma A.1 describes the behavior of the function $\alpha$ where $t$ is small and where $|z|$ is small relative to $t$. This lemma elaborates on the observation at the end of Part 3 in Section 9b. The second lemma describes the small $t$ part of the $\alpha = 0$ locus in $(0,\infty) \times \mathbb{R}^2 \times \mathbb{R}$. Keep in mind Lemma 9.2 which says in effect that $t\alpha$ can differ significantly from $t\alpha^{(k)}$ (the $(A^{(k)}, \mathfrak{a}^{(k)})$ version of $\alpha$) at small $t$ only where $\tfrac{|z|}{t}$ is very small.

**Lemma A.1**: *Given positive numbers $(\mathfrak{z}, \delta, r, \Xi)$, there exists $\kappa > 1$ with the following significance: Let $(A, \mathfrak{a})$ denote a solution to (1.4) that is described by the $(\mathfrak{z}, \delta, r, \Xi)$ version of (6.1) (and assume that $(A, \mathfrak{a})$ is not a model solution from Section 1c). Fix a sufficiently small positive time $t_* \in (0, 1)$, then $t \in (0, \frac{1}{\kappa} t_*]$, and then $z \in \mathbb{R}^2$ with $|z| \leq \frac{1}{\kappa} t$.*

- *If $\alpha > -\frac{(k+1 - \frac{1}{1000})}{2t}$ at $(t, z, 0)$, then $\alpha$ increases monotonically on $[t, \kappa t]$ (with $z$ fixed) and $\alpha$ differs from $\frac{k+1}{2s}$ by at most $\frac{1}{2000 s}$ at $(s, z, 0)$ if $s \in [\kappa t, t_*]$.*
- *If $(t, z) \in (0, t_*] \times \mathbb{R}^2$ and if $\alpha < \frac{(k+1 - \frac{1}{1000})}{2t}$ at $(t, z, 0)$, then $\frac{\partial}{\partial t} \alpha > 0$ on $[\frac{1}{\kappa} t, t]$ (with $z$ fixed) and $\alpha$ differs from $-\frac{k+1}{2s}$ by at most $\frac{1}{2000 s}$ at $(s, z, 0)$ if $s < \frac{1}{\kappa} t$.*

This lemma is proved momentarily. What follows directly is the promised lemma about the small t zero locus of $\alpha$.

**Lemma A.2** *Let $(A, \mathfrak{a})$ denote a solution to (1.4) that is described by (6.1). Assume that $(A, \mathfrak{a})$ is not a model solution from Section 1c. If $t_*$ is positive but sufficiently small, then the $\alpha = 0$ locus in $(0, t_*) \times \mathbb{R}^2 \times \{0\}$ is a non-empty, smooth, 2-dimensional submanifold with the projection map to $\{t_*\} \times \mathbb{R}^2 \times \{0\}$ being a diffeomorphism onto an open set. Said differently, the $\alpha = 0$ locus in $(0, t_*) \times \mathbb{R}^2 \times \{0\}$ is a graph over its image in $\{t_*\} \times \mathbb{R}^2 \times \{0\}$. In addition, there is an integer $k > m$ and a positive number $\rho_*$ such that*

- $\alpha > \frac{(k+1 - \frac{1}{1000})}{2t_*}$ *on the $|z| \leq \rho_* t_*$ part of $\{t_*\} \times \mathbb{R}^2 \times \{0\}$.*
- *The $0 < |z| \leq \rho_* t$ part of $\{t_*\} \times \mathbb{R}^2 \times \{0\}$ is in the projected image of the $(0, t_*) \times \mathbb{R}^2 \times \{0\}$ part of the $\alpha = 0$ locus.*

The rest of this subsection of the Appendix has the proofs of these lemmas.

***Proof of Lemma A.1***: Fix small but positive $\varepsilon$ and then a time $t_\varepsilon > 0$ so that (9.5) holds where $t \in (0, t_\varepsilon]$ and $|z| < \varepsilon t$. Write $\alpha$ as $\frac{y}{t}$ and, having fixed some small, positive time $t_* \leq t_\varepsilon$ and then $t_\ddagger < t_*$, let $\tau$ denote $\ln \frac{t}{t_\ddagger}$. The equation (9.5) when written using $y$ and $\tau$ as

$$\frac{\partial}{\partial \tau} y + 2y^2 - \frac{(k+1)^2}{2} = \mathfrak{e}_4 \ .$$
(A.1)

Let $\mu$ denote the least upper bound for $|\mathfrak{e}_4|$ where $t \in (0, t_*]$ and $|z| \leq \varepsilon t$. If $y^2 < \frac{(k+1)^2}{4} - \mu$ at time $\tau = 0$, then the equation in (A.1) has $y$ increasing until at least the first positive $\tau$ where $y^2 > \frac{(k+1)^2}{4} - c_0 \mu$. It is a straightforward exercise to prove that the time it takes for $y$ to come within any given specified distance of $\frac{(k+1)}{2} - c_0 \mu$ can be bounded from below knowing only $k$ and an upper bound for $\frac{(k+1)^2}{4} - \mu - x^2$ at the initial time $\tau = 0$ (thus, $t=t_\ddagger$). (The exercise can be solved as follows: Introduce by way of notation a positive number, $\lambda$, by setting $\lambda^2 = \frac{(k+1)^2}{4} - \mu$. Then (A.1) leads to the inequality

$$\tfrac{\partial}{\partial \tau} y \geq 2(\lambda^2 - y^2)$$

(A.2)

This can be solved using partial fractions to see that

$$y(\tau) \geq \lambda \frac{ce^{4\lambda\tau} - 1}{ce^{4\lambda\tau} + 1}$$

(A.3)

with $c$ a positive number that is determined by the size of $\lambda^2 - y^2$ at $\tau = 0$ (thus $t = t_\ddagger$.).)

To continue: If $y$ at $\tau = 0$ were to obey $y < -\frac{(k+1)}{2} - c_0\mu$, then (A.1) says that the derivative of $y$ is negative with norm getting ever larger. It is a straightforward exercise to prove that the maximal interval where the solution exists is determined by $y + \lambda$ at $\tau = 0$ with $\lambda$ now denoting the positive square root of $\frac{(k+1)^2}{4} + \mu$. (The solution in this case is such that $|y|(\tau) \geq \lambda \frac{ce^{4\lambda\tau} - 1}{ce^{4\lambda\tau} + 1}$ with $c$ again determined by the value of $y$ at $\tau = 0$.)

Finally, if $y$ at $\tau = 0$ were to obey $y > \frac{(k+1)}{2} + c_0\mu$, then (A.1) says that $y$ has negative derivative so it is driven back towards $\frac{(k+1)}{2}$ until such time as it is less than $\frac{(k+1)}{2} + c_0$. One can again obtain an upper bound for how long this takes in terms of the $\tau = 0$ value of $y - (\frac{(k+1)}{2} + c_0\mu)$.

The preceding properties of $x$ when reinterpreted as properties of $\alpha$ restate the bullets of the lemma.

**Proof of Lemma A.2**: The proof has two parts. Keep in mind for both parts that the $\alpha = 0$ locus in $\{t\} \times \mathbb{R}^2 \times \mathbb{R}$ is not empty when $t$ is small (but it sits necessarily where $\frac{|z|}{t} \ll 1$).

*Part 1*:. Fix $\varepsilon > 0$ but very small and then a time $t_\varepsilon > 0$ so that (9.6) holds where $t < t_\varepsilon$. Because the t-derivative of $\alpha$ is not zero on the $\alpha = 0$ locus, the differential of $\alpha$ is not zero there, which implies in turn that the $t < t_\Diamond$ part of the $\alpha = 0$ locus is a smooth manifold. The equation in (9.6) also implies the vector field $\frac{\partial}{\partial t}$ is never tangent to the $\alpha = 0$ locus where $t$ is small; and that implies in turn that the projection from the $\alpha = 0$ locus in $(0, t_*) \times \mathbb{R}^2 \times \{0\}$ to $\{t_*\} \times \mathbb{R}^2 \times \{0\}$ is 1-1 onto its image with non-zero differential. The positivity of $\frac{\partial}{\partial t} \alpha$ on the $\alpha = 0$ locus in the $t < t_*$ part of $(0, \infty) \times \mathbb{R}^2 \times \{0\}$ implies (via the inverse function theorem) that this locus can be depicted as a smooth graph over its projected image in $\{t_*\} \times \mathbb{R}^2 \times \{0\}$.

*Part 2*: The first and second bullets of the lemma follow from the third and fourth bullets because of what is said by Lemma 9.2. The specific fact from Lemma 9.2 is that its top bullet holds where $|z| > \varepsilon t$ for any given positive $\varepsilon$ if t is sufficiently small given $\varepsilon$. This implies that $\alpha < -\frac{(k+1 - \frac{1}{1000})}{2t}$ if $z \neq 0$ and t is sufficiently small given $|z|$. The latter fact (with an appeal to the third and fourth bullets of Lemma A.1) implies that the closure of the projected image of the $(0, t_*)$ part of the $\alpha = 0$ locus in $\{t_*\} \times \mathbb{R}^2 \times \{0\}$ must contain the $z = 0$ point and that $\alpha$ at this point must be greater than $\frac{(k+1 - \frac{1}{1000})}{2t_*}$. (The projected image must contain a sequence of points limiting to the $z = 0$ point if $t\alpha$ is not bounded away from zero at small t.) Because $\alpha$ is smooth, there is an open disk in $\{t_*\} \times \mathbb{R}^2 \times \{0\}$ containing the $z = 0$ point where $\alpha$ has this same lower bound. Any non-zero point in this disk must be in the projected image of the $t < t_*$ part of the $z = 0$ locus (which is also a consequence of the lemma's last two bullets and the top bullet of Lemma 9.2).

**b) More small t behavior of $\alpha$; and the behavior of $|\beta|, |\varphi|$ and $|\mathcal{b}|$**

Fix $\varepsilon > 0$. Then fix $t_*$ to be positive but much less than the corresponding $t_\varepsilon$ and invoke Lemma A.2 to define the number $\rho_*$. Let $\mu$ denote the least upper bound for the norm of the function $\mathfrak{e}_4$ in (A.1) when $t \in (0, t_*]$ and $|z| < \rho_* t_*$. If $\varepsilon$ is small and if $t_*$ is sufficiently small (given $\varepsilon$), then $\mu$ will be less than $\varepsilon^2$. Assume that this is the case.

Having fixed $z \in \mathbb{R}^2$ with $|z| < \rho_* t_*$, let $t_z \in (0, t_*)$ denote the (unique) time where $\alpha(\cdot, z, 0)$ is zero. The analysis that led to (A.3) can be repeated to see that the function $\alpha$ can be written as

$$\alpha = \frac{k+1}{2t} \left( \frac{t^{2(k+1)} - t_z^{2(k+1)}}{t^{2(k+1)} + t_z^{2(k+1)}} \right) + \mathfrak{e}\frac{1}{t}$$

(A.5)

with $\mathfrak{e}$ here and in what follows denoting a function of t and z that obeys $|\mathfrak{e}| < c_0 \varepsilon$. (The precise function that is depicted by $\mathfrak{e}$ can be different in each of $\mathfrak{e}$'s appearances.)

The preceding approximation for $\alpha$ can be used in conjunction with top bullet of (9.1) to obtain an approximation for $|\beta|$:

$$|\beta| = \frac{k+1}{t} \frac{t_z^{2(k+1)}}{t^{2(k+1)} + t_z^{2(k+1)}} + \mathfrak{e}\frac{1}{t} .$$

(A.6)

An approximation for the function $|\varphi|$ can also be obtained from (A.5) by integrating the equation in the top bullet of (3.9). In particular $|\varphi|$ at times $t < t_*$ for such z can be written as

$$|\varphi| = |\varphi|(t_z, z, 0) \left( \left(\tfrac{t_z}{t}\right)^{(k+1)} + \left(\tfrac{t}{t_z}\right)^{(k+1)} \right)(1 + \mathfrak{e})$$

(A.7)

To write this another way, note that if $|z|$ is sufficiently small, then $|\varphi|(t_*, z, 0)$ can be written as $\kappa |z|^m (1+\mathfrak{e})$ with $\kappa$ being positive and independent of z. Since $t_z$ has to be much less than $t_*$ when $|z|$ is small, this depiction of $|\varphi|$ at $(t_*, z, 0)$ can be used to write (A.7) as

$$|\varphi| = \kappa |z|^m \left(\tfrac{t_z}{t_*}\right)^{(k+1)} \left( \left(\tfrac{t_z}{t}\right)^{(k+1)} + \left(\tfrac{t}{t_z}\right)^{(k+1)} \right)(1 + \mathfrak{e})$$

(A.8)

Finally, an approximation for the function $|\mathfrak{b}|$ can be had by integrating (6.6) using the formula for $\alpha$ in (A.5) and the fact that $|E_A| \leq c_{\mathfrak{z}} \tfrac{1}{t^2}$:

$$|\mathfrak{b}| = |\mathfrak{b}|(t_z, z, 0) \left( \left(\tfrac{t_z}{t}\right)^{(k+1)} + \left(\tfrac{t}{t_z}\right)^{(k+1)} \right)^{-1} .$$

(A.9)

By way of a parenthetical remark, I expect that $\mathfrak{b}$ at very small t can be written with respect to a particular product structure on $\mathcal{L}^+$ as

$$\mathfrak{b} = \tfrac{p}{\sqrt{2}} \frac{t^{-(k+1)} \bar{z}^{p-1}}{(1 + t^{-2(k+1)} |z|^{2p})} + \mathfrak{e}$$

(A.10)

where $p = \tfrac{1}{2}(k - m)$.

c) **The size of $t_z$**

The depiction of $|\varphi|$ in (A.8) can be used to obtain an approximation for $t_z$. To do this, note that $|\varphi|$ where the size of t is on the order of the $|z|$ must be very close to $|\varphi^{(k)}|$ which has the form $(\tfrac{1}{t})^{(k+1)} |z|^k (1 + \mathcal{O}(\tfrac{|z|}{t}))$. Taking $t < \varsigma^{-1} |z|$ but with $\tfrac{|z|}{t}$ being none-the-less of size $\mathcal{O}(1)$ leads from (A.8) to the conclusion that

$$|z|^{k-m} \sim t_z^{2(k+1)}$$

(A.11)

with $\sim$ indicating that their ratio is bounded from above and below by positive, z independent numbers when z is very small. Writing $k = m + 2p$ as done in Theorem 2, this says that says $t_z \sim |z|^{p/(m+2p+1)}$.